%% file: Blow_up_.tex
\documentclass[11pt,a4paper]{article}
\usepackage[utf8]{inputenc}

\usepackage{setspace} 
\usepackage{shellesc}
\usepackage{mathrsfs,amsthm,xcolor,verbatim,bbm,amsmath,amsfonts,amssymb,nicefrac,enumitem,hyperref,bm,mathtools,xparse,etoolbox,textcomp,float}
\usepackage[margin=1in]{geometry}
\usepackage[sort,capitalise]{cleveref} 
\usepackage{xurl}
\usepackage{frontespizio}
\usepackage{breakurl}
\usepackage{subfigure}

\numberwithin{equation}{section}

\crefname{enumi}{item}{items}

\crefname{equation}{}{}
\crefname{subsection}{Subsection}{Subsections}

\hypersetup{
    colorlinks,
    linkcolor={red!50!black},
    citecolor={green},
    urlcolor={blue!80!black}
}
    
\theoremstyle{plain}
\newtheorem{theorem}{Theorem} [section]
\newtheorem{lemma}[theorem]{Lemma}
\newtheorem{prop}[theorem]{Proposition}
\newtheorem{cor}[theorem]{Corollary}

\newtheorem{setting}[theorem]{Setting}

\theoremstyle{remark}

\newcommand\scalemath[2]{\scalebox{#1}{\mbox{\ensuremath{\displaystyle #2}}}}

\theoremstyle{definition}

\usepackage{tikz}
\usetikzlibrary{matrix,chains,positioning,decorations.pathreplacing,arrows}
\usetikzlibrary{shapes,arrows}
\tikzset{
	font={\fontsize{9pt}{12}\selectfont}}
\usepackage{adjustbox}

\input{0_commands.tex}

\begin{document}

\title{Blow up phenomena for gradient descent optimization \\
methods in the training of artificial neural networks}

\author{Davide Gallon$^{1}$, Arnulf Jentzen$^{2,3}$, and Felix Lindner$^4$
\bigskip
    \\
	\small{$^1$University of Padua, Italy, e-mail: davide.gallon@studenti.unipd.it} 
	\smallskip
	\\
	\small{$^2$School of Data Science and Shenzhen Research Institute of Big Data,}
	\\
	\small{The Chinese University of Hong Kong, Shenzhen, China, e-mail: ajentzen@cuhk.edu.cn} 
	\smallskip
	\\
	\small{$^3$Applied Mathematics: Institute for Analysis and Numerics,}
	\\
	\small{University of M\"unster, Germany, e-mail: ajentzen@uni-muenster.de} 
	\smallskip
	\\
	\small{$^4$Faculty of Mathematics and Natural Sciences, University of Kassel,}
	\\
	\small{Germany, e-mail: lindner@mathematik.uni-kassel.de}
	\smallskip
	\\
}

\date{\today}

\maketitle

\begin{abstract}
\noindent
In this article we investigate blow up phenomena for gradient descent optimization methods in the training of artificial neural networks (ANNs). 
Our theoretical analysis is focused on shallow ANNs with one neuron on the input layer, one neuron on the output layer, and one hidden layer. 
For ANNs with ReLU activation and at least two neurons on the hidden layer we establish the existence of a 
target function such that there exists a lower bound for the risk values of the critical points of the associated risk function which is strictly greater than the infimum of the image of the risk function. 
This allows us to demonstrate that every gradient flow trajectory with an initial risk smaller than this lower bound diverges. 
Furthermore, we analyze and compare various popular types of activation functions with regard to 
the divergence of gradient flow trajectories and gradient descent trajectories 
in the training of ANNs and with regard to the closely related question concerning the existence of global minimum points of the risk function.
\end{abstract}
\tableofcontents

\section{Introduction}
%

While artificial neural networks (ANNs) are widely used and increasingly popular in a large variety of scientific and industrial applications, training methods for ANNs are still far from being well-understood from an analytical perspective.

In the training of ANNs one is ultimately interested in mimimizing the true risk, i.e., the expected loss of the realization function associated to the ANN. 
A natural direction of research aiming at a better theoretical understanding of 
such 
optimization problems concerns the analysis of the associated gradient flow (GF) differential equations. 
Loosely speaking, each GF trajectory represents a 
path of steepest descent in the risk landscape.
In order to 
render the concept of 
GFs useful for the 
practical 
training of ANNs, at least two 
types
of approximation have to be taken into account. 
First, the unknown true risk function has to be approximated by an empirical risk function based on the training data at hand. 
Second, the continuous-time GF has to be approximated by a discrete-time gradient descent (GD) 
optimization scheme. 
Discretization parameters associated to these 
types of approximation 
are the size of the training data set and the learning rate respectively.
A possible further approximation 
regarding  
the gradient of the empirical risk 
function 
by means of Monte Carlo estimation 
leads to the class of stochastic GD 
optimization methods, involving the batch size as an additional approximation parameter. 

In order to ensure that a GD optimization scheme produces trajectories which lie suitably close to the corresponding GF trajectories associated to the true risk, the discretization parameters, say, the learning rate and the size of the training data set, have to be chosen sufficiently small and large respectively. 
However, just how small or large the discretization parameters need to be chosen is generally relative to the object 
to be approximated, i.e., relative to the GF. 
In particular, if a GF trajectory is such that the norm of the ANN parameter vector specifying the realization function diverges to infinity, problems concerning an adequate choice of the 
discretization parameters specifying the approximation of the GF 
may arise. 
As a matter of fact, available results in the research literature 
pertaining to 
the convergence analysis of GFs and GD type optimization algorithms 
are typically based, 
either explicitly or implicitly, 
on suitable
boundedness assumptions 
concerning 
the GF and GD trajectories; 
see 
the literature overview in Subsection~\ref{subsec:intro_literature} below. 

It is the contribution of this article to 
uncover and analyze situations where 
the standard boundedness assumptions on GF and GD trajectories fail 
to hold and blow up phenomena occur instead. 
In our theoretical analysis we focus 
on shallow ANNs with one neuron in the input layer, one neuron in the output layer, and one hidden layer made up of  
an arbitrary number of 
neurons. 
In this introductory section we present four selected key results 
of this article 
to elucidate the scope of the study. 
The first key result, 
see \cref{theorem:blowup2} 
below, 
concerns ANNs with ReLU activation function and establishes the existence of a 
target function such that 
every GF trajectory with an initial risk
below a certain threshhold 
diverges, 
assuming that the hidden layer consists of at least two neurons. 
In the second 
key
result  of this article, 
see \cref{theorem:blowup} 
below, 
we show for 
various popular types of activation functions 
that there exist target functions such that GF and GD trajectories diverge whenever the associated risk values satisfy certain asymptotic optimality conditions
and the hidden layer consists of at least two neurons. 
From an analytical perspective, blow up phenomena 
in the training of ANNs are closely connected to the question whether there exist global minimum points of the risk function. 
Further key results of this article, see \cref{theorem:globalminima} and \cref{theorem:superiority} below, therefore concern the existence and non-existence of global minima of the risk function depending on the choice of the activation function and a superior role played  in this regard by piecewise affine activation functions.

In 
the remainder of 
this introductory section we provide the precise statements of 
the selected results mentioned above
together with some additional comments 
in Subsections~\ref{subsec:intro_blowup}--\ref{subsec:intro_superiority}, 
we give a short overview of related research literature in Subsection~\ref{subsec:intro_literature}, 
and we outline
the overall structure of this article in Subsection~\ref{subsec:intro_structure}.

\begin{figure}
	\centering
	\begin{adjustbox}{width=\textwidth}
		\begin{tikzpicture}[shorten >=1pt,->,draw=black!50, node distance=\layersep]
		\tikzstyle{every pin edge}=[<-,shorten <=1pt]
		\tikzstyle{input neuron}=[very thick, circle,draw=red, fill=red!30, minimum size=15pt,inner sep=0pt] 
		\tikzstyle{output neuron}=[very thick, circle,draw=green, fill=green!30, minimum size=20pt,inner sep=0pt]
		\tikzstyle{hidden neuron}=[very thick, circle, draw=blue, fill=blue!30, minimum size=15pt,inner sep=0pt]
		\tikzstyle{annot} = [text width=9em, text centered]
		\tikzstyle{annot2} = [text width=4em, text centered]
		\node[input neuron] (I) at (0,-2.5 cm) {$x$};
		\foreach \name / \y in {1,...,5}
		\path[yshift=0.5cm]
		node[hidden neuron] (H-\name) at (5 cm,-\y cm) {};      
		\path[yshift=-1cm]
		node[output neuron](O) at (10 cm,-1.5 cm) {$\realization{\theta}(x)$};   
		\foreach \dest in {1,...,5}
		\path[line width = 0.8] (I) edge (H-\dest);     
		\foreach \source in {1,...,5}
		\path[line width = 0.8] (H-\source) edge (O);     
		\node[annot,above of=H-1, node distance=1cm, align=center] (hl) {Hidden layer\\(2nd layer)};
		\node[annot, above of=I, node distance = 3cm, align=center] {Input layer\\ (1st layer)};
		\node[annot, above of=O, node distance=3cm, align=center] {Output layer\\(3rd layer)};		
		\node[annot2,below of=H-5, node distance=1cm, align=center] (sl) { $\width=5$};
		\end{tikzpicture}
	\end{adjustbox}
	\caption{Graphical illustration of the shallow ANN architecture considered in Theorems~\ref{theorem:blowup2}--\ref{theorem:superiority} in the special case of $h=5$ neurons on the hidden layer. In this situation every ANN parameter vector $\theta\in\R^{3h+1}=\R^{16}$ specifies a realization function $\realization{\theta}\colon\R\to\R$ which depends on the choice of the activation function $\bbA\colon\R\mapsto\R$ and maps the scalar input $x\in[\scra,\scrb]$ to the scalar output $\realization{\theta}(x)=\theta_{3\width+1} + \sum_{i=1}^{\width} \theta_{2\width + i } \bbA( \theta_{\width  + i}  +\theta_{i } x )\in\R $.}
	\label{figure_shallow}
\end{figure}
\subsection{Blow up phenomena}
\label{subsec:intro_blowup}
Before stating
the key results
we shortly comment on 
the mathematical setting and the employed notation. 
We consider ANNs with one neuron on the input layer, one neuron on the output layer, and one hidden layer made up of $\width \in \N$ neurons; compare the illustration in Figure~\ref{figure_shallow}. 
The specification of such an ANN involves $2\width$ real weight parameters and $\width+1$ real bias parameters, so that the overall number of ANN parameters amounts to $3\width+1$.
The real numbers $\scra \in \R$, $\scrb \in (\scra,\infty)$ define the domain of the target function $v \colon [\scra,\scrb]\to \R$ and thus the interval on which the realization function associated to the ANN should approximate $v$. 
 
Our first main result on blow up phenomena in the training of ANNs 
formulated
in \cref{theorem:blowup2} below 
focuses on ANNs with ReLU activation function $\R \ni x \mapsto \max\{x,0\}\in \R$ and 
establishes the existence of a non-decreasing target function $v \colon [\scra,\scrb]\to \R$ such that 
for every choice 
$\width \in \N\backslash\{1\}$ 
of the number of hidden neurons 
there exists a positive threshold value
such that every GF trajectory $\Theta \colon [0, \infty) \to \R^{3\width+1}$ with a initial risk smaller than this threshold value 
diverges to infinity in the sense that $\liminf_{t\to\infty}\|\Theta_t\|=\infty$.
Here we denote for every 
non-decreasing 
$v \colon [\scra,\scrb]\to \R$ and every 
$\width \in \N$ by $\cL^{v,\width} \colon  \R^{3\width+1} \to \R$ the risk function associated to the ANN measuring how well the realization function approximates the target function. 
\cref{theorem:blowup2} is a slightly simplified version of a more general result in \cref{H:theorem:nozero} in Subsection~\ref{1div} below. 
Note that the lack of differentiability of the ReLU activation function at zero entails a lack of differentiability of the risk function associated to the ANN, 
so that we need to work with an appropriate generalized gradient of the risk function 
to be able to specify the GF. 
While the general statement in \cref{H:theorem:nozero} 
is based on a generalized gradient defined in terms of 
continuously differentiable 
approximations of the ReLU activation function,
we simplify the exposition 
in \cref{theorem:blowup2} below 
by employing the less involved concept of the left gradient of the risk function, 
denoted for every non-decreasing 
$v \colon [\scra,\scrb]\to \R$ and every 
$\width \in \N$ by $\cG^{v,\width} \colon  \R^{3\width+1} \to \R^{3\width+1}$. 
In order to ensure that both concepts of generalized gradients coincide, 
we additionally assume in \cref{theorem:blowup2} 
that $\scra\ge0$.  
This additional assumption is not needed in the general statement in \cref{H:theorem:nozero} .

\begin{theorem}\label{theorem:blowup2}
Let $\scra\in[0,\infty)$, $\scrb \in (\scra,\infty)$
and
for every non-decreasing $v \colon [\scra,\scrb] \to \R$  and every $\width \in \N$ 
let $\cL^{v,\width} \colon  \R^{3 \width +1} \to \R$
satisfy for all $\theta = (\theta_1, \ldots, \theta_{3 \width + 1}) \in \R^{3\width+1}$ that
\begin{equation}
    \cL^{v,\width} ( \theta ) = \int_{\scra}^{\scrb} \rbr[\big]{ v(x) -  \theta_{3\width+1} - \smallsum_{i=1}^{\width} \theta_{2\width + i } \br[\big]{ \max\{ \theta_{\width  + i}  +\theta_{i } x ,0\} } }^2  \, \d x
\end{equation}
and let $\cG^{v,\width}\colon \R^{3\width+1} \to \R^{3\width+1}$ be the left gradient\footnote{For every $h\in\N$ let $e_1^{(h)},\ldots,e_{3h+1}^{(h)}\in\R^{3h+1}$ satisfy $e_1^{(h)}=(1,0,\ldots,0), \,\ldots\,, e_{3h+1}^{(h)}=(0,\ldots,0,1)$. \Nobs that for all non-decreasing $v\colon[\scra,\scrb]\to\R$ and all $\width\in\N$, $\theta\in\R^{3\width+1}$ it holds that $\cG^{v,h}(\theta)=(\lim_{\varepsilon \nearrow 0}(\cL^{v,\width}(\theta+\varepsilon e^{(h)}_1)-\cL^{v,\width}(\theta))\varepsilon^{-1}, \ldots, \lim_{\varepsilon \nearrow 0}(\cL^{v,\width}(\theta+\varepsilon e^{(h)}_{3\width+1})-\cL^{v,\width}(\theta))\varepsilon^{-1}).$
}
of $\cL^{v,\width}$.
 Then\footnote{\Nobs that the function 
$ \norm{ \cdot } \colon ( \cup_{ n \in \N } \R^n ) \to \R $ 
satisfies for all $ n \in \N $,
$ x = ( x_1, \dots, x_n ) \in \R^n $
that $ \norm{x} = ( \sum_{ i = 1 }^n | x_i |^2 )^{ 1 / 2 } $.}  there exists a non-decreasing $v\colon [\scra,\scrb] \to \R$ such that for all $\width \in \N \backslash\{1\}$ there exists $\varepsilon \in (0,\infty)$ such that for all 
    $\Theta \in C([0, \infty) , \R^{3\width+1})$ with  $\forall \, t \in [0, \infty) \colon \Theta_t = \Theta_0 - \int_0^t \cG^{v,\width} ( \Theta_s ) \, \d s$ and $\cL^{v,\width}(\Theta_0)<\varepsilon+\inf_{ \theta \in \R^{ 3 \width + 1 } } \cL^{v,\width}( \theta )$
    it holds that  $\liminf_{ t \to \infty } \| \Theta_t \|= \infty$.
\end{theorem}
\cref{theorem:blowup2} is a direct consequence of \cref{H:theorem:nozero} 
in Subsection~\ref{1div} below. 
\cref{H:theorem:nozero}, in turn, follows from combining \cref{theo:intro:convergence}
in Subsection~\ref{1div} below, 
which is a slight modification of \cite[Theorem 1.3]{JentzenRiekert2021aa}, 
and \cref{H:theorem:zerograd} 
in Subsection~\ref{1div} below,  
which is one of the main results of this paper. 
Roughly speaking, \cref{theo:intro:convergence} states that 
every GF trajectory which does not diverge to infinity converges to a critical point of the risk function and the risk values associated to the GF trajectory converge to the risk of the critical point.
\cref{H:theorem:zerograd}, rather, asserts  
for the case of the indicator function 
$\mathbbm{1}_{(\nicefrac{(\scra+\scrb)}2,\scrb]}\colon [\scra,\scrb]\to \R$ 
being the target function that 
there exists a lower bound $\varepsilon \in (0,\infty)$ for the risk of critical points.  
In this situation, \cref{H:theorem:nozero} infers the
divergence of every GF trajectory with an initial risk smaller than $\varepsilon$, 
which readily implies \cref{theorem:blowup2} in view of the fact that the different concepts of generalized gradients used in these results are compatible and in view of the fact that the infimum of the risk values associated to the 
considered 
target function f
equals zero in the case of at least two neurons in the hidden layer; compare, e.g., \cref{prop:relu:limzero}.

Our second main result on blow up phenomena in the training of ANNs formulated in \cref{theorem:blowup} below treats various types of activation functions and establishes for each type the existence of a target function such that all GF and GD trajectories which fulfill certain asymptotic optimality conditions 
w.r.t. the associated risk values 
diverge to infinity. 
Observe that the family of functions $A_{k,\gamma}\colon\R\to\R$, $k\in\Z$, $\gamma\in\R$, appearing in 
\cref{theorem:blowup} 
is such that depending on the choice of $k\in\Z$, $\gamma\in\R$ we have that $A_{k,\gamma}\colon\R\to\R$ refers to 
the softsign activation function in the case $k<-5$, 
the arctangent activation function in the case $k=-5$, 
the inverse square root unit activation function with parameter $\xi\in(0,3)$ in the case $k=-4$, 
the exponential linear unit activation function in the case $k=-3$, 
the hyperbolic tangent activation function in the case $k=-2$, 
the logistic activation function in the case $k=-1$, 
the softplus activation function in the case $k=0$, 
the ReLU activation function in the case $k=1$, $\gamma=0$, 
the leaky ReLU activation function with parameter $\gamma$ in the case $k=1$, $\gamma\in(0,1)$, 
and the rectified power unit activation function with exponent $k$ in the case $k>1$. 
Here we denote for every Lebesgue square integrable target function $v\colon[\scra,\scrb]\to\R$ and every $\width\in\N$, $k\in\Z$, $\gamma\in\R$ the associated risk function and its left gradient by $\cL_{k,\gamma}^{v,\width} \colon  \R^{3\width+1} \to \R$ and $\cG_{k,\gamma}^{v,\width} \colon  \R^{3\width+1} \to \R^{3\width+1}$ respectively. 
Similar to \cref{theorem:blowup2} above, the 
assertions in \cref{theorem:blowup} are slightly simplified versions of more general results employing a generalized gradient of the risk function defined in terms of continuously differentiable approximations of the activation function in \cref{generalframe} below.

\begin{theorem}\label{theorem:blowup}
Let $\scra\in[0,\infty)$, $\scrb \in (\scra,\infty)$, $\xii \in (0,3)$,
for every $k \in \Z $, $\gamma \in \R$
let $A_{k, \gamma} \colon \R \to \R$, satisfy for all $x \in \R$ that
\begin{equation}\label{eq:activation_family}
    A_{k,\gamma}(x)=
    \begin{cases}
           x(1+|x|)^{-1}, & \colon k <-5\\
           \arctan (x) & \colon  k=-5  \\
           x (1+\xii x^2)^{-\nicefrac12}  & \colon k=-4 \\
            x \indicator{(0,\infty)}(x)+(\exp(x)-1)\indicator{(-\infty,0]}(x) & \colon k =-3 \\
           (\exp(x)-\exp(-x))(\exp(x)+\exp(-x))^{-1} & \colon k=-2  \\
          (1+\exp(-x))^{-1} & \colon k=-1 \\
          \ln (1+\exp(x))   & \colon k =0 \\
          (\max\{x,0\})^k+\min\{\gamma x,0\} & \colon k >0,
    \end{cases} 
\end{equation}
and for every Lebesgue square integrable $v \colon [\scra,\scrb] \to \R$ and every $\width \in \N$,
 $k \in \Z$, $\gamma \in \R$ 
let $\cL_{k,\gamma}^{v,\width} \colon  \R^{3\width+1} \to \R$
satisfy for all $\theta= (\theta_1, \ldots, \theta_{3 \width + 1}) \in \R^{3\width+1}$ that
\begin{equation}
    \cL_{k,\gamma}^{v,\width} ( \theta ) = \int_{\scra}^{\scrb} \rbr[\big]{ v(x) -  \theta_{3\width+1} - \smallsum_{i=1}^{\width} \theta_{2\width + i } \br[\big]{ A _{k,\gamma} ( \theta_{\width  + i}  +\theta_{i } x ) } }^2  \, \d x
\end{equation}
and let $\cG^{v,\width}_{k,\gamma} \colon \R^{3\width+1} \to \R^{3\width+1}$ be the left gradient of $\cL^{v,\width}_{k,\gamma}$.
Then
\begin{enumerate} [label=(\roman*)]
    \item \label{theorem:blowupit1}
    there exists a polynomial $v\colon [\scra,\scrb] \to \R$ such that
    for all $\width \in \N\backslash\{1\}$, $k \in \Z \backslash \N$,
    $\Theta 
    \in C([0,\infty), \R^{3\width+1})$ with $\liminf_{ t \to \infty } \cL_{k,0}^{v,\width}( \Theta_t ) = \inf_{ \theta \in  \R^{ 3\width+1 } } \cL_{k,0}^{v,\width}( \theta )$ and $\forall \, t \in [0, \infty ) \colon \Theta_t = \Theta_0 - \int_0^t \cG_{k,0}^{v,\width} ( \Theta_s ) \, \d s$ 
    it holds that  $\liminf_{ t \to \infty } \| \Theta_t \| = \infty$,
    \item \label{theorem:blowupit1b}
    there exists a polynomial $v\colon [\scra,\scrb] \to \R$ such that
    for all $\width \in \N\backslash\{1\}$, $k \in \Z \backslash \N$ and all
    $\Theta \colon \N_0 \to \R^{3\width+1}$ with $\limsup_{ n \to \infty } \cL_{k,0}^{v,\width}( \Theta_n ) = \inf_{ \theta \in  \R^{ 3\width+1 } } \cL_{k,0}^{v,\width}( \theta )$ 
    it holds that  $\liminf_{ n \to \infty } \| \Theta_n \|\allowbreak = \infty$,
    \item \label{theorem:blowupit2}
     for all $k \in \Z \backslash\{1\}$ there exists a Lipschitz continuous $v \colon [\scra,\scrb] \to \R$ such that for all $\width \in \N\backslash\{1\}$,
    $\Theta 
    \in C([0,\infty), \R^{3\width+1})$ with $\liminf_{ t \to \infty } \cL_{k,0}^{v,\width}( \Theta_t ) = \inf_{ \theta \in  \R^{ 3\width+1 } } \cL_{k,0}^{v,\width}( \theta )$ and $\forall \, t \in [0, \infty ) \colon \Theta_t = \Theta_0 - \int_0^t \cG_{k,0}^{v,\width} ( \Theta_s ) \, \d s$ 
    it holds that  $\liminf_{ t \to \infty } \| \Theta_t \| = \infty$,
    \item \label{theorem:blowupit2b}
    for all $k \in \Z \backslash\{1\}$ there exists a Lipschitz continuous $v \colon [\scra,\scrb] \to \R$ such that for all $\width \in \N\backslash\{1\}$ and all
    $\Theta \colon \N_0 \to \R^{3\width+1}$ with $\limsup_{ n \to \infty } \cL_{k,0}^{v,\width}( \Theta_n ) = \inf_{ \theta \in  \R^{ 3\width+1 } } \cL_{k,0}^{v,\width}( \theta )$  
    it holds that  $\liminf_{ n \to \infty } \| \Theta_n \| = \infty$,
    \item \label{theorem:blowupit3}
    there exists a non-decreasing $v\colon [\scra,\scrb] \to \R$ such that for all $\width \in \N \backslash\{1\}$, $\gamma \in \R\backslash\{1\}$,
    $\Theta 
    \in C([0,\infty), \R^{3\width+1})$ with $\liminf_{ t \to \infty } \cL_{1,\gamma}^{v,\width}( \Theta_t ) = \inf_{ \theta \in  \R^{ 3\width+1 } } \cL_{1,\gamma}^{v,\width}( \theta )$ 
     and $\forall \, t \in [0, \infty ) \colon \Theta_t = \Theta_0 - \int_0^t \cG_{1,\gamma}^{v,\width} ( \Theta_s ) \, \d s$ 
    it holds that  $\liminf_{ t \to \infty } \| \Theta_t \| = \infty$, and
    \item \label{theorem:blowupit3b}
    there exists a non-decreasing $v\colon [\scra,\scrb] \to \R$ such that for all $\width \in \N \backslash\{1\}$, $\gamma \in \R\backslash\{1\}$ and all
    $\Theta \colon \N_0 \to \R^{3\width+1}$
    with $\limsup_{ n \to \infty } \cL_{1,\gamma}^{v,\width}( \Theta_n ) = \inf_{ \theta \in  \R^{ 3\width+1 } } \cL_{1,\gamma}^{v,\width}( \theta )$ 
    it holds that  $\liminf_{ n \to \infty } \| \Theta_n \| = \infty$.
\end{enumerate}
\end{theorem}

\Cref{theorem:blowupit1},  \cref{theorem:blowupit2}, and  \cref{theorem:blowupit3} in \cref{theorem:blowup} are direct consequences of \cref{theorem:globalminima:proof}, \cref{theorem:globalminima:proof2}, and \cref{newtechnique:risk} in \cref{subdivergence} below.  \Cref{theorem:globalminima:proof}, \cref{theorem:globalminima:proof2}, and \cref{newtechnique:risk}, in turn, are based on non-existence results concerning global minima of the risk function, compare \cref{theorem:globalminima} below, and an abstract divergence result for GF trajectories in \cref{lemma:new} in \cref{subdivergence} below.
 \Cref{theorem:blowupit1b},  \cref{theorem:blowupit2b}, and  \cref{theorem:blowupit3b} in \cref{theorem:blowup} are direct consequences of \cref{theorem:globalminima:proofb}, \cref{theorem:globalminima:proof2b}, and \cref{newtechnique:riskb} in \cref{subdivergenceGD} below.  \Cref{theorem:globalminima:proofb}, \cref{theorem:globalminima:proof2b}, and \cref{newtechnique:riskb}, in turn, are based on non-existence results concerning global minima of the risk function, compare \cref{theorem:globalminima} below, and an abstract divergence result for GD trajectories in \cref{lemma:new2} in \cref{subdivergenceGD} below. 
Related results can be found in \cite[Proposition 3.6]{petersen}.
%
%
%
\subsection{Existence and non-existence of global minima}
\label{subsec:intro_existence}
The analysis of blow up phenomena for GFs and GD optimization methods in the training of ANNs  is closely related to the question whether there exist global minimum points of the risk function associated to the ANN. In fact, the divergence results in \cref{theorem:blowup} above heavily rely on the non-existence of global minimum points of the risk function for certain target functions. 
In our third key result formulated in \cref{theorem:globalminima} below we consider the activation functions introduced in \cref{theorem:blowup} and establish the existence of several target functions $v\colon [\scra,\scrb]\to \R$ such that for every choice $\width \in \N\backslash\{1\}$ of the number of hidden neurons and for specific choices of $k \in \Z$, $\gamma \in \R$ the set of global minimum points
    $\mathscr{M}_{ k, \gamma }^{v,\width}=\{ \theta \in \R^{3\width+1 }\colon \cL^{v,\width}_{k,\gamma}(\theta ) = 
    \inf\nolimits_{ \vartheta \in \R^{3\width+1} } \cL^{v,\width}_{k,\gamma}(\vartheta )
    \}$
is empty. 
In particular, in the case of 
softsign, arctangent, inverse square root unit, exponential linear unit, hyperbolic tangent, standard logistic, and softplus activation we employ a polynomial target function, 
in the case of rectified power unit activation we employ a Lipschitz continuous target function, and in the case of ReLU and leaky ReLU activation we employ a non-decreasing target function. 

\begin{theorem}\label{theorem:globalminima}
Let $\scra\in \R$, $\scrb \in (\scra,\infty)$, $\xii \in (0,3)$,
for every $k \in \Z$, $\gamma \in \R$
let 
$A_{k, \gamma} \colon \R \to \R$ 
satisfy for all $x \in \R$ that
\begin{equation}
    A_{k,\gamma}(x)=
    \begin{cases}
           x(1+|x|)^{-1}, & \colon k <-5  \\
          \arctan (x) & \colon k=-5\\
          x (1+\xii x^2)^{-\nicefrac12}  & \colon k=-4 \\
          x \indicator{(0,\infty)}(x)+(\exp(x)-1)\indicator{(-\infty,0]}(x) & \colon k=-3 \\
           (\exp(x)-\exp(-x))(\exp(x)+\exp(-x))^{-1} & \colon k=-2  \\
          (1+\exp(-x))^{-1} & \colon k=-1  \\
          \ln (1+\exp(x))   & \colon k =0 \\
          (\max\{x,0\})^k+\min\{\gamma x,0\} & \colon k >0,
    \end{cases} 
\end{equation}
and for every measurable $v \colon [\scra,\scrb] \to \R$ and every $\width \in \N$,
$k \in \Z$, $\gamma \in \R$ 
let $\cL^{v,\width}_{k,\gamma} \colon  \R^{3\width+1} \to [0,\infty]$
satisfy for all $\theta= (\theta_1, \ldots, \theta_{3\width+1}) \in \R^{3\width+1}$ that
\begin{equation}
    \cL^{v,\width}_{k,\gamma} ( \theta ) = \int_{\scra}^{\scrb} \rbr[\big]{ v(x) -  \theta_{3\width+1} - \smallsum_{i=1}^{\width} \theta_{2\width + i } \br[\big]{ A_{k,\gamma} ( \theta_{\width  + i}  +\theta_{i } x ) } }^2  \, \d x
\end{equation}
and let $\mathscr{M}_{ k, \gamma }^{v,\width}\subseteq \R^{3\width+1}$ satisfy  $\mathscr{M}_{ k, \gamma }^{v,\width}= \{ \theta \in \R^{3\width+1 }\colon \cL^{v,\width}_{k,\gamma}(\theta ) = 
    \inf\nolimits_{ \vartheta \in \R^{3\width+1} } \cL^{v,\width}_{k,\gamma}(\vartheta )
    \} $.
Then
\begin{enumerate} [label=(\roman*)]
    \item \label{theorem:globalminima:negk}
    there exists a polynomial $v\colon [\scra,\scrb] \to \R$ such that 
    $\cup_{\width \in \N\backslash\{1\},\, k \in \Z \backslash \N } \mathscr{M}^{v,\width}_{k,0} = \emptyset$,
    \item \label{theorem:globalminima:lips}
    it holds for all $k \in \Z \backslash\{1\}$ that there exists a Lipschitz continuous $v \colon [\scra,\scrb] \to \R$ such that
    $\cup_{\width \in \N\backslash\{1\}}\mathscr{M}^{v,\width}_{k,0} = \emptyset$,
    \item \label{theorem:globalminima:nondecr}
    there exists a non-decreasing $v\colon [\scra,\scrb] \to \R$ such that 
    $\cup_{\width \in \N\backslash\{1\},\, \gamma\in \R\backslash\{1\}}\mathscr{M}^{v,\width}_{1,\gamma}
    = \emptyset,$ and
    %
    \item \label{theorem:globalminima:other} 
    it holds for all Lipschitz continuous $v \colon [\scra,\scrb] \to \R$ and all $\width\in \N$  that
    $\mathscr{M}^{v,\width}_{1,0} \neq \emptyset$.
\end{enumerate}
\end{theorem}

\Cref{theorem:globalminima:negk} and \cref{theorem:globalminima:lips} in \cref{theorem:globalminima} follow directly from combining \cref{H:noinf:soft:k2} in \cref{soft}, \cref{H:noinf:log2} in \cref{some}, \cref{H:noinf:repu} in \cref{repu}, \cref{H:noinf:elu} in \cref{elu}, and \cref{H:noinf:softsign} in \cref{softsign} below. 
\Cref{theorem:globalminima:nondecr} in \cref{theorem:globalminima} is a direct consequence of \cref{H:noinf} in \cref{ReLU} below.
The strategy in the proofs of these results is to identify the infimum of the image of the risk function and to consequently show that the set of global minimum points is empty. 
\Cref{theorem:globalminima:other}  in \cref{theorem:globalminima} has been proven in \cite[Theorem 1.1]{JentzenRiekert2021aa}. 
Related results can be found in \cite[Theorem 3.3]{petersen}.
\subsection{Superiority of piecewise affine activation functions with respect to the existence of minimum points}
\label{subsec:intro_superiority}
In practical applications of ANNs the choice of the activation functions is typically guided by heuristic arguments and numerical experiments. 
\Cref{theorem:globalminima} above suggest from an analytic perspective a superiority of piecewise affine activation functions with respect to the existence of minimum points of the risk function. 
Indeed, the non-existence results for global minima of the risk function in \cref{theorem:globalminima:negk} and \cref{theorem:globalminima:lips} in \cref{theorem:globalminima} are based on polynomial and Lip-schitz continuous target functions and exclusively involve continuously differentiable activation functions, whereas the existence result for global minima of the risk function in \cref{theorem:globalminima:other}  in \cref{theorem:globalminima} holds for all all Lipschitz continuous target functions and involves the piecewise affine ReLU activation function. 
This aspect is further highlighted in \cref{theorem:superiority} below. 

\begin{theorem}\label{theorem:superiority}
Let $\width\in\N\backslash\{1\}$, $\scra\in \R$, $\scrb \in (\scra,\infty)$, $\xii \in (0,3)$,
for every $k \in \Z$ 
let $A_{k} \in C(\R,\R)$ satisfy for all $x \in \R$ that
\begin{equation}
    A_{k}(x)=
    \begin{cases}
           x(1+|x|)^{-1}, & \colon k <-5  \\
          \arctan (x) & \colon k=-5\\
          x (1+\xii x^2)^{-\nicefrac12}  & \colon k=-4 \\
          x \indicator{(0,\infty)}(x)+(\exp(x)-1)\indicator{(-\infty,0]}(x) & \colon k=-3 \\
           (\exp(x)-\exp(-x))(\exp(x)+\exp(-x))^{-1} & \colon k=-2  \\
          (1+\exp(-x))^{-1} & \colon k=-1  \\
          \ln (1+\exp(x))   & \colon k =0 \\
          (\max\{x,0\})^k& \colon k >0,
    \end{cases} 
\end{equation}
let $k\in\Z$, 
and for every measurable $v \colon [\scra,\scrb] \to \R$ 
let $\cL^{v} \colon  \R^{3\width+1} \to [0,\infty]$
satisfy for all $\theta= (\theta_1, \ldots, \theta_{3\width+1}) \in \R^{3\width+1}$ that
\begin{equation}
    \cL^{v} ( \theta ) = \int_{\scra}^{\scrb} \rbr[\big]{ v(x) -  \theta_{3\width+1} - \smallsum_{i=1}^{\width} \theta_{2\width + i } \br[\big]{ A_{k} ( \theta_{\width  + i}  +\theta_{i } x ) } }^2  \, \d x.
\end{equation}
Then the following three statements are equivalent:
\begin{enumerate} [label=(\roman*)]
    %
    %
    \item \label{theorem:superiority:lip}
    It holds for every Lipschitz continuous $v\colon [\scra,\scrb] \to \R$ that there exists $\theta\in\R^{3\width+1}$ such that $\cL^v(\theta)=\inf_{\vartheta\in\R^{3\width+1}}\cL^v(\vartheta)$.
    \item \label{theorem:superiority:notC1}
    It holds that $A_k\notin C^1(\R,\R)$.
    \item \label{theorem:superiority:k=1}
    It holds that $k=1$. 
    \end{enumerate}
\end{theorem}

\Cref{theorem:superiority} is a direct consequence of \cref{theorem:globalminima} and the elementary observation that the ReLU activation function is the only activation function appearing in \cref{theorem:superiority} which is not continuously differentiable.
\subsection{Literature overview}
\label{subsec:intro_literature}
Let us complement the presentation of the findings of our work by a short review of related research literature. 
Despite the lack of a full-fledged convergence analysis for GFs and GD optimization schemes in the training of ANNs in literature, there are several promising mathematical approaches. 
For the convergence of GFs and GD type methods in the case of convex target functions 
we refer, e.g, to \cite{JentzenKuckuckNeufeldVonWurstemberger2021, BachMoulines2013,BachMoulines2011} and the references mentioned therein.
More complications are encountered in the case of non-convex problems: in principle there could be many local minima. For more details
on abstract convergence results for GD and GF optimization methods
we refer, e.g., to \cite{AkyildizSabanis2021, FehrmanGessJentzen2020,LeiHuLiTang2020,LovasSabanis2020, BertsekasTsitsiklis2000, DereichMuller_Gronbach2019}.

Another promising direction of research considers the overparametrized regime, where the number of parameters of the model exceeds the number of training points; see, e.g., \cite{JentzenKroeger2021,AllenzhuLiLiang2019,AllenzhuLiSong2019,DuLeeLiWangZhai2019,DuZhaiPoczosSingh2018arXiv,EMaWu2020,LiLiang2019}.
Under Lojasiewicz type assumptions convergence results for GD and GF type optimization schemes can be found, e.g, in \cite{Jentzen2,LeiHuLiTang2020,DereichKassing2021,AbsilMahonyAndrews2005,AttouchBolte2009}.
A further interesting method is to consider only special target functions; see, e.g.,\cite{JentzenRiekert2021,CheriditoJentzenRiekert2021}  for a convergence analysis for GF and GD processes in the case of constant
target functions and 
\cite{JentzenRiekert2021a}
for a convergence analysis for GF and GD processes in the training
of ANNs with piecewise linear target functions. 

For lower bounds and divergence results for GD and GF optimization methods we refer, e.g., to \cite{JentzenvonWurstemberger2020,CheriditoJentzenRossmannek2020,LuShinSuKarniadakis2020}. 
Results related to the findings of the present article can be found in \cite[Section~3]{petersen}.
\subsection{Structure of the article}
\label{subsec:intro_structure}
The remainder of this article is structured as follows.
In \cref{Annsrelu} we present the mathematical framework used to prove \cref{theorem:blowup2}, we establish several properties for critical points of the risk function in order to find a strictly positive lower bound for the risk values of critical points, and 
we demonstrate that the considered GF trajectories blow up. In \cref{generalframe} we introduce the mathematical framework needed for the proof of \cref{theorem:blowup} and \cref{theorem:globalminima}, we establish the non-existence of global minima employing various target functions, and we finally verify the assiciated blow up phenomena, thus proving \cref{theorem:globalminima} and \cref{theorem:blowup}.
In \cref{discretemeasure} we complement our findings be investigating the non-existence of global minima in the case in which the risk is defined using a discrete measure, the activation function is the standard logistic function, and there is one neuron in the hidden layer.

\section{Blow up phenomena for gradient flows (GFs) in the training of artificial neural networks (ANNs) with ReLU activation}\label{Annsrelu}
%
In this section we investigate blow up phenomena for GFs in the training of shallow ANNs with ReLU activation function in the case where the target function is given by the indicator function 
$\mathbbm{1}_{(\nicefrac{(\scra+\scrb)}2,\scrb]}\colon [\scra,\scrb]\to \R$.  
In particular, in \cref{H:theorem:nozero} in \cref{1div} below we demonstrate that every GF trajectory $\Theta\colon[0,\infty)\to\R^{3h+1}$ with an initial risk 
smaller than 
a certain threshold diverges to infinity in the sense that $\liminf_{t\to\infty}\|\Theta_t\|=\infty$. \cref{theorem:blowup2} in the introduction is a direct consequence of \cref{H:theorem:nozero}.

The two main ingredients in our proof of \cref{H:theorem:nozero} are \cref{theo:intro:convergence} and \cref{H:theorem:zerograd} 
in Subsection~\ref{1div} below. 
\cref{theo:intro:convergence} is a slight modification of \cite[Theorem 1.3]{JentzenRiekert2021aa} and states that for every GF trajectory $\Theta\colon[0,\infty)\to\R^{3h+1}$ with $\liminf_{t\to\infty}\|\Theta_t\|<\infty$ there exists $\beta\in (0,\infty)$ such that  
the GF trajectory converges with order $\beta$ to a critical point of the risk function and 
the risk values associated to the GF trajectory converge with order $1$ to the risk of the critical point. 
\cref{H:theorem:zerograd} is one of the main results of this article and establishes
a positive lower bound $\varepsilon \in (0,\infty)$ for the risk of critical points. 
In the proof of \cref{H:theorem:nozero} we combine \cref{theo:intro:convergence},  \cref{H:theorem:zerograd}, and the well-known fact that the risk values associated to a GF trajectory are non-increasing (see, e.g., \cite[Lemma 3.1]{Jentzen}) to conclude that every GF trajectory with an initial risk smaller than $\varepsilon$ 
diverges to infinity. 
The positive lower bound for the risk of critical points in \cref{H:theorem:zerograd} is based on an analogous result for the specific case $[\scra,\scrb]=[0,1]$ in \cref{lemma:H:bound} in Subsection~\ref{1div} in combination with an affine coordinate transformation. 
\cref{lemma:H:bound}, in turn, is proved by induction w.r.t.\ the number of hidden neurons $h\in\N$ and relies on a series of auxiliary results in Subsections~\ref{Generalresults1}--\ref{Annsrelu_esti} below. 

In \cref{Annsrelu_esti} we provide in \cref{lemma:H:1} and \cref{lemma:H:active} the base step and the induction step for our proof of the positive lower bound for the risk of critical points in \cref{lemma:H:bound}. 
More precisely, in \cref{lemma:H:active} we establish a positive lower bound for the risk of critical points in the case in which all hidden neurons are active and no combination of parameters allows for a representation of the same realization function using less neurons. 
\cref{lemma:H:active} is built on a detailed analysis of all possible parameter constellations of such critical points in \cref{lemma:H:q0},   \cref{lemma:H:qq}, \cref{lemma:H:qqneg}, \cref{lemma:H:case4sx}, \cref{lemma:H:case4dx}, \cref{lemma:H:updown}, \cref{lemma:H:N0:finN}, and \cref{lemma:H:N0:finN2} in \cref{Annsrelu_esti}. 
These results, in turn, employ 
properties of critical points 
of the risk function 
derived in \cref{proporties_crit} below and several elementary and well-known estimates and conclusions regarding specific integrals associated to the risk function provided in \cref{Generalresults1} and \cref{Generalresults2} below.

In \cref{subsec:descrANNReLU} below we specify in \cref{H:setting:snn} the mathematical framework regarding the training of shallow ANNs with ReLU activation used repeatedly throughout this section. 
Note that here for every ANN parameter vector $\theta\in\R^{3h+1}$ and every index $j\in\{1,2,\ldots,\width\}$ associated to a hidden neuron with weight parameter $\w{\theta}_{j}\neq0$ we have that the real number $\q{\theta}_{j}\in\R$ represents a possible breakpoint of the piecewise affine realization function $\realization{\width,\theta}_{\infty}\colon\R\to\R$ associated to the ANN.

Finally, in \cref{subsectionestimate} and \cref{subsec:upper_bounds} below we complement our findings in this section by providing an explicit lower bound for the risk of critical points in the specific case of $\width=2$ hidden neurons in \cref{lemma:2} in \cref{subsectionestimate} and by providing a general upper bound for the norms of GF trajectories in \cref{prop:upperbound} in \cref{subsec:upper_bounds}.

\subsection{Mathematical description of ANNs}
\label{subsec:descrANNReLU}
\begin{setting} \label{H:setting:snn}
For every $ \width \in \N$ let $\fd_\width \in \N$ satisfy $\fd_\width = 3 \width + 1$,
for every $\width\in \N$, $\theta  = ( \theta_1 ,  \ldots, \theta_{\fd_\width}) \in \R^{\fd_\width}$
let $\w{\theta} _ {0}, \w{\theta} _ {1}, \ldots,\w{\theta}_{\width+1},  \b{\theta} _ 1,  \b{\theta} _ 2 \ldots, \b{\theta}_{\width}, \v{\theta} _ 1, \v{\theta} _ 2, \ldots, \v{\theta}_{\width}, \c{\theta}, \q{\theta}_0, \q{\theta}_1, \ldots,\q{\theta}_{\width+1} \in (- \infty , \infty]$
 satisfy for all $j \in \{1,2, \ldots, \width \}$ that $\w{\theta}_{ 0}=-\w{\theta}_{ \width+1}=-1$, $\q{\theta}_{ 0}=1-\q{\theta}_{ \width+1}=0$, $\w{\theta}_{ j} = \theta_{ j}$, $\b{\theta}_j = \theta_{\width + j}$, 
$\v{\theta}_j = \theta_{2 \width + j}$, 
$\c{\theta} = \theta_{\fd_\width}$, and
\begin{equation} \label{H:setting:eq:defqi}
    \q{\theta}_j = \begin{cases}
    - \nicefrac{\b{\theta}_j}{\w{\theta}_j} & \colon \w{\theta}_j \not= 0 \\
    \infty & \colon \w{\theta}_j = 0,
    \end{cases}
\end{equation}
for every $\width\in \N$, $\theta \in \R^{\fd_\width }$ let $\M^{\theta}_{v} \subseteq \N$, $v \in \{0,1\}$, satisfy $\M^{\theta}_{0}=\{k \in \{1,2,\ldots,\width\} \colon 0\leq\q{\theta}_k\leq \nicefrac12\}$ and $\M^{\theta}_{1}=\{k \in \{1,2,\ldots,\width\} \colon \nicefrac12\leq\q{\theta}_k\leq 1\}$,
for every $\width\in \N$, $\theta \in \R^{\fd_\width }$, $v\in\{0,1\}$ let
$\m^{\theta}_{v,1},\m^{\theta}_{v,2}\in \R$ satisfy
\begin{equation} \label{H:setting:eq:defm}
    (\m^{\theta}_{v,1},\m^{\theta}_{v,2}) = \begin{cases}
   \big(\min(\M^{\theta}_{v}),\max(\M^{\theta}_{v})\big) & \colon \M^{\theta}_{v} \neq \emptyset \\
    (0,\width+1) & \colon \M^{\theta}_{v} = \emptyset,
    \end{cases} 
\end{equation}
let $\act_r \in C ( \R , \R )$, $r \in \N \cup \cu{ \infty } $,
satisfy for all $x \in \R$ that $( \bigcup_{r \in \N} \cu{ \act_r } ) \subseteq C^1( \R , \R)$, $\act_\infty ( x ) = \max \cu{ x , 0 }$,
 $\sup_{r \in \N} \sup_{y \in [- \abs{x}, \abs{x} ] }  \abs{ ( \act_r)'(y)}< \infty$, and
\begin{equation} \label{H:theo:intro:general:eq1}
    \limsup\nolimits_{r \to \infty}  \rbr*{ \abs { \act_r ( x ) - \act _\infty ( x ) } + \abs { (\act_r)' ( x ) - \indicator{(0, \infty)} ( x ) } } = 0,
\end{equation}
for every $\width\in \N$, $r \in \N \cup \{\infty\}$, $\theta \in \R^{\fd_\width}$
let $\realization{\width,\theta}_r \colon \R \to \R$
satisfy for all $x \in \R$ that
\begin{equation} \label{H:setting:snn:width1:eq:realization}
    \realization{\width,\theta}_r (x) = \c{\theta} +  \smallsum_{i=1}^{\width} \v{\theta}_i \br[\big]{ \act_r ( \w{\theta}_i x + \b{\theta}_i )},
\end{equation}
for every $\width\in \N$, $r \in \N \cup \{\infty\}$ let $\cL^\width_r \colon \R^{\fd_\width  } \to \R$  satisfy for all  $\theta \in \R^{\fd_\width}$ that
 \begin{equation} \label{H:theo:intro:general:L}
    \cL_r^\width ( \theta ) =  \int_0^1 \rbr[\big]{ \realization{\width,\theta}_r (x) -\mathbbm{1}_{(\nicefrac12,\infty)}(x) } ^2 \, \d x,
\end{equation}
for every $\width\in \N$, $\theta \in \R^{\fd_\width }$, $i\in \{1,2, \ldots, \width \}$ let $I^{\theta}_i \subseteq \R$  satisfy  $I^{\theta}_i = \cu{ x \in [0 , 1 ] \colon \w{\theta}_i x + \b{\theta}_i  > 0 }$,
 for every $\width\in \N$,  $\theta \in \R^{\fd_\width }$ with $\m^{\theta}_{0,1}\neq0$ let $\alpha^{\theta} \in \R$ satisfy for all $x\in \big[0,\q{\theta}_{\m^{\theta}_{0,1}}\big]$ that $\realization{\width,\theta}_\infty(x)=\alpha^{\theta} x + \realization{\width,\theta}_\infty(0)$,
 for every $\width\in \N$,  $\theta \in \R^{\fd_\width }$ with $\m^{\theta}_{1,2}\neq0$ let $\beta^{\theta} \in \R$ satisfy for all $x\in \big[\q{\theta}_{\m^{\theta}_{1,2}},1\big]$ that $\realization{\width,\theta}_\infty(x)=\beta^{\theta} (x-1) + \realization{\width,\theta}_\infty(1)$,
and for every $\width\in \N$ let $\cG^\width = ( \cG_1^\width , \ldots, \cG_{\fd_\width}^\width ) \colon \R^{\fd_\width} \to \R^{\fd_\width}$ satisfy for all
$\theta \in \cu{ \vartheta \in \R^{\fd_\width} \colon ( ( \nabla \cL_r^\width ) ( \vartheta ) ) _{r \in \N} \text{ is convergent} }$
that $\cG^\width ( \theta ) = \lim_{r \to \infty} (\nabla \cL_r^\width) ( \theta )$. 
\end{setting}
\subsection{Estimates of integrals}\label{Generalresults1}
\begin{prop} \label{lemma:line}
Let $\alpha, \beta \in \R$. Then
\begin{equation} \label{lemma:min:line}
    \int_0^1 (\alpha x+ \beta -\mathbbm{1}_{(\nicefrac12,\infty)}(x))^2 \, \d x
    \geq \int_0^1 \left(\frac{3x}2 -\frac14 -\mathbbm{1}_{(\nicefrac12,\infty)}(x)\right)^2 \, \d x=\frac1{16}.
\end{equation}
\end{prop}
\begin{cproof}{lemma:line}
Throughout this proof let $\Phi \colon \R^2 \to \R$ satisfy for all $\bfa, \bfb \in \R$ that
\begin{equation}\label{lemma:line:phi}
    \Phi(\bfa,\bfb)=\int_0^1 (\bfa x+ \bfb -\mathbbm{1}_{(\nicefrac12,\infty)}(x))^2 \, \d x.
\end{equation}
\Nobs that \cref{lemma:line:phi} ensures that for all $\bfa, \bfb \in \R$ it holds that
\begin{equation}\label{lemma:min:line:equation}
\begin{split}
      \Phi(\bfa,\bfb) &=\int_0^1 (\bfa x+\bfb -\mathbbm{1}_{(\nicefrac12,\infty)}(x))^2 \, \d x \\
      &=\int_0^1 [(\bfa x+\bfb)^2-2(\bfa x+\bfb)\mathbbm{1}_{(\nicefrac12,\infty)}(x) +\mathbbm{1}_{(\nicefrac12,\infty)}(x)] \, \d x \\
      &=\int_0^1 (\bfa^2 x^2+2\bfa \bfb x +\bfb^2) \, \d x -2 \int_{\frac12}^1 (\bfa x+\bfb) \, \d x +\frac12
      \\
      &= \frac{\bfa^2}{3}+\bfa\bfb+\bfb^2-\frac{3 \bfa}{4}-\bfb+\frac12.
\end{split}
\end{equation}
\Hence for all $\bfa,\bfb \in \R$ that
\begin{equation}
    (\nabla \Phi)(\bfa,\bfb)=\left(\frac{2\bfa}{3}+\bfb-\frac34,\bfa+2\bfb-1\right).
\end{equation}
This implies that 
\begin{equation}
\begin{split}
    &\{(\bfa,\bfb)\in \R^2 \colon (\nabla \Phi)(\bfa,\bfb)=0\} \\
    &= \{(\bfa,\bfb)\in \R^2 \colon [(\bfa=1-2\bfb) \wedge (8\bfa+12 \bfb-9=0)]\}\\
    &= \{(\bfa,\bfb)\in \R^2 \colon [(\bfa=1-2\bfb) \wedge (8(1-2\bfb)+12 \bfb-9=0)]\}
    \\
    &= \{(\bfa,\bfb)\in \R^2 \colon [(\bfa=1-2\bfb) \wedge (\bfb=-\nicefrac14)]\}
    = \{(\nicefrac32,-\nicefrac14)\}.
\end{split}
\end{equation}
Combining this with the fact that
\begin{equation}
   \Phi(\nicefrac32,-\nicefrac14)=\frac34-\frac38+\frac1{16}-\frac98+\frac14+\frac12=\frac1{16}
\end{equation}
and the fact that
\begin{equation}
   \lim\nolimits_{\|(\bfa,\bfb)\|\to \infty}\Phi(\bfa,\bfb)=\infty
\end{equation}
establishes that for all $\bfa,\bfb \in \R$ it holds that
\begin{equation}
    \Phi(\bfa,\bfb) \geq \Phi(\nicefrac32,-\nicefrac14)=\int_0^1 \left(\frac{3x}2 -\frac14 -\mathbbm{1}_{(\nicefrac12,\infty)}(x)\right)^2 \, \d x=\frac1{16}.
\end{equation}
\end{cproof}
\begin{prop} \label{lemma:half}
Let $\alpha, \beta  \in \R$, $\scra\in [0,\nicefrac12]$ satisfy $\alpha \scra + \beta =0$. Then
\begin{equation} \label{lemma:half:line}
    \int_\scra^{\frac34} (\alpha x+\beta -\mathbbm{1}_{(\nicefrac12,\infty)}(x))^2 \, \d x \geq\int_{\frac38}^{\frac34} \left( \frac{32x}9 -\frac43 -\mathbbm{1}_{(\nicefrac12,\infty)}(x)\right)^2 \, \d x = \frac1{36}.
\end{equation}
\end{prop}
\begin{cproof}{lemma:half}
Throughout this proof let $\Phi \colon \R^2\times[0, \nicefrac12] \to \R$ satisfy for all $\bfa, \bfb \in \R$, $\fa \in [0, \nicefrac12]$ that
\begin{equation}\label{lemma:half:phi}
    \Phi(\bfa,\bfb,\fa)=\int_\fa^{\frac34} (\bfa x+ \bfb -\mathbbm{1}_{(\nicefrac12,\infty)}(x))^2
\end{equation}
and let $L \colon \R^3\times[0, \nicefrac12]$ satisfy for all $\bfa, \bfb, \lambda \in \R$, $\fa \in [0, \nicefrac12]$ that 
\begin{equation}
    L(\bfa, \bfb, \lambda,\fa)=\Phi(\bfa,\bfb,\fa)-\lambda(\bfa \fa + \bfb).
\end{equation}
\Nobs that \cref{lemma:half:phi} ensures that
\begin{equation}\label{lemma:half:equation}
\begin{split}
      \Phi(\bfa,\bfb,\fa) &=\int_\fa^{\frac34} (\bfa x+\bfb -\mathbbm{1}_{(\nicefrac12,\infty)}(x))^2 \, \d x\\ &=\int_\fa^{\frac34} [(\bfa x+\bfb)^2-2(\bfa x+\bfb)\mathbbm{1}_{(\nicefrac12,\infty)}(x) +\mathbbm{1}_{(\nicefrac12,\infty)}(x)] \, \d x\\ &=\int_\fa^{\frac34} (\bfa^2 x^2+2\bfa \bfb x +\bfb^2) \, \d x -2 \int_{\frac12}^{\frac34} (\bfa x+\bfb) \, \d x +\frac14
      \\
      &=\frac{\bfa^2}{3}\Big(\frac{27}{64}-\fa^3\Big)+\bfa \bfb \Big(\frac9{16}-\fa^2\Big)+\bfb^2\Big(\frac34-\fa\Big)-\frac{5\bfa}{16}-\frac\bfb2+\frac14.
\end{split}
\end{equation}
This implies for all $\bfa,\bfb, \lambda \in \R$, $\fa \in [0, \nicefrac12]$ that
\begin{equation}\label{lemma:half:grad}
\begin{split}
    \frac{\partial}{\partial \bfa} L(\bfa, \bfb, \lambda,\fa)&=\frac{\partial}{\partial \bfa} \Phi(\bfa,\bfb,\fa)- \lambda\fa=
    \frac{2\bfa}{3}\Big(\frac{27}{64}-\fa^3\Big)+ \bfb \Big(\frac9{16}-\fa^2\Big)-\frac{5}{16}- \lambda\fa,\\
    \frac{\partial}{\partial \bfb} L(\bfa, \bfb, \lambda,\fa)&=\frac{\partial}{\partial \bfb} \Phi(\bfa,\bfb,\fa)- \lambda=
    \bfa \Big(\frac9{16}-\fa^2\Big)+2\bfb\Big(\frac34-\fa\Big)-\frac12- \lambda,\\
    \frac{\partial}{\partial \lambda} L(\bfa, \bfb, \lambda, \fa)&=-\bfa\fa-\bfb, \qand \\
    \frac{\partial}{\partial \fa} L(\bfa, \bfb, \lambda, \fa)&=\frac{\partial}{\partial \fa} \Phi(\bfa,\bfb,\fa)-\lambda \bfa=-\bfa^2 \fa^2-2\bfa \bfb \fa -\bfb^2-\lambda \bfa.
\end{split}
\end{equation}
\Hence that
\begin{equation}\label{lemma:half:zero}
\begin{split}
    \begin{cases}
     \frac{\partial}{\partial \bfa} L(\bfa, \bfb, \lambda,\fa)=0 \\
     \frac{\partial}{\partial \bfb} L(\bfa, \bfb, \lambda,\fa)=0 \\
     \frac{\partial}{\partial \lambda} L(\bfa, \bfb, \lambda,\fa)=0 \\
     \frac{\partial}{\partial \fa} L(\bfa, \bfb, \lambda,\fa)=0
     \end{cases}
     &=
     \begin{cases}
     \frac{2\bfa}{3}\left(\frac{27}{64}-\fa^3\right)+ \bfb \left(\frac9{16}-\fa^2\right)-\frac{5}{16}- \lambda\fa=0 \\
      \bfa \left(\frac9{16}-\fa^2\right)+2\bfb\left(\frac34-\fa\right)-\frac12- \lambda=0 \\
      -\bfa\fa-\bfb=0 \\
      -\bfa^2 \fa^2-2\bfa \bfb \fa -\bfb^2-\lambda \bfa=0
     \end{cases}\\
     &  =
     \begin{cases}
     \frac{2\bfa}{3}\left(\frac{27}{64}-\fa^3\right)-\bfa\fa \left(\frac9{16}-\fa^2\right)-\frac{5}{16}- \lambda\fa=0 \\
      \bfa \left(\frac9{16}-\fa^2\right)-2\bfa\fa\left(\frac34-\fa\right)-\frac12- \lambda=0 \\
      \bfb=-\bfa\fa \\
      \lambda \bfa=0
     \end{cases}
    \\
     &  =\begin{cases}
     \frac{9\bfa}{32}+\frac{\bfa \fa^3}{3}-\frac{9\bfa\fa}{16} -\frac{5}{16}- \lambda\fa=0 \\
      \frac{9\bfa}{16}+\bfa \fa^2-\frac{3\bfa\fa}{2}-\frac12- \lambda=0 \\
      \bfb=-\bfa\fa \\
      \lambda \bfa=0
     \end{cases}\\
     & =\begin{cases}
     \bfa =\left(\frac{5}{16}+ \lambda\fa\right)\left(\frac{9}{32}-\frac{9\fa}{16}+\frac{ \fa^3}{3}\right)^{-1}  \\
     \left(\frac{5}{16}+ \lambda\fa\right)\left(\frac{9}{32}-\frac{9\fa}{16}+\frac{ \fa^3}{3}\right)^{-1}
      \left(\frac{9}{16}+ \fa^2-\frac{3\fa}{2}\right)-\frac12 = \lambda \\
      \bfb=-\bfa\fa \\
      \lambda \bfa=0.
     \end{cases}
\end{split}
\end{equation}
This implies in the case $\nabla=L(\bfa, \bfb, \lambda,\fa)=\bfa=0$ that it holds that 
$\bfa=\bfb=0$, $\lambda=-\nicefrac12$, and $\fa=\nicefrac58$ which is not in the domain of $L$.
In the following we distinguish between the case $\lambda=0$, the case $\fa=0$, and the case $\fa=\nicefrac12$. We first show \cref{lemma:half:line} in the case
\begin{equation}\label{lemma:half:case2}
    \lambda=0.
\end{equation} 
\Nobs that \cref{lemma:half:zero} and \cref{lemma:half:case2} ensure that 
\begin{equation}
\begin{split}
    0&=\left(\frac{5}{16}\right)\left(\frac{9}{32}-\frac{9\fa}{16}+\frac{ \fa^3}{3}\right)^{-1}
      \left(\frac{9}{16}+ \fa^2-\frac{3\fa}{2}\right)-\frac12
      =
      5\left(\frac{9}{16}+ \fa^2-\frac{3\fa}{2}\right)
      \\ &\quad-8\left(\frac{9}{32}-\frac{9\fa}{16}+\frac{ \fa^3}{3}\right)=\frac9{16}+5\fa^2-3\fa-\frac{ 8\fa^3}{3}=-\frac1{48}(3-4\fa)^2(8\fa-3).
\end{split}
\end{equation}
\Hence that $\fa=\nicefrac38$. Combining this with 
\cref{lemma:half:zero} shows that 
\begin{equation}
    \nabla L(\nicefrac{32}9,-\nicefrac43,0,\nicefrac38)=0.
\end{equation}
This,  Lagrange multiplier theorem, and the fact that for all $\fa \in (0,\nicefrac12)$ it holds that \begin{equation}
    \lim\nolimits_{\|(\bfa,\bfb)\|\to \infty} \Phi(\bfa,\bfb,\fa)=\infty
\end{equation}
imply for all $\bfa,\bfb \in \R$, $\fa \in (0,\nicefrac12)$ with $\bfa \fa + \bfb =0$ that 
\begin{equation}\label{lemma:half:case}
    \Phi(\bfa,\bfb,\fa)\geq \Phi(\nicefrac{32}9,-\nicefrac43,\nicefrac38)=\frac1{36}.
\end{equation}
This establishes \cref{lemma:half:line} in the case $\lambda=0$.
In the next step we prove \cref{lemma:half:line} in the case
\begin{equation}\label{lemma:half:case3}
    \fa=0.
\end{equation}
\Nobs that  \cref{lemma:half:equation} and \cref{lemma:half:case3} assure for all $\bfa,\bfb \in \R$ with $\bfa \fa + \bfb =0$ that \begin{equation} \label{lemma:half:case0}
    \Phi(\bfa,\bfb,0)=\Phi(\bfa,0,0)=\frac{9\bfa^2}{64}-\frac{5\bfa}{16}+\frac14\geq \frac{9}{64}\left(\frac{10}9\right)^2-\frac{5}{16}\left(\frac{10}9\right)+\frac14= \frac{11}{144}.
\end{equation}
This establishes \cref{lemma:half:line} in the case $\fa=0.$
Finally we demonstrate \cref{lemma:half:line} in the case
\begin{equation}\label{lemma:half:case4}
    \fa=\frac12.
\end{equation}
\Nobs that  \cref{lemma:half:equation} and \cref{lemma:half:case4} assure for all $\bfa,\bfb \in \R$ with $\bfa \fa + \bfb =0$ that \begin{equation}
\begin{split}
     \Phi(\bfa,\bfb,\nicefrac12)&=\Phi(\bfa,-\nicefrac\bfa2,\nicefrac12) =
   \frac{19\bfa^2}{192}-  \frac{5\bfa^2}{32}+\frac{\bfa^2}{16}-\frac{5\bfa}{16}+\frac{\bfa}4+\frac14 \\
   &=\frac{\bfa^2}{192}-\frac{\bfa}{16}+\frac14\geq\frac{6^2}{192}-\frac{6}{16}+\frac14= \frac1{16}.
\end{split}
\end{equation}
This establishes \cref{lemma:half:line} in the case $\fa=\nicefrac12.$
\end{cproof}
\begin{prop} \label{lemma:halfdx}
Let $\alpha, \beta \in \R$, $\scrb\in [\nicefrac12,1]$ satisfy $\alpha \scrb + \beta =1$. Then
\begin{equation} \label{lemma:halfdx:line}
    \int_{\frac14}^\scrb (\alpha x+\beta -\mathbbm{1}_{(\nicefrac12,\infty)}(x))^2 \, \d x\geq\int_{\frac14}^{\frac58} \left( \frac{32x}9 -\frac{11}
    9 -\mathbbm{1}_{(\nicefrac12,\infty)}(x)\right)^2 \, \d x = \frac1{36}.
\end{equation}
\end{prop}
\begin{cproof}{lemma:halfdx}
Throughout this proof let $\Phi \colon \R^2\times[ \nicefrac12,1] \to \R$ satisfy for all $\bfa, \bfb \in \R$, $\fb \in [\nicefrac12,1]$ that
\begin{equation}\label{lemma:halfdx:phi}
    \Phi(\bfa,\bfb,\fb)=\int_{\frac14}^\fb (\bfa x+ \bfb -\mathbbm{1}_{(\nicefrac12,\infty)}(x))^2
\end{equation}
and let $L \colon \R^3\times[ \nicefrac12,1]$ satisfy for all $\bfa, \bfb, \lambda \in \R$, $\fb \in [\nicefrac12,1]$ that 
\begin{equation}
    L(\bfa, \bfb, \lambda,\fb)=\Phi(\bfa,\bfb,\fb)-\lambda(\bfa \fb + \bfb-1).
\end{equation}
\Nobs that \cref{lemma:halfdx:phi} ensures for all $\bfa, \bfb, \lambda \in \R$, $\fb \in [\nicefrac12,1]$ that 
\begin{equation}\label{lemma:halfdx:equation}
\begin{split}
      \Phi(\bfa,\bfb,\fb) &=\int_{\frac14}^\fb (\bfa x+\bfb -\mathbbm{1}_{(\nicefrac12,\infty)}(x))^2 \, \d x\\
      & =\int_{\frac14}^\fb [(\bfa x+\bfb)^2-2(\bfa x+\bfb)\mathbbm{1}_{(\nicefrac12,\infty)}(x) +\mathbbm{1}_{(\nicefrac12,\infty)}(x)] \, \d x \\ &=\int_{\frac14}^\fb (\bfa^2 x^2+2\bfa \bfb x +\bfb^2) \, \d x -2 \int_{\frac12}^{\fb} (\bfa x+\bfb) \, \d x +\fb-\frac12
      \\
      &=\frac{\bfa^2}{3}\Big(\fb^3-\frac{1}{64}\Big)+\bfa \bfb \Big(\fb^2-\frac1{16}\Big)+\bfb^2\Big(\fb-\frac14\Big)-\bfa\left(\fb^2-\frac14\right)\\
      &\quad -2\bfb\left(\fb-\frac12\right)+\fb-\frac12.
\end{split}
\end{equation}
This implies for all $\bfa,\bfb, \lambda \in \R$, $\fb \in [ \nicefrac12,1]$ that
\begin{equation}\label{lemma:halfdx:grad}
\begin{split}
    \frac{\partial}{\partial \bfa} L(\bfa, \bfb, \lambda,\fb)&=\frac{\partial}{\partial \bfa} \Phi(\bfa,\bfb,\fb)- \lambda\fb=
    \frac{2\bfa}{3}\Big(\fb^3-\frac{1}{64}\Big)+ \bfb \Big(\fb^2-\frac1{16}\Big)-\fb^2+\frac14- \lambda\fb,\\
    \frac{\partial}{\partial \bfb} L(\bfa, \bfb, \lambda,\fb)&=\frac{\partial}{\partial \bfb} \Phi(\bfa,\bfb,\fb)- \lambda=
    \bfa \Big(\fb^2-\frac1{16}\Big)+2\bfb\Big(\fb-\frac14\Big)-2\fb+1- \lambda,\\
    \frac{\partial}{\partial \lambda} L(\bfa, \bfb, \lambda, \fb)&=-\bfa\fb-\bfb+1, \qand \\
    \frac{\partial}{\partial \fb} L(\bfa, \bfb, \lambda, \fb)&=\frac{\partial}{\partial \fb} \Phi(\bfa,\bfb,\fb)-\lambda \bfa=\bfa^2 \fb^2+2\bfa \bfb \fb + \bfb^2-2\bfa\fb-2\bfb+1-\lambda \bfa.
\end{split}
\end{equation}
\Hence that
\begin{equation}\label{lemma:halfdx:zero}
\begin{split}
    \begin{cases}
     \frac{\partial}{\partial \bfa} L(\bfa, \bfb, \lambda,\fb)=0 \\
     \frac{\partial}{\partial \bfb} L(\bfa, \bfb, \lambda,\fb)=0 \\
     \frac{\partial}{\partial \lambda} L(\bfa, \bfb, \lambda,\fb)=0 \\
     \frac{\partial}{\partial \fb} L(\bfa, \bfb, \lambda,\fb)=0
     \end{cases}
     &=
     \begin{cases}
     \frac{2\bfa}{3}\left(\fb^3-\frac{1}{64}\right)+ \bfb \left(\fb^2-\frac1{16}\right)-\fb^2+\frac14- \lambda\fb=0 \\
      \bfa \left(\fb^2-\frac1{16}\right)+2\bfb\left(\fb-\frac14\right)-2\fb+1- \lambda=0 \\
     -\bfa\fb-\bfb+1=0 \\
      \bfa^2 \fb^2+2\bfa \bfb \fb + \bfb^2-2\bfa\fb-2\bfb+1-\lambda \bfa=0
     \end{cases}\\
     &  =
     \begin{cases}
     \frac{2\bfa}{3}\left(\fb^3-\frac{1}{64}\right)+ \bfb \left(\fb^2-\frac1{16}\right)-\fb^2+\frac14- \lambda\fb=0 \\
      \bfa \left(\fb^2-\frac1{16}\right)+2\bfb\left(\fb-\frac14\right)-2\fb+1- \lambda=0 \\
      \bfb=-\bfa\fb+1 \\
      \lambda \bfa=0
     \end{cases}
    \\
     &  =\begin{cases}
      \frac{2\bfa\fb^3}{3} -\frac{\bfa}{96}-\bfa \fb^3+\fb^2+\frac{\bfa\fb}{16}-\frac1{16}
      -\fb^2+\frac14- \lambda\fb=0 \\
     \bfa \fb^2 - \frac{\bfa}{16} - 2\bfa\fb^2+2\fb + \frac{\bfa \fb}{2} -\frac12
     -2\fb+1- \lambda=0 \\
      \bfb=-\bfa\fb+1 \\
      \lambda \bfa=0
     \end{cases}\\
      &  =\begin{cases}
      -\frac{\bfa\fb^3}{3} -\frac{\bfa}{96}+\frac{\bfa\fb}{16}+\frac3{16}
      - \lambda\fb=0 \\
     -\bfa \fb^2 - \frac{\bfa}{16}  + \frac{\bfa \fb}{2} +\frac12- \lambda=0 \\
      \bfb=-\bfa\fb+1 \\
      \lambda \bfa=0
     \end{cases}\\
     & =\begin{cases}
     \bfa =\left(-\frac{3}{16}+ \lambda\fb\right)\big(-\frac{1}{96}+\frac{\fb}{16}-\frac{ \fb^3}{3}\big)^{-1}  \\
     \left(\frac\fb2-\fb^2 - \frac1{16} \right)\left(-\frac{3}{16}+ \lambda\fb\right)\big(-\frac{1}{96}+\frac{\fb}{16}-\frac{ \fb^3}{3}\big)^{-1}+\frac12 = \lambda \\
      \bfb=-\bfa\fb+1 \\
      \lambda \bfa=0.
     \end{cases}
\end{split}
\end{equation}
This implies that in the case $\nabla L(\bfa, \bfb, \lambda,\fb)=\lambda-\nicefrac12=0$ it holds that $\bfa=0$, $\bfb=1$, and $\fb=\nicefrac38$ which is not in the domain of $L$.
In the following we distinguish between the case $\lambda=0$,  the case $\fb=\nicefrac12$, and the case $\fb=1$. We first demonstrate \cref{lemma:halfdx:line} in the case
\begin{equation}\label{lemma:halfdx:case2}
    \lambda=0.
\end{equation} 
\Nobs that \cref{lemma:halfdx:zero} and \cref{lemma:halfdx:case2} ensure that for all $\bfa, \bfb \in \R$, $\fb \in [\nicefrac12,1]$ such that $\nabla L(\bfa, \bfb, \lambda,\fb)=0$ it holds that
\begin{equation}
    0= \left(\frac\fb2-\fb^2 - \frac1{16} \right)\left(-\frac{3}{16}\right)\left(-\frac{1}{96}+\frac{\fb}{16}-\frac{ \fb^3}{3}\right)^{-1}+\frac12.
\end{equation}
\Hence that for all $\bfa, \bfb \in \R$, $\fb \in [\nicefrac12,1]$ such that $\nabla L(\bfa, \bfb, \lambda,\fb)=0$ it holds that
\begin{equation}
     0= \left(48\fb-96\fb^2 -6 \right)-\frac83\left(-1+6\fb-32 \fb^3\right)=\frac23(4\fb-1)^2(8\fb-5).
\end{equation}
Combining this with 
\cref{lemma:halfdx:zero} shows that for all $\bfa, \bfb \in \R$, $\fb \in [\nicefrac12,1]$ such that $\nabla L(\bfa, \bfb, \lambda,\fb)=0$ it holds that $\fb=\nicefrac58$, $\bfa=\nicefrac{32}9$, and $\bfb=-\nicefrac{11}9$. 
This,  Lagrange multiplier theorem, and the fact that for all $\fb \in (\nicefrac12,1)$ it holds that \begin{equation}
    \lim\nolimits_{\|(\bfa,\bfb)\|\to \infty} \Phi(\bfa,\bfb,\fb)=\infty
\end{equation}
imply for all $\bfa,\bfb \in \R$, $\fb \in (\nicefrac12,1)$ with $\bfa \fb + \bfb =1$ that 
\begin{equation}\label{lemma:halfdx:case}
    \Phi(\bfa,\bfb,\fb)\geq \Phi(\nicefrac{32}9,-\nicefrac{11}9,\nicefrac58)=\frac1{36}.
\end{equation}
This establishes \cref{lemma:halfdx:line} in the case $\lambda=0.$
In the next step we prove \cref{lemma:halfdx:line} in the case
\begin{equation}\label{lemma:halfdx:case3}
    \fb=1.
\end{equation}
\Nobs that \cref{lemma:halfdx:equation} and \cref{lemma:halfdx:case3} assure for all $\bfa,\bfb \in \R$ with $\bfa \fb + \bfb =1$ that \begin{equation} \label{lemma:halfdx:case0}
\begin{split}
     \Phi(\bfa,\bfb,1)&=\Phi(\bfa,1-\bfa,1)=\frac{21\bfa^2}{64}+\frac{15\bfa(1-\bfa)}{16}+\frac{3(1-\bfa)^2}4
    -\frac{3\bfa}{4}-(1-\bfa)+1-\frac12 \\
    &=\frac1{64}(9 \bfa^2-20\bfa+16)\geq\frac1{64}\left(9 \bigg(\frac{10}9\right)^2-20\left(\frac{10}9\right)+16\bigg)=\frac{11}{144}.
\end{split}
\end{equation}
This establishes \cref{lemma:halfdx:line} in the case $\fb=1$.
Finally we show \cref{lemma:halfdx:line} in the case \begin{equation}\label{lemma:halfdx:case4}
    \fb=\frac12.
\end{equation}
\Nobs that \cref{lemma:halfdx:equation} and \cref{lemma:halfdx:case4} assure for all $\bfa,\bfb \in \R$ with $\bfa \fb + \bfb =1$ that 
\begin{equation}\label{lemma:halfdx:fin}
\begin{split}
     \Phi(\bfa,\bfb,\nicefrac12)&=\Phi(\bfa,1-\nicefrac\bfa2,\nicefrac12) =
   \frac{\bfa^2}{64}+\frac{3\bfa(1-\frac\bfa2)}{16}+\frac{(1-\frac\bfa2)^2}4 \\
   &=\frac{\bfa^2}{192}-\frac{\bfa}{16}+\frac14\geq\frac{6^2}{192}-\frac{6}{16}+\frac14= \frac1{16}.
\end{split}
\end{equation}
This establishes \cref{lemma:halfdx:line} in the case $\fb=\nicefrac12$.
\end{cproof}
\begin{prop} \label{lemma:halfc}
Let $\alpha, \beta \in \R$. Then
\begin{equation} \label{lemma:halfc:line}
    \int_{\frac14}^{\frac34} (\alpha x+\beta -\mathbbm{1}_{(\nicefrac12,\infty)}(x))^2 \, \d x \geq \int_{\frac14}^{\frac34} \left(3 x-1 -\mathbbm{1}_{(\nicefrac12,\infty)}(x)\right)^2 \, \d x=\frac1{32}.
\end{equation}
\end{prop}
\begin{cproof}{lemma:halfc}
Throughout this proof let $\Phi \colon \R^2 \to \R$ satisfy for all $\bfa, \bfb \in \R$ that
\begin{equation}\label{lemma:halfc:phi}
    \Phi(\bfa,\bfb)=\int_{\frac14}^{\frac34} (\bfa x+ \bfb -\mathbbm{1}_{(\nicefrac12,\infty)}(x))^2 \, \d x.
\end{equation}
\Nobs that \cref{lemma:halfc:phi} ensures that
\begin{equation}\label{lemma:min:halfc:equation}
\begin{split}
      \Phi(\bfa,\bfb) &=\int_{\frac14}^{\frac34} (\bfa x+\bfb -\mathbbm{1}_{(\nicefrac12,\infty)}(x))^2 \, \d x\\ &=\int_{\frac14}^{\frac34} [(\bfa x+\bfb)^2-2(\bfa x+\bfb)\mathbbm{1}_{(\nicefrac12,\infty)}(x) +\mathbbm{1}_{(\nicefrac12,\infty)}(x)] \, \d x \\ &=\int_{\frac14}^{\frac34} (\bfa^2 x^2+2\bfa \bfb x +\bfb^2) \, \d x -2 \int_{\frac12}^{\frac34} (\bfa x+\bfb) \, \d x +\frac14
      \\
      &=\frac{26\bfa^2}{192}+\frac{\bfa\bfb}2 +\frac{\bfb^2}2 -\frac{5\bfa}{16}-\frac{\bfb}2+\frac14.
\end{split}
\end{equation}
\Hence for all $\bfa,\bfb \in \R$ that
\begin{equation}
    (\nabla \Phi)(\bfa,\bfb)=\left(\frac{13\bfa}{48}+\frac\bfb2-\frac5{16},\frac\bfa2 +\bfb-\frac12\right).
\end{equation}
This implies that 
\begin{equation}
\begin{split}
    &\{(\bfa,\bfb)\in \R^2 \colon (\nabla \Phi)(\bfa,\bfb)=0\} \\
    &= \{(\bfa,\bfb)\in \R^2 \colon [(\bfa=1-2\bfb) \wedge (13\bfa+24 \bfb-15=0)]\}\\
    &= \{(\bfa,\bfb)\in \R^2 \colon [(\bfa=1-2\bfb) \wedge (13(1-2\bfb)+24 \bfb-15=0)]\}
    \\&= \{(\bfa,\bfb)\in \R^2 \colon [(\bfa=1-2\bfb) \wedge (\bfb=-1)]\}
    = \{(3,-1)\}.
\end{split}
\end{equation}
Combining this with the fact that
\begin{equation}
   \Phi(3,-1)=\frac{39}{32}-\frac32+\frac1{2}-\frac{15}{16}+\frac12+\frac14=\frac{39}{32}-\frac{30}{32}-\frac{1}{4}=\frac1{32}
\end{equation}
and the fact that
\begin{equation}
   \lim\nolimits_{\|(\bfa,\bfb)\|\to \infty}\Phi(\bfa,\bfb)=\infty
\end{equation}
establishes that for all $\bfa,\bfb \in \R$ it holds that
\begin{equation}
    \Phi(\bfa,\bfb) \geq \Phi(3,-1)=\int_{\frac14}^{\frac34} \left(3 x-1 -\mathbbm{1}_{(\nicefrac12,\infty)}(x)\right)^2 \, \d x=\frac1{32}.
\end{equation}
\end{cproof}
\begin{lemma} \label{lemma:intmin}
Let $\alpha, \beta, \scra\in \R$, $\scrb\in (\scra,\infty)$. Then
\begin{equation}\label{lemma:intmin:thesis}
    \int_\scra^{\scrb} (\alpha x+\beta)^2 \, \d x \geq  \int_\scra^{\scrb} \left[\alpha x-\frac{\alpha(\scrb+\scra)}2 \right]^2 \, \d x=
    \frac{\alpha^2(\scrb-\scra)^3}{12}.
\end{equation}
\end{lemma}
\begin{cproof}{lemma:intmin}
\Nobs that, e.g., \cite[Lemma 5.1]{Jentzen} establishes \cref{lemma:intmin:thesis}.
\end{cproof}
\subsection{Properties of integrands}\label{Generalresults2}
\begin{lemma} \label{lem:affine:integral:zero}
Let $\scra \in \R$, $\scrb \in (\scra , \infty)$, $\alpha, \beta \in \R$ satisfy $\nicefrac12 \notin (\scra,\scrb)$ and
\begin{equation} \label{lem:affine:integral:zero:assumption}
    \int_\scra^\scrb x ( \alpha x + \beta - \mathbbm{1}_{(\nicefrac12,\infty)}(x)) \, \d x = \int_\scra^\scrb (\alpha x + \beta - \mathbbm{1}_{(\nicefrac12,\infty)}(x)) \, \d x = 0.
\end{equation}
Then $\alpha = 0$ and $\beta = \mathbbm{1}_{ (-\infty, \scra] }( \nicefrac{ 1 }{ 2 } )$.
\end{lemma}
 \begin{cproof}{lem:affine:integral:zero}
\Nobs that the assumption that $\nicefrac12 \notin (\scra,\scrb)$ and, e.g., \cite[Lemma 6.6]{Jentzen} demonstrate that $\alpha = 0$ and $\beta = \mathbbm{1}_{ (-\infty, \scra] }( \nicefrac{ 1 }{ 2 } )$.
 \end{cproof}
 \begin{lemma} \label{lem:affine1:integral:zero}
Let $\scra \in [0,\nicefrac12),$ $\scrb \in (\nicefrac12,1]$, $\alpha, \beta \in \R$ satisfy 
\begin{equation} \label{lem:affine1:integral:zero:assumption}
    \int_{\scra}^{\scrb} x ( \alpha x + \beta - \mathbbm{1}_{(\nicefrac12,\infty)}(x)) \, \d x = \int_{\scra}^{\scrb} (\alpha x + \beta - \mathbbm{1}_{(\nicefrac12,\infty)}(x)) \, \d x = 0.
\end{equation}
Then \begin{equation}
    \alpha=\frac{3 (2 \scra - 1) (2 \scrb - 1)}{2 (\scra - \scrb)^3} \qandq \beta= -\frac{(2 \scrb - 1) (8 \scra^2 + \scra (2 \scrb - 3) + \scrb (2 \scrb - 3))}{4 (\scra - \scrb)^3}.
\end{equation}
\end{lemma}
 \begin{cproof}{lem:affine1:integral:zero}
 \Nobs that \cref{lem:affine1:integral:zero:assumption} implies that
 \begin{equation}
    0= \int_{\scra}^{\scrb} (\alpha x + \beta - \mathbbm{1}_{(\nicefrac12,\infty)}(x)) \, \d x=
    \frac{\alpha}{2}(\scrb^2-\scra^2)+ \beta (\scrb -\scra)-(\scrb-\nicefrac12).
 \end{equation}
 \Hence that 
 \begin{equation}\label{lem:affine1:integral:zero:io}
     \beta=(\scrb-\nicefrac12) (\scrb -\scra)^{-1}-\frac{\alpha}{2}(\scrb+\scra).
 \end{equation}
 \Moreover \cref{lem:affine1:integral:zero:assumption} assures that
 \begin{equation}
     0=\int_{\scra}^{\scrb}  [\alpha x^2 + \beta x - x\mathbbm{1}_{(\nicefrac12,\infty)}(x) ] \, \d x= \frac{\alpha}{3}(\scrb^3-\scra^3)+ \frac{\beta}{2} (\scrb^2 -\scra^2)-\frac12(\scrb^2-\nicefrac14).
 \end{equation}
 This and \cref{lem:affine1:integral:zero:io} demonstrate that 
 \begin{equation}
 \begin{split}
     0&= \frac{\alpha}{3}(\scrb^3-\scra^3)+ \frac{\beta}{2} (\scrb^2 -\scra^2)-\frac12(\scrb^2-\nicefrac14) \\ &=\frac{\alpha}{3}(\scrb^3-\scra^3)+ \frac{1}{2}(\scrb-\nicefrac12) (\scrb +\scra)-\frac{\alpha}{4}(\scrb^2-\scra^2)(\scrb +\scra)-\frac12(\scrb^2-\nicefrac14)\\
     &=\frac{\alpha}{12}(\scrb-\scra)[4(\scrb^2+\scra\scrb+\scra^2)-3(\scrb+\scra)^2]+ \frac{1}{2}(\scrb-\nicefrac12)[ (\scrb+\scra)-(\scrb +\nicefrac12)]\\
     &=\frac{\alpha}{12}(\scrb-\scra)^3+ \frac{1}{2}(\scrb-\nicefrac12)(\scra-\nicefrac12).
 \end{split}
 \end{equation}
 \Hence that 
 \begin{equation}
      \alpha= \frac{3(2\scra-1)(2\scrb-1)}{2(\scra-\scrb)^3}.
 \end{equation}
 This and \cref{lem:affine1:integral:zero:io} establish that
 \begin{equation}
 \begin{split}
      \beta&=(\scrb-\frac12) (\scrb -\scra)^{-1}-\frac{\alpha}{2}(\scrb+\scra)
      =(\scrb-\frac12) (\scrb -\scra)^{-1}-\frac{3(2\scrb-1)(2\scra-1)(\scrb+\scra)}{4(\scra-\scrb)^3}\\
      &=-\frac{(2\scrb-1)[2(\scra-\scrb)^2+3(2\scra-1)(\scrb+\scra)]}{4(\scra-\scrb)^3}\\
      &= -\frac{(2 \scrb - 1) (8 \scra^2 + \scra (2 \scrb - 3) + \scrb (2 \scrb - 3))}{4 (\scra - \scrb)^3}.
 \end{split}
 \end{equation}
 \end{cproof}
\begin{cor} \label{cor:affine1:integral:zero}
Let $\scra \in [0,\nicefrac12)$, $\scrb \in (\nicefrac12,1]$, $\alpha, \beta \in \R$ satisfy $\alpha\scra+\beta=0$ and
\begin{equation} 
    \int_{\scra}^{\scrb} x ( (\alpha x + \beta) - \mathbbm{1}_{(\nicefrac12,\infty)}(x)) \, \d x = \int_{\scra}^{\scrb} ((\alpha x + \beta) - \mathbbm{1}_{(\nicefrac12,\infty)}(x)) \, \d x = 0.
\end{equation}
Then  \begin{equation}
    \alpha=\frac{16}{9(2\scrb-1)} \qandq \beta=\frac{4(2\scrb-3)}{9(2\scrb-1)}.
\end{equation}
\end{cor}
\begin{cproof}{cor:affine1:integral:zero}
\Nobs that \cref{lem:affine1:integral:zero} and the assumption that $\alpha\scra+\beta=0$ prove that
\begin{equation}
    0=\frac{3\scra(2\scra-1)(2\scrb-1)}{2(\scra-\scrb)^3}-\frac{(2 \scrb - 1) (8 \scra^2 + \scra (2 \scrb - 3) + \scrb (2 \scrb - 3))}{4 (\scra - \scrb)^3}.
\end{equation}
\Hence that
\begin{equation}
\begin{split}
0&=6\scra(2\scra-1)-(8 \scra^2 + \scra (2 \scrb - 3) + \scrb (2 \scrb - 3)) \\
    &= 4 \scra^2-\scra  (2 \scrb +3)- \scrb (2 \scrb - 3)=
    4 \scra^2-2\scra \scrb -3\scra - 2\scrb^2 + 3\scrb\\ &=
    4\scra(\scra-\scrb)+2\scra\scrb-2\scrb^2-3(\scra-\scrb)=(\scra-\scrb)(4\scra+2\scrb-3).
\end{split}
\end{equation}
This, the assumption that $\scra \in [0,\nicefrac12)$, and the assumption that $\scrb \in (\nicefrac12,1]$ imply that $\scra=-\nicefrac{\scrb}2+\nicefrac34$.
Combining this with \cref{lem:affine1:integral:zero} 
demonstrates that
\begin{equation}
\begin{split}
    \alpha &= \frac{3(-\scrb+\frac12)(2\scrb-1)}{2(-\frac{3\scrb}2+\frac34)^3}=\frac{-\frac32(2\scrb-1)^2}{\frac{27}{32}(-2\scrb+1)^3}=\frac{-\frac32}{-\frac{27}{32}(2\scrb-1)}
=\frac{16}{9(2\scrb-1)} \qand \\ 
\beta &= -\frac{(2 \scrb - 1) [8 (-\frac\scrb2+\frac34)^2 + (-\frac\scrb2+\frac34) (2 \scrb - 3) + \scrb (2 \scrb - 3)]}{4 (-\frac{3\scrb}2+\frac34)^3}\\
&=\frac{\frac12 (-2\scrb+3)^2 + \frac14(2\scrb+3) (2 \scrb - 3)}{\frac{27}{16} (-2\scrb+1)^2}=
\frac{4 (2\scrb-3) [2(2\scrb-3)+ (2\scrb+3)]}{27 (-2\scrb+1)^2}\\
&=\frac{4 (2\scrb-3) [4\scrb-6+ 2\scrb+3]}{27 (2\scrb-1)^2}
=\frac{4 (2\scrb-3) [6\scrb-3]}{27 (2\scrb-1)^2}
=\frac{4(2\scrb-3)}{9(2\scrb-1)}.
\end{split}
\end{equation}
\end{cproof}
\begin{cor} \label{cor:affine1:integral:zero2}
Let $\scra \in [0,\nicefrac12)$, $\scrb \in (\nicefrac12,1]$, $\alpha, \beta \in \R$ satisfy $\alpha\scrb+\beta=1$ and
\begin{equation} 
    \int_{\scra}^{\scrb} x ( (\alpha x + \beta) - \mathbbm{1}_{(\nicefrac12,\infty)}(x)) \, \d x = \int_{\scra}^{\scrb} ((\alpha x + \beta) - \mathbbm{1}_{(\nicefrac12,\infty)}(x)) \, \d x = 0.
\end{equation}
Then  \begin{equation}
    \alpha=-\frac{16}{9(2\scra-1)} \qandq \beta=\frac{10\scra+3}{9(2\scra-1)}.
\end{equation}
\end{cor}
\begin{cproof}{cor:affine1:integral:zero2}
\Nobs that \cref{lem:affine1:integral:zero} and the assumption that $\alpha\scrb+\beta=1$ prove that
\begin{equation}
    0=-1+\frac{3\scrb(2\scrb-1)(2\scra-1)}{2(\scra-\scrb)^3}-\frac{(2 \scrb - 1) (8 \scra^2 + \scra (2 \scrb - 3) + \scrb (2 \scrb - 3))}{4 (\scra - \scrb)^3}.
\end{equation}
\Hence that
\begin{equation}
\begin{split}    
    0&=-4 (\scra - \scrb)^3+6\scrb(2\scrb-1)(2\scra-1)-(2\scrb-1)(8 \scra^2 + \scra (2 \scrb - 3) + \scrb (2 \scrb - 3)) \\
    &=-4 (\scra - \scrb)^3 +(2\scrb-1)(12 \scra\scrb-6\scrb   -8 \scra^2 -2 \scra\scrb + 3\scra - 2\scrb^2 +3\scrb)\\
    &=-4 \scra^3+4\scrb^3 +12 \scra^2\scrb-12 \scra \scrb^2+(2\scrb-1)(10\scra\scrb-3\scrb-8\scra^2+3\scra-2\scrb^2)\\
    &=-4 \scra^3 -4 \scra^2\scrb+8 \scra \scrb^2
    - 6 \scrb^2 +6\scra \scrb 
    -10\scra\scrb+3\scrb+8\scra^2-3\scra+2\scrb^2\\
    &=-4 \scra^3 -4 \scra^2\scrb+8 \scra \scrb^2-4\scra \scrb + 8 \scra^2 -4\scrb^2-3\scra+3\scrb \\
    &= \scra(-4\scra^2 -8 \scra\scrb+8\scra+4\scrb-3)+\scrb(4\scra^2+8 \scra\scrb-8\scra-4\scrb+3)\\
    &=(\scra-\scrb)(4\scrb-8\scra\scrb+8\scra-4\scra^2-3)\\
    &=
    (\scra-\scrb)(4\scrb(1-2\scra)+2\scra(1-2\scra)+6\scra-3)=(\scra-\scrb)(1-2\scra)(2\scra+4\scrb-3).
\end{split}
\end{equation}
This, the assumption that $\scra \in [0,\nicefrac12)$, and the assumption that $\scrb \in (\nicefrac12,1]$ imply that $\scrb=-\nicefrac{\scra}2+\nicefrac34$.
Combining this with \cref{lem:affine1:integral:zero} 
demonstrates that
\begin{equation}
\begin{split}
    \alpha &=\frac{3 (2 \scra - 1) (-\scra+\frac12)}{2 (\frac{3\scra}2 - \frac34)^3}=
    -\frac{\frac32(2\scra-1)^2}{\frac{27}{32}(2\scra-1)^3}
=-\frac{16}{9(2\scra-1)} \qand \\ 
\beta &=  -\frac{(-\scra +\frac12) [8 \scra^2 + \scra (-\scra - \frac32) + (-\frac\scra2+\frac34) (-\scra-\frac32)]}{4 (\frac{3\scra}2 -\frac34)^3}\\
&=-\frac{\frac12(-2\scra +1) (\frac{15\scra^2}{2}  -\frac{3\scra}2  -\frac98)}{\frac{27}{16} (2\scra -1)^3}=
\frac{8(2\scra-1) (\frac{15\scra^2}{2}  -\frac{3\scra}2  -\frac98)}{27 (2\scra -1)^3}=
\frac{8 (\frac{5\scra^2}{2}  -\frac{\scra}2  -\frac38)}{9 (2\scra -1)^3}\\
&=\frac{20\scra^2-4\scra-3}{9(2\scra-1)^2}=\frac{(2\scra-1)(10\scra+3)}{9(2\scra-1)^2}
=\frac{10\scra+3}{9(2\scra-1)}.
\end{split}
\end{equation}
\end{cproof}
\begin{lemma}\label{corollary:integralzero}
Let $N \in \N $, $\fx_0, \fx_1, \ldots, \fx_N, \alpha_1, \alpha_2, \ldots, \alpha_N, \beta_1, \beta_2, \ldots, \beta_N, \scra, \scrb \in \R$,  $f \in  C( [\scra,\scrb],\allowbreak \R )$ 
satisfy for all $i \in \{1, 2, \ldots, N\}$, $x \in [\fx_{i-1}, \fx_i ]$ that $\scra = \fx_0 < \fx_1 < \cdots < \fx_N = \scrb$, $f(x) = \alpha_i x + \beta_i$, and 
\begin{equation} \label{cor:integralzero:function}
    \int_{\fx_{i-1}}^{ \fx_i}  f(x) \, dx=0.
\end{equation}
Then it holds for all $i \in \{1, 2, \ldots, N\}$ that $f(\fx_0)=(-1)^i f(\fx_i)$.
\end{lemma}
\begin{cproof}{corollary:integralzero}
\Nobs that \cref{cor:integralzero:function} implies for all $i \in \{1, 2, \ldots, N\}$ that
\begin{equation}
    \beta_i=-\frac{\alpha_i}{2}(\fx_i+\fx_{i-1}).
\end{equation}
This proves for all $i \in \{1, 2, \ldots, N\}$, $x \in [\fx_{i-1}, \fx_i ]$ that $f(x) = \alpha_i x -\nicefrac{\alpha_i}{2}(\fx_i+\fx_{i-1})$. 
\Hence for all $i \in \{1, 2, \ldots, N\}$ that $f(\fx_{i-1})=-f(\fx_{i})$.
This establishes for all $i \in \{1, 2, \ldots, N\}$ that $f(\fx_0)=(-1)^i f(\fx_i)$.
\end{cproof}
\begin{lemma} \label{lemma:doubleintegral:zero}
Let $\scra \in \R$, $\scrb \in (\scra , \infty)$, $\scrc \in (\scrb , \infty)$, $\alpha_1, \alpha_2 \in \R\backslash\{0\},$  $\beta_1, \beta_2 \in \R$ satisfy 
\begin{equation} \label{lemma:doubleintegral:zero:assumption}
    \int_\scra^\scrb (\alpha_1 x + \beta_1) \, \d x =
    \int_\scrb^\scrc (\alpha_2 x + \beta_2)  \, \d x =
     \int_\scra^\scrb x(\alpha_1 x + \beta_1) \, \d x +\int_\scrb^\scrc x(\alpha_2 x + \beta_2)  \, \d x 
    = 0
\end{equation}
and $\alpha_1 \scrb +\beta_1= \alpha_2 \scrb +\beta_2$.
Then $ \scrb-\scra= \scrc -\scrb $, $\beta_1 =- \nicefrac{\alpha_1(\scra+\scrb)}2  $, and $\beta_2=- \nicefrac{\alpha_2(\scrb+\scrc)}2  $.
\end{lemma}
 \begin{cproof}{lemma:doubleintegral:zero}
\Nobs that the assumption that \begin{equation}
    \int_\scra^\scrb (\alpha_1 x + \beta_1) \, \d x =
    \int_\scrb^\scrc (\alpha_2 x + \beta_2)  \, \d x 
    = 0
\end{equation} 
implies that 
\begin{equation}\label{lemma:doubleintegral:zero:beta}
    \beta_1 =- \frac{ \alpha_1 (\scra+\scrb)}2 \qandq \beta_2=- \frac{ \alpha_2 (\scrb+\scrc)}2.
\end{equation}
Combining this with the assumption that $\alpha_1 \scrb +\beta_1= \alpha_2 \scrb +\beta_2$ demonstrates that
\begin{equation}
    \frac{ \alpha_1 (\scrb-\scra)}2=\alpha_1 \scrb - \frac{ \alpha_1 (\scra+\scrb)}2=\alpha_2 \scrb - \frac{ \alpha_2 (\scrb+\scrc)}2=\frac{ \alpha_2 (\scrb-\scrc)}2.
\end{equation}
\Hence that
\begin{equation}
    \alpha_2= \frac{\alpha_1(\scrb-\scra)}{\scrb-\scrc}.
\end{equation}
This, \cref{lemma:doubleintegral:zero:beta}, and the assumption that
\begin{equation} 
     \int_\scra^\scrb x(\alpha_1 x + \beta_1) \, \d x +\int_\scrb^\scrc x(\alpha_2 x + \beta_2)  \, \d x 
    = 0
\end{equation}
ensure that
\begin{equation}
\begin{split}
    0&=\int_\scra^\scrb x\left(\alpha_1 x - \frac{ \alpha_1 (\scra+\scrb)}2\right) \, \d x +\int_\scrb^\scrc x\left(\alpha_2 x - \frac{ \alpha_2 (\scrb+\scrc)}2\right)  \, \d x \\
    &=\int_\scra^\scrb \alpha_1x\left( x - \frac{ \scra+\scrb}2\right) \, \d x +\int_\scrb^\scrc \alpha_2x\left( x - \frac{ \scrb+\scrc}2\right)  \, \d x \\
    &=\int_\scra^\scrb \alpha_1x\left( x - \frac{ \scra+\scrb}2\right) \, \d x +\frac{\scrb-\scra}{\scrb-\scrc}\int_\scrb^\scrc \alpha_1 x\left( x - \frac{ \scrb+\scrc}2\right) \, \d x.
\end{split}
\end{equation}
\Hence that
\begin{equation}
\begin{split}
    0&=\int_\scra^\scrb x\left( x - \frac{ \scra+\scrb}2\right) \, \d x -\frac{\scrb-\scra}{\scrc-\scrb}\int_\scrb^\scrc  x\left( x - \frac{ \scrb+\scrc}2\right) \, \d x\\
    &=\int_\scra^\scrb (x-\scra)\left( x - \frac{ \scra+\scrb}2\right) \, \d x -\frac{\scrb-\scra}{\scrc-\scrb}\int_\scrb^\scrc  (x-\scrb)\left( x - \frac{ \scrb+\scrc}2\right) \, \d x\\
    &=\int_0^{\scrb-\scra} x\left( x - \frac{ \scrb-\scra}2\right) \, \d x -\frac{\scrb-\scra}{\scrc-\scrb}\int_0^{\scrc-\scrb}  x\left( x - \frac{ \scrc-\scrb}2\right) \, \d x \\
    &=(\scrb-\scra)^3\int_0^1 x\left( x - \frac{ 1}2\right) \, \d x -(\scrb-\scra)(\scrc-\scrb)^2\int_0^1  x\left( x - \frac{1}2\right) \, \d x.
\end{split}
\end{equation}
\Hence that
\begin{equation}
    0=[(\scrb-\scra)^2-(\scrc-\scrb)^2]\int_0^1  x\left( x - \frac{1}2\right) \, \d x =\frac1{12}[(\scrb-\scra)^2-(\scrc-\scrb)^2].
\end{equation}
\Hence that $ \scrb-\scra= \scrc -\scrb $.
Combining this with \cref{lemma:doubleintegral:zero:beta} establishes that $ \scrb-\scra= \scrc -\scrb $, $\beta_1 =- \nicefrac{\alpha_1(\scra+\scrb)}2  $, and $\beta_2=- \nicefrac{\alpha_2(\scrb+\scrc)}2  $.
 \end{cproof}
\begin{lemma}\label{lemma:f:case1}
Let $f\colon \R^3 \to \R$ satisfy for all $\bfa,\bfn \in \R$, $\bfq \in [0,\nicefrac12]$ that
satisfy \begin{equation}\label{lemma:f:case1:assumption}
    f(\bfa,\bfn,\bfq)=\frac{\bfq\bfn^2}{3}+
      \frac{\bfa^2(1-\bfq^3)}{3}+\frac12-\frac{3\bfa}4
      +(-\bfa \bfq+\bfn)^2(1-\bfq)+(\bfa\bfq-\bfn)
      +\bfa(-\bfa\bfq+\bfn)(1-\bfq^2).
\end{equation}
Then \begin{equation}\label{lemma:f:case1:thesis}
    f(\bfa,\bfn,\bfq)\geq f(\nicefrac{16}{9},0,\nicefrac14)=\nicefrac1{18}.
\end{equation} 
\end{lemma}
\begin{cproof}{lemma:f:case1}
\Nobs that for all $\bfq \in [0,\nicefrac12]$ it holds that \begin{equation}\label{lemma:f:case1:infty}
    \liminf_{\norm{(\bfa,\bfn)}\to \infty}f(\bfa,\bfn,\bfq)=\infty.
\end{equation}
This implies that the minimum of $f$ occurs in the case $\nabla f(\bfa,\bfn,\bfq)=0$, in the case $\bfq=0$, or in the case $\bfq=\nicefrac12$.
\Moreover \cref{lemma:f:case1:assumption} assures for all $\bfa,\bfn \in \R$, $\bfq \in [0,\nicefrac12]$ that 
\begin{equation}\label{lemma:f:case1:grad}
\begin{split}
    \frac{\partial }{\partial \bfa}f(\bfa,\bfn,\bfq)&=\frac{2\bfa (1-\bfq^3)}{3} -\frac{3}4-2\bfq(-\bfa \bfq+\bfn)(1-\bfq)+\bfq+(-\bfa\bfq+\bfn)(1-\bfq^2)\\
    &\quad +\bfa(-\bfq)(1-\bfq^2)= \frac{2\bfa }{3}-\frac{2\bfa\bfq^3}{3}-\frac{3}4 +2 \bfa \bfq^2-2\bfn \bfq -2\bfa \bfq^3+2\bfn\bfq^2\\ & \quad+\bfq-\bfa\bfq+\bfn  +\bfa\bfq^3-\bfn\bfq^2-\bfa\bfq  + \bfa \bfq^3 \\
    &= \frac{2\bfa }{3}-\frac{2\bfa\bfq^3}{3}-\frac{3}4 +2\bfa\bfq^2-2\bfn\bfq+\bfn\bfq^2+\bfq-2\bfa\bfq+\bfn,  \\
    \frac{\partial }{\partial \bfn}f(\bfa,\bfn,\bfq)&=\frac{2\bfn\bfq}{3}+2(-\bfa \bfq+\bfn)(1-\bfq)-1+\bfa(1-\bfq^2)\\
    &=  \frac{2\bfn\bfq}{3}-2\bfa \bfq+2\bfn+2\bfa \bfq^2 -2\bfn \bfq -1 +\bfa -\bfa \bfq^2 \\
    &=  -\frac{4\bfn\bfq}{3}-2\bfa \bfq+2\bfn+\bfa \bfq^2 -1 +\bfa, \qand \\
    \frac{\partial }{\partial \bfq}f(\bfa,\bfn,\bfq)&=\frac{\bfn^2}{3}-
      \bfa^2\bfq^2-2\bfa(-\bfa \bfq+\bfn)(1-\bfq)-(-\bfa \bfq+\bfn)^2+\bfa-\bfa^2(1-\bfq^2)\\
    &\quad -2\bfa(-\bfa\bfq+\bfn)\bfq = \frac{\bfn^2}{3}-
      \bfa^2\bfq^2
    +2\bfa^2 \bfq-2\bfa\bfn  -2\bfa^2 \bfq^2 +2\bfa\bfn \bfq
    -\bfa^2\bfq^2 \\
    & \quad -\bfn^2+2\bfa \bfn \bfq+\bfa-\bfa^2+\bfa^2\bfq^2+2\bfa^2 \bfq^2 -2 \bfa \bfn \bfq\\
    &= -\frac{2\bfn^2}{3}-\bfa^2 \bfq^2
    +2\bfa^2 \bfq -2\bfa\bfn   +2\bfa\bfn \bfq
     +\bfa-\bfa^2.
\end{split}
\end{equation}
This implies that for all $\bfa,\bfn \in \R$, $\bfq \in [0,\nicefrac12]$ such that $\frac{\partial }{\partial \bfn}f(\bfa,\bfn,\bfq)=0$ it holds that 
\begin{equation}
    0=-\frac{4\bfq\bfn}{3}-2\bfa \bfq+2\bfn+\bfa \bfq^2  -1 +\bfa=2\bfn\left(-\frac{2\bfq}{3}+1\right)
    -2\bfa \bfq+\bfa \bfq^2 -1+\bfa.
\end{equation}
\Hence that for all $\bfa,\bfn \in \R$, $\bfq \in [0,\nicefrac12]$ such that $\frac{\partial }{\partial \bfn}f(\bfa,\bfn,\bfq)=0$ it holds that 
\begin{equation}\label{lemma:f:case1:n}
    \bfn =\frac12(2\bfa \bfq-\bfa \bfq^2 +1-\bfa)\left(-\frac{2\bfq}{3}+1\right)^{-1}.
\end{equation}
This and \cref{lemma:f:case1:grad} assure that for all $\bfa,\bfn \in \R$, $\bfq \in [0,\nicefrac12]$ such that $\frac{\partial }{\partial \bfn}f(\bfa,\bfn,\bfq)=\frac{\partial }{\partial \bfq}f(\bfa,\bfn,\bfq)=0$ it holds that \begin{equation}
\begin{split}
    0&=-\frac1{6}(2\bfa \bfq-\bfa \bfq^2 +1-\bfa)^2\left(-\frac{2\bfq}{3}+1\right)^{-2}
    -\bfa^2 \bfq^2
    +2\bfa^2 \bfq  +\bfa-\bfa^2 \\
    & \quad -\bfa(2\bfa \bfq-\bfa \bfq^2 +1-\bfa)\left(-\frac{2\bfq}{3}+1\right)^{-1}  +\bfa \bfq (2\bfa \bfq-\bfa \bfq^2 +1-\bfa)\left(-\frac{2\bfq}{3}+1\right)^{-1} 
    \\
     &=-\frac1{6}(2\bfa \bfq-\bfa \bfq^2 +1-\bfa)^2\left(-\frac{2\bfq}{3}+1\right)^{-2}+( -\bfa^2 \bfq^2 +2\bfa^2 \bfq +\bfa-\bfa^2)\\
     & \quad +(2\bfa \bfq-\bfa \bfq^2 +1-\bfa)\left(-\frac{2\bfq}{3}+1\right)^{-1} (-\bfa+\bfa \bfq).
\end{split}
\end{equation}
\Hence that  for all $\bfa,\bfn \in \R$, $\bfq \in [0,\nicefrac12]$ such that $\frac{\partial }{\partial \bfn}f(\bfa,\bfn,\bfq)=\frac{\partial }{\partial \bfq}f(\bfa,\bfn,\bfq)=0$ it holds that
\begin{equation}
\begin{split}
    0&=(2\bfa \bfq-\bfa \bfq^2 +1-\bfa)  \bigg(-2\bfa \bfq+\bfa \bfq^2 -1+\bfa+6\bfa\left(-\frac{2\bfq}{3}+1\right)^{2}\\
    & \quad +6(-\bfa+\bfa\bfq)\left(-\frac{2\bfq}{3}+1\right)\bigg)= (2\bfa \bfq-\bfa \bfq^2 +1-\bfa)  \bigg(-2\bfa \bfq+\bfa \bfq^2 -1+\bfa \\
    & \quad + \frac{8\bfa\bfq^2}{3}+6\bfa-8\bfa\bfq
    +4\bfa \bfq-6\bfa-4\bfa\bfq^2+6\bfa\bfq\bigg)\\
    &=(2\bfa \bfq-\bfa \bfq^2 +1-\bfa)\bigg(-\frac{\bfa \bfq^2}3-1+\bfa \bigg)\\
    &=\bigg(\bfq-1-\frac{\sqrt{\bfa}}\bfa\bigg)\bigg(\bfq-1+\frac{\sqrt{\bfa}}\bfa\bigg)\bigg(\bfq-\sqrt{3-\frac{3}{\bfa}}\bigg)\bigg(\bfq+\sqrt{3-\frac{3}{\bfa}}\bigg).
\end{split}
\end{equation}
\Nobs that  \cref{lemma:f:case1:n} assures that in the case $\bfq=1+\nicefrac{\sqrt{\bfa}}{\bfa}$ for all $\bfa,\bfn \in \R$ such that $\frac{\partial }{\partial \bfn}f(\bfa,\bfn,\bfq)=0$ it holds that $\bfn=0$ and 
\begin{equation}
\begin{split}
    \frac{\partial }{\partial \bfa}f(\bfa,\bfn,\bfq)&= \frac{2\bfa }{3}-\frac{2\bfa\bfq^3}{3}-\frac{3}4 +2\bfa\bfq^2-2\bfn\bfq+\bfn\bfq^2+\bfq-2\bfa\bfq+\bfn \\
    &=\frac{2\bfa }{3}-\frac{2\bfa(1+\nicefrac{\sqrt{\bfa}}{\bfa})^3}{3}-\frac{3}4 +2\bfa(1+\nicefrac{\sqrt{\bfa}}{\bfa})^2+1+\nicefrac{\sqrt{\bfa}}{\bfa}-2\bfa(1+\nicefrac{\sqrt{\bfa}}{\bfa})\\
    &=\frac{3\sqrt{\bfa}+4}{12\sqrt{\bfa}}\geq\frac{3}{12}.
\end{split}
\end{equation} 
In the following we distinguish between the case $\bfq=1-\nicefrac{\sqrt{\bfa}}{\bfa}$, the case $\bfq=(3-\nicefrac{3}{\bfa})^{\nicefrac12}$, the case $\bfq=0$, and the case $\bfq=\nicefrac12$.
We first prove \cref{lemma:f:case1:thesis} in the case \begin{equation}\label{lemma:f:case1:q1}
\bfq=1-\frac{\sqrt{\bfa}}{\bfa}.
\end{equation}
\Nobs that \cref{lemma:f:case1:grad}, \cref{lemma:f:case1:n}, and \cref{lemma:f:case1:q1} imply that for all $\bfa,\bfn \in \R$ such that $\nabla f(\bfa,\bfn,\bfq)=0$ it holds that $\bfn=0$ and 
\begin{equation}
\begin{split}
    0=\frac{\partial }{\partial \bfa}f(\bfa,\bfn,\bfq)&= \frac{2\bfa }{3}-\frac{2\bfa\bfq^3}{3}-\frac{3}4 +2\bfa\bfq^2-2\bfn\bfq+\bfn\bfq^2+\bfq-2\bfa\bfq+\bfn \\
    &=\frac{2\bfa }{3}-\frac{2\bfa(1-\nicefrac{\sqrt{\bfa}}{\bfa})^3}{3}-\frac{3}4 +2\bfa(1-\nicefrac{\sqrt{\bfa}}{\bfa})^2+1-\nicefrac{\sqrt{\bfa}}{\bfa}-2\bfa(1-\nicefrac{\sqrt{\bfa}}{\bfa})\\
    &=\frac{3\sqrt{\bfa}-4}{12\sqrt{\bfa}}.
\end{split}
\end{equation}
\Hence that for all $\bfa,\bfn \in \R$ such that $\nabla f(\bfa,\bfn,\bfq)=0$ it holds that $\bfa=\nicefrac{16}{9}$. This and the fact that $\bfn=0$ ensure that 
\begin{equation}
    f\bigg(\nicefrac{16}{9},0,1-\frac{\nicefrac{4}{3}}{\nicefrac{16}{9}}\bigg)=f(\nicefrac{16}{9},0,\nicefrac14)=\frac1{18}.
\end{equation}
This implies \cref{lemma:f:case1:thesis} in the case
$\bfq=1-\frac{\sqrt{\bfa}}{\bfa}$.
We show \cref{lemma:f:case1:thesis} in the case
\begin{equation}\label{lemma:f:case1:q2}
    \bfq=(3-\nicefrac{3}{\bfa})^{\nicefrac12}.
\end{equation}
\Nobs that \cref{lemma:f:case1:n} and \cref{lemma:f:case1:q2} show that for all $\bfa,\bfn \in \R$ such that $\frac{\partial }{\partial \bfn}f(\bfa,\bfn,\bfq)=0$ it holds that
\begin{equation}
    \bfn=\frac{3\sqrt{\bfa}(-2\bfa+\sqrt{3\bfa^2-3\bfa}+2)}{2\sqrt{3\bfa-3}-3\sqrt{\bfa}}.
\end{equation}
This and \cref{lemma:f:case1:grad} ensure that for all $\bfa,\bfn \in \R$ such that $\nabla f(\bfa,\bfn,\bfq)=0$ it holds that
\begin{equation}
    \bfa=\frac98, \qquad \bfn=\frac{3\sqrt{3}}{8}, \qandq \bfq=\frac1{\sqrt{3}}.
\end{equation}
\Hence that 
\begin{equation}
    f(\nicefrac98,\nicefrac{3\sqrt{3}}{8},\nicefrac1{\sqrt{3}})=\frac5{64}.
\end{equation}
This demonstrates \cref{lemma:f:case1:thesis} in the case $\bfq=(3-\nicefrac{3}{\bfa})^{\nicefrac12}.$
We now prove \cref{lemma:f:case1:thesis} in the case
\begin{equation}\label{lemma:f:case1:q3}
    \bfq=0.
\end{equation}
\Nobs that \cref{lemma:f:case1:n} and \cref{lemma:f:case1:q3} ensure that for all $\bfa,\bfn \in \R$ such that $\frac{\partial }{\partial \bfn}f(\bfa,\bfn,0)=0$ it holds that 
\begin{equation}\label{lemma:f:case1:q3:n}
    \bfn =\frac12(1-\bfa).
\end{equation}
Combining this and \cref{lemma:f:case1:grad} shows that for all $\bfa,\bfn \in \R$ such that $\frac{\partial }{\partial \bfn}f(\bfa,\bfn,0)=\frac{\partial }{\partial \bfa}f(\bfa,\bfn,0)=0$ it holds that 
\begin{equation}
    0=\frac{\partial }{\partial \bfa}f(\bfa,\bfn,0)= \frac{2\bfa }{3}-\frac{3}4+\bfn = \frac{2\bfa}{3}-\frac{3}4+\frac12(1-\bfa)=\frac\bfa6-\frac14.
\end{equation}
This and \cref{lemma:f:case1:infty} imply that for all $\bfa,\bfn \in \R$ it holds that 
\begin{equation}
    f(\bfa,\bfn,0)\geq f(\nicefrac32,-\nicefrac14,0)=\frac1{16}.
\end{equation}
This assures \cref{lemma:f:case1:thesis} in the case $\bfq=0$.
We demonstrate \cref{lemma:f:case1:thesis} in the case
\begin{equation}\label{lemma:f:case1:q4}
    \bfq=\frac12.
\end{equation}
\Nobs that \cref{lemma:f:case1:n} and \cref{lemma:f:case1:q4} ensure that for all $\bfa,\bfn \in \R$ such that $\frac{\partial }{\partial \bfn}f(\bfa,\bfn,\nicefrac12)=0$ it holds that 
\begin{equation}\label{lemma:f:case1:q4:n}
    \bfn =\frac12\left(\bfa-\frac\bfa4 +1-\bfa\right)\left(-\frac{2}{6}+1\right)^{-1}=-\frac{3\bfa}{16}+\frac34.
\end{equation}
Combining this and \cref{lemma:f:case1:grad} shows that for all $\bfa,\bfn \in \R$ such that $\frac{\partial }{\partial \bfn}f(\bfa,\bfn,\nicefrac12)=\frac{\partial }{\partial \bfa}f(\bfa,\bfn,\nicefrac12)=0$ it holds that 
\begin{equation}
\begin{split}
      0=\frac{\partial }{\partial \bfa}f(\bfa,\bfn,\nicefrac12)&= \frac{2\bfa }{3}-\frac{\bfa}{12}-\frac{3}4 +\frac\bfa2-\bfn+\frac\bfn4+\frac12-\bfa+\bfn \\
      &=\frac{\bfa}{12}+\frac\bfn4-\frac14=\frac{\bfa}{12}  -\frac{3\bfa}{64}+\frac3{16} -\frac14=\frac{7\bfa}{192} -\frac1{16}.
\end{split}
\end{equation}
This and \cref{lemma:f:case1:infty} imply that for all $\bfa,\bfn \in \R$ it holds that 
\begin{equation}
    f(\bfa,\bfn,\nicefrac12)\geq f(\nicefrac{12}{7},\nicefrac37,\nicefrac12)=\frac1{14}.
\end{equation}
This ensures \cref{lemma:f:case1:thesis} in the case $\bfq=0$.
\end{cproof}
\begin{lemma}\label{lemma:f:case2}
Let $f\colon \R^3 \to \R$ satisfy for all $\bfa,\bfn \in \R$, $\bfq \in [\nicefrac12,1]$ that
satisfy \begin{equation}\label{lemma:f:case2:assumption}
\begin{split}
    f(\bfa,\bfn,\bfq)&=\frac{\bfa^2\bfq^3}{3}+(-\bfa\bfq+\bfn)^2\bfq-\bfa \left(\bfq^2-\frac14\right) +\bfa (-\bfa \bfq+\bfn)\bfq^2+  \frac13(1-\bfq) (\bfn-1)^2 \\ 
    & \quad +(1+2\bfa\bfq-2\bfn)\left(\bfq-\frac12\right).
\end{split}
\end{equation}
Then \begin{equation}\label{lemma:f:case2:thesis}
    f(\bfa,\bfn,\bfq)\geq f(\nicefrac{16}{9},1,\nicefrac34)=\nicefrac1{18}.
\end{equation} 
\end{lemma}
\begin{cproof}{lemma:f:case2}
\Nobs that for all $\bfq \in [\nicefrac12,1]$ it holds that \begin{equation}\label{lemma:f:case2:infty}
    \liminf_{\norm{(\bfa,\bfn)}\to \infty}f(\bfa,\bfn,\bfq)=\infty.
\end{equation}
This implies that the minimum of $f$ occurs in the case $\nabla f(\bfa,\bfn,\bfq)=0$, in the case $\bfq=\nicefrac12$, or in the case $\bfq=1$.
\Moreover \cref{lemma:f:case2:assumption} assures for all $\bfa,\bfn \in \R$, $\bfq \in [0,\nicefrac12]$ that 
\begin{equation}\label{lemma:f:case2:grad}
\begin{split}
    \frac{\partial }{\partial \bfa}f(\bfa,\bfn,\bfq)&=\frac{2\bfa \bfq^3}{3} -2(-\bfa \bfq+\bfn)\bfq^2-\bigg(\bfq^2-\frac14\bigg)+(-\bfa\bfq+\bfn)\bfq^2-\bfa\bfq^3\\
    &\quad +2\bfq\bigg(\bfq-\frac12\bigg)= \frac{2\bfa\bfq^3}{3} +2 \bfa \bfq^3-2\bfn \bfq^2 - \bfq^2+\frac14-\bfa\bfq^3+\bfn\bfq^2\\ & \quad -\bfa\bfq^3+2\bfq^2-\bfq \\
    &= \frac{2\bfa\bfq^3}{3}-\bfn \bfq^2+\bfq^2+\frac14 -\bfq,  \\
    \frac{\partial }{\partial \bfn}f(\bfa,\bfn,\bfq)&=2(-\bfa \bfq+\bfn)\bfq+\bfa\bfq^2+\frac23 (1-\bfq)(\bfn-1)-2\bigg(\bfq-\frac12\bigg)\\
    &= -2\bfa \bfq^2+2\bfn\bfq+\bfa \bfq^2+\frac{2\bfn}3 -\frac23-\frac{2\bfn\bfq}3+\frac{2\bfq}{3} -2 \bfq +1 \\
    &=  -\bfa\bfq^2+\frac{4\bfn\bfq}{3}+\frac{2\bfn}3+\frac13-\frac{4\bfq}{3}, \qand \\
    \frac{\partial }{\partial \bfq}f(\bfa,\bfn,\bfq)&=\bfa^2\bfq^2-2\bfa\bfq(-\bfa \bfq+\bfn)+(-\bfa \bfq+\bfn)^2-2\bfa\bfq-\bfa^2\bfq^2+2\bfa\bfq(-\bfa \bfq+\bfn)\\
    &\quad -\frac13(\bfn-1)^2+2\bfa\bigg(\bfq-\frac12\bigg)+1+2\bfa\bfq-2\bfn\\
    &= \bfa^2\bfq^2+2\bfa^2\bfq^2-2\bfa \bfq \bfn+\bfa^2 \bfq^2+\bfn^2-2\bfa\bfq\bfn-2\bfa\bfq-\bfa^2\bfq^2-2\bfa^2\bfq^2 \\
    & \quad +2\bfa \bfq\bfn -\frac{\bfn^2}3-\frac13+\frac{2\bfn}{3}+2\bfa\bfq-\bfa+1+2\bfa\bfq-2\bfn\\
    &= \bfa^2 \bfq^2-2\bfa\bfn \bfq+\frac{2\bfn^2}3+2\bfa\bfq+\frac23-\frac{4\bfn}{3}-\bfa.
\end{split}
\end{equation}
This implies that for all $\bfa,\bfn \in \R$, $\bfq \in [\nicefrac12,1]$ such that $\frac{\partial }{\partial \bfn}f(\bfa,\bfn,\bfq)=0$ it holds that 
\begin{equation}
    0=-\bfa\bfq^2+\frac{4\bfn\bfq}{3}+\frac{2\bfn}3+\frac13-\frac{4\bfq}{3}
    =\frac{2\bfn}3\left(2\bfq+1\right)
    -\bfa \bfq^2 +\frac13-\frac{4\bfq}{3}.
\end{equation}
\Hence that for all $\bfa,\bfn \in \R$, $\bfq \in [\nicefrac12,1]$ such that $\frac{\partial }{\partial \bfn}f(\bfa,\bfn,\bfq)=0$ it holds that 
\begin{equation}\label{lemma:f:case2:n}
    \bfn =\frac32\left(\bfa \bfq^2 -\frac13+\frac{4\bfq}{3}\right)(2\bfq+1)^{-1}.
\end{equation}
This and \cref{lemma:f:case2:grad} assure that for all $\bfa,\bfn \in \R$, $\bfq \in [0,\nicefrac12]$ such that $\frac{\partial }{\partial \bfn}f(\bfa,\bfn,\bfq)=\frac{\partial }{\partial \bfq}f(\bfa,\bfn,\bfq)=0$ it holds that \begin{equation}
\begin{split}
    0&=\bfa^2\bfq^2-3\bfa \bfq\left(\bfa \bfq^2 -\frac13+\frac{4\bfq}{3}\right)(2\bfq+1)^{-1}  +\frac32\left(\bfa \bfq^2 -\frac13+\frac{4\bfq}{3}\right)^2(2\bfq+1)^{-2}\\
    &\quad +2\bfa\bfq+\frac23-2\left(\bfa \bfq^2 -\frac13+\frac{4\bfq}{3}\right)(2\bfq+1)^{-1}-\bfa\\
     &=\bfa^2\bfq^2+2\bfa\bfq-\bfa+\frac23+\left(-3\bfa^2 \bfq^3 +\bfa\bfq-4\bfa\bfq^2-2\bfa\bfq^2+\frac32-\frac{8\bfq}{3}\right)(2\bfq+1)^{-1}\\
     &\quad + \frac32\left(\bfa^2 \bfq^4 +\frac19+\frac{16\bfq^2}{9}-\frac{2\bfa\bfq^2}3+\frac{8\bfa\bfq^3}{3}-\frac{8\bfq}{9}\right)(2\bfq+1)^{-2}.
\end{split}
\end{equation}
\Hence that  for all $\bfa,\bfn \in \R$, $\bfq \in [\nicefrac12,1]$ such that $\frac{\partial }{\partial \bfn}f(\bfa,\bfn,\bfq)=\frac{\partial }{\partial \bfq}f(\bfa,\bfn,\bfq)=0$ it holds that
\begin{equation}
\begin{split}
    0&=\bigg(\bfa^2\bfq^2+2\bfa\bfq-\bfa+\frac23\bigg)(2\bfq+1)^2  
    +\bigg(-3\bfa^2 \bfq^3 +\bfa\bfq-4\bfa\bfq^2-2\bfa\bfq^2+\frac23-\frac{8\bfq}{3}\bigg) \\
    & \quad(2\bfq+1) + \frac32\left(\bfa^2 \bfq^4 +\frac19+\frac{16\bfq^2}{9}-\frac{2\bfa\bfq^2}3+\frac{8\bfa\bfq^3}{3}-\frac{8\bfq}{9}\right)\\
    &=\bigg(\bfa^2\bfq^2+2\bfa\bfq-\bfa+\frac23\bigg)(4\bfq^2+1+4\bfq) -6\bfa^2 \bfq^4 +2\bfa\bfq^2-12\bfa\bfq^3+\frac{4\bfq}{3} \\
    & \quad -\frac{16\bfq^2}{3} -3\bfa^2 \bfq^3 +\bfa\bfq-4\bfa\bfq^2-2\bfa\bfq^2+\frac23-\frac{8\bfq}{3} 
    + \frac{3\bfa^2 \bfq^4}2 +\frac16+\frac{8\bfq^2}{3}-\bfa\bfq^2\\
    & \quad +4\bfa\bfq^3-\frac{4\bfq}{3}=4\bfa^2\bfq^4+8\bfa\bfq^3-4\bfa\bfq^2+\frac{8\bfq^2}3+\bfa^2\bfq^2+2\bfa\bfq-\bfa+\frac23+4\bfa^2\bfq^3\\
    &\quad +8\bfa\bfq^2-4\bfa\bfq+\frac{8\bfq}3  -\frac{9\bfa^2\bfq^4}{2}-5\bfa\bfq^2-8\bfa\bfq^3+\frac{5}6-\frac{8\bfq^2}{3}-3\bfa^2\bfq^3+\bfa\bfq-\frac{8\bfq}{3}\\
    &=-\frac{\bfa^2\bfq^4}{2} -\bfa\bfq^2+\bfa^2\bfq^2 -\bfa\bfq
    -\bfa+\frac32+\bfa^2\bfq^3\\
    & = \frac{\bfa\bfq^2}{2}(-\bfa\bfq^2+2\bfa\bfq+2\bfa-3)+\frac{\bfa\bfq^2}{2}-\bfa\bfq-\bfa+\frac32\\
    &=\frac12(\bfa\bfq^2-1)(-\bfa\bfq^2+2\bfa\bfq+2\bfa-3)\\
    &=\frac12(\bfa\bfq^2-1)\bigg(\bfq-1-\sqrt{3-\frac{3}{\bfa}}\bigg)\bigg(\bfq-1+\sqrt{3-\frac{3}{\bfa}}\bigg).
\end{split}
\end{equation}
In the following we distinguish between the case $\bfq=\nicefrac{\sqrt{\bfa}}{\bfa}$, the case $\bfq=1-(3-\nicefrac{3}{\bfa})^{\nicefrac12}$, the case $\bfq=\nicefrac12$, and the case $\bfq=1$.
We first prove \cref{lemma:f:case2:thesis} in the case \begin{equation}\label{lemma:f:case2:q1}
\bfq=\frac{\sqrt{\bfa}}{\bfa}.
\end{equation}
\Nobs that \cref{lemma:f:case2:grad}, \cref{lemma:f:case2:n}, and \cref{lemma:f:case2:q1} imply that for all $\bfa,\bfn \in \R$ such that $\nabla f(\bfa,\bfn,\bfq)=0$ it holds that 
\begin{equation}
    \bfn =\frac32\left(\frac{\bfa^2}{\bfa^2} -\frac13+\frac{4\sqrt{\bfa}}{3\bfa}\right)\left(\frac{2\sqrt{\bfa}}{\bfa}+1\right)^{-1}=1.
\end{equation} and 
\begin{equation}
\begin{split}
    0=\frac{\partial }{\partial \bfa}f(\bfa,\bfn,\bfq)&= \frac{2\bfa\bfq^3}{3}-\bfn \bfq^2+\bfq^2+\frac14 -\bfq =\frac{2\bfa(\nicefrac{\sqrt{\bfa}}{\bfa})^3}{3}+\frac{1}4-\nicefrac{\sqrt{\bfa}}{\bfa}\\
    &=\frac{3\sqrt{\bfa}-4}{12\sqrt{\bfa}}.
\end{split}
\end{equation}
\Hence that for all $\bfa,\bfn \in \R$ such that $\nabla f(\bfa,\bfn,\bfq)=0$ it holds that $\bfa=\nicefrac{16}{9}$. This and the fact that $\bfn=1$ ensure that 
\begin{equation}
    f\bigg(\nicefrac{16}{9},1,\frac{\nicefrac{4}{3}}{\nicefrac{16}{9}}\bigg)=f(\nicefrac{16}{9},0,\nicefrac34)=\frac1{18}.
\end{equation}
This implies \cref{lemma:f:case2:thesis} in the case
$\bfq=\frac{\sqrt{\bfa}}{\bfa}$.
We show \cref{lemma:f:case2:thesis} in the case
\begin{equation}\label{lemma:f:case2:q2}
    \bfq=1-(3-\nicefrac{3}{\bfa})^{\nicefrac12}.
\end{equation}
\Nobs that \cref{lemma:f:case2:n} and \cref{lemma:f:case2:q2} show that for all $\bfa,\bfn \in \R$ such that $\frac{\partial }{\partial \bfn}f(\bfa,\bfn,\bfq)=0$ it holds that
\begin{equation}
\bfn = \frac{6\bfa^2-3\bfa\sqrt{3\bfa^2-3\bfa}-3\bfa-2\sqrt{3\bfa^2-3\bfa}}{3\bfa-2\sqrt{3\bfa^2-3\bfa}}
\end{equation}
This and \cref{lemma:f:case2:grad} ensure that for all $\bfa,\bfn \in \R$ such that $\nabla f(\bfa,\bfn,\bfq)=0$ it holds that
\begin{equation}
    \bfa=\frac98, \qquad \bfn=1-\frac{3\sqrt{3}}{8}, \qandq \bfq=1-\frac1{\sqrt{3}}.
\end{equation}
\Hence that 
\begin{equation}
    f(\nicefrac98,1-\nicefrac{3\sqrt{3}}{8},1-\nicefrac1{\sqrt{3}})=\frac5{64}.
\end{equation}
This demonstrates \cref{lemma:f:case2:thesis} in the case $\bfq=(3-\nicefrac{3}{\bfa})^{\nicefrac12}.$
We now prove \cref{lemma:f:case2:thesis} in the case
\begin{equation}\label{lemma:f:case2:q3}
    \bfq=\nicefrac12.
\end{equation}
\Nobs that \cref{lemma:f:case2:n} and \cref{lemma:f:case2:q3} ensure that for all $\bfa,\bfn \in \R$ such that $\frac{\partial }{\partial \bfn}f(\bfa,\bfn,0)=0$ it holds that 
\begin{equation}\label{lemma:f:case2:q3:n}
    \bfn =\frac32\bigg(\frac\bfa4-\frac13+\frac23\bigg)(1+1)^{-1}=\frac{3\bfa}{16}+\frac14.
\end{equation}
Combining this and \cref{lemma:f:case2:grad} shows that for all $\bfa,\bfn \in \R$ such that $\frac{\partial }{\partial \bfn}f(\bfa,\bfn,\nicefrac12)=\frac{\partial }{\partial \bfa}f(\bfa,\bfn,\nicefrac12)=0$ it holds that 
\begin{equation}
    0=\frac{\partial }{\partial \bfa}f(\bfa,\bfn,\nicefrac12)= \frac{2\bfa }{24}-\frac{\bfn}{4}+\frac14+\frac14-\frac12 = \frac{\bfa}{12}-\frac{3\bfa}{64}-\frac{1}{16}=\frac{7\bfa}{192}-\frac1{16}.
\end{equation}
This and \cref{lemma:f:case2:infty} imply that for all $\bfa,\bfn \in \R$ it holds that 
\begin{equation}
    f(\bfa,\bfn,\nicefrac12)\geq f(\nicefrac{12}7,\nicefrac47,\nicefrac12)=\frac1{14}.
\end{equation}
This assures \cref{lemma:f:case2:thesis} in the case $\bfq=\nicefrac12$.
We demonstrate \cref{lemma:f:case2:thesis} in the case
\begin{equation}\label{lemma:f:case2:q4}
    \bfq=1.
\end{equation}
\Nobs that \cref{lemma:f:case2:n} and \cref{lemma:f:case2:q4} ensure that for all $\bfa,\bfn \in \R$ such that $\frac{\partial }{\partial \bfn}f(\bfa,\bfn,1)=0$ it holds that 
\begin{equation}\label{lemma:f:case2:q4:n}
    \bfn =\frac32\left(\bfa-\frac13 +\frac43\right)\left(2+1\right)^{-1}=\frac{\bfa}{2}+\frac12.
\end{equation}
Combining this and \cref{lemma:f:case2:grad} shows that for all $\bfa,\bfn \in \R$ such that $\frac{\partial }{\partial \bfn}f(\bfa,\bfn,1)=\frac{\partial }{\partial \bfa}f(\bfa,\bfn,1)=0$ it holds that 
\begin{equation}
\begin{split}
      0=\frac{\partial }{\partial \bfa}f(\bfa,\bfn,\nicefrac12)&= \frac{2\bfa }{3}-\bfn+1+\frac14-1 =
      \frac{2\bfa }{3}-\frac{\bfa}{2}-\frac12+\frac14
      =\frac{\bfa}{6} -\frac1{4}.
\end{split}
\end{equation}
This and \cref{lemma:f:case2:infty} imply that for all $\bfa,\bfn \in \R$ it holds that 
\begin{equation}
    f(\bfa,\bfn,1)\geq f(\nicefrac32,\nicefrac54,1)=\frac1{16}.
\end{equation}
This ensures \cref{lemma:f:case2:thesis} in the case $\bfq=0$.
\end{cproof}
\begin{lemma}\label{cor:H:updown}
Let $H \in \N$, $c\in (0,\infty)$ and let $(\bfq_{0,n})_{n\in \N}\subseteq [0,\nicefrac12]$, $(\bfq_{1,n})_{n\in \N}\subseteq [\nicefrac12,1]$, $(\bfm_{0,n})_{n\in \N},(\bfm_{1,n})_{n\in \N} \subseteq [1,H]$, $(\bfa_{n})_{n\in \N}, (\bfb_{n})_{n\in \N} \subseteq \R$ satisfy for all $n \in \N$ that
\begin{equation}\label{cor:H:updown:max}
\max\{ |\bfa_n |, |\bfb_n| \}\geq c
\end{equation}
and
\begin{equation}\label{cor:H:updown:limits}
    \begin{split}
    0&=\lim\nolimits_{n\to \infty} \left(\bfq_{1,n}-\frac12\right)^3(\bfm_{0, n} \bfa_{n}+(H-\bfm_{1,n}+1) \bfb_{n})^2\\
    &=\lim\nolimits_{n\to \infty} \left(\frac12-\bfq_{0,n}\right)^3(\bfm_{0, n} \bfa_{n}+(H-\bfm_{1,n}+1) \bfb_{n})^2 \\
      &=\lim\nolimits_{n\to \infty} (\bfb_{n})^2 (1-\bfq_{1,n})^3\\
       &=\lim\nolimits_{n\to \infty} (\bfa_{n})^2
   (\bfq_{0,n})^3.
\end{split}
\end{equation}
Then there exists a strictly increasing $n \colon \N \to \N$ such that  \begin{equation}\label{cor:H:updown:thesis}
\lim\nolimits_{k\to \infty}\bfq_{1,n(k)}=1 \qandq \lim\nolimits_{k\to \infty}\bfq_{0,n(k)}=0.
\end{equation}
\end{lemma}
\begin{cproof}{cor:H:updown}
\Nobs that \cref{cor:H:updown:max} demonstrates that there exists $C \in [c,\infty]$ which satisfies
$C=\max\{ \limsup\nolimits_{n \to \infty}|\bfa_n|, \limsup\nolimits_{n \to \infty}|\bfb_n|\}$.
In the following we distinguish between the case $\limsup\nolimits_{n \to \infty}|\bfa_n|=C$ and the case $\limsup\nolimits_{n \to \infty}|\bfb_n|=C$.
We first prove \cref{cor:H:updown:thesis} in the case \begin{equation}\label{cor:H:updown:case1}
    \limsup\nolimits_{n \to \infty}|\bfa_n|=C.
\end{equation}
\Nobs that \cref{cor:H:updown:case1} shows that there exists a strictly increasing $n \colon \N \to \N$ such that $\lim\nolimits_{k \to \infty}|\bfa_{n(k)}|=C$.
Combining this with \cref{cor:H:updown:limits} implies that 
\begin{equation}\label{cor:H:updown:zero}
    0=\lim\nolimits_{k\to \infty} (\bfa_{n(k)})^2
   (\bfq_{0,n(k)})^3=
     \lim\nolimits_{k\to \infty} C^2
     (\bfq_{0,n(k)})^3.
\end{equation}
\Hence that $\lim\nolimits_{k\to \infty}
     \bfq_{0,n(k)}=0$.
This and \cref{cor:H:updown:limits} assure that 
\begin{equation}
\begin{split}
   0& =\lim\nolimits_{k\to \infty} (\bfm_{0, n(k)} \bfa_{n(k)}+(H-\bfm_{1,n(k)}+1) \bfb_{n(k)})^2\\
    &=\lim\nolimits_{k\to \infty} \bfm_{0, n(k)} \bfa_{n(k)}+(H-\bfm_{1,n(k)}+1) \bfb_{n(k)}.
\end{split}
\end{equation}
\Hence that 
\begin{equation}
\begin{split}
    \lim\nolimits_{k\to \infty}|\bfb_{n(k)}|&=\lim\nolimits_{k\to \infty} |\bfm_{0,n(k)}  (H-\bfm_{1,n(k)}+1)^{-1}\bfa_{n(k)}|\\
    &>\lim\nolimits_{k\to \infty} |  (H+1)^{-1}\bfa_{n(k)}|=(H+1)^{-1} C.
\end{split}
\end{equation}
Combining this with \cref{cor:H:updown:limits} demonstrates that 
\begin{equation}
    0=\lim\nolimits_{k\to \infty}(\bfb_{n(k)})^2 (1-\bfq_{1,n(k)})^3=\lim\nolimits_{k\to \infty} (1-\bfq_{1,n(k)})^3.
\end{equation}
This and \cref{cor:H:updown:zero} show that 
$\lim\nolimits_{k\to \infty}\bfq_{1,n(k)}=1$ and $\lim\nolimits_{k\to \infty}\bfq_{0,n(k)}=0.$
This implies \cref{cor:H:updown:thesis} in the case  $\limsup\nolimits_{n \to \infty}|\bfa_n|=C$.
In the next step we prove \cref{cor:H:updown:thesis} in the case
\begin{equation}\label{cor:H:updown:case2}
    \limsup\nolimits_{n \to \infty}|\bfb_n|=C.
\end{equation}
\Nobs that \cref{cor:H:updown:case2} shows that there exists a strictly increasing $n \colon \N \to \N$ such that $\lim\nolimits_{k \to \infty}|\bfb_{n(k)}|=C$.
Combining this with \cref{cor:H:updown:limits} implies that 
\begin{equation}\label{cor:H:updown:zero2}
    0=\lim\nolimits_{k\to \infty} (\bfb_{n(k)})^2
   (1-\bfq_{1,n(k)})^3=
     \lim\nolimits_{k\to \infty} C^2
     (1-\bfq_{1,n(k)})^3.
\end{equation}
\Hence that $\lim\nolimits_{k\to \infty}
     \bfq_{1,n(k)}=1$.
This and \cref{cor:H:updown:limits} assure that 
\begin{equation}
\begin{split}
   0& =\lim\nolimits_{k\to \infty} (\bfm_{0, n(k)} \bfa_{n(k)}+(H-\bfm_{1,n(k)}+1) \bfb_{n(k)})^2\\
    &=\lim\nolimits_{k\to \infty} \bfm_{0, n(k)} \bfa_{n(k)}+(H-\bfm_{1,n(k)}+1) \bfb_{n(k)}.
\end{split}
\end{equation}
\Hence that 
\begin{equation}
\begin{split}
    \lim\nolimits_{k\to \infty}|\bfa_{n(k)}|&=\lim\nolimits_{k\to \infty} |(\bfm_{0,n(k)} )^{-1} (H-\bfm_{1,n(k)}+1)\bfb_{n(k)}|\\
    &>\lim\nolimits_{k\to \infty} |(H+1)^{-1}  \bfb_{n(k)}|=(H+1)^{-1} C.
\end{split}
\end{equation}
Combining this with \cref{cor:H:updown:limits} demonstrates that 
\begin{equation}
    0=\lim\nolimits_{k\to \infty}(\bfa_{n(k)})^2 (\bfq_{0,n(k)})^3=\lim\nolimits_{k\to \infty} (\bfq_{0,n(k)})^3.
\end{equation}
This and \cref{cor:H:updown:zero2} show that $\lim\nolimits_{k\to \infty}\bfq_{1,n(k)}=1$ and $\lim\nolimits_{k\to \infty}\bfq_{0,n(k)}=0.$
This implies \cref{cor:H:updown:thesis} in the case  $\limsup\nolimits_{n \to \infty}|\bfb_n|=C$.
\end{cproof}
\subsection{Properties of critical points of the risk function}\label{proporties_crit}
\begin{prop} \label{prop:H:approximate:gradient}
Assume \cref{H:setting:snn} and let $\width \in\N$. Then it holds for all $\theta \in \R^{\fd_\width}$, $i \in \{1, 2, \ldots, \width \}$ that
\begin{equation} \label{H:eq:loss:gradient}
\begin{split}
        \cG^\width_{ i} ( \theta) &= 2 \v{\theta}_i \int_{I_i^{\theta}} x  ( \realization{\width,\theta}_\infty (x) - \mathbbm{1}_{(\nicefrac12,\infty)}(x)) \, \d x , \\
        \cG^\width_{\width + i} ( \theta) &= 2 \v{\theta}_i \int_{I_i^{\theta}} (\realization{\width,\theta}_\infty (x) - \mathbbm{1}_{(\nicefrac12,\infty)}(x)) \,  \d x , \\
        \cG^\width_{2\width + i} ( \theta) &= 2 \int_{0}^1 \br[\big]{\act_\infty \rbr{ \w{\theta}_{i} x + \b{\theta}_i } } ( \realization{\width,\theta}_\infty (x) - \mathbbm{1}_{(\nicefrac12,\infty)}(x) ) \, \d x , \\
        \text{and} \qquad \cG^\width_{\fd_\width} ( \theta) &= 2 \int_{0}^1 (\realization{\width,\theta}_\infty (x) - \mathbbm{1}_{(\nicefrac12,\infty)}(x) ) \, \d x .
        \end{split}
\end{equation}
\end{prop}
\begin{cproof}{prop:H:approximate:gradient}
\Nobs that \cref{H:theo:intro:general:eq1}--\cref{H:theo:intro:general:L} imply \cref{H:eq:loss:gradient} (cf., e.g., 
\cite[Proposition~2.3]{JentzenRiekert2021}).
\end{cproof}
\begin{cor}\label{cor:H:m1}
Assume \cref{H:setting:snn}, let $\width\in \N$, $\theta\in (\cG^\width)^{-1}(\{0\})$ satisfy
for all $x\in [0,\q{\theta}_{\m^{\theta}_{0,1}}]$ that $\realization{\width,\theta}_\infty(x)=0$, $\M^{\theta}_{0}\neq\emptyset$, and $\prod_{k=1}^\width \v{\theta}_k\neq0$.
Then it holds for all $x\in [0,\q{\theta}_{\m^{\theta}_{0,2}}]$ that $\realization{\width,\theta}_\infty(x)=0$.
\end{cor}
\begin{cproof}{cor:H:m1}
Assume without loss of generality that $\q{\theta}_1 \leq \ldots \leq \q{\theta}_\width$.
\Nobs that the assumption that $\theta\in(\cG^\width)^{-1}(\{0\})$,  the assumption that $\prod_{k=1}^\width\v{\theta}_k\neq0$, and \cref{prop:H:approximate:gradient} imply that for all $j \in \M^{\theta}_{0}$ it holds that
\begin{equation}
   \int_{0}^{\q{\theta}_{\m^{\theta}_{0,1}}}  ( \realization{\width,\theta}_\infty (x) - \mathbbm{1}_{(\nicefrac12,\infty)}(x)) \, \d x = \int_{\q{\theta}_{j-1}}^{\q{\theta}_{j}}  ( \realization{\width,\theta}_\infty (x) - \mathbbm{1}_{(\nicefrac12,\infty)}(x)) \, \d x =0.
\end{equation}
 Combining this with  \cref{corollary:integralzero} demonstrates that for all $j \in \M^{\theta}_{0}$ it holds that $0=\realization{\width,\theta}_\infty(0)=\realization{\width,\theta}_\infty(\q{\theta}_j)$.
 This and piecewise linearity of $\realization{\width,\theta}_\infty$ assure that for all $x\in [0,\q{\theta}_{\m^{\theta}_{0,2}}]$ it holds that $\realization{\width,\theta}_\infty(x)=0$.
\end{cproof}
\begin{cor}\label{cor:H:l2}
Assume \cref{H:setting:snn}, let $\width \in \N$, $\theta\in (\cG^\width)^{-1}(\{0\})$ satisfy
for all $x\in [\q{\theta}_{\m^{\theta}_{1,2}},1]$ that $\realization{\width,\theta}_\infty(x)=1$, $\M^{\theta}_{1}\neq\emptyset$,
and $\prod_{k=1}^\width \v{\theta}_k\neq0$.
Then it holds for all $x\in [\q{\theta}_{\m^{\theta}_{1,1}},1]$ that $\realization{\width,\theta}_\infty(x)=1$.
\end{cor}
\begin{cproof}{cor:H:l2}
Assume without loss of generality that $\q{\theta}_1 \leq \ldots \leq \q{\theta}_\width$.
\Nobs that the assumption that $\theta\in (\cG^\width)^{-1}(\{0\})$,  the assumption that $\prod_{k=1}^\width\v{\theta}_k\neq0$, and \cref{prop:H:approximate:gradient} imply that for all $j \in \M^{\theta}_{1}$ it holds that
\begin{equation}
   \int_{\q{\theta}_{\m^{\theta}_{1,2}}}^1  ( \realization{\width,\theta}_\infty (x) - \mathbbm{1}_{(\nicefrac12,\infty)}(x)) \, \d x = \int_{\q{\theta}_{j}}^{\q{\theta}_{j+1}}  ( \realization{\width,\theta}_\infty (x) - \mathbbm{1}_{(\nicefrac12,\infty)}(x)) \, \d x =0.
\end{equation}
 Combining this with  \cref{corollary:integralzero} demonstrates that for all $j \in \M^{\theta}_{1}$ it holds that $1=\realization{\width,\theta}_\infty(1)=\realization{\width,\theta}_\infty(\q{\theta}_j)$.
 This and the fact that $\realization{\width,\theta}_\infty$ is piecewise affine linear assure that for all $x\in [\q{\theta}_{\m^{\theta}_{1,1}},1]$ it holds that $\realization{\width,\theta}_\infty(x)=1$.
\end{cproof}
\begin{lemma}\label{cor:H:w0sx}
Assume \cref{H:setting:snn}, let $\width \in \N$, $\theta\in (\cG^\width)^{-1}(\{0\})$ satisfy for all $j \in \{0,1,\ldots,\width\}$ that
$\q{\theta}_j <  \q{\theta}_{j+1}$ and $\prod_{k=1}^\width\v{\theta}_k\neq0$,
and let $i \in \{0,1,\ldots,\width\}$ satisfy
 $0<\w{\theta}_i \w{\theta}_{i+1}$ and $\q{\theta}_{i+1}\leq \nicefrac12$. Then it holds for all $x\in [0,\q{\theta}_{\m^{\theta}_{0,2}}]$ that $\realization{\width,\theta}_\infty(x)=0$.
\end{lemma}
\begin{cproof}{cor:H:w0sx}
\Nobs that the assumption that $\theta\in (\cG^\width)^{-1}(\{0\})$, the fact that $\prod_{k=1}^\width\v{\theta}_k\allowbreak \neq0<\w{\theta}_i \w{\theta}_{i+1}$, and \cref{prop:H:approximate:gradient} imply that 
\begin{equation}
    \int_{\q{\theta}_i}^{\q{\theta}_{i+1}} x  ( \realization{\width,\theta}_\infty (x) - \mathbbm{1}_{(\nicefrac12,\infty)}(x)) \, \d x = \int_{\q{\theta}_i}^{\q{\theta}_{i+1}}   ( \realization{\width,\theta}_\infty (x) - \mathbbm{1}_{(\nicefrac12,\infty)}(x)) \, \d x=0.
\end{equation}
This, \cref{lem:affine:integral:zero}, and the assumption that $\q{\theta}_{i+1}\leq \nicefrac12$ assure that for all $x \in [\q{\theta}_i,\q{\theta}_{i+1}]$ it holds that
\begin{equation}\label{cor:H:w0sx:eq0}
    \realization{\width,\theta}_\infty(x)=0.
\end{equation}
\Moreover the assumption that $\theta\in (\cG^\width)^{-1}(\{0\})$, the assumption that $\prod_{k=1}^\width\v{\theta}_k\neq0$, and \cref{prop:H:approximate:gradient} imply that for all $j \in \M^{\theta}_{0}$ it holds that
\begin{equation}
    \int_{\q{\theta}_{j-1}}^{\q{\theta}_{j}}  ( \realization{\width,\theta}_\infty (x) - \mathbbm{1}_{(\nicefrac12,\infty)}(x)) \, \d x =0.
\end{equation}
 Combining this with \cref{cor:H:w0sx:eq0} and \cref{corollary:integralzero} demonstrates that for all $j \in \M^{\theta}_{0}$ it holds that $\realization{\width,\theta}_\infty(0)=\realization{\width,\theta}_\infty(\q{\theta}_j)=0$.
 This and the fact that $\realization{\width,\theta}_\infty$ is piecewise affine linear ensure that for all $x\in [0,\q{\theta}_{\m^{\theta}_{0,2}}]$ it holds that $\realization{\width,\theta}_\infty(x)=0$.
\end{cproof}
\begin{lemma}\label{cor:H:w0sdx}
Assume \cref{H:setting:snn}, let $\width \in \N$, $\theta\in (\cG^\width)^{-1}(\{0\})$ satisfy for all $j \in \{0,1,\ldots,\width\}$ that
$\q{\theta}_j <  \q{\theta}_{j+1}$ and $\prod_{k=1}^\width\v{\theta}_k\neq0$, and let $i \in \{0,1,\ldots,\width\}$ satisfy $0<\w{\theta}_i \w{\theta}_{i+1}$ and $\nicefrac12\leq\q{\theta}_i $. Then it holds for all $x\in [\q{\theta}_{\m^{\theta}_{1,1}},1]$ that $\realization{\width,\theta}_\infty(x)=1$.
\end{lemma}
\begin{cproof}{cor:H:w0sdx}
\Nobs that the assumption that $\theta\in (\cG^\width)^{-1}(\{0\})$, the fact that $\prod_{k=1}^\width\v{\theta}_k\allowbreak \neq0<\w{\theta}_i \w{\theta}_{i+1}$, and \cref{prop:H:approximate:gradient} imply that
\begin{equation}
    \int_{\q{\theta}_i}^{\q{\theta}_{i+1}} x  ( \realization{\width,\theta}_\infty (x) - \mathbbm{1}_{(\nicefrac12,\infty)}(x)) \, \d x = \int_{\q{\theta}_i}^{\q{\theta}_{i+1}}   ( \realization{\width,\theta}_\infty (x) - \mathbbm{1}_{(\nicefrac12,\infty)}(x)) \, \d x=0.
\end{equation}
This, \cref{lem:affine:integral:zero}, and the assumption that $\nicefrac12\leq \q{\theta}_{i}$  assure that for all $x \in [\q{\theta}_i,\q{\theta}_{i+1}]$ it holds that \begin{equation}\label{cor:H:w0dx:eq0}
    \realization{\width,\theta}_\infty(x)=1.
\end{equation}
\Moreover the assumption that $\theta\in (\cG^\width)^{-1}(\{0\})$,  the assumption that $\prod_{k=1}^\width\v{\theta}_k\neq0$, and \cref{prop:H:approximate:gradient} imply that for all $j \in \M^{\theta}_{1}$ it holds that
\begin{equation}
    \int_{\q{\theta}_j}^{\q{\theta}_{j+1}}  ( \realization{\width,\theta}_\infty (x) - \mathbbm{1}_{(\nicefrac12,\infty)}(x)) \, \d x =0.
\end{equation}
 Combining this with \cref{cor:H:w0dx:eq0} and \cref{corollary:integralzero} demonstrates that for all $j \in \M^{\theta}_{1}$ it holds that $\realization{\width,\theta}_\infty(1)=\realization{\width,\theta}_\infty(\q{\theta}_j)=1$.
  This  and the fact that $\realization{\width,\theta}_\infty$ is piecewise affine linear ensure that for all $x\in [\q{\theta}_{\m^{\theta}_{1,1}},1]$ it holds that $\realization{\width,\theta}_\infty(x)=1$.
\end{cproof}
\begin{lemma}\label{lemma:H:N0}
Assume \cref{H:setting:snn} and let $\width\in \N$, $\theta\in \R^{\fd_\width}$ satisfy $\M^{\theta}_{0}\neq\emptyset\neq\M^{\theta}_{1}$,
$\prod_{k=1}^\width\v{\theta}_k\neq0$,
for all $x \in [0,\q{\theta}_{\m^{\theta}_{0,2}}]$ 
that $\realization{\width,\theta}_\infty(x)=0$, and
for all $x \in [\q{\theta}_{\m^{\theta}_{1,1}},1]$ 
that $\realization{\width,\theta}_\infty(x)=1$. Then \begin{equation}\label{lemma:H:N0:thesis}
    \cG^\width(\theta)\neq 0.
\end{equation}
\end{lemma}
\begin{cproof}{lemma:H:N0}
\Nobs that continuity of $\realization{\width,\theta}_\infty$ implies for all $x \in [\q{\theta}_{\m^{\theta}_{0,2}},\q{\theta}_{\m^{\theta}_{1,1}}]$ that $\q{\theta}_{\m^{\theta}_{0,2}}<\q{\theta}_{\m^{\theta}_{1,1}}$ and  \begin{equation}
    \realization{\width,\theta}_\infty(x)=\frac{x}{\q{\theta}_{\m^{\theta}_{1,1}}-\q{\theta}_{\m^{\theta}_{0,2}}}-\frac{\q{\theta}_{\m^{\theta}_{0,2}}}{\q{\theta}_{\m^{\theta}_{1,1}}-\q{\theta}_{\m^{\theta}_{0,2}}}.
\end{equation}
This shows that there exists $i\in\{1,2,\ldots,\width\}$ which satisfies that
\begin{equation}\label{lemma:H:N0:I}
    I_i^{\theta}\cap(\q{\theta}_{\m^{\theta}_{0,2}},\q{\theta}_{\m^{\theta}_{1,1}})=(\q{\theta}_{\m^{\theta}_{0,2}},\q{\theta}_{\m^{\theta}_{1,1}}).
\end{equation}
We prove \cref{lemma:H:N0:thesis} by contradiction. Assume that \begin{equation}\label{lemma:H:N0:zero}
    \cG^\width(\theta)= 0.
\end{equation}
\Nobs that \cref{lemma:H:N0:I}, \cref{lemma:H:N0:zero}, and \cref{prop:H:approximate:gradient} imply that 
\begin{equation}
    \int_{\q{\theta}_{\m^{\theta}_{0,2}}}^{\q{\theta}_{\m^{\theta}_{1,1}}} x  ( \realization{\width,\theta}_\infty (x) - \mathbbm{1}_{(\nicefrac12,\infty)}(x)) \, \d x = \int_{\q{\theta}_{\m^{\theta}_{0,2}}}^{\q{\theta}_{\m^{\theta}_{1,1}}}   ( \realization{\width,\theta}_\infty (x) - \mathbbm{1}_{(\nicefrac12,\infty)}(x)) \, \d x=0.
\end{equation}
Combining this and the fact that $\realization{\width,\theta}_\infty(\q{\theta}_{\m^{\theta}_{0,2}})=0$ with \cref{cor:affine1:integral:zero} demonstrates that
 for all $x \in [\q{\theta}_{\m^{\theta}_{0,2}},\q{\theta}_{\m^{\theta}_{1,1}}]$ it holds that \begin{equation}
    \realization{\width,\theta}_\infty(x)=\frac{16x}{9(2\q{\theta}_{\m^{\theta}_{1,1}}-1)}-\frac{4(2\q{\theta}_{\m^{\theta}_{1,1}}-3)}{9(2\q{\theta}_{\m^{\theta}_{1,1}}-1)}.
\end{equation} 
This establishes that $\realization{\width,\theta}_\infty(\q{\theta}_{\m^{\theta}_{1,1}})=\nicefrac43$ which is a contradiction.
\end{cproof}
\begin{lemma}\label{cor:H:x0}
Assume \cref{H:setting:snn} and let $\width \in \N$, $\theta\in (\cG^\width)^{-1}(\{0\})$ satisfy $\prod_{k=1}^\width \v{\theta}_k\neq0$ and 
\begin{equation}\label{cor:H:x0:hp}
   \int_{0}^{1}  x( \realization{\width,\theta}_\infty (x) - \mathbbm{1}_{(\nicefrac12,\infty)}(x)) \, \d x =0.
\end{equation}
Then \begin{equation}\label{cor:H:x0:thesis}
    \cL^\width_\infty(\theta)\geq\nicefrac1{36}.
\end{equation}
\end{lemma}
\begin{cproof}{cor:H:x0}
Assume without loss of generality that $\q{\theta}_1\leq\ldots\leq\q{\theta}_\width$.
In the following we distinguish between the case $\M^{\theta}_{0}=\emptyset = \M^{\theta}_{1}$, the case $\M^{\theta}_{0}\neq\emptyset $, and the case $\M^{\theta}_{0}=\emptyset \neq \M^{\theta}_{1}$.
We first prove \cref{cor:H:x0:thesis} in the case 
\begin{equation}\label{cor:H:x0:case1}
    \M^{\theta}_{0}=\emptyset = \M^{\theta}_{1}.
\end{equation}
\Nobs that \cref{cor:H:x0:case1} implies that there exist $a,b \in \R$ wich satisfy for all $x \in [0,1]$ that $\realization{\width,\theta}_\infty(x)=a x +b$. This and
\cref{lemma:line} establish that 
\begin{equation}\label{cor:H:x0:1}
    \cL^\width_\infty(\theta)\geq\frac{1}{16}.
\end{equation}
This establishes \cref{cor:H:x0:thesis} in the case $\M^{\theta}_{0}=\emptyset = \M^{\theta}_{1}$. Next we prove \cref{cor:H:x0:thesis} in the case
\begin{equation}\label{cor:H:x0:case2}
    \M^{\theta}_{0}\neq\emptyset.
\end{equation}
\Nobs that the assumption that $\theta\in (\cG^\width)^{-1}(\{0\})$, the assumption that $\prod_{k=1}^\width \v{\theta}_k\neq0$, \cref{cor:H:x0:hp}, \cref{cor:H:x0:case2}, and \cref{prop:H:approximate:gradient} demonstrate that
\begin{equation}
     \int_{0}^{\q{\theta}_{\m^{\theta}_{0,1}}}  ( \realization{\width,\theta}_\infty (x) - \mathbbm{1}_{(\nicefrac12,\infty)}(x)) \, \d x =\int_{0}^{\q{\theta}_{\m^{\theta}_{0,1}}}  x( \realization{\width,\theta}_\infty (x) - \mathbbm{1}_{(\nicefrac12,\infty)}(x)) \, \d x =0.
\end{equation}
Combining this with \cref{lem:affine:integral:zero} and \cref{cor:H:m1} proves that for all $x\in [0,\q{\theta}_{\m^{\theta}_{0,2}}]$ it holds that $\realization{\theta}_\infty(x)=0$.
This and \cref{lemma:H:N0} imply that $\M^{\theta}_{1}=\emptyset$. \Hence that there exists $a \in \R$ which satisfy for all $x\in [\q{\theta}_{\m^{\theta}_{0,2}},1]$ that $\realization{\width,\theta}_\infty(x)=a(x-\q{\theta}_{\m^{\theta}_{0,2}})$.
This and \cref{lemma:half} establish that 
\begin{equation}\label{cor:H:x0:2}
    \cL^\width_\infty(\theta)\geq\frac1{36}.
\end{equation}
This demonstrates \cref{cor:H:x0:thesis} in the case $\M^{\theta}_{0}\neq\emptyset$. Next we prove \cref{cor:H:x0:thesis} in the case \begin{equation}\label{cor:H:x0:case3}
    \M^{\theta}_{0}=\emptyset \neq \M^{\theta}_{1}.
\end{equation}
\Nobs that the assumption that $\theta\in (\cG^\width)^{-1}(\{0\})$, the assumption that $\prod_{k=1}^\width \v{\theta}_k\neq0$, \cref{cor:H:x0:hp}, \cref{cor:H:x0:case3}, and \cref{prop:H:approximate:gradient} show that
\begin{equation}
     \int_{\q{\theta}_{\m^{\theta}_{1,2}}}^1  ( \realization{\width,\theta}_\infty (x) - \mathbbm{1}_{(\nicefrac12,\infty)}(x)) \, \d x =\int_{\q{\theta}_{\m^{\theta}_{1,2}}}^1  x( \realization{\width,\theta}_\infty (x) - \mathbbm{1}_{(\nicefrac12,\infty)}(x)) \, \d x =0.
\end{equation}
Combining this with \cref{lem:affine:integral:zero} and \cref{cor:H:m1} proves that for all $x\in [\q{\theta}_{\m^{\theta}_{1,1}},1]$ it holds that $\realization{\width,\theta}_\infty(x)=1$.
\Hence that there exists $a \in \R$ which satisfies for all $x\in [0,\q{\theta}_{\m^{\theta}_{1,1}}]$ that $\realization{\width,\theta}_\infty(x)=a(x-\q{\theta}_{\m^{\theta}_{1,1}})+1$.
This and \cref{lemma:halfdx} establish that 
\begin{equation}
    \cL^\width_\infty(\theta)\geq\frac1{36}.
\end{equation}
This demonstrates \cref{cor:H:x0:thesis} in the case $\M^{\theta}_{0}=\emptyset \neq \M^{\theta}_{1}.$
\end{cproof}
\begin{lemma}\label{lemma:H:updown:sx}
Assume \cref{H:setting:snn}, let $\width\in \N$, $\theta\in (\cG^\width)^{-1}(\{0\})$ satisfy for all $j\in \{0,1,\ldots,\width\}$, $i \in \{1,2,\ldots, \m^{\theta}_{0,2} \}$ that
$\q{\theta}_j < \q{\theta}_{j+1}$, $\w{\theta}_i \w{\theta}_{i-1}<0\neq\v{\theta}_i$, $\M^{\theta}_{0}\neq\emptyset$, and $\alpha^{\theta} \in \R\backslash\{0\}$. 
Then 
\begin{enumerate} [label=(\roman*)]
    \item \label{lemma:H:updown:sx:item1} 
    for all $j \in \{1,2,\ldots, \m^{\theta}_{0,2}\}$ it holds that $\q{\theta}_j=j\q{\theta}_1$,
    \item \label{lemma:H:updown:sx:item2} 
    for all $j \in \{0,1,\ldots, \m^{\theta}_{0,2}\}$ it holds that \begin{equation}
    - \frac{ \alpha^{\theta}\q{\theta}_1}2=(-1)^{j} \realization{\width,\theta}_\infty(\q{\theta}_j),
    \end{equation} 
    \item \label{lemma:H:updown:sx:item3} for all $j \in \{1,2,\ldots,\m^{\theta}_{0,2}\}$, $x \in [\q{\theta}_{j-1},\q{\theta}_{j}]$
it holds that
\begin{equation}
    \realization{\width,\theta}_\infty(x)=(-1)^{j+1}\alpha^{\theta}x +(-1)^{j}\bigg( j-\frac12\bigg) \alpha^{\theta}\q{\theta}_1,
\end{equation}
    and 
    \item \label{lemma:H:updown:sx:item4}
    it holds that \begin{equation}
\int_0^{\q{\theta}_{\m^{\theta}_{0,1}}}   ( \realization{\width,\theta}_\infty (x) - \mathbbm{1}_{(\nicefrac12,\infty)}(x))^2 \, \d x=\frac{\m^{\theta}_{0,2}}{12}(\alpha^{\theta})^2(\q{\theta}_1)^3.
\end{equation} 
\end{enumerate}
\end{lemma}
\begin{cproof}{lemma:H:updown:sx}
\Nobs that the assumption that $\M^{\theta}_{0}\neq\emptyset$ ensures that $\q{\theta}_1\leq\nicefrac12$.
This, the assumption that $\theta\in (\cG^\width)^{-1}(\{0\})$, the fact that $\w{\theta}_0 \w{\theta}_{1}<0\neq\v{\theta}_1$,
and \cref{prop:H:approximate:gradient} assure that 
\begin{equation}
     0=\int_0^{\q{\theta}_{1}}   ( \realization{\width,\theta}_\infty (x) - \mathbbm{1}_{(\nicefrac12,\infty)}(x)) \, \d x=  \int_0^{\q{\theta}_{1}}    ( \realization{\width,\theta}_\infty (x)) \, \d x 
     =\frac{\alpha^{\theta}}{2}(\q{\theta}_{1})^2+\realization{\width,\theta}_\infty(0)\q{\theta}_{1}.
\end{equation}
This demonstrates that for all $x \in [0,\q{\theta}_1]$ it holds that
\begin{equation}\label{lemma:H:updown:sx:initcase1}
    \realization{\width,\theta}_\infty(x)=\alpha^{\theta}x
    -\frac{ \alpha^{\theta}\q{\theta}_1}2.
\end{equation}
\Hence that 
\begin{equation}\label{lemma:H:updown:sx:qnote}
    \realization{\width,\theta}_\infty(0)=-\frac{ \alpha^{\theta}\q{\theta}_1}2=-\realization{\width,\theta}_\infty(\q{\theta}_1)
\end{equation}
and 
\begin{equation}\label{lemma:H:updown:sx:initcase3}
    \int_0^{\q{\theta}_1}   ( \realization{\width,\theta}_\infty (x) - \mathbbm{1}_{(\nicefrac12,\infty)}(x))^2 \, \d x
    = \frac{1}{12}(\alpha^{\theta})^2(\q{\theta}_1)^3.
\end{equation}
\Nobs that \cref{lemma:H:updown:sx:initcase1}, \cref{lemma:H:updown:sx:qnote}, and \cref{lemma:H:updown:sx:initcase3} establish items \ref{lemma:H:updown:sx:item1},  \ref{lemma:H:updown:sx:item2},  \ref{lemma:H:updown:sx:item3}, and  \ref{lemma:H:updown:sx:item4}
in the case  $\m^{\theta}_{0,2}=1$.
Assume now that 
\begin{equation}\label{lemma:H:updown:sx:case2}
    \m^{\theta}_{0,2}>1.
\end{equation}
\Nobs that the assumption that $\theta\in (\cG^\width)^{-1}(\{0\})$, the assumption that for all $i \in \{1,2,\ldots, \m^{\theta}_{0,2} \}$ it holds that $\w{\theta}_i \w{\theta}_{i-1}<0\neq\v{\theta}_i$, \cref{lemma:H:updown:sx:case2}, and  \cref{prop:H:approximate:gradient} imply that for all $i \in \{1,2,\ldots,\m^{\theta}_{0,2}\}$ it holds that
\begin{equation}\label{lemma:H:updown:sx:zerogradient}
      \int_{\q{\theta}_{i-1}}^{\q{\theta}_{i}}   ( \realization{\width,\theta}_\infty (x) - \mathbbm{1}_{(\nicefrac12,\infty)}(x)) \, \d x=0.
\end{equation}
This, \cref{corollary:integralzero}, and \cref{lemma:H:updown:sx:qnote} show that for all $i \in \{0,1,\ldots, \m^{\theta}_{0,2}\}$ it holds that 
\begin{equation} \label{lemma:H:updown:sx:real0}
    - \frac{ \alpha^{\theta}\q{\theta}_1}2=(-1)^{i} \realization{\width,\theta}_\infty(\q{\theta}_j).
\end{equation}
This establishes that for all $i \in \{1,2,\ldots,\m^{\theta}_{0,2} \}$, $x \in [\q{\theta}_{i-1},\q{\theta}_{i}]$ it holds that 
\begin{equation}\label{lemma:real2}
    \realization{\width,\theta}_\infty(x)=\frac{(-1)^{i+1}\alpha^{\theta}\q{\theta}_1 x}{\q{\theta}_{i}-\q{\theta}_{i-1}} +(-1)^{i}\bigg( \frac{\alpha^{\theta}\q{\theta}_1}{\q{\theta}_{i}-\q{\theta}_{i-1}}\q{\theta}_{i-1}+\frac{ \alpha^{\theta}\q{\theta}_1}2\bigg).
\end{equation}
\Moreover the assumption that $\theta\in (\cG^\width)^{-1}(\{0\})$, the assumption that for all $i \in \{1,2,\ldots, \m^{\theta}_{0,2} \}$ it holds that $\w{\theta}_i \w{\theta}_{i-1}<0\neq\v{\theta}_i$, \cref{lemma:H:updown:sx:case2}, and  \cref{prop:H:approximate:gradient} imply that for all $i \in \{1,2,\ldots,\m^{\theta}_{0,2}-1\}$ it holds that
\begin{equation}\label{lemma:H:updown:sx:zerogradient2}
      \int_{\q{\theta}_{i-1}}^{\q{\theta}_{i+1}}  x ( \realization{\width,\theta}_\infty (x) - \mathbbm{1}_{(\nicefrac12,\infty)}(x)) \, \d x=0.
\end{equation}
Combining this and \cref{lemma:H:updown:sx:zerogradient} with \cref{lemma:doubleintegral:zero} demonstrates that for all $i \in \{1,2,\ldots,\m^{\theta}_{0,2}-1\}$ it holds that $ \q{\theta}_i-\q{\theta}_{i-1}= \q{\theta}_{i+1}-\q{\theta}_{i}$. This and the fact that $\q{\theta}_{0}=0$
show that for all $i \in \{1,2,\ldots,\m^{\theta}_{0,2}\}$ it holds that
\begin{equation}\label{lemma:H:updown:sx:real1}
   \q{\theta}_i=i \q{\theta}_1.
\end{equation}
Combining this with \cref{lemma:real2} assures that for all $i \in \{1,2,\ldots,\m^{\theta}_{0,2}\}$, $x \in [\q{\theta}_{i-1},\q{\theta}_{i}]$
it holds that
\begin{equation}
    \realization{\width,\theta}_\infty(x)=(-1)^{i+1}\alpha^{\theta}x +(-1)^{i}\bigg( i-\frac12\bigg) \alpha^{\theta}\q{\theta}_1.
\end{equation}
This ensures that
\begin{equation}\label{lemma:H:updown:sx:real2}
\begin{split}
     \int_0^{\q{\theta}_{\m^{\theta}_{0,2}}}   ( \realization{\width,\theta}_\infty (x) - \mathbbm{1}_{(\nicefrac12,\infty)}(x))^2 \, \d x&=
    \m^{\theta}_{0,2} \int_0^{\q{\theta}_1}   ( \realization{\width,\theta}_\infty (x) - \mathbbm{1}_{(\nicefrac12,\infty)}(x))^2 \, \d x\\
    &= \frac{\m^{\theta}_{0,2}}{12}(\alpha^{\theta})^2(\q{\theta}_1)^3.
\end{split}
\end{equation}
This, \cref{lemma:H:updown:sx:real0}, \cref{lemma:H:updown:sx:real1}, and \cref{lemma:H:updown:sx:real2} prove items \ref{lemma:H:updown:sx:item1},  \ref{lemma:H:updown:sx:item2},  \ref{lemma:H:updown:sx:item3}, and  \ref{lemma:H:updown:sx:item4} in the case $\m^{\theta}_{0,2}>1$.
\end{cproof}
\begin{lemma}\label{lemma:H:updown:dx}
Assume \cref{H:setting:snn}, let $\width\in \N$, $\theta\in (\cG^\width)^{-1}(\{0\})$ satisfy for all $j\in \{0,1,\ldots,\width\}$, $i \in \{\m^{\theta}_{1,1},\m^{\theta}_{1,1}+1,\ldots, \width \}$ that
$\q{\theta}_j < \q{\theta}_{j+1}$, $\w{\theta}_i \w{\theta}_{i+1}<0\neq \v{\theta}_i$, $\M^{\theta}_{1}\neq\emptyset$, and $\beta^{\theta} \in \R\backslash\{0\} $.
Then 
\begin{enumerate} [label=(\roman*)]
    \item \label{lemma:H:updown:dx:item1}
    for all $j \in \{\m^{\theta}_{1,1},\m^{\theta}_{1,1}+1,\ldots, \width\}$ it holds that $\q{\theta}_j=1-(\width+1-j)( 1-\q{\theta}_\width)$,
    \item \label{lemma:H:updown:dx:item2} for all $j \in \{\m^{\theta}_{1,1},\m^{\theta}_{1,1}+1,\ldots, \width+1\}$ it holds that \begin{equation}
    \frac{\beta^{\theta}}{2}(1-\q{\theta}_\width)=(-1)^{ \width +1-j} (\realization{\width,\theta}_\infty(\q{\theta}_j)-1),
    \end{equation} 
    \item \label{lemma:H:updown:dx:item3} for all $j \in \{\m^{\theta}_{1,1},\m^{\theta}_{1,1}+1,\ldots, \width\}$, $x \in [\q{\theta}_{j},\q{\theta}_{j+1}]$ it holds that
    \begin{equation}
    \realization{\width,\theta}_\infty(x)=(-1)^{\width-j}\beta^{\theta}x+1
    +(-1)^{\width+1-j}\beta^{\theta}\left(1-\left(\width+\frac12-j\right)(1-\q{\theta}_\width)\right),
    \end{equation}
    and 
    \item \label{lemma:H:updown:dx:item4}
    it holds that \begin{equation}
\int_{\q{\theta}_{\m^{\theta}_{1,1}}}^{1}   ( \realization{\width,\theta}_\infty (x) - \mathbbm{1}_{(\nicefrac12,\infty)}(x))^2 \, \d x=\frac{1}{12}(\width+1-\m^{\theta}_{1,1})(\beta^{\theta})^2(1-\q{\theta}_\width)^3.
\end{equation}
\end{enumerate}
\end{lemma}
\begin{cproof}{lemma:H:updown:dx}
\Nobs that the assumption that $\M^{\theta}_{1}\neq \emptyset$ ensures that $\q{\theta}_\width\geq\nicefrac12$.
This, the assumption that $\theta\in (\cG^\width)^{-1}(\{0\})$, 
the assumption that $\w{\theta}_\width \w{\theta}_{\width+1}<0\neq\v{\theta}_1$,
and \cref{prop:H:approximate:gradient} assure that 
\begin{equation}
\begin{split}
     0&=\int_{\q{\theta}_{\width}}^{1}   ( \realization{\width,\theta}_\infty (x) - \mathbbm{1}_{(\nicefrac12,\infty)}(x)) \, \d x=  \int_{\q{\theta}_{\width}}^{1}   ( \realization{\width,\theta}_\infty (x) - 1) \, \d x \\ 
     &=\frac{\beta^{\theta}}{2}(1-(\q{\theta}_{\width})^2)+(\realization{\width,\theta}_\infty(1)-1-\beta^{\theta})(1-\q{\theta}_{\width}).
\end{split}
\end{equation}
This demonstrates that for all $x \in [\q{\theta}_{\width},1]$ it holds that
\begin{equation}\label{lemma:H:updown:dx:initcase1}
    \realization{\width,\theta}_\infty(x)=\beta^{\theta}x+1
    -\frac{ \beta^{\theta}}2 (1+\q{\theta}_\width).
\end{equation}
\Hence that 
\begin{equation}\label{lemma:H:updown:dx:qnote}
    \frac{\beta^{\theta}}{2}(1-\q{\theta}_\width)=\realization{\width,\theta}_\infty(1)-1=-\realization{\width,\theta}_\infty(\q{\theta}_\width)+1
\end{equation}
and 
\begin{equation}\label{lemma:H:updown:dx:initcase3}
    \int_{\q{\theta}_{\width}}^{1}   ( \realization{\width,\theta}_\infty (x) - \mathbbm{1}_{(\nicefrac12,\infty)}(x))^2 \, \d x
    = \frac{1}{12}(\beta^{\theta})^2(1-\q{\theta}_\width)^3.
\end{equation}
\Nobs that \cref{lemma:H:updown:dx:initcase1}, \cref{lemma:H:updown:dx:qnote}, and \cref{lemma:H:updown:dx:initcase3} establish items \ref{lemma:H:updown:dx:item1},  \ref{lemma:H:updown:dx:item2},  \ref{lemma:H:updown:dx:item3}, and  \ref{lemma:H:updown:dx:item4}
in the case  $\m^{\theta}_{1,1}=\width$.
Assume now \begin{equation}\label{lemma:H:updown:dx:case2}
    \m^{\theta}_{1,1}<\width.
\end{equation}
\Nobs that the assumption that $\theta\in (\cG^\width)^{-1}(\{0\})$, the assumption that for all $i \in \{\m^{\theta}_{1,1},\m^{\theta}_{1,1}+1,\ldots, \width \}$ it holds that $\w{\theta}_i \w{\theta}_{i+1}<0\neq \v{\theta}_i$, \cref{lemma:H:updown:dx:case2}, and  \cref{prop:H:approximate:gradient} imply that for all $i \in \{\m^{\theta}_{1,1},\m^{\theta}_{1,1}+1,\ldots, \width \}$ it holds that
\begin{equation}\label{lemma:H:updown:dx:zerogradient}
      \int_{\q{\theta}_{i}}^{\q{\theta}_{i+1}}   ( \realization{\width,\theta}_\infty (x) - \mathbbm{1}_{(\nicefrac12,\infty)}(x)) \, \d x=0.
\end{equation}
This, \cref{corollary:integralzero}, and \cref{lemma:H:updown:dx:qnote} show that for all $i \in \{\m^{\theta}_{1,1},\m^{\theta}_{1,1}+1,\ldots, \width+1 \}$ it holds that 
\begin{equation} \label{lemma:H:updown:dx:real0}
    \frac{\beta^{\theta}}{2}(1-\q{\theta}_\width)=\realization{\width,\theta}_\infty(1)-1=(-1)^{ \width +1-i} (\realization{\width,\theta}_\infty(\q{\theta}_i)-1).
\end{equation}
This establishes that for all $i \in \{\m^{\theta}_{1,1},\m^{\theta}_{1,1}+1,\ldots, \width \}$, $x \in [\q{\theta}_i,\q{\theta}_{i+1}]$ it holds that 
\begin{equation}\label{lemma:real:dx}
    \realization{\width,\theta}_\infty(x)=\frac{(-1)^{\width-i}\beta^{\theta}(1-\q{\theta}_\width) x}{\q{\theta}_{i+1}-\q{\theta}_i }+1
    +(-1)^{\width+1-i}\bigg(\frac{\beta^{\theta}(1-\q{\theta}_\width)}{\q{\theta}_{i+1}-\q{\theta}_i}\q{\theta}_i +
    \frac{\beta^{\theta}}{2}(1-\q{\theta}_\width)\bigg).
\end{equation}
\Moreover the assumption that $\theta\in (\cG^\width)^{-1}(\{0\})$, the assumption that for all $i \in \{\m^{\theta}_{1,1},\m^{\theta}_{1,1}+1,\ldots, \width \}$ it holds that $\w{\theta}_i \w{\theta}_{i+1}<0\neq \v{\theta}_i$, \cref{lemma:H:updown:dx:case2}, and  \cref{prop:H:approximate:gradient} imply that 
for all $i \in \{\m^{\theta}_{1,1},\m^{\theta}_{1,1}+1,\ldots, \width -1\}$ it holds that
\begin{equation}
          \int_{\q{\theta}_{i}}^{\q{\theta}_{i+2}}  x ( \realization{\width,\theta}_\infty (x) - \mathbbm{1}_{(\nicefrac12,\infty)}(x)) \, \d x=0.
\end{equation}
Combining this and \cref{lemma:H:updown:dx:zerogradient} with \cref{lemma:doubleintegral:zero} demonstrates that for all $i \in \{\m^{\theta}_{1,1}+1,\m^{\theta}_{1,1}+2,\ldots,\width\}$ it holds that
$\q{\theta}_i-\q{\theta}_{i-1}= \q{\theta}_{i+1}-\q{\theta}_{i}$.
This and the fact that $\q{\theta}_{\width+1}=1$
show that for all $i \in \{\m^{\theta}_{1,1},\m^{\theta}_{1,1}+1,\ldots,\width\}$ it holds that
\begin{equation}\label{lemma:H:updown:dx:real1}
    \q{\theta}_i=1-(\width+1-i) (1-\q{\theta}_\width).
\end{equation}
 Combining this with \cref{lemma:real:dx} assures that for all $i \in \{\m^{\theta}_{1,1},\m^{\theta}_{1,1}+1,\ldots,\width\}$, $x \in [\q{\theta}_{i},\q{\theta}_{i+1}]$ it holds that
\begin{equation}\label{lemma:H:updown:dx:real2}
    \realization{\width,\theta}_\infty(x)=(-1)^{\width-i}\beta^{\theta}x+1
    +(-1)^{\width+1-i}\beta^{\theta}\left(1-\left(\width+\frac12-i\right)(1-\q{\theta}_\width)\right).
\end{equation}
This ensures that
\begin{equation}
\begin{split}
     \int_{\q{\theta}_{\m^{\theta}_{1,1}}}^{1}   ( \realization{\width,\theta}_\infty (x) - \mathbbm{1}_{(\nicefrac12,\infty)}(x))^2 \, \d x&=
    (\width+1-\m^{\theta}_{1,1})\int_{\q{\theta}_{\width}}^{1}   ( \realization{\width,\theta}_\infty (x) - \mathbbm{1}_{(\nicefrac12,\infty)}(x))^2 \, \d x\\
    &= \frac{1}{12}(\width+1-\m^{\theta}_{1,1})(\beta^{\theta})^2(1-\q{\theta}_\width)^3.
\end{split}
\end{equation}
This, \cref{lemma:H:updown:dx:real0}, \cref{lemma:H:updown:dx:real1}, and \cref{lemma:H:updown:dx:real2} prove items \ref{lemma:H:updown:dx:item1}, \ref{lemma:H:updown:dx:item2}, \ref{lemma:H:updown:dx:item3}, and \ref{lemma:H:updown:dx:item4} in the case $\m^{\theta}_{1,1}<\width$.
\end{cproof}
\begin{lemma}\label{lemma:H:updownsx}
Assume \cref{H:setting:snn}, let $\width \in \N $, $\theta\in (\cG^\width)^{-1}(\{0\})$ satisfy for all $j \in \{0,1,\ldots,\width\}$, $i \in \{1,2,\ldots,\m^{\theta}_{0,2}\}$, $x \in [\q{\theta}_{i-1},\q{\theta}_{i}]$ that 
$\q{\theta}_j <  \q{\theta}_{j+1}$,  $\w{\theta}_i\w{\theta}_{i-1}<0\neq\prod_{k=1}^\width\v{\theta}_k$, $\M^{\theta}_{0}\neq\emptyset$, and
    \begin{equation}\label{lemma:H:updownsx:real}
     \realization{\width,\theta}_\infty(x)=(-1)^{i+1}\alpha^{\theta}x +(-1)^{i}\bigg( i-\frac12\bigg) \alpha^{\theta}\q{\theta}_1.
    \end{equation}
Then it holds for all $j\in \{1,2,\ldots, \m^{\theta}_{0,2}\}\backslash\{\m^{\theta}_{0,2}\}$ that \begin{equation}
        \w{\theta}_j\v{\theta}_j=- 2\alpha^{\theta}.
    \end{equation}
\end{lemma}
\begin{cproof}{lemma:H:updownsx}
\Nobs that the assumption that for all $i \in \{1,2,\ldots,\m^{\theta}_{0,2}\}$ it holds that $\w{\theta}_i\w{\theta}_{i-1}<0$ and the fact that $\w{\theta}_0=-1$ implies that for all $j \in \{1,2,\ldots,\m^{\theta}_{0,2}\}\backslash\{\m^{\theta}_{0,2}\}$ with 
$\{n\in \N\colon j= 2n-1\}\neq\emptyset$
 it holds that
\begin{equation}
 I^{\theta}_j=(\q{\theta}_j,1] \qandq  \w{\theta}_j>0.
\end{equation}
This and \cref{lemma:H:updownsx:real} show that for all  $j \in \{1,2,\ldots,\m^{\theta}_{0,2}\}\backslash\{\m^{\theta}_{0,2}\}$ with 
$\{n\in \N\colon j= 2n-1\}\neq\emptyset$ it holds that
 \begin{equation}\label{lemma:H:updownsx:case1}
        \w{\theta}_j\v{\theta}_j=- 2\alpha^{\theta}.
    \end{equation}
\Moreover the assumption that for all $i \in \{1,2,\ldots,\m^{\theta}_{0,2}\}$ it holds that $\w{\theta}_i\w{\theta}_{i-1}<0$ and the fact that $\w{\theta}_0=-1$ ensures that for all $j \in \{1,2,\ldots,\m^{\theta}_{0,2}\}\backslash\{\m^{\theta}_{0,2}\}$ with 
$\{n\in \N\colon j= 2n\}\neq\emptyset$ it holds that
\begin{equation}
   I^{\theta}_j=[0,\q{\theta}_j) \qandq \w{\theta}_j<0.
\end{equation}
This and \cref{lemma:H:updownsx:real} show that for all  $j \in \{1,2,\ldots,\m^{\theta}_{0,2}\}\backslash\{\m^{\theta}_{0,2}\}$  with 
$\{n\in \N\colon j= 2n\}\neq\emptyset$ it holds that
 \begin{equation}
        \w{\theta}_j\v{\theta}_j=- 2\alpha^{\theta}.
    \end{equation}
Combining this and \cref{lemma:H:updownsx:case1} demonstrates that for all $j \in \{1,2,\ldots,\m^{\theta}_{0,2}\}\backslash\{\m^{\theta}_{0,2}\}$ it holds that
 \begin{equation}
        \w{\theta}_j\v{\theta}_j=- 2\alpha^{\theta}.
    \end{equation}
\end{cproof}
\begin{lemma}\label{lemma:H:updowndx}
Assume \cref{H:setting:snn}, let $\width \in \N $, $\theta\in (\cG^\width)^{-1}(\{0\})$ satisfy for all $j \in \{0,1,\ldots,\width\}$, $i \in \{\m^{\theta}_{1,1},\m^{\theta}_{1,1}+1,\ldots,\width\}$, $x \in [\q{\theta}_{i},\q{\theta}_{i+1}]$ that  
$\q{\theta}_j <  \q{\theta}_{j+1}$,  $\w{\theta}_i\w{\theta}_{i+1}<0\neq\prod_{k=1}^\width\v{\theta}_k$, $\M^{\theta}_{1}\neq     \emptyset$, and
    \begin{equation}\label{lemma:H:updowndx:real}
   \realization{\width,\theta}_\infty(x)=(-1)^{\width-i}\beta^{\theta}x+1
    +(-1)^{\width+1-i}\beta^{\theta}\left(1-\left(\width+\frac12-i\right)(1-\q{\theta}_\width)\right).
    \end{equation}
Then it holds
for all $j\in \{\m^{\theta}_{1,1},\m^{\theta}_{1,1}+1,\ldots, \width\}\backslash\{\m^{\theta}_{1,1}\}$ that \begin{equation}
        \w{\theta}_j\v{\theta}_j=- 2\beta^{\theta}.
    \end{equation}
\end{lemma}
\begin{cproof}{lemma:H:updowndx}
\Nobs that the assumption that for all $i \in \{\m^{\theta}_{1,1},\m^{\theta}_{1,1}+1,\ldots,\width\}$ it holds that $\w{\theta}_i\w{\theta}_{i+1}<0$ and the fact that $\w{\theta}_{\width+1}=1$ implies that for all $j \in \{\m^{\theta}_{1,1},\m^{\theta}_{1,1}+1,\ldots,\width\}\backslash\{\m^{\theta}_{1,1}\}$  
with $\{n\in \N\colon j= \width-2n+1\}\neq\emptyset$ it holds that
\begin{equation}
  I^{\theta}_j=(\q{\theta}_j,1] \qandq  \w{\theta}_j>0.
\end{equation}
This and \cref{lemma:H:updowndx:real} show that for all $j \in \{\m^{\theta}_{1,1},\m^{\theta}_{1,1}+1,\ldots,\width\}\backslash\{\m^{\theta}_{1,1}\}$  with $\{n\in \N\colon j= \width-2n+1\}\neq\emptyset$ it holds that
 \begin{equation}\label{lemma:H:updowndx:case1}
        \w{\theta}_j\v{\theta}_j=- 2\beta^{\theta}.
    \end{equation}
\Moreover the assumption that for all $i \in \{\m^{\theta}_{1,1},\m^{\theta}_{1,1}+1,\ldots,\width\}$ it holds that $\w{\theta}_i\w{\theta}_{i+1}<0$ and the fact that $\w{\theta}_{\width+1}=-1$ ensures that for all $j \in \{\m^{\theta}_{1,1},\m^{\theta}_{1,1}+1,\ldots,\width\}\backslash\{\m^{\theta}_{1,1}\}$  with $\{n\in \N\colon j= \width-2n+2\}\neq\emptyset$ it holds that
\begin{equation}
   I^{\theta}_j=[0,\q{\theta}_j) \qandq \w{\theta}_j<0.
\end{equation}
This and \cref{lemma:H:updowndx:real} establish that for all $j \in \{\m^{\theta}_{1,1},\m^{\theta}_{1,1}+1,\ldots,\width\}\backslash\{\m^{\theta}_{1,1}\}$ with $\{n\in \N\colon j= \width-2n+2\}\neq\emptyset$ it holds that
 \begin{equation}
        \w{\theta}_j\v{\theta}_j=- 2\beta^{\theta}.
    \end{equation}
Combining this and \cref{lemma:H:updowndx:case1} demonstrates that for all $j \in \{\m^{\theta}_{1,1},\m^{\theta}_{1,1}+1,\ldots,\width\}\backslash\{\m^{\theta}_{1,1}\}$ it holds that
 \begin{equation}
        \w{\theta}_j\v{\theta}_j=- 2\beta^{\theta}.
    \end{equation}
\end{cproof}
\begin{lemma}\label{lemma:H:N0:init}
Assume \cref{H:setting:snn}, let $\width \in \N \cap(1,\infty)$, $\theta\in (\cG^\width)^{-1}(\{0\})$ satisfy for all $j \in \{0,1,\ldots,\width\}$, $x \in [0,\q{\theta}_{\m^{\theta}_{0,2}}]$ that
$\q{\theta}_j <  \q{\theta}_{j+1}$, $\realization{\width,\theta}_\infty(x)=0$, $\M^{\theta}_{0}\neq\emptyset\neq\M^{\theta}_{1}$, and $0\neq\prod_{k=1}^\width\v{\theta}_k$.
Then 
\begin{equation}\label{lemma:H:N0:init:thesis}
   \cL^\width_\infty(\theta)\geq\frac{1}{384(1+\width)^2}.
\end{equation}
\end{lemma}
\begin{cproof}{lemma:H:N0:init}
\Nobs that the assumption that  $M^{\theta}_{0}\neq\emptyset\neq\M^{\theta}_{1}$, the assumption that for all  $x \in [0,\q{\theta}_{\m^{\theta}_{0,2}}]$ it holds that $\realization{\width,\theta}_\infty(x)=0$,  \cref{cor:H:m1}, \cref{cor:H:w0sdx}, and  \cref{lemma:H:N0} assure that for all $j \in \{\m^{\theta}_{1,1},\m^{\theta}_{1,1}+1,\ldots,\width\}$ it holds that $\w{\theta}_j\w{\theta}_{j+1}<0$ and $\beta^{\theta}\neq0$. This and \cref{lemma:H:updown:dx} demonstrate that 
\begin{enumerate}[label=(\roman*)]
\item for all $j \in \{\m^{\theta}_{1,1},\m^{\theta}_{1,1}+1,\ldots,\width\}$ it holds that $\q{\theta}_j=1-(\width+1-j) (1-\q{\theta}_\width)$,
\item
for all $j \in \{\m^{\theta}_{1,1},\m^{\theta}_{1,1}+1,\ldots, \width+1\}$ it holds that 
\begin{equation}\label{lemma:H:N0:init:max}
  \frac{\beta^{\theta}}{2}(1-\q{\theta}_\width)=(-1)^{ \width +1-j} (\realization{\width,\theta}_\infty(\q{\theta}_j)-1),
\end{equation}
\item for all $j \in \{\m^{\theta}_{1,1},\m^{\theta}_{1,1}+1,\ldots, \width\}$, $x \in [\q{\theta}_{j},\q{\theta}_{j+1}]$ it holds that
    \begin{equation}\label{lemma:H:N0:init:real}
    \realization{\width,\theta}_\infty(x)=(-1)^{\width-j}\beta^{\theta}x+1 +(-1)^{\width+1-j}\beta^{\theta}\left(1-\left(\width+\frac12-j\right)(1-\q{\theta}_\width)\right),
    \end{equation}
and 
\item
it holds that
\begin{equation}\label{lemma:H:N0:init:risk}
 \int_{\q{\theta}_{\m^{\theta}_{1,1}}}^{1}   ( \realization{\width,\theta}_\infty (x) - \mathbbm{1}_{(\nicefrac12,\infty)}(x))^2 \, \d x=\frac{1}{12}(\width+1-\m^{\theta}_{1,1})(\beta^{\theta})^2(1-\q{\theta}_\width)^3.
\end{equation}
\end{enumerate}
\Nobs that \cref{lemma:H:N0:init:max}  assures that in the case $\q{\theta}_{\m^{\theta}_{0,2}}=\q{\theta}_{\m^{\theta}_{1,1}}=\nicefrac12$ it holds that
\begin{equation} \label{lemma:H:N0:init:hpbeta}
    1=|\realization{\width,\theta}_\infty(\q{\theta}_{\m^{\theta}_{0,2}})-1|=|\realization{\width,\theta}_\infty(\q{\theta}_{\m^{\theta}_{1,1}})-1|=|\nicefrac12\beta^{\theta}(1-\q{\theta}_\width)|.
\end{equation} 
\Moreover in the case $\q{\theta}_{\m^{\theta}_{0,2}}=\q{\theta}_{\m^{\theta}_{1,1}}=\nicefrac12$ it holds that $(\width+1-\m^{\theta}_{1,1})(1-\q{\theta}_\width)=\nicefrac12$.
Combining this and \cref{lemma:H:N0:init:hpbeta} with \cref{lemma:H:N0:init:risk} implies that in the case $\q{\theta}_{\m^{\theta}_{0,2}}=\q{\theta}_{\m^{\theta}_{1,1}}=\nicefrac12$ it holds that 
\begin{equation}\label{lemma:H:N0:init:l1m2}
    \cL^\width_\infty(\theta)=\frac16.
\end{equation}
This establishes \cref{lemma:H:N0:init:thesis} in the case $\q{\theta}_{\m^{\theta}_{0,2}}=\q{\theta}_{\m^{\theta}_{1,1}}=\nicefrac12$.
Assume now that $\q{\theta}_{\m^{\theta}_{0,2}}<\nicefrac12<\q{\theta}_{\m^{\theta}_{1,2}}$.
In the following we distinguish between the case $\w{\theta}_{\m^{\theta}_{0,2}}<0< \w{\theta}_{\m^{\theta}_{1,1}}$ and  the case $\max\{\w{\theta}_{\m^{\theta}_{0,2}}\w{\theta}_{\m^{\theta}_{1,1}},-\w{\theta}_{\m^{\theta}_{1,1}}\}>0$.
We first prove \cref{lemma:H:N0:init:thesis} in the case
\begin{equation}\label{lemma:H:N0:init:case1}
    \max\{\w{\theta}_{\m^{\theta}_{0,2}}\w{\theta}_{\m^{\theta}_{1,1}},-\w{\theta}_{\m^{\theta}_{1,1}}\}>0.
\end{equation}
\Nobs that \cref{lemma:H:N0:init:case1}, the assumption that $\theta\in (\cG^\width)^{-1}(\{0\})$, and \cref{prop:H:approximate:gradient} assure that 
\begin{equation}\label{lemma:H:N0:init:m2l1}
      \int_{\q{\theta}_{\m^{\theta}_{0,2}}}^{\q{\theta}_{\m^{\theta}_{1,1}}}   x( \realization{\width,\theta}_\infty (x) - \mathbbm{1}_{(\nicefrac12,\infty)}(x)) \, \d x=
       \int_{\q{\theta}_{\m^{\theta}_{0,2}}}^{\q{\theta}_{\m^{\theta}_{1,1}}}   ( \realization{\width,\theta}_\infty (x) - \mathbbm{1}_{(\nicefrac12,\infty)}(x)) \, \d x
      =0.
\end{equation}
Combining this with \cref{cor:affine1:integral:zero} demonstrates that for all $x \in [\q{\theta}_{\m^{\theta}_{0,2}},\q{\theta}_{\m^{\theta}_{1,1}}]$ it holds that \begin{equation}
    \realization{\width,\theta}_\infty(x)= \frac{16 x}{9(2\q{\theta}_{\m^{\theta}_{1,1}}-1)}  +\frac{4(2\q{\theta}_{\m^{\theta}_{1,1}}-3)}{9(2\q{\theta}_{\m^{\theta}_{1,1}}-1)}.
\end{equation}
\Hence that $ \realization{\width,\theta}_\infty(\q{\theta}_{\m^{\theta}_{1,1}})=\nicefrac43$ and $\q{\theta}_{\m^{\theta}_{0,2}}=\nicefrac34-\nicefrac{\q{\theta}_{\m^{\theta}_{1,1}}}2$.
Combining this with \cref{lemma:H:N0:init:max}  and \cref{lemma:H:N0:init:risk} shows that $|\nicefrac{\beta^{\theta}}2(1-\q{\theta}_\width)|=|\realization{\width,\theta}_\infty(\q{\theta}_{\m^{\theta}_{1,1}})-1|=\nicefrac13$ and
\begin{equation}\label{lemma:H:N0:init:prove1}
\begin{split}
     \cL^\width_\infty(\theta)&=\int_{\q{\theta}_{\m^{\theta}_{0,2}}}^{\q{\theta}_{\m^{\theta}_{1,1}}}  \bigg(\frac{16x}{9(2\q{\theta}_{\m^{\theta}_{1,1}}-1)}  +\frac{4(2\q{\theta}_{\m^{\theta}_{1,1}}-3)}{9(2\q{\theta}_{\m^{\theta}_{1,1}}-1)} - \mathbbm{1}_{(\nicefrac12,\infty)}(x)\bigg)^2 \, \d x\\
     & \quad +\frac{\width+1-\m^{\theta}_{1,1}}{12} (\beta^{\theta})^2(1-\q{\theta}_\width)^3= \frac7{36}(2\q{\theta}_{\m^{\theta}_{1,1}}-1)+\frac1{27} (1-\q{\theta}_{\m^{\theta}_{1,1}}) \\
&= \frac{19 \q{\theta}_{\m^{\theta}_{1,1}}}{54} -\frac{17}{108}\geq \frac{19}{108} -\frac{17}{108}=\frac1{54}.
\end{split}
\end{equation}
This implies \cref{lemma:H:N0:init:thesis} in the case $\max\{\w{\theta}_{\m^{\theta}_{0,2}}\w{\theta}_{\m^{\theta}_{1,1}},-\w{\theta}_{\m^{\theta}_{1,1}}\}>0$.
Next we demonstrate \cref{lemma:H:N0:init:thesis} in the case
\begin{equation}\label{lemma:H:N0:init:case2}
    \w{\theta}_{\m^{\theta}_{0,2}}<0< \w{\theta}_{\m^{\theta}_{1,1}}.
\end{equation}
\Nobs that \cref{lemma:H:N0:init:case2} and the fact that for all $j \in \{\m^{\theta}_{1,1},\m^{\theta}_{1,1}+1,\ldots,\width\}$ it holds that $\w{\theta}_j\w{\theta}_{j+1}<0$ prove that there exists $k^* \in \N$ which satisfies that $\width=\m^{\theta}_{1,1}+2k^*-1$.
Combining this with the fact that for all $x \in [0,\q{\theta}_{\m^{\theta}_{0,2}}]$ it holds that $\realization{\width,\theta}_\infty(x)=0$ and the fact that for all $j \in \{\m^{\theta}_{1,1},\m^{\theta}_{1,1}+1,\ldots,\width\}$ it holds that $\w{\theta}_j\w{\theta}_{j+1}<0$ establishes that
\begin{equation}\label{lemma:H:N0:init:active}
    \{k\in\{1,2,\ldots,\width\}\colon I^{\theta}_k\cup (\q{\theta}_{\m^{\theta}_{0,2}},\q{\theta}_{\m^{\theta}_{1,1}})\neq \emptyset\}=
\cup_{k=1}^{k^*}\{\m^{\theta}_{1,1}+2k-1\}.
\end{equation}
\Nobs that \cref{lemma:H:updowndx} ensures that for all $j \in \cup_{k=1}^{k^*}\{\m^{\theta}_{1,1}+2k-1\}$ it holds that $\w{\theta}_j\v{\theta}_j=-2 \beta^{\theta}$.
Combining this with \cref{lemma:H:N0:init:active} and the fact that $\realization{\width,\theta}_\infty(\q{\theta}_{\m^{\theta}_{0,2}})=0$ assures that for all $x \in [\q{\theta}_{\m^{\theta}_{0,2}},\q{\theta}_{\m^{\theta}_{1,1}}]$
 it holds that 
\begin{equation}\label{lemma:H:N0:init:real2}
    \realization{\width,\theta}_\infty(x)=-2k^*\beta^{\theta}x  +2k^*\beta^{\theta} \q{\theta}_{\m^{\theta}_{0,2}}.
\end{equation}
This and \cref{lemma:H:N0:init:max} demonstrate that 
\begin{equation}
    -2k^*\beta^{\theta}\q{\theta}_{\m^{\theta}_{1,1}}  +2k^*\beta^{\theta} \q{\theta}_{\m^{\theta}_{0,2}}=1+\frac{\beta^{\theta}}2(1-\q{\theta}_\width).
\end{equation}
\Hence that
\begin{equation}
    |\beta^{\theta}|=\frac{1}{|\frac12(1-\q{\theta}_\width)+2k^*(\q{\theta}_{\m^{\theta}_{1,1}} -\q{\theta}_{\m^{\theta}_{0,2}} )|}\geq \frac{1}{1+2k^*}.
\end{equation}
Combining this, \cref{lemma:half}, \cref{lemma:H:N0:init:risk}, and \cref{lemma:H:N0:init:real2} shows that 
\begin{equation}
\begin{split}
     \cL^{\width}_\infty(\theta) &\geq 
     \begin{cases}
     \int_{\q{\theta}_{\m^{\theta}_{0,2}}}^{\q{\theta}_{\m^{\theta}_{1,1}}}  (-2k^*\beta^{\theta}x  +2k^*\beta^{\theta} \q{\theta}_{\m^{\theta}_{0,2}} - \mathbbm{1}_{(\nicefrac12,\infty)}(x))^2 \, \d x
      & \colon \q{\theta}_{\m^{\theta}_{1,1}}\geq\frac34 \\
    \int_{\q{\theta}_{\m^{\theta}_{1,1}}}^{1}   ( \realization{\width,\theta}_\infty (x) - \mathbbm{1}_{(\nicefrac12,\infty)}(x))^2 \, \d x= \frac{2k^*}{12} (\beta^{\theta})^2(1-\q{\theta}_\width)^3 &\colon \q{\theta}_{\m^{\theta}_{1,1}}\leq\frac34
     \end{cases}
     \\ 
     & \geq
    \begin{cases}
    \frac{1}{36} & \colon \q{\theta}_{\m^{\theta}_{1,1}}\geq\frac34\\
    \frac{k^*}{6} \frac{1}{(1+2k^*)^2}\frac1{64} \geq
    \frac{1}{384(1+\width)^2} &\colon \q{\theta}_{\m^{\theta}_{1,1}}\leq\frac34.
     \end{cases}
     \end{split}
\end{equation}
This, \cref{lemma:H:N0:init:l1m2}, and \cref{lemma:H:N0:init:prove1} assure that
\begin{equation}
    \cL^{\width}_\infty(\theta)\geq\frac{1}{384(1+\width)^2}.
\end{equation}
This shows \cref{lemma:H:N0:init:thesis} in the case $    \w{\theta}_{\m^{\theta}_{0,2}}<0< \w{\theta}_{\m^{\theta}_{1,1}}$.
\end{cproof}
\begin{lemma}\label{lemma:H:N0:fin}
Assume \cref{H:setting:snn}, let $\width \in \N\cap(1,\infty)$, $\theta\in (\cG^\width)^{-1}(\{0\})$ satisfy for all $j \in \{0,1,\ldots,\width\}$, $x \in [\q{\theta}_{\m^{\theta}_{1,2}},1]$ that
$\q{\theta}_j <  \q{\theta}_{j+1}$, $\realization{\width,\theta}_\infty(x)=1$, $M^{\theta}_{0}\neq\emptyset\neq\M^{\theta}_{1}$, and $0\neq \prod_{k=1}^\width\v{\theta}_k$. Then 
\begin{equation}\label{lemma:H:N0:fin:thesis}
    \cL^{\width}_\infty(\theta)\geq\frac{1}{384(1+\width)^2}.
\end{equation}
\end{lemma}
\begin{cproof}{lemma:H:N0:fin}
\Nobs that the assumption that $\theta\in (\cG^\width)^{-1}(\{0\})$, the assumption that  $M^{\theta}_{0}\neq\emptyset\neq\M^{\theta}_{1}$, the assumption that for all  $x \in [\q{\theta}_{\m^{\theta}_{1,1}},1]$ it holds that $\realization{\theta}_\infty(x)=1$,  \cref{cor:H:l2}, \cref{cor:H:w0sx}, and  \cref{lemma:H:N0} assure that for all $j \in \{1,2,\ldots,\m^{\theta}_{0,2}\}$ it holds that $\w{\theta}_j\w{\theta}_{j-1}<0$ and $\alpha^{\theta}\neq0$. This and \cref{lemma:H:updown:sx} demonstrate that 
\begin{enumerate}[label=(\roman*)]
\item 
for all $j \in \{1,2,\ldots,\m^{\theta}_{0,2}\}$ it holds that $\q{\theta}_j=j\q{\theta}_1$,
\item
for all $j \in \{0,1,\ldots, \m^{\theta}_{0,2}\}$ it holds that
\begin{equation}\label{lemma:H:N0:fin:max}
   - \frac{ \alpha^{\theta}\q{\theta}_1}2=(-1)^{j} \realization{\width,\theta}_\infty(\q{\theta}_j),
\end{equation}
\item
for all $j \in \{1,2,\ldots,\m^{\theta}_{0,2}\}$, $x \in [\q{\theta}_{j-1},\q{\theta}_{j}]$ it holds that
    \begin{equation}\label{lemma:H:N0:fin:real}
     \realization{\width,\theta}_\infty(x)=(-1)^{j+1}\alpha^{\theta}x +(-1)^{j}\bigg( j-\frac12\bigg) \alpha^{\theta}\q{\theta}_1,
    \end{equation}
and
\item it holds that
\begin{equation}\label{lemma:H:N0:fin:risk}
\int_0^{\q{\theta}_{\m^{\theta}_{0,1}}}   ( \realization{\width,\theta}_\infty (x) - \mathbbm{1}_{(\nicefrac12,\infty)}(x))^2 \, \d x=\frac{\m^{\theta}_{0,2}}{12}(\alpha^{\theta})^2(\q{\theta}_1)^3.
\end{equation}
\end{enumerate}
\Nobs that \cref{lemma:H:N0:fin:max} assures that in the case $\q{\theta}_{\m^{\theta}_{0,2}}=\q{\theta}_{\m^{\theta}_{1,1}}=\nicefrac12$ it holds that
\begin{equation}\label{lemma:H:N0:fin:hpbeta}
    1=|\realization{\width,\theta}_\infty(\q{\theta}_{\m^{\theta}_{1,1}})|=|\realization{\width,\theta}_\infty(\q{\theta}_{\m^{\theta}_{0,2}})|=\left|\frac{\alpha^{\theta}\q{\theta}_1}2\right|.
\end{equation} 
\Moreover in the case $\q{\theta}_{\m^{\theta}_{0,2}}=\q{\theta}_{\m^{\theta}_{1,1}}=\nicefrac12$ it holds that $\m^{\theta}_{0,2} \q{\theta}_1=\nicefrac12$.
Combining this and \cref{lemma:H:N0:fin:hpbeta} with \cref{lemma:H:N0:fin:risk} implies that in the case $\q{\theta}_{\m^{\theta}_{0,2}}=\q{\theta}_{\m^{\theta}_{1,1}}=\nicefrac12$
it holds that \begin{equation}\label{lemma:H:N0:fin:l1m2}
    \cL^{\width}_\infty(\theta)=\frac16.
\end{equation}
This demonstrates \cref{lemma:H:N0:fin:thesis} in the case $\q{\theta}_{\m^{\theta}_{0,2}}=\q{\theta}_{\m^{\theta}_{1,1}}=\nicefrac12$.
Assume now that $\q{\theta}_{\m^{\theta}_{0,2}}<\nicefrac12<\q{\theta}_{\m^{\theta}_{1,2}}$.
In the following we distinguish between the case $\w{\theta}_{\m^{\theta}_{0,2}}<0< \w{\theta}_{\m^{\theta}_{1,1}}$ and the case $\max\{\w{\theta}_{\m^{\theta}_{0,2}}\w{\theta}_{\m^{\theta}_{1,1}},\w{\theta}_{\m^{\theta}_{0,2}}\}>0$.
We first demonstrate \cref{lemma:H:N0:fin:thesis} in the case
\begin{equation}\label{lemma:H:N0:fin:case1}
    \max\{\w{\theta}_{\m^{\theta}_{0,2}}\w{\theta}_{\m^{\theta}_{1,1}},\w{\theta}_{\m^{\theta}_{0,2}}\}>0.
\end{equation}
\Nobs that \cref{lemma:H:N0:fin:case1}, the assumption that $\theta\in (\cG^\width)^{-1}(\{0\})$, and \cref{prop:H:approximate:gradient} assure that 
\begin{equation}\label{lemma:H:N0:fin:m2l1}
      \int_{\q{\theta}_{\m^{\theta}_{0,2}}}^{\q{\theta}_{\m^{\theta}_{1,1}}}   x( \realization{\width,\theta}_\infty (x) - \mathbbm{1}_{(\nicefrac12,\infty)}(x)) \, \d x=
       \int_{\q{\theta}_{\m^{\theta}_{0,2}}}^{\q{\theta}_{\m^{\theta}_{1,1}}}   ( \realization{\width,\theta}_\infty (x) - \mathbbm{1}_{(\nicefrac12,\infty)}(x)) \, \d x
      =0.
\end{equation}
Combining this with \cref{cor:affine1:integral:zero2} demonstrates that for all $x \in [\q{\theta}_{\m^{\theta}_{0,2}},\q{\theta}_{\m^{\theta}_{1,1}}]$ it holds that \begin{equation}
    \realization{\width,\theta}_\infty(x)= -\frac{16 x}{9(2\q{\theta}_{\m^{\theta}_{0,2}}-1)}  +\frac{10\q{\theta}_{\m^{\theta}_{0,2}}+3}{9(2\q{\theta}_{\m^{\theta}_{0,2}}-1)}.
\end{equation}
\Hence that $ \realization{\width,\theta}_\infty(\q{\theta}_{\m^{\theta}_{0,2}})=-\nicefrac13$ and $\q{\theta}_{\m^{\theta}_{1,1}}=\nicefrac34-\nicefrac12\q{\theta}_{\m^{\theta}_{0,2}}$.
Combining this with \cref{lemma:H:N0:fin:max}  and \cref{lemma:H:N0:fin:risk} shows that $|\nicefrac12\alpha^{\theta}\q{\theta}_1|=|\realization{\width,\theta}_\infty(\q{\theta}_{\m^{\theta}_{0,2}})|=\nicefrac13$ and
\begin{equation}\label{lemma:H:N0:fin:prove1}
\begin{split}
     \cL^{\width}_\infty(\theta)&=\int_{\q{\theta}_{\m^{\theta}_{0,2}}}^{\q{\theta}_{\m^{\theta}_{1,1}}}  \bigg(-\frac{16 x}{9(2\q{\theta}_{\m^{\theta}_{0,2}}-1)} +\frac{10\q{\theta}_{\m^{\theta}_{0,2}}+3}{9(2\q{\theta}_{\m^{\theta}_{0,2}}-1)}- \mathbbm{1}_{(\nicefrac12,\infty)}(x)\bigg)^2 \, \d x \\ &\quad +\frac{\m^{\theta}_{0,2}}{12}(\alpha^{\theta})^2(\q{\theta}_1)^3
     =\frac1{18}(1-2\q{\theta}_{\m^{\theta}_{0,2}})+\frac1{27} \q{\theta}_{\m^{\theta}_{0,2}} \\
&= -\frac{2 \q{\theta}_{\m^{\theta}_{0,2}}}{27}+ \frac1{18}\geq -\frac{1}{27}+\frac1{18}=\frac1{54}.
\end{split}
\end{equation}
This implies \cref{lemma:H:N0:fin:thesis} in the case $ \max\{\w{\theta}_{\m^{\theta}_{0,2}}\w{\theta}_{\m^{\theta}_{1,1}},\w{\theta}_{\m^{\theta}_{0,2}}\}>0.$
Next we demonstrate \cref{lemma:H:N0:fin:thesis} in the case \begin{equation}\label{lemma:H:N0:fin:case2}
    \w{\theta}_{\m^{\theta}_{0,2}}<0< \w{\theta}_{\m^{\theta}_{1,1}}.
\end{equation}
\Nobs that \cref{lemma:H:N0:fin:case2} and the fact that for all $i \in \{0,1,\ldots,\m^{\theta}_{0,2}-1\}$ it holds that $\w{\theta}_i\w{\theta}_{i+1}<0$ prove that there exists $k^* \in \N$ which satisfies that $\m^{\theta}_{0,2}=2k^*$.
Combining this with the fact that for all $x \in [\q{\theta}_{\m^{\theta}_{1,1}},1]$ it holds that $\realization{\width,\theta}_\infty(x)=1$ and the fact that for all $j \in \{0,1,\ldots,\m^{\theta}_{0,2}-1\}$ it holds that $\w{\theta}_j\w{\theta}_{j+1}<0$ establishes that
\begin{equation}\label{lemma:H:N0:fin:active}
    \{k\in\{1,2,\ldots,\width\}\colon I^{\theta}_k\cap (\q{\theta}_{\m^{\theta}_{0,2}},\q{\theta}_{\m^{\theta}_{1,1}})\neq \emptyset\}= \cup_{k =1}^{k^*}\{2k-1\}.
\end{equation}
\Nobs that \cref{lemma:H:updowndx}  ensures that for all $j \in\cup_{k =1}^{k^*}\{2k-1\}$ it holds that $\w{\theta}_j\v{\theta}_j=-2 \alpha^{\theta}$.
Combining this with \cref{lemma:H:N0:fin:active} and the fact that $\realization{\width,\theta}_\infty(\q{\theta}_{\m^{\theta}_{1,1}})=1$ assures that for all $x \in [\q{\theta}_{\m^{\theta}_{0,2}},\q{\theta}_{\m^{\theta}_{1,1}}]$
 it holds that 
\begin{equation}\label{lemma:H:N0:fin:real2}
    \realization{\width,\theta}_\infty(x)=-2 k^*\alpha^{\theta}x  +2k^*\alpha^{\theta} \q{\theta}_{\m^{\theta}_{1,1}}+1.
\end{equation}
This and \cref{lemma:H:N0:fin:max} demonstrate that 
\begin{equation}
    -2k^*\alpha^{\theta}\q{\theta}_{\m^{\theta}_{0,2}}  +2k^*\alpha^{\theta} \q{\theta}_{\m^{\theta}_{1,1}}+1=-\frac{\alpha^{\theta}\q{\theta}_1}2.
\end{equation}
\Hence that 
\begin{equation}
    |\alpha^{\theta}|=\frac{1}{|\frac{\q{\theta}_1}2+2k^*(\q{\theta}_{\m^{\theta}_{1,1}} -\q{\theta}_{\m^{\theta}_{0,2}} )|}\geq \frac{1}{1+2k^*}.
\end{equation}
Combining this, \cref{lemma:halfdx}, \cref{lemma:H:N0:fin:risk}, and \cref{lemma:H:N0:fin:real2} shows that 
\begin{equation}
\begin{split}
     \cL^{\width}_\infty(\theta) &\geq 
     \begin{cases}
     \int_{\q{\theta}_{\m^{\theta}_{0,2}}}^{\q{\theta}_{\m^{\theta}_{1,1}}}  (-2 k^*\alpha^{\theta}x  +2k^*\alpha^{\theta} \q{\theta}_{\m^{\theta}_{1,1}}+1 - \mathbbm{1}_{(\nicefrac12,\infty)}(x))^2 \, \d x
      & \colon \q{\theta}_{\m^{\theta}_{0,2}}\leq\frac14 \\
    \int_{\q{\theta}_{\m^{\theta}_{1,1}}}^{1}   ( \realization{\width,\theta}_\infty (x) - \mathbbm{1}_{(\nicefrac12,\infty)}(x))^2 \, \d x= \frac{2k^*}{12} (\alpha^{\theta})^2(\q{\theta}_1)^3 &\colon \q{\theta}_{\m^{\theta}_{0,2}}\geq\frac14
     \end{cases}
     \\ 
     & \geq
    \begin{cases}
    \frac{1}{36} & \colon \q{\theta}_{\m^{\theta}_{0,2}}\leq\frac14\\
    \frac{k^*}{6} \frac{1}{(1+2k^*)^2}\frac1{64} \geq
    \frac{1}{384(1+\width)^2} &\colon \q{\theta}_{\m^{\theta}_{0,2}} \geq\frac14.
     \end{cases}
     \end{split}
\end{equation}
This, \cref{lemma:H:N0:fin:l1m2}, and \cref{lemma:H:N0:fin:prove1} assure that
\begin{equation}
    \cL^{\width}_\infty(\theta)\geq\frac{1}{384(1+\width)^2}.
\end{equation}
This shows \cref{lemma:H:N0:fin:thesis} in the case $\w{\theta}_{\m^{\theta}_{0,2}}<0< \w{\theta}_{\m^{\theta}_{1,1}}.$
\end{cproof}
\begin{cor}\label{lemma:H:updown:leqm}
Assume \cref{H:setting:snn} and let $\width\in \N\cap(1,\infty)$, $\theta\in (\cG^\width)^{-1}(\{0\})$ satisfy  for all $i \in \{0,\ldots,\width\}$ that 
$\w{\theta}_i\w{\theta}_{i+1}<0\neq\alpha^{\theta}\beta^{\theta}$,
$\prod_{k=1}^\width\v{\theta}_k\neq0<\q{\theta}_1  <\q{\theta}_2 < \ldots < \q{\theta}_\width<1$, and $\q{\theta}_{\m^{\theta}_{0,2}}=\q{\theta}_{\m^{\theta}_{1,1}}=\nicefrac12$.
 Then $\cL^{\width}_\infty(\theta)\geq\nicefrac1{12}$.
\end{cor}
\begin{cproof}{lemma:H:updown:leqm}
\Nobs that 
\cref{lemma:H:updown:sx}, \cref{lemma:H:updown:dx}, and the assumption that $\q{\theta}_{\m^{\theta}_{0,2}}=\q{\theta}_{\m^{\theta}_{1,1}}=\nicefrac12$ assure that $|1+\nicefrac12(-1)^{\width+1-\m^{\theta}_{1,1}}\beta^{\theta}(1-\q{\theta}_\width)|=|\realization{\width,\theta}_\infty (\q{\theta}_{\m^{\theta}_{1,1}})|=|\realization{\width,\theta}_\infty (\q{\theta}_{\m^{\theta}_{0,2}})|=|\nicefrac{\alpha^{\theta} \q{\theta}_1}2 |$ and
\begin{equation}\label{lemma:H:updown:ml:int}
     \int_0^1 ( \realization{\width,\theta}_\infty (x) - \mathbbm{1}_{(\nicefrac12,\infty)}(x))^2 \, \d x
     =\frac{1}{24}(\alpha^{\theta})^2(\q{\theta}_1)^2+\frac{1}{24}(\beta^{\theta})^2(1-\q{\theta}_\width)^2.
\end{equation}
\Hence that
\begin{equation}
\begin{split}
     \cL^{\width}_\infty(\theta)&=\frac{1}{6}\left(1+\frac12(-1)^{\width+1-\m^{\theta}_{1,1}}\beta^{\theta}(1-\q{\theta}_\width)\right)^2 +\frac{1}{24}(\beta^{\theta})^2(1-\q{\theta}_\width)^2\\
     &=\frac1{12}(\beta^{\theta}(1-\q{\theta}_\width))^2+\frac16(-1)^{\width+1-\m^{\theta}_{1,1}} \beta^{\theta}(1-\q{\theta}_\width) + \frac16 \\
     & \geq \frac1{12}(-1)^2+\frac16(-1) + \frac16=\frac1{12}.
\end{split}
\end{equation}
This demonstrates that $\cL^{\width}_\infty(\theta)\geq \nicefrac1{12}$.
\end{cproof}
\begin{lemma}\label{lemma:H:updown:weq}
Assume \cref{H:setting:snn} and let $\width\in \N \cap(1,\infty)$, $\theta\in \R^{\fd_\width}$ satisfy 
for all $i \in \{0,1,\ldots,\width\}\allowbreak\backslash \{\m^{\theta}_{0,2}\}$ that $\w{\theta}_i\w{\theta}_{i+1}<0$,
$\prod_{k=1}^\width\v{\theta}_k\neq0<\q{\theta}_1  < \q{\theta}_2  <\ldots < \q{\theta}_\width<1$, $\alpha^{\theta}\beta^{\theta}\neq 0<\m^{\theta}_{0,2}<\m^{\theta}_{1,1}<\width+1$, and $0<\w{\theta}_{\m^{\theta}_{0,2}}\w{\theta}_{\m^{\theta}_{1,1}}$.
 Then \begin{equation}\label{lemma:H:updown:weq:thesis}
     \cG^\width(\theta)\neq0.
 \end{equation}
\end{lemma}
\begin{cproof}{lemma:H:updown:weq}
We prove \cref{lemma:H:updown:weq:thesis} by contradiction. Assume that 
\begin{equation}\label{lemma:H:updown:weq:thesisab}
    \cG^\width(\theta)=0.
\end{equation}
\Nobs that \cref{lemma:H:updown:weq:thesisab},
\cref{lemma:H:updown:sx}, and \cref{lemma:H:updown:dx} ensure that
 \begin{enumerate}[label=(\roman*)]
    \item for all $j \in \{1,2\ldots, \m^{\theta}_{0,2}\}$ it holds that \begin{equation}\label{lemma:H:updown:weq:alpha}
    \q{\theta}_j=j\q{\theta}_1,
    \end{equation} 
     \item  for all $j \in \{\m^{\theta}_{1,1},\m^{\theta}_{1,1}+1,\ldots, \width\}$ it holds that
     \begin{equation}\label{lemma:H:updown:weq:beta}
    \q{\theta}_j=1-(\width+1-j)( 1-\q{\theta}_\width),
    \end{equation} 
    and
   \item 
   for all $j \in \{0,1,\ldots,\width\}\backslash\{\m^{\theta}_{0,2}\}$, $x \in [\q{\theta}_{j},\q{\theta}_{j+1}]$
   it holds that
   \begin{equation}\label{lemma:H:updown:weq:real}
   \begin{split}
     &\realization{\width,\theta}_\infty(x)\\
     & =    \begin{cases}
    (-1)^{j}\alpha^{\theta}x +(-1)^{j+1}( j+\nicefrac12) \alpha^{\theta}\q{\theta}_1 & \colon x\leq \q{\theta}_{\m^{\theta}_{0,2}} \\
    (-1)^{\width-j}\beta^{\theta}x+1
    +(-1)^{\width+1-j}\beta^{\theta}(1-(\width+\nicefrac12-j) (1-\q{\theta}_\width)) & \colon x\geq \q{\theta}_{\m^{\theta}_{1,1}}.
    \end{cases}
   \end{split}
\end{equation}
\end{enumerate}
This, \cref{lemma:H:updownsx}, and \cref{lemma:H:updowndx} assure that for all $j \in \{1,2,\ldots,\width\}\backslash\{\m^{\theta}_{0,2},\m^{\theta}_{1,1}\}$ it holds that
\begin{equation}\label{lemma:H:updown:weq:wv}
    \w{\theta}_{j}\v{\theta}_{j}=
    \begin{cases}
    -2\alpha^{\theta} & \colon j < \m^{\theta}_{0,2}\\
    -2\beta^{\theta} & \colon j > \m^{\theta}_{1,1}.
    \end{cases}
\end{equation}
In the following we distinguish between the case $\max\{\w{\theta}_{\m^{\theta}_{0,2}},\w{\theta}_{\m^{\theta}_{1,1}}\} <0$ and the case $\min\{\w{\theta}_{\m^{\theta}_{0,2}},\allowbreak \w{\theta}_{\m^{\theta}_{1,1}}\}>0$. 
We first establish the contradiction in the case 
\begin{equation}\label{lemma:H:updown:weq:case1}
    \max\{\w{\theta}_{\m^{\theta}_{0,2}},\w{\theta}_{\m^{\theta}_{1,1}}\}<0.
\end{equation}  
\Nobs that \cref{lemma:H:updown:weq:case1} and the assumption that
for all $i \in \{0,1,\ldots,\width\}\backslash\{\m^{\theta}_{0,2}\}$ it holds that 
 $\w{\theta}_i\w{\theta}_{i+1}<0$ 
 prove that there exist $k_1,k_2 \in \N$ which satisfy that $ \m^{\theta}_{0,2}=2k_1$, $\width-\m^{\theta}_{1,1}=2(k_2-1)$, and
\begin{equation}
    \{k\in\{1,2,\ldots,\width\}\colon I^{\theta}_k\cap (\q{\theta}_{\m^{\theta}_{1,1}},\q{\theta}_{\m^{\theta}_{1,1}+1})\neq \emptyset\}=
(\cup_{k=0}^{\nicefrac12(\width-1)}\{1+2k\})\backslash\{\m^{\theta}_{1,1}\}.
\end{equation}
This, \cref{lemma:H:updown:weq:wv}, and the fact that for all $x \in [\q{\theta}_{\m^{\theta}_{1,1}},\q{\theta}_{\m^{\theta}_{1,1}+1}]$ it holds that
 $\realization{\width,\theta}_\infty(x)=
    \beta^{\theta}x+1
    -\beta^{\theta}\left(1-(\width+\nicefrac12-\m^{\theta}_{1,1}) (1-\q{\theta}_\width)\right) $
show that
\begin{equation}
   \beta^{\theta}=-\m^{\theta}_{0,2} \alpha^{\theta}-(\width-\m^{\theta}_{1,1})\beta^{\theta}.
\end{equation}
Combining this, the fact that $\max\{\w{\theta}_{\m^{\theta}_{0,2}},\w{\theta}_{\m^{\theta}_{1,1}}\}<0<\q{\theta}_{\m^{\theta}_{0,2}}<\q{\theta}_{\m^{\theta}_{1,1}}<1$, \cref{lemma:H:updown:weq:alpha}, \cref{lemma:H:updown:weq:beta}, \cref{lemma:H:updown:weq:real}, the assumption that $\theta\in \cG^{-1}(\{0\})$, and \cref{prop:H:approximate:gradient} assures that 
\begin{equation}
\begin{split}
       0 &= \int_{\q{\theta}_{\m^{\theta}_{0,2}-1}}^{\q{\theta}_{\m^{\theta}_{1,1}+1}}   x( \realization{\width,\theta}_\infty (x) - \mathbbm{1}_{(\nicefrac12,\infty)}(x)) \, \d x \\
       &=
       \int_{\q{\theta}_{\m^{\theta}_{0,2}-1}}^{\q{\theta}_{\m^{\theta}_{0,2}}}   x( \realization{\width,\theta}_\infty (x) - \mathbbm{1}_{(\nicefrac12,\infty)}(x)) \, \d x +
        \int_{\q{\theta}_{\m^{\theta}_{1,1}}}^{\q{\theta}_{\m^{\theta}_{1,1}+1}}   x( \realization{\width,\theta}_\infty (x) - \mathbbm{1}_{(\nicefrac12,\infty)}(x)) \, \d x\\
        &= -\frac{\alpha^{\theta}}{12} (\q{\theta}_1)^3+ \frac{\beta^{\theta}}{12} (1-\q{\theta}_\width)^3
        = -\frac{\alpha^{\theta}}{12}  \bigg((\q{\theta}_1)^3+\frac{\m^{\theta}_{0,2}}{\width-\m^{\theta}_{1,1}+1}(1-\q{\theta}_\width)^3\bigg).
\end{split}
\end{equation}
This implies that \begin{equation}
    \q{\theta}_1=-\sqrt[3]{\frac{\m^{\theta}_{0,2}}{\width-\m^{\theta}_{1,1}+1}}(1-\q{\theta}_\width)<0
\end{equation} which is a contradiction.
In the next step we establish the contradiction in the case
\begin{equation}\label{lemma:H:updown:weq:case2}
     \min\{\w{\theta}_{\m^{\theta}_{0,2}},\w{\theta}_{\m^{\theta}_{1,1}}\}>0.
\end{equation}
\Nobs that \cref{lemma:H:updown:weq:case2} and the assumption that
for all $i \in \{0,1,\ldots,\width\}\backslash\{\m^{\theta}_{0,2}\}$ it holds that 
 $\w{\theta}_i\w{\theta}_{i+1}<0$ 
 prove that there exist $k_1,k_2 \in \N$ which satisfy that $ \m^{\theta}_{0,2}=2k_1-1$, $\width-\m^{\theta}_{1,1}=2k_2-1$,  and
\begin{equation}
    \{k\in\{1,2,\ldots,\width\}\colon I^{\theta}_k\cap (\q{\theta}_{\m^{\theta}_{0,2}-1},\q{\theta}_{\m^{\theta}_{0,2}})\neq \emptyset\}=
(\cup_{k=0}^{\nicefrac12(\width-1)}\{1+2k\})\backslash\{\m^{\theta}_{0,2}\}.
\end{equation} 
This, \cref{lemma:H:updown:weq:wv}, and the fact that for all $x \in [\q{\theta}_{\m^{\theta}_{0,2}-1},\q{\theta}_{\m^{\theta}_{0,2}}]$ it holds that
 $\realization{\width,\theta}_\infty(x)=
    \alpha^{\theta}x -(\m^{\theta}_{0,2}-\nicefrac12) \alpha^{\theta}\q{\theta}_1$
show that
\begin{equation}
   \alpha^{\theta}=-(\m^{\theta}_{0,2}-1) \alpha^{\theta}-(\width-\m^{\theta}_{1,1}+1)\beta^{\theta}.
\end{equation}
Combining this, the fact that $-\min\{\w{\theta}_{\m^{\theta}_{0,2}},\w{\theta}_{\m^{\theta}_{1,1}}\}<0<\q{\theta}_{\m^{\theta}_{0,2}}<\q{\theta}_{\m^{\theta}_{1,1}}<1$, \cref{lemma:H:updown:weq:alpha}, \cref{lemma:H:updown:weq:beta}, \cref{lemma:H:updown:weq:real}, the assumption that $\theta\in \cG^{-1}(\{0\})$, and \cref{prop:H:approximate:gradient} assures that 
\begin{equation}
\begin{split}
       0 &= \int_{\q{\theta}_{\m^{\theta}_{0,2}-1}}^{\q{\theta}_{\m^{\theta}_{1,1}+1}}   x( \realization{\width,\theta}_\infty (x) - \mathbbm{1}_{(\nicefrac12,\infty)}(x)) \, \d x \\
       &=
       \int_{\q{\theta}_{\m^{\theta}_{0,2}-1}}^{\q{\theta}_{\m^{\theta}_{0,2}}}   x( \realization{\width,\theta}_\infty (x) - \mathbbm{1}_{(\nicefrac12,\infty)}(x)) \, \d x +
        \int_{\q{\theta}_{\m^{\theta}_{1,1}}}^{\q{\theta}_{\m^{\theta}_{1,1}+1}}   x( \realization{\width,\theta}_\infty (x) - \mathbbm{1}_{(\nicefrac12,\infty)}(x)) \, \d x\\
        &= \frac{\alpha^{\theta}}{12} (\q{\theta}_1)^3 - \frac{\beta^{\theta} }{12}(1-\q{\theta}_\width)^3
        = \frac{\alpha^{\theta}}{12}  \bigg((\q{\theta}_1)^3+\frac{\m^{\theta}_{0,2}}{\width-\m^{\theta}_{1,1}+1}(1-\q{\theta}_\width)^3\bigg).
\end{split}
\end{equation}
This implies that \begin{equation}
    \q{\theta}_1=-\sqrt[3]{\frac{\m^{\theta}_{0,2}}{\width-\m^{\theta}_{1,1}+1}}(1-\q{\theta}_\width)<0
\end{equation} which is a contradiction.
\end{cproof}
\begin{lemma}\label{lemma:alphabeta}
Assume \cref{H:setting:snn}, let $\width\in \N \cap(1,\infty)$ and let $(\theta_n)_{n \in \N}\subseteq \R^{\fd_\width}$ satisfy for all $n \in \N$ that
 $0<\m^{\theta_n}_{0,2}<\m^{\theta_n}_{1,1}<\width+1$ and
\begin{equation}\label{lemma:alphabeta:assumption}
    (\m^{\theta_n}_{0,2} \alpha^{\theta_n}+(\width-\m^{\theta_n}_{1,1}+1) \beta^{\theta_n})(\q{\theta_n}_{\m^{\theta_n}_{1,1}}-\q{\theta_n}_{\m^{\theta_n}_{0,2}})+ \frac12\alpha^{\theta_n}\q{\theta_n}_1=1-\frac12\beta^{\theta_n}(1-\q{\theta_n}_\width).
\end{equation}
Then there exists $c \in (0,\infty)$ such that for all $n \in \N$ it holds that
\begin{equation}\label{lemma:alphabeta:thesis}
    \max\{ |\alpha^{\theta_n}|, |\beta^{\theta_n}|\}\geq c.
\end{equation} 
\end{lemma}
\begin{cproof}{lemma:alphabeta}
We prove \cref{lemma:alphabeta:thesis} by contradiction. We thus assume that for every $n \in \N$ it holds that \begin{equation}\label{lemma:alphabeta:ab}
    \max\{ |\alpha^{\theta_n}|, |\beta^{\theta_n}|\}\leq \frac1n.
\end{equation} 
\Nobs that \cref{lemma:alphabeta:ab} implies that $\lim_{n\to \infty}\alpha^{\theta_n}=\lim_{n\to \infty}\beta^{\theta_n}=0$. 
Combining this with \cref{lemma:alphabeta:assumption} demonstrates that
\begin{equation}
\begin{split}
     0&=\lim\nolimits_{n\to \infty}(\m^{\theta_n}_{0,2} \alpha^{\theta_n}+(\width-\m^{\theta_n}_{1,1}+1) \beta^{\theta_n})(\q{\theta_n}_{\m^{\theta_n}_{1,1}}-\q{\theta_n}_{\m^{\theta_n}_{0,2}})+ \frac12\alpha^{\theta_n}\q{\theta_n}_1\\
     &=\lim\nolimits_{n\to \infty}1-\frac12\beta^{\theta_n}(1-\q{\theta_n}_\width)=1.
\end{split}
\end{equation}
This is a contradiction.
\end{cproof}
\begin{lemma}\label{lemma:H:updown:ml}
Assume \cref{H:setting:snn} and let $\width \in \N\cap(1,\infty) $.
Then there exists $\varepsilon \in (0,\infty)$ such that for all 
$\theta\in (\cG^\width)^{-1}(\{0\})$ with $ \forall \, i \in \{0,1,\ldots,\width\}\backslash\{\m^{\theta}_{0,2}\} \colon \w{\theta}_i\w{\theta}_{i+1}<0$,
$\prod_{k=1}^\width\v{\theta}_k\neq0<\q{\theta}_1  < \q{\theta}_2  < \ldots < \q{\theta}_\width<1$, 
$\alpha^{\theta}\beta^{\theta}\neq 0<\m^{\theta}_{0,2}<\m^{\theta}_{1,1}<\width+1$, and
$\w{\theta}_{\m^{\theta}_{0,2}}<0<\w{\theta}_{\m^{\theta}_{1,1}}$ it holds that
 \begin{equation}\label{lemma:H:updown:ml:thesis}
     \cL^{\width}_\infty(\theta)\geq \varepsilon.
 \end{equation}
\end{lemma}
\begin{cproof}{lemma:H:updown:ml}
Throughout this proof let $ R\subseteq \R^{\fd_\width} $ satisfy  \begin{equation}
\begin{split}
    R&=\{\theta\in (\cG^\width)^{-1}(\{0\}) \colon [\forall \, i \in \{0,1,\ldots,\width\}\backslash\{\m^{\theta}_{0,2}\} \colon \w{\theta}_i\w{\theta}_{i+1}<0], \w{\theta}_{\m^{\theta}_{0,2}}<0<\w{\theta}_{\m^{\theta}_{1,1}}, \\
     & \quad \textstyle{\prod_{k=1}^\width}\v{\theta}_k\neq0<\q{\theta}_1  < \q{\theta}_2  < \ldots < \q{\theta}_\width<1, 
 \alpha^{\theta}\beta^{\theta}\neq 0<\m^{\theta}_{0,2}<\m^{\theta}_{1,1}<\width+1 \}.
\end{split}
\end{equation} 
\Nobs that 
\cref{lemma:H:updown:sx} and \cref{lemma:H:updown:dx} ensure
for all $\theta \in R$ that
\begin{enumerate}[label=(\roman*)]
    \item for all $j \in \{0,1,\ldots, \m^{\theta}_{0,2}\}$ it holds that \begin{equation}\label{lemma:H:updown:ml:alpha}
    - \frac{ \alpha^{\theta}\q{\theta}_1}2=(-1)^{j} \realization{\width,\theta}_\infty(\q{\theta}_j),
    \end{equation} 
     \item for all $j \in \{\m^{\theta}_{1,1},\m^{\theta}_{1,1}+1,\ldots, \width+1\}$ it holds that \begin{equation}\label{lemma:H:updown:ml:beta}
    \frac{\beta^{\theta}}{2}(1-\q{\theta}_\width)=(-1)^{ \width +1-j} (\realization{\width,\theta}_\infty(\q{\theta}_j)-1),
    \end{equation} 
    \item it holds that \begin{equation}\label{lemma:H:updown:ml:int2}
    \begin{split}
    & \int_{[0,1]\backslash [\q{\theta}_{\m^{\theta}_{0,2}},\q{\theta}_{\m^{\theta}_{1,1}}] }  ( \realization{\width,\theta}_\infty (x) - \mathbbm{1}_{(\nicefrac12,\infty)}(x))^2 \, \d x \\
     & =\frac{\m^{\theta}_{0,2}}{12} (\alpha^{\theta})^2(\q{\theta}_1)^3+\frac{1}{12}(\width+1-\m^{\theta}_{1,1})(\beta^{\theta})^2(1-\q{\theta}_\width)^3,
     \end{split}
     \end{equation}
     and
     \item  for all $j \in \{0,1,\ldots,\width\}\backslash\{\m^{\theta}_{0,2}\}$, $x \in [\q{\theta}_{j},\q{\theta}_{j+1}]$
it holds that
\begin{equation}\label{lemma:H:updown:ml:real}
\begin{split}
&\realization{\width,\theta}_\infty(x)\\
& =    \begin{cases}
(-1)^{j}\alpha^{\theta}x +(-1)^{j+1}( j+\nicefrac12) \alpha^{\theta}\q{\theta}_1 & \colon x\leq \q{\theta}_{\m^{\theta}_{0,2}} \\
   \scalemath{0.97}{ (-1)^{\width-j}\beta^{\theta}x+1
    +(-1)^{\width+1-j}\beta^{\theta}\left(1-(\width+\nicefrac12-j) (1-\q{\theta}_\width)\right)} & \colon x\geq \q{\theta}_{\m^{\theta}_{1,1}}.
    \end{cases}
\end{split}
\end{equation}
\end{enumerate}
This, \cref{lemma:H:updownsx}, and \cref{lemma:H:updowndx} assure that for all $\theta \in R$, $j \in \{1,2,\ldots,\width\}\backslash\{\m^{\theta}_{0,2},\m^{\theta}_{1,1}\}$ it holds that
\begin{equation}\label{lemma:H:updown:ml:wv}
    \w{\theta}_{j}\v{\theta}_{j}=
    \begin{cases}
    -2\alpha^{\theta} & \colon j < \m^{\theta}_{0,2}\\
    -2\beta^{\theta} & \colon j > \m^{\theta}_{1,1}.
    \end{cases}
\end{equation}
\Nobs that the fact that for all $\theta \in R$ it holds that $\w{\theta}_{\m^{\theta}_{0,2}}<0<\w{\theta}_{\m^{\theta}_{1,1}}$ and the fact that
 for all $\theta \in R$, $i \in \{0,1,\ldots,\width\}\backslash\{\m^{\theta}_{0,2}\}$ it holds that
 $\w{\theta}_i\w{\theta}_{i+1}<0$ prove that
 there exist $k_1,k_2 \in \N$ such that for all $\theta \in R$ it holds that $ \m^{\theta}_{0,2}=2k_1$, $\width-\m^{\theta}_{1,1}=2k_2-1$, and
\begin{equation}
    \{k\in\{1,2,\ldots,\width\}\colon I^{\theta}_k\cap (\q{\theta}_{\m^{\theta}_{0,2}},\q{\theta}_{\m^{\theta}_{1,1}})\neq \emptyset\}=
(\cup_{k=0}^{k_1-1}\{2k+1\})\cup (\cup_{k=0}^{k_2-1}\{\m^{\theta}_{1,1}+1+2k\}).
\end{equation}
This, \cref{lemma:H:updown:ml:alpha}, and \cref{lemma:H:updown:ml:wv}
assure that for all $\theta \in R$, $x \in [\q{\theta}_{\m^{\theta}_{0,2}},\q{\theta}_{\m^{\theta}_{1,1}}]$ it holds that
\begin{equation}\label{lemma:H:updown:ml:realization}
   \realization{\width,\theta}_\infty(x)=-(\m^{\theta}_{0,2} \alpha^{\theta}+(\width-\m^{\theta}_{1,1}+1) \beta^{\theta})(x-\q{\theta}_{\m^{\theta}_{0,2}})- \frac{\alpha^{\theta}\q{\theta}_1}2. 
\end{equation}
Combining this, \cref{lemma:H:updown:ml:beta}, and the fact that for all $\theta \in R$ the function $\realization{\width,\theta}_\infty$ is continuous proves for all $\theta \in R$ that
\begin{equation}\label{lemma:H:updown:ml:continuity}
    -(\m^{\theta}_{0,2} \alpha^{\theta}+(\width-\m^{\theta}_{1,1}+1) \beta^{\theta})(\q{\theta}_{\m^{\theta}_{1,1}}-\q{\theta}_{\m^{\theta}_{0,2}})- \frac{\alpha^{\theta}\q{\theta}_1}2=1+\frac{\beta^{\theta}}2(1-\q{\theta}_\width).
\end{equation}
\Moreover \cref{lemma:H:updown:ml:realization} and \cref{lemma:intmin} imply for all $\theta \in R$ that 
\begin{equation}
      \int_{\q{\theta}_{\m^{\theta}_{0,2}}}^{\frac12}   ( \realization{\width,\theta}_\infty (x) - \mathbbm{1}_{(\nicefrac12,\infty)}(x))^2 \, \d x  \geq 
     \frac1{12}\bigg(\frac12-\q{\theta}_{\m^{\theta}_{0,2}}\bigg)^3(\m^{\theta}_{0,2} \alpha^{\theta}+(\width-\m^{\theta}_{1,1}+1) \beta^{\theta})^2
\end{equation}
and
\begin{equation}
     \int_{\frac12}^{\q{\theta}_{\m^{\theta}_{1,1}}}   ( \realization{\width,\theta}_\infty (x) - \mathbbm{1}_{(\nicefrac12,\infty)}(x))^2 \, \d x  \geq 
     \frac1{12}\bigg(\q{\theta}_{\m^{\theta}_{1,1}}-\frac12\bigg)^3(\m^{\theta}_{0,2} \alpha^{\theta}+(\width-\m^{\theta}_{1,1}+1) \beta^{\theta})^2.
\end{equation}
This and \cref{lemma:H:updown:ml:int2} establish for all $\theta \in R$ that
\begin{equation}\label{lemma:H:updown:ml:lowerbound}
\begin{split}
     \cL^{\width}_\infty(\theta) &\geq  \frac1{12}\bigg(\bigg(\q{\theta}_{\m^{\theta}_{1,1}}-\frac12\bigg)^3+\bigg(\frac12-\q{\theta}_{\m^{\theta}_{0,2}}\bigg)^3\bigg)(\m^{\theta}_{0,2} \alpha^{\theta}+(\width-\m^{\theta}_{1,1}+1) \beta^{\theta})^2
     \\
     & \quad +\frac{1}{12}(\beta^{\theta})^2(1-\q{\theta}_{\m^{\theta}_{1,1}})^3
     +\frac{1}{12}(\alpha^{\theta})^2(\q{\theta}_{\m^{\theta}_{0,2}})^3.
\end{split}
\end{equation}
We prove \cref{lemma:H:updown:ml:thesis} by contradiction.
Assume that for every $n \in \N$ there exists $\theta_{n}\in (\cG^\width)^{-1}(\{0\})$ with
$0<\q{\theta_{n}}_1  < \q{\theta_{n}}_2< \ldots < \q{\theta_{n}}_\width<1$, $\alpha^{\theta_{n}}\beta^{\theta_{n}}\neq 0<\m^{\theta_n}_{0,2}<\m^{\theta_n}_{1,1}<\width+1$, 
$\w{\theta_{n}}_{\m^{\theta_n}_{0,2}}<0<\w{\theta_{n}}_{\m^{\theta_n}_{1,1}}$, and
 $\forall \, i \in \{0,1,\ldots,\width\}\backslash\{\m^{\theta_n}_{0,2}\}  \colon \w{\theta_{n}}_i\w{\theta_{n}}_{i+1}<0$
which satisfies that
 \begin{equation}\label{lemma:H:updown:ml:ab}
    \cL^{\width}_\infty(\theta_n)\leq\frac1n. 
\end{equation}
\Nobs that \cref{lemma:H:updown:ml:lowerbound} and \cref{lemma:H:updown:ml:ab} assure that
\begin{equation}
    \begin{split}
    \lim\nolimits_{n\to \infty}& \frac1{12}\bigg(\bigg(\q{\theta_n}_{\m^{\theta_n}_{1,1}}-\frac12\bigg)^3+\bigg(\frac12-\q{\theta_n}_{\m^{\theta_n}_{0,2}}\bigg)^3\bigg)(\m^{\theta_n}_{0,2} \alpha^{\theta_n}+(\width-\m^{\theta_n}_{1,1}+1) \beta^{\theta_n})^2
     \\
     & \quad +\frac{1}{12}(\beta^{\theta_n})^2(1-\q{\theta_n}_{\m^{\theta_n}_{1,1}})^3
     +\frac{1}{12}(\alpha^{\theta_n})^2(\q{\theta_n}_{\m^{\theta_n}_{0,2}})^3=0.
\end{split}
\end{equation}
\Hence that
\begin{equation}\label{lemma:H:updown:ml:limits}
    \begin{split}
    \lim\nolimits_{n\to \infty}& \left(\q{\theta_n}_{\m^{\theta_n}_{1,1}}-\frac12\right)^3(\m^{\theta_n}_{0,2} \alpha^{\theta_n}+(\width-\m^{\theta_n}_{1,1}+1) \beta^{\theta_n})^2=0,\\
    \lim\nolimits_{n\to \infty}& \left(\frac12-\q{\theta_n}_{\m^{\theta_n}_{0,2}}\right)^3(\m^{\theta_n}_{0,2} \alpha^{\theta_n}+(\width-\m^{\theta_n}_{1,1}+1) \beta^{\theta_n})^2=0, \\
      \lim\nolimits_{n\to \infty}& (\beta^{\theta_n})^2(1-\q{\theta_n}_{\m^{\theta_n}_{1,1}})^3=0, \qand\\
       \lim\nolimits_{n\to \infty}& (\alpha^{\theta_n})^2(\q{\theta_n}_{\m^{\theta_n}_{0,2}})^3=0.
\end{split}
\end{equation}
\Nobs that \cref{lemma:H:updown:ml:continuity} and \cref{lemma:alphabeta} demonstrate that there exists $c \in (0,\infty)$ such that for all $n \in \N$ it holds that $\max\{ |\alpha^{\theta_n}|, |\beta^{\theta_n}|\}\geq c$. 
Combining this and \cref{lemma:H:updown:ml:limits} with \cref{cor:H:updown} assure that there exists a strictly increasing $n \colon \N \to \N$ such that  \begin{equation}
    \lim\nolimits_{k\to \infty}\q{\theta_{n(k)}}_{\m^{\theta_{n(k)}}_{1,1}}=1 \qandq \lim\nolimits_{k\to \infty}\q{\theta_{n(k)}}_{\m^{\theta_{n(k)}}_{0,2}}=0.
\end{equation}
This shows that there exists $k^* \in \N$ such that for all $k\in \N \cap[k^*,\infty)$  it holds that $\q{\theta_{n(k)}}_{\m^{\theta}_{1,1}}\geq \nicefrac34$ and $\q{\theta_{n(k)}}_{\m^{\theta}_{0,2}}\leq \nicefrac14$.
Combining this and \cref{lemma:halfc} implies that $\lim\nolimits_{k\to \infty}\cL^{\width}_\infty(\theta_{n(k)})\allowbreak \geq \nicefrac1{32}$ which is a contradiction.
\end{cproof}
\begin{lemma}\label{lemma:H:updown:lm}
Assume \cref{H:setting:snn} and let $\width \in \N\cap(1,\infty)$.
Then there exists $\varepsilon \in (0,\infty)$ such that for all
$\theta\in (\cG^\width)^{-1}(\{0\})$ with $\forall\, i \in \{0,\ldots,\width\}\backslash\{\m^{\theta}_{0,2}\} \colon \w{\theta}_i\w{\theta}_{i+1}<0$, $\prod_{k=1}^\width\v{\theta}_k\neq0<\q{\theta}_1  <\q{\theta}_2  < \ldots < \q{\theta}_\width<1$, $\alpha^{\theta}\beta^{\theta}\neq 0<\m^{\theta}_{0,2}<\m^{\theta}_{1,1}<\width+1$, and
$\w{\theta}_{\m^{\theta}_{1,1}}<0<\w{\theta}_{\m^{\theta}_{0,2}}$ it holds that
  \begin{equation}\label{lemma:H:updown:lm:thesis}
      \cL^{\width}_\infty(\theta)\geq \varepsilon.
 \end{equation}
\end{lemma}
\begin{cproof}{lemma:H:updown:lm}
Throughout this proof let $ R\subseteq \R^{\fd_\width} $ satisfy  \begin{equation}
\begin{split}
    R&=\{\theta\in (\cG^\width)^{-1}(\{0\}) \colon [\forall\, i \in \{0,1,\ldots,\width\}\backslash\{\m^{\theta}_{0,2}\} \colon \w{\theta}_i\w{\theta}_{i+1}<0], \w{\theta}_{\m^{\theta}_{1,1}}<0<\w{\theta}_{\m^{\theta}_{0,2}}, \\
     & \quad  \textstyle{\prod_{k=1}^\width}\v{\theta}_k\neq0<\q{\theta}_1  < \q{\theta}_2  < \ldots < \q{\theta}_\width<1, 
 \alpha^{\theta}\beta^{\theta}\neq 0<\m^{\theta}_{0,2}<\m^{\theta}_{1,1}<\width+1 \}.
\end{split}
\end{equation}
\Nobs that 
\cref{lemma:H:updown:sx} and \cref{lemma:H:updown:dx} ensure for all $\theta  \in R$ that
\begin{enumerate}[label=(\roman*)]
    \item for all $j \in \{0,1,\ldots, \m^{\theta}_{0,2}\}$ it holds that \begin{equation}\label{lemma:H:updown:lm:alpha}
    - \frac{ \alpha^{\theta}\q{\theta}_1}2=(-1)^{j} \realization{\width,\theta}_\infty(\q{\theta}_j),
    \end{equation} 
     \item for all $j \in \{\m^{\theta}_{1,1},\m^{\theta}_{1,1}+1,\ldots, \width+1\}$ it holds that \begin{equation}\label{lemma:H:updown:lm:beta}
    \frac{\beta^{\theta}}{2}(1-\q{\theta}_\width)=(-1)^{ \width +1-j} (\realization{\width,\theta}_\infty(\q{\theta}_j)-1),
    \end{equation} 
    \item it holds that \begin{equation}\label{lemma:H:updown:lm:int}
    \begin{split}
    & \int_{[0,1]\backslash [\q{\theta}_{\m^{\theta}_{0,2}},\q{\theta}_{\m^{\theta}_{1,1}}] }  ( \realization{\width,\theta}_\infty (x) - \mathbbm{1}_{(\nicefrac12,\infty)}(x))^2 \, \d x \\
     & =\frac{ \m^{\theta}_{0,2}}{12}(\alpha^{\theta})^2(\q{\theta}_1)^3+\frac{1}{12}(\width+1-\m^{\theta}_{1,1})(\beta^{\theta})^2(1-\q{\theta}_\width)^3,
     \end{split}
     \end{equation}
     and 
     \item for all $j \in \{0,1,\ldots,\width\}\backslash\{\m^{\theta}_{0,2}\}$, $x \in [\q{\theta}_{j},\q{\theta}_{j+1}]$
it holds that
\begin{equation}\label{lemma:H:updown:lm:real}
\begin{split}
     &\realization{\width,\theta}_\infty(x)\\
     &  =    \begin{cases}
    (-1)^{j}\alpha^{\theta}x +(-1)^{j+1}( j+\nicefrac12) \alpha^{\theta}\q{\theta}_1 & \colon x\leq \q{\theta}_{\m^{\theta}_{0,2}} \\
    \scalemath{0.97}{(-1)^{\width-j}\beta^{\theta}x+1
    +(-1)^{\width+1-j}\beta^{\theta}\left(1-(\width+\nicefrac12-j) (1-\q{\theta}_\width)\right)} & \colon x\geq \q{\theta}_{\m^{\theta}_{1,1}}.
    \end{cases}
\end{split}
\end{equation}
\end{enumerate}
This, \cref{lemma:H:updownsx}, and \cref{lemma:H:updowndx} assure that for all  $\theta  \in R$, $j \in \{1,2,\ldots,\width\}\backslash\{\m^{\theta}_{0,2},\m^{\theta}_{1,1}\}$ it holds that
\begin{equation}\label{lemma:H:updown:lm:wv}
    \w{\theta}_{j}\v{\theta}_{j}=
    \begin{cases}
    -2\alpha^{\theta} & \colon j < \m^{\theta}_{0,2}\\
    -2\beta^{\theta} & \colon j > \m^{\theta}_{1,1}.
    \end{cases}
\end{equation}
\Moreover the fact that for all  $\theta  \in R$ it holds that $\w{\theta}_{\m^{\theta}_{1,1}}<0<\w{\theta}_{\m^{\theta}_{0,2}}$ and the fact that
 for all  $\theta  \in R$, $i \in \{0,1,\ldots,\width\}\backslash\{\m^{\theta}_{0,2}\}$ it holds that 
 $\w{\theta}_i\w{\theta}_{i+1}<0$ prove that
 there exist $k_1,k_2 \in \N$ such that for all $\theta  \in R$ it holds that $ \m^{\theta}_{0,2}=2k_1-1$, $\width-\m^{\theta}_{1,1}=2k_2-2$,
 \begin{equation}\label{lemma:H:updown:lm:set}
 \begin{split}
 \scalemath{0.97}{
 \{k\in\{1,\ldots,\width\}\colon I^{\theta}_k\cap (\q{\theta}_{\m^{\theta}_{0,2}},\q{\theta}_{\m^{\theta}_{1,1}})\neq \emptyset\}}&=(\cup_{k=0}^{k_1-1}\{2k+1\}) \cup (\cup_{k=0}^{k_2-1}\{\m^{\theta}_{1,1}+2k\}), \\
\scalemath{0.97}{\{k\in\{1,\ldots,\width\}\colon I^{\theta}_k\cap (\q{\theta}_{\m^{\theta}_{0,2}-1},\q{\theta}_{\m^{\theta}_{0,2}})\neq \emptyset\}}&=\big((\cup_{k=0}^{k_1-1}\{2k+1\}) \backslash\{\m^{\theta}_{0,2}\}\\
 &\quad \cup (\cup_{k=0}^{k_2-1}\{\m^{\theta}_{1,1}+2k\})\big), \qand\\
 \scalemath{0.97}{\{k\in\{1,\ldots,\width\}\colon I^{\theta}_k\cap (\q{\theta}_{\m^{\theta}_{1,1}},\q{\theta}_{\m^{\theta}_{1,1}+1})\neq \emptyset\}}&=\big((\cup_{k=0}^{k_1-1}\{2k+1\}) \\
 &\quad \cup (\cup_{k=0}^{k_2-1}\{\m^{\theta}_{1,1}+2k\})\backslash\{\m^{\theta}_{1,1}\}\big).
 \end{split}
 \end{equation}
 This, \cref{lemma:H:updown:lm:wv}, the fact that for all $\theta  \in R$, $x \in [\q{\theta}_{\m^{\theta}_{0,2}-1},\q{\theta}_{\m^{\theta}_{0,2}}]$ it holds that
$\realization{\width,\theta}_\infty(x)=
    \alpha^{\theta}x -(\m^{\theta}_{0,2}-\nicefrac12) \alpha^{\theta}\q{\theta}_1$,
    and the fact that for all $\theta  \in R$, $x \in [\q{\theta}_{\m^{\theta}_{1,1}},\q{\theta}_{\m^{\theta}_{1,1}+1}]$ it holds that
$\realization{\width,\theta}_\infty(x)=
    \beta^{\theta}x+1
    -\beta^{\theta}\left(1-(\width+\nicefrac12-\m^{\theta}_{1,1}) (1-\q{\theta}_\width)\right) $
demonstrate for all $\theta  \in R$ that 
\begin{equation}
\begin{split}
     &\alpha^{\theta}=-((\m^{\theta}_{0,2}-1) \alpha^{\theta}+(\width-\m^{\theta}_{1,1}) \beta^{\theta})+\v{\theta}_{\m^{\theta}_{1,1}}\w{\theta}_{\m^{\theta}_{1,1}} \\
     \andq 
     &\beta^{\theta}=-((\m^{\theta}_{0,2}-1) \alpha^{\theta}+(\width-\m^{\theta}_{1,1}) \beta^{\theta})+\v{\theta}_{\m^{\theta}_{0,2}}\w{\theta}_{\m^{\theta}_{0,2}}.
\end{split}
\end{equation}
Combining this, \cref{lemma:H:updown:lm:alpha}, and \cref{lemma:H:updown:lm:set} ensures that for all $\theta  \in R$, $x \in [\q{\theta}_{\m^{\theta}_{0,2}},\q{\theta}_{\m^{\theta}_{1,1}}]$ it holds that
\begin{equation}\label{lemma:H:updown:lm:realization}
\begin{split}
    \realization{\width,\theta}_\infty(x)&=
    \big(-((\m^{\theta}_{0,2}-1) \alpha^{\theta}+(\width-\m^{\theta}_{1,1}) \beta^{\theta})+\v{\theta}_{\m^{\theta}_{0,2}}\w{\theta}_{\m^{\theta}_{0,2}}\\
    & \quad +\v{\theta}_{\m^{\theta}_{1,1}}\w{\theta}_{\m^{\theta}_{1,1}}\big)(x-\q{\theta}_{\m^{\theta}_{0,2}}) + \frac{\alpha^{\theta}\q{\theta}_1}2
    \\
    &=(\m^{\theta}_{0,2} \alpha^{\theta}+(\width-\m^{\theta}_{1,1}+1) \beta^{\theta})(x-\q{\theta}_{\m^{\theta}_{0,2}})+ \frac{\alpha^{\theta}\q{\theta}_1}2. 
\end{split}
\end{equation}
This, \cref{lemma:H:updown:lm:beta}, and the fact that for all $\theta  \in R$ it holds that the function $\realization{\width,\theta}_\infty$ is continuous prove for all $\theta  \in R$ that 
\begin{equation}\label{lemma:H:updown:lm:continuity}
    (\m^{\theta}_{0,2} \alpha^{\theta}+(\width-\m^{\theta}_{1,1}+1) \beta^{\theta})(\q{\theta}_{\m^{\theta}_{1,1}}-\q{\theta}_{\m^{\theta}_{0,2}})+ \frac{\alpha^{\theta}\q{\theta}_1}2=1-\frac{\beta^{\theta}}2(1-\q{\theta}_\width).
\end{equation}
\Moreover \cref{lemma:H:updown:lm:realization} and \cref{lemma:intmin} imply for all $\theta  \in R$ that 
\begin{equation}
      \int_{\q{\theta}_{\m^{\theta}_{0,2}}}^{\frac12}   ( \realization{\width,\theta}_\infty (x) - \mathbbm{1}_{(\nicefrac12,\infty)}(x))^2 \, \d x  \geq 
     \frac1{12}\bigg(\frac12-\q{\theta}_{\m^{\theta}_{0,2}}\bigg)^3(\m^{\theta}_{0,2} \alpha^{\theta}+(\width-\m^{\theta}_{1,1}+1) \beta^{\theta})^2
\end{equation}
and
\begin{equation}
      \int_{\frac12}^{\q{\theta}_{\m^{\theta}_{1,1}}}   ( \realization{\width,\theta}_\infty (x) - \mathbbm{1}_{(\nicefrac12,\infty)}(x))^2 \, \d x  \geq 
     \frac1{12}\bigg(\q{\theta}_{\m^{\theta}_{1,1}}-\frac12\bigg)^3(\m^{\theta}_{0,2} \alpha^{\theta}+(\width-\m^{\theta}_{1,1}+1) \beta^{\theta})^2.
\end{equation}
This and \cref{lemma:H:updown:lm:int} establish for all $\theta  \in R$ that
\begin{equation}\label{lemma:H:updown:lm:lowerbound}
\begin{split}
     \cL^{\width}_\infty(\theta) &\geq  \frac1{12}\bigg(\bigg(\q{\theta}_{\m^{\theta}_{1,1}}-\frac12\bigg)^3+\bigg(\frac12-\q{\theta}_{\m^{\theta}_{0,2}}\bigg)^3\bigg)(\m^{\theta}_{0,2} \alpha^{\theta}+(\width-\m^{\theta}_{1,1}+1) \beta^{\theta})^2
     \\
     & \quad +\frac{1}{12}(\beta^{\theta})^2(1-\q{\theta}_{\m^{\theta}_{1,1}})^3
     +\frac{1}{12}(\alpha^{\theta})^2(\q{\theta}_{\m^{\theta}_{0,2}})^3.
\end{split}
\end{equation}
We prove \cref{lemma:H:updown:lm:thesis} by contradiction.
Assume that for every $n \in \N$ there exists $\theta_{n}\in (\cG^\width)^{-1}(\{0\})$ with
$0<\q{\theta_{n}}_1  <\q{\theta_{n}}_2<   \ldots < \q{\theta_{n}}_\width<1$,
$\alpha^{\theta_{n}}\beta^{\theta_{n}}\neq 0<\m^{\theta_n}_{0,2}<\m^{\theta_n}_{1,1}<\width+1$, 
$\w{\theta_{n}}_{\m^{\theta_n}_{1,1}}<0<\w{\theta_{n}}_{\m^{\theta_n}_{0,2}}$, and $\forall \, i \in \{0,1,\ldots,\width\}\backslash\{\m^{\theta_n}_{0,2}\} \colon \w{\theta_{n}}_i\w{\theta_{n}}_{i+1}<0$
 which satisfies that
 \begin{equation}\label{lemma:H:updown:lm:ab}
      \cL^{\width}_\infty(\theta_n)\leq \frac1n. 
  \end{equation}
\Nobs that \cref{lemma:H:updown:lm:lowerbound} and \cref{lemma:H:updown:lm:ab} assure that
\begin{equation}
    \begin{split}
    \lim\nolimits_{n\to \infty}& \frac1{12}\left(\left(\q{\theta_n}_{\m^{\theta_n}_{1,1}}-\frac12\right)^3+\left(\frac12-\q{\theta_n}_{\m^{\theta_n}_{0,2}}\right)^3\right)(\m^{\theta_n}_{0,2} \alpha^{\theta_n}+(\width-\m^{\theta_n}_{1,1}+1) \beta^{\theta_n})^2
     \\
     & \quad +\frac{1}{12}(\beta^{\theta_n})^2(1-\q{\theta_n}_{\m^{\theta_n}_{1,1}})^3
     +\frac{1}{12}(\alpha^{\theta_n})^2(\q{\theta_n}_{\m^{\theta_n}_{0,2}})^3=0.
\end{split}
\end{equation}
\Hence that
\begin{equation}\label{lemma:H:updown:lm:limits}
    \begin{split}
    \lim\nolimits_{n\to \infty}& \left(\q{\theta_n}_{\m^{\theta_n}_{1,1}}-\frac12\right)^3(\m^{\theta_n}_{0,2} \alpha^{\theta_n}+(\width-\m^{\theta_n}_{1,1}+1) \beta^{\theta_n})^2=0,\\
    \lim\nolimits_{n\to \infty}& \left(\frac12-\q{\theta_n}_{\m^{\theta_n}_{0,2}}\right)^3(\m^{\theta}_{0,2} \alpha^{\theta_n}+(\width-\m^{\theta_n}_{1,1}+1) \beta^{\theta_n})^2=0, \\
      \lim\nolimits_{n\to \infty}& (\beta^{\theta_n})^2(1-\q{\theta_n}_{\m^{\theta_n}_{1,1}})^3=0, \qand\\
       \lim\nolimits_{n\to \infty}& (\alpha^{\theta_n})^2(\q{\theta_n}_{\m^{\theta_n}_{0,2}})^3=0.
\end{split}
\end{equation}
\Nobs that \cref{lemma:H:updown:lm:continuity} and \cref{lemma:alphabeta} demonstrate that there exists $c \in (0,\infty)$ such that for all $n \in \N$ it holds that $\max\{ |\alpha^{\theta_n}|, |\beta^{\theta_n}|\}\geq c$. 
Combining this and \cref{lemma:H:updown:lm:limits} with \cref{cor:H:updown} assure that there exists a strictly increasing $n \colon \N \to \N$ such that  \begin{equation}
    \lim\nolimits_{k\to \infty}\q{\theta_{n(k)}}_{\m^{\theta_{n(k)}}_{1,1}}=1 \qandq \lim\nolimits_{k\to \infty}\q{\theta_{n(k)}}_{\m^{\theta_{n(k)}}_{0,2}}=0.
\end{equation}
This  shows that there exists $k^* \in \N$ such that for all $k\in \N \cap[k^*,\infty)$  it holds that $\q{\theta_{n(k)}}_{\m^{\theta}_{1,1}}\geq \nicefrac34$ and $\q{\theta_{n(k)}}_{\m^{\theta}_{0,2}}\leq \nicefrac14$.
Combining this and \cref{lemma:halfc} assures that $\lim\nolimits_{k\to \infty}\cL^{\width}_\infty(\theta_{n(k)})\allowbreak \geq \nicefrac1{32}$ which is a contradiction.
\end{cproof}
\subsection{Estimates for the risk of critical points}\label{Annsrelu_esti}
\begin{cor}\label{lemma:H:q0}
Assume \cref{H:setting:snn} and let $\width\in \N$, $j\in \{1,2,\ldots,\width\}$, $\theta\in (\cG^\width)^{-1}(\{0\})$ satisfy 
$(0,1)\subseteq I^{\theta}_j$ and $\prod_{k=1}^\width \v{\theta}_k\neq0$. Then 
 $\cL_\infty^\width(\theta)\geq\nicefrac1{36}$.
\end{cor}
\begin{cproof}{lemma:H:q0}
\Nobs that the assumption that $\theta\in (\cG^\width)^{-1}(\{0\})$, the assumption that $\prod_{k=1}^\width \v{\theta}_k\neq0$, the assumption that $(0,1)\subseteq I^{\theta}_j$, and \cref{prop:H:approximate:gradient} demonstrate that
\begin{equation}
     \int_{0}^{1}  x( \realization{\width,\theta}_\infty (x) - \mathbbm{1}_{(\nicefrac12,\infty)}(x)) \, \d x =0.
\end{equation}
Combining this with \cref{cor:H:x0} assures that $\cL_\infty^\width(\theta)\geq\nicefrac1{36}$.
\end{cproof}
\begin{cor}\label{lemma:H:qq}
Assume \cref{H:setting:snn} and let $\width\in \N \cap (1,\infty)$, $i,j\in \{1,2,\ldots,\width\}$, $\theta\in (\cG^\width)^{-1}(\{0\})$ satisfy  $i\neq j$, $\q{\theta}_i=\q{\theta}_j \in (0,1)$, and $\w{\theta}_j\w{\theta}_i<0\neq \prod_{k=1}^\width \v{\theta}_k$. Then 
 $\cL_\infty^\width(\theta)\geq\nicefrac1{36}$.
\end{cor}
\begin{cproof}{lemma:H:qq}
\Nobs that the assumption that $\theta\in (\cG^\width)^{-1}(\{0\})$, the assumption that $\w{\theta}_j\w{\theta}_i<0\neq \prod_{k=1}^\width \v{\theta}_k$, and \cref{prop:H:approximate:gradient} assure that 
\begin{equation}
    \int_{0}^{\q{\theta}_i}  x( \realization{\width,\theta}_\infty (x) - \mathbbm{1}_{(\nicefrac12,\infty)}(x)) \, \d x =
    \int_{\q{\theta}_i}^{1}  x( \realization{\width,\theta}_\infty (x) - \mathbbm{1}_{(\nicefrac12,\infty)}(x)) \, \d x =0.
\end{equation}
This implies that
\begin{equation}
    \int_{0}^{1}  x( \realization{\width,\theta}_\infty (x) - \mathbbm{1}_{(\nicefrac12,\infty)}(x)) \, \d x =0.
\end{equation}
Combing this with \cref{cor:H:x0} demonstrates that $\cL_\infty^\width(\theta)\geq\nicefrac1{36}$.
\end{cproof}
%
%
\begin{prop}\label{lemma:H:qqneg}
Assume \cref{H:setting:snn} and let $\width\in \N\cap(1,\infty)$, $i,j\in \{1,2,\ldots,\width\}$, $\theta\in (\cG^\width)^{-1}(\{0\})$ satisfy  $i\neq j$, $\q{\theta}_i=\q{\theta}_j \in (0,1)$, and $ \prod_{k=1}^\width \v{\theta}_k\neq0<\w{\theta}_j\w{\theta}_i$. Then there exist $\vartheta \in (\cG^{\width})^{-1}(\{0\})$, $k\in \{1,2,\ldots,\width\}$ such that $I^{\vartheta}_k=\emptyset$ and ${\realization{\width,\theta}_\infty}|_{[0,1]}={\realization{\width,\vartheta}_\infty}|_{[0,1]}$.
\end{prop}
\begin{cproof}{lemma:H:qqneg}
\Nobs that the assumption that $\q{\theta}_i=\q{\theta}_j$ and the assumption that $ 0<\w{\theta}_j\w{\theta}_i$ demonstrate that
\begin{equation}\label{lemma:H:qqneg:eq}
   \frac{\b{\theta}_i}{|\w{\theta}_i|}=\frac{\b{\theta}_j}{|\w{\theta}_j|}.
\end{equation}
Assume without loss of generality that $i=1$ and $j=2$. This and \cref{lemma:H:qqneg:eq} ensure for all $x\in[0,1]$ that
\begin{equation}\label{lemma:H:qqneg:eq2}
\begin{split}
    \v{\theta}_1 \br[\big]{ \act_\infty ( \w{\theta}_1 x + \b{\theta}_1 )}+\v{\theta}_2 \br[\big]{ \act_\infty ( \w{\theta}_2 x + \b{\theta}_2 )} &=
     \v{\theta}_1 |\w{\theta}_1| \act_\infty \bigg( \frac{\w{\theta}_1}{|\w{\theta}_1|} x + \frac{\b{\theta}_1 }{|\w{\theta}_1|}\bigg)\\
     & \quad +\v{\theta}_2 |\w{\theta}_2|  \act_\infty 
      \bigg( \frac{\w{\theta}_2}{|\w{\theta}_2|} x + \frac{\b{\theta}_2 }{|\w{\theta}_2|}\bigg) \\
      &=\big(\v{\theta}_1|\w{\theta}_1|+\v{\theta}_2|\w{\theta}_2|\big)  \act_\infty \bigg( \frac{\w{\theta}_1}{|\w{\theta}_1|} x + \frac{\b{\theta}_1 }{|\w{\theta}_1|}\bigg).
\end{split}
\end{equation}
Let $\vartheta \in \R^{\fd_\width}$ satisfy for all $m \in \{1,2,\ldots,\width\}\backslash\{1,2\}$ that $\v{\vartheta}_1=\v{\theta}_1|\w{\theta}_1|+\v{\theta}_2|\w{\theta}_2|$,
$\w{\vartheta}_1=\nicefrac{\w{\theta}_1}{|\w{\theta}_1|}$,
$\b{\vartheta}_1=\nicefrac{\b{\theta}_1 }{|\w{\theta}_1|}$, 
$\v{\vartheta}_2=\w{\vartheta}_2=\b{\vartheta}_2=0$, $\c{\vartheta}=\c{\theta}$, $\v{\vartheta}_m=\v{\theta}_m$, $\w{\vartheta}_m=\w{\theta}_m$, and $\b{\vartheta}_m=\b{\theta}_m$.
This and \cref{lemma:H:qqneg:eq2} imply for all $x \in [0,1]$ that $I^{\vartheta}_2=\emptyset$,  $\cG^\width(\vartheta)=0$, and 
\begin{equation}
     \realization{\width,\theta}_   \infty(x)=\realization{\width,\vartheta}_\infty(x).
\end{equation}
\end{cproof}
\begin{lemma}\label{lemma:H:case4sx}
Assume \cref{H:setting:snn} and let $\width \in \N$, $\theta\in (\cG^\width)^{-1}(\{0\})$ satisfy 
 for all $i \in \{1,2,\ldots,\width\}$ that $\q{\theta}_{i-1}<\q{\theta}_i\leq\nicefrac12$ and
 $\prod_{k=1}^\width{\v{\theta}_k}\neq0$.
Then \begin{equation}\label{lemma:H:case4sx:thesis}
    \cL^{\width}_\infty(\theta)\geq\nicefrac1{36}.
\end{equation}
\end{lemma}
\begin{cproof}{lemma:H:case4sx}
In the following we distinguish between the case $\alpha^{\theta}=0$ and the case $\alpha^{\theta}\neq0$.
We first show \cref{lemma:H:case4sx:thesis} in the case \begin{equation}\label{lemma:H:case4sx:case1}
    \alpha^{\theta}=0.
\end{equation} 
\Nobs that \cref{lemma:H:case4sx:case1}, the assumption that $\theta\in (\cG^\width)^{-1}(\{0\})$, and \cref{cor:H:m1} establish that for all $x \in [0,\q{\theta}_{\m^{\theta}_{0,2}}]$ it holds that $\realization{\width,\theta}_\infty(x)=0$. 
\Hence that there exists $a^{\theta} \in \R$ such that for all $x \in [\q{\theta}_{\m^{\theta}_{0,2}},1]$ it holds that $\realization{\width,\theta}_\infty(x)=a^{\theta}(x-\q{\theta}_{\m^{\theta}_{0,2}})$.
Combining this with \cref{lemma:half} ensures that
\begin{equation}\label{lemma:H:case4sx:res1}
    \cL^{\width}_\infty(\theta)\geq\frac1{36}.
\end{equation}
This establishes \cref{lemma:H:case4sx:thesis} in the case $ \alpha^{\theta}=0.$
In the next step we demonstrate \cref{lemma:H:case4sx:thesis} in the case 
\begin{equation}\label{lemma:H:case4sx:case2}
    \alpha^{\theta}\neq0.
\end{equation}
\Nobs that \cref{lemma:H:case4sx:case2} and \cref{cor:H:w0sx} prove that for all $i \in \{1,2,\ldots, \m^{\theta}_{0,2} \}$ it holds that $\w{\theta}_i \w{\theta}_{i-1}<0$.
Combining this with \cref{lemma:H:updown:sx} demonstrates that
\begin{enumerate}[label=(\roman*)]
    \item for all $j \in \{1,2,\ldots, \m^{\theta}_{0,2}\}$ it holds that $\q{\theta}_j=j\q{\theta}_1$,
    \item for all $j \in \{0,1,\ldots, \m^{\theta}_{0,2}\}$ it holds that 
  \begin{equation}\label{lemma:H:case4sx:nq}
    - \frac{ \alpha^{\theta}\q{\theta}_1}2=(-1)^{j} \realization{\width,\theta}_\infty(\q{\theta}_j),
    \end{equation} 
    and 
   \item it holds that
 \begin{equation}\label{lemma:H:case4sx:int}
     \int_0^{\q{\theta}_{\m^{\theta}_{0,2}}}   ( \realization{\width,\theta}_\infty (x) - \mathbbm{1}_{(\nicefrac12,\infty)}(x))^2 \, \d x=\frac{\m^{\theta}_{0,2}}{12}(\alpha^{\theta})^2(\q{\theta}_1)^3.
 \end{equation}
\end{enumerate}
\Moreover the assumption that for all $i \in \{1,2,\ldots,\width\}$ it holds that $\q{\theta}_i\leq \nicefrac12$ implies that there exists $a^{\theta}\in \R$ such that for all $x \in [\q{\theta}_{\m^{\theta}_{0,2}},1]$ it holds that $\realization{\width,\theta}_\infty(x)= a^{\theta}(x-\q{\theta}_{\m^{\theta}_{0,2}})+\realization{\width,\theta}_\infty(\q{\theta}_{\m^{\theta}_{0,2}})$.
This, \cref{lemma:H:case4sx:nq}, and \cref{lemma:H:case4sx:int} establish that 
\begin{equation}
\label{lemma:four:case4sx:zerogradient:q0left:eq1}
\begin{split}
     \cL^{\width}_\infty(\theta) & = \frac{\q{\theta}_{\m^{\theta}_{0,2}}}{3}(\realization{\width,\theta}_\infty(\q{\theta}_{\m^{\theta}_{0,2}}))^2+
      \int_{\q{\theta}_{\m^{\theta}_{0,2}}}^1 (a^{\theta}(x-\q{\theta}_{\m^{\theta}_{0,2}})+\realization{\width,\theta}_\infty(\q{\theta}_{\m^{\theta}_{0,2}}) - \mathbbm{1}_{(\nicefrac12,\infty)}(x) )^2 \, \d x \\
      & =\frac{\q{\theta}_{\m^{\theta}_{0,2}}}{3}(\realization{\width,\theta}_\infty(\q{\theta}_{\m^{\theta}_{0,2}}))^2+
      \frac{{(a^{\theta}})^2}{3}(1-(\q{\theta}_{\m^{\theta}_{0,2}})^3)+\frac12-\frac{3a^{\theta}}4\\
      & \quad +(-a^{\theta}\q{\theta}_{\m^{\theta}_{0,2}}+\realization{\width,\theta}_\infty (\q{\theta}_{\m^{\theta}_{0,2}}))^2(1-\q{\theta}_{\m^{\theta}_{0,2}})+(a^{\theta}\q{\theta}_{\m^{\theta}_{0,2}}-\realization{\width,\theta}_\infty(\q{\theta}_{\m^{\theta}_{0,2}})) \\
      & \quad  
      +a^{\theta}(-a^{\theta}\q{\theta}_{\m^{\theta}_{0,2}}+\realization{\width,\theta}_\infty (\q{\theta}_{\m^{\theta}_{0,2}}))(1-(\q{\theta}_{\m^{\theta}_{0,2}})^2).
\end{split}
\end{equation}
Combining this and \cref{lemma:f:case1} ensures that $\cL^{\width}_\infty(\theta)\geq \nicefrac1{18}$
This shows \cref{lemma:H:case4sx:thesis} in the case $\alpha^{\theta}\neq0$.
\end{cproof}
\begin{lemma}\label{lemma:H:case4dx}
Assume \cref{H:setting:snn} and let $\width \in \N$, $\theta\in (\cG^\width)^{-1}(\{0\})$ satisfy 
for all $i \in \{1,2,\ldots,\width\}$ that $\nicefrac12\leq \q{\theta}_{i}<\q{\theta}_{i+1}$ and
 $\prod_{k=1}^\width{\v{\theta}_k}\neq0$.
Then \begin{equation}\label{lemma:H:case4dx:thesis}
    \cL^{\width}_\infty(\theta)\geq\nicefrac1{36}.
\end{equation}
\end{lemma}
\begin{cproof}{lemma:H:case4dx}
In the following we distinguish between the case $\beta^{\theta}=0$ and the case $\beta^{\theta}\neq0$.
We first demonstrate \cref{lemma:H:case4dx:thesis} in the case
\begin{equation}\label{lemma:H:case4dx:case1}
    \beta^{\theta}=0.
\end{equation}
\Nobs that \cref{lemma:H:case4dx:case1}, the assumption that $\theta\in (\cG^\width)^{-1}(\{0\})$, and \cref{cor:H:l2} establish that for all $x \in [\q{\theta}_{\m^{\theta}_{1,1}},1]$ it holds that $\realization{\width,\theta}_\infty(x)=1$.
\Hence that there exists $a^{\theta} \in \R$ such that for all $x \in [0,\q{\theta}_{\m^{\theta}_{1,1}}]$ it holds that $\realization{\width,\theta}_\infty(x)=a^{\theta}(x-\q{\theta}_{\m^{\theta}_{1,1}})+1$.
Combining this with \cref{lemma:halfdx} ensures that
\begin{equation}\label{lemma:H:case4dx:case}
    \cL^{\width}_\infty(\theta)\geq\frac1{36}.
\end{equation}
This establishes \cref{lemma:H:case4dx:thesis} in the case $ \beta^{\theta}=0.$
In the next step we demonstrate \cref{lemma:H:case4dx:thesis} in the case 
\begin{equation}\label{lemma:H:case4dx:case2}
    \beta^{\theta}\neq0.
\end{equation}
\Nobs that \cref{lemma:H:case4dx:case2} and \cref{cor:H:w0sdx} prove that for all $i \in \{\m^{\theta}_{1,1},\m^{\theta}_{1,1}+1,\ldots, \width \}$ it holds that $\w{\theta}_i \w{\theta}_{i+1}<0$.
Combining this with \cref{lemma:H:updown:dx} demonstrates that 
\begin{enumerate}[label=(\roman*)]
    \item for all $j \in \{\m^{\theta}_{1,1},\m^{\theta}_{1,1}+1,\ldots, \width\}$ it holds that $\q{\theta}_j=1-(\width+1-j)( 1-\q{\theta}_\width)$,
    \item
    for all $j \in \{\m^{\theta}_{1,1},\m^{\theta}_{1,1}+1,\ldots, \width+1\}$ it holds that 
\begin{equation}\label{lemma:H:case4dx:nq}
    \frac{\beta^{\theta}}{2}(1-\q{\theta}_\width)=(-1)^{ \width +1-j} (\realization{\width,\theta}_\infty(\q{\theta}_j)-1),
\end{equation} 
    \item
    and for all $j \in \{\m^{\theta}_{1,1},\m^{\theta}_{1,1}+1,\ldots, \width\}$, $x \in [\q{\theta}_{j},\q{\theta}_{j+1}]$ it holds that
 \begin{equation}\label{lemma:H:case4dx:int}
\int_{\q{\theta}_{\m^{\theta}_{1,1}}}^{1}   ( \realization{\width,\theta}_\infty (x) - \mathbbm{1}_{(\nicefrac12,\infty)}(x))^2 \, \d x=\frac{1}{12}(\width+1-\m^{\theta}_{1,1})(\beta^{\theta})^2(1-\q{\theta}_\width)^3.
\end{equation}
\end{enumerate}
\Moreover the assumption that for all $i \in \{1,2,\ldots,\width\}$ it holds that $\nicefrac12\leq\q{\theta}_i$ implies that there exists $a^{\theta}\in \R$ such that for all $x \in [0,\q{\theta}_{\m^{\theta}_{1,1}}]$ it holds that $\realization{\width,\theta}_\infty(x)= a^{\theta}(x-\q{\theta}_{\m^{\theta}_{1,1}})+\realization{\width,\theta}_\infty(\q{\theta}_{\m^{\theta}_{1,1}})$.
This, \cref{lemma:H:case4dx:nq}, and \cref{lemma:H:case4dx:int} establish that 
\begin{equation}
\label{lemma:four:case4dx:zerogradient:q0left:eq1}
\begin{split}
     \cL^{\width}_\infty(\theta) & =\int_0^{\q{\theta}_{\m^{\theta}_{1,1}}} (a^{\theta}(x-\q{\theta}_{\m^{\theta}_{1,1}})+\realization{\width,\theta}_\infty(\q{\theta}_{\m^{\theta}_{1,1}}) - \mathbbm{1}_{(\nicefrac12,\infty)}(x) )^2 \, \d x\\
     & \quad +\frac1{3}(1-\q{\theta}_{\m^{\theta}_{1,1}})(\realization{\width,\theta}_\infty(\q{\theta}_{\m^{\theta}_{1,1}})-1)^2
       \\
      & =\frac{(a^{\theta})^2}{3}(\q{\theta}_{\m^{\theta}_{1,1}})^3+(-a^{\theta}\q{\theta}_{\m^{\theta}_{1,1}}+\realization{\width,\theta}_\infty (\q{\theta}_{\m^{\theta}_{1,1}}))^2\q{\theta}_{\m^{\theta}_{1,1}}-a^{\theta}\left((\q{\theta}_{\m^{\theta}_{1,1}})^2-\frac14\right)
       \\
      & \quad +a^{\theta}(-a^{\theta}\q{\theta}_{\m^{\theta}_{1,1}}+\realization{\width,\theta}_\infty (\q{\theta}_{\m^{\theta}_{1,1}}))(\q{\theta}_{\m^{\theta}_{1,1}})^2+ \frac13(1-\q{\theta}_{\m^{\theta}_{1,1}}) (\realization{\width,\theta}_\infty (\q{\theta}_{\m^{\theta}_{1,1}})-1)^2
      \\
      & \quad  +(1+2a^{\theta}\q{\theta}_{\m^{\theta}_{1,1}}-2\realization{\width,\theta}_\infty (\q{\theta}_{\m^{\theta}_{1,1}}))\left(\q{\theta}_{\m^{\theta}_{1,1}}-\frac12\right).
\end{split}
\end{equation}
Combining this and \cref{lemma:f:case2} ensures that $\cL^{\width}_\infty(\theta)\geq \nicefrac1{18}$.
This shows \cref{lemma:H:case4dx:thesis} in the case $ \beta^{\theta}\neq0.$
\end{cproof}
\begin{cor}\label{lemma:H:updown}
Assume \cref{H:setting:snn} and let $\width \in \N\cap(1,\infty)$. Then there exists $\varepsilon \in (0,\infty)$ such that for all $\theta\in (\cG^\width)^{-1}(\{0\})$ with $\forall \, i \in \{0,1,\ldots,\width\}\backslash\{\m^{\theta}_{0,2}\} \colon \w{\theta}_i\w{\theta}_{i+1}<0\neq \alpha^{\theta}\beta^{\theta}$, $\prod_{k=1}^\width\v{\theta}_k\neq0<\q{\theta}_1 <\q{\theta}_2< \ldots < \q{\theta}_\width<1$, and $\M^{\theta}_{0}\neq\emptyset\neq\M^{\theta}_{1}$ it holds that
$\cL^{\width}_\infty(\theta)\geq \varepsilon$.
\end{cor}
\begin{cproof}{lemma:H:updown}
\Nobs that \cref{lemma:H:updown:leqm} implies that 
for all $\theta\in (\cG^\width)^{-1}(\{0\})$ with $\forall \, i \in \{0,1,\ldots,\width\} \colon \w{\theta}_i\w{\theta}_{i+1}<0\neq\alpha^{\theta}\beta^{\theta}$,
$\prod_{k=1}^\width\v{\theta}_k\neq0<\q{\theta}_1  <\q{\theta}_2 < \ldots < \q{\theta}_\width<1$, and $\q{\theta}_{\m^{\theta}_{0,2}}=\q{\theta}_{\m^{\theta}_{1,1}}=\nicefrac12$ it holds that
\begin{equation}\label{lemma:H:updown:1}
    \cL^{\width}_\infty(\theta)\geq\frac1{12}.
\end{equation}
\Moreover \cref{lemma:H:updown:weq} assures that for all $\theta\in \R^{\fd_\width}$ with $\forall \, i \in \{0,1,\ldots,\width\}\backslash \allowbreak \{\m^{\theta}_{0,2}\} \colon \w{\theta}_i\w{\theta}_{i+1}<0$,
$\prod_{k=1}^\width\v{\theta}_k\neq0<\q{\theta}_1  < \q{\theta}_2  <\ldots < \q{\theta}_\width<1$, $\alpha^{\theta}\beta^{\theta}\neq 0<\m^{\theta}_{0,2}<\m^{\theta}_{1,1}<\width+1$, and $0<\w{\theta}_{\m^{\theta}_{0,2}}\w{\theta}_{\m^{\theta}_{1,1}}$ it holds that 
\begin{equation}\label{lemma:H:updown:2}
    \cG^\width(\theta)\neq0.
\end{equation}
\Moreover \cref{lemma:H:updown:ml} and \cref{lemma:H:updown:lm} demonstrate that there exists $\delta \in (0,\infty)$ such that
for all $\theta\in (\cG^\width)^{-1}(\{0\})$ with $\forall \, i \in \{0,1,\ldots,\width\}\backslash \allowbreak \{\m^{\theta}_{0,2}\} \colon \w{\theta}_i\w{\theta}_{i+1}<0$, 
$\prod_{k=1}^\width\v{\theta}_k\neq0<\q{\theta}_1  < \q{\theta}_2  < \ldots < \q{\theta}_\width<1$, 
$\alpha^{\theta}\beta^{\theta}\neq 0<\m^{\theta}_{0,2}<\m^{\theta}_{1,1}<\width+1$, and $\w{\theta}_{\m^{\theta}_{0,2}}\w{\theta}_{\m^{\theta}_{1,1}}<0$ 
it holds that
$\cL^{\width}_\infty(\theta)\geq\delta$.
Combining this with \cref{lemma:H:updown:1} and \cref{lemma:H:updown:2}
shows that
for all $\theta\in (\cG^\width)^{-1}(\{0\})$ with $\forall \, i \in \{0,1,\ldots,\width\}\backslash\{\m^{\theta}_{0,2}\} \colon \w{\theta}_i\w{\theta}_{i+1}<0\neq \alpha^{\theta}\beta^{\theta}$, $\prod_{k=1}^\width\v{\theta}_k\neq0<\q{\theta}_1 <\q{\theta}_2< \ldots < \q{\theta}_\width<1$, and $\M^{\theta}_{0}\neq\emptyset\neq\M^{\theta}_{1}$ it holds that
$\cL^{\width}_\infty(\theta)\geq\min\{\nicefrac1{12},\delta\}$.
\end{cproof}
\begin{lemma}\label{lemma:H:N0:finN}
Assume \cref{H:setting:snn}, let $\width \in \N\cap(1,\infty)$, $\theta\in (\cG^\width)^{-1}(\{0\})$, $i \in \{0,1,\ldots,\width\}\backslash\{\m^{\theta}_{0,2}\}$  satisfy $\prod_{k=1}^\width\v{\theta}_k\neq0<\q{\theta}_1 <\q{\theta}_2 <   \ldots < \q{\theta}_\width<1$, $-\w{\theta}_i \w{\theta}_{i+1}<0$, and $\M^{\theta}_{0}\neq\emptyset\neq\M^{\theta}_{1}$. Then 
\begin{equation}\label{lemma:H:N0:finN:thesis}
    \cL^{\width}_\infty(\theta)\geq\frac{1}{384(1+\width)^2}.
\end{equation}
\end{lemma}
\begin{cproof}{lemma:H:N0:finN}
\Nobs that the fact that  $i \in \{0,1,\ldots,\width\}\backslash\{\m^{\theta}_{0,2}\}$ demonstrates that $\nicefrac12 \notin (\q{\theta}_i,\q{\theta}_{i+1})$. In the following we distinguish between the case $\q{\theta}_{i+1}\leq\nicefrac12$ and the case $\q{\theta}_{i}\geq \nicefrac12$.
We first prove \cref{lemma:H:N0:finN:thesis} in the case
\begin{equation}\label{lemma:H:N0:finN:case1}
    \q{\theta}_{i+1}\leq\frac12.
\end{equation}
\Nobs that \cref{lemma:H:N0:finN:case1} and \cref{cor:H:w0sx} ensure that for all $x \in [0,\q{\theta}_{\m^{\theta}_{0,2}}]$ it holds that $\realization{\width,\theta}_\infty(x)=0$.
This and \cref{lemma:H:N0:init} imply that 
\begin{equation}
    \cL^{\width}_\infty(\theta)\geq\frac{1}{384(1+\width)^2}.
\end{equation}
This establishes \cref{lemma:H:N0:finN:thesis} in the case $\q{\theta}_{i+1}\leq\nicefrac12$.
In the next step we demonstrate \cref{lemma:H:N0:finN:thesis} in the case
\begin{equation}\label{lemma:H:N0:finN:case2}
    \q{\theta}_{i}\geq \frac12.
\end{equation}
\Nobs that \cref{lemma:H:N0:finN:case2} and \cref{cor:H:w0sdx} show that for all $x \in [\q{\theta}_{\m^{\theta}_{1,1}},1]$ it holds that $\realization{\width,\theta}_\infty(x)=1$. This and \cref{lemma:H:N0:fin} establish that 
\begin{equation}
    \cL^{\width}_\infty(\theta)\geq\frac{1}{384(1+\width)^2}.
\end{equation}
This establishes \cref{lemma:H:N0:finN:thesis} in the case $\q{\theta}_{i}\geq \nicefrac12.$
\end{cproof}
\begin{lemma}\label{lemma:H:N0:finN2}
Assume \cref{H:setting:snn}, let $\width \in \N\cap(1,\infty)$, $\theta\in (\cG^\width)^{-1}(\{0\})$ satisfy 
$\prod_{k=1}^\width\v{\theta}_k\neq0<\q{\theta}_1 <\q{\theta}_2<   \ldots < \q{\theta}_\width<1$, $\alpha^{\theta}\beta^{\theta}=0$, and $\M^{\theta}_{0}\neq\emptyset\neq\M^{\theta}_{1}$. Then 
\begin{equation}\label{lemma:H:N0:finN2:thesis}
    \cL^{\width}_\infty(\theta)\geq\frac{1}{384(1+\width)^2}.
\end{equation}
\end{lemma}
\begin{cproof}{lemma:H:N0:finN2}
\Nobs that \cref{cor:H:m1}, \cref{cor:H:l2}, \cref{lemma:H:N0}, and the assumption that $\theta\in (\cG^\width)^{-1}(\{0\})$ prove that $\max\{|\alpha^{\theta}|,|\beta^{\theta}|\}\neq0$.
In the following we distinguish between the case $\alpha^{\theta}=0\neq\beta^{\theta}$ and the case $\alpha^{\theta}\neq0=\beta^{\theta}$.
We first demonstrate \cref{lemma:H:N0:finN2:thesis} in the case
\begin{equation}\label{lemma:H:N0:finN2:case1}
    \alpha^{\theta}=0\neq\beta^{\theta}.
\end{equation}
\Nobs that \cref{lemma:H:N0:finN2:case1}, the assumption that $\theta\in (\cG^\width)^{-1}(\{0\})$, and \cref{cor:H:m1} ensure that for all $x \in [0,\q{\theta}_{\m^{\theta}_{0,2}}]$ it holds that $\realization{\width,\theta}_\infty(x)=0$. This and \cref{lemma:H:N0:init} imply that 
\begin{equation}
    \cL^{\width}_\infty(\theta)\geq\frac{1}{384(1+\width)^2}.
\end{equation}
This establishes \cref{lemma:H:N0:finN2:thesis} in the case $\alpha^{\theta}=0\neq\beta^{\theta}.$
In the next step we prove \cref{lemma:H:N0:finN2:thesis} in the case
\begin{equation}\label{lemma:H:N0:finN2:case2}
    \alpha^{\theta}\neq0=\beta^{\theta}.
\end{equation}
\Nobs that \cref{lemma:H:N0:finN2:case2}, the assumption that $\theta\in (\cG^\width)^{-1}(\{0\})$, and \cref{cor:H:l2} show that for all $x \in [\q{\theta}_{\m^{\theta}_{1,1}},1]$ it holds that $\realization{\width,\theta}_\infty(x)=1$. This and \cref{lemma:H:N0:fin} establish that 
\begin{equation}
    \cL^{\width}_\infty(\theta)\geq\frac{1}{384(1+\width)^2}.
\end{equation}
This demonstrates \cref{lemma:H:N0:finN2:thesis} in the case $\alpha^{\theta}\neq0=\beta^{\theta}.$
\end{cproof}
\begin{prop}\label{lemma:H:1}
Assume \cref{H:setting:snn} and let $\theta\in (\cG^1)^{-1}(\{0\})$. Then 
\begin{equation}\label{lemma:H:1:thesis}
    \cL^{1}_\infty(\theta)\geq\nicefrac1{36}.
\end{equation}
\end{prop}
\begin{cproof}{lemma:H:1}
\Nobs that in the case $\v{\theta}_1=0$ there exists $b \in \R$ such that for all $x \in \R$ it holds that $\realization{1,\theta}_\infty(x)=b$. This and \cref{lemma:line} establish that in the case $\v{\theta}_1=0$ it holds that
\begin{equation}\label{lemma:H:1:eq1}
    \cL^1_\infty(\theta)\geq\frac1{16}.
\end{equation}
Assume now that $\v{\theta}_1 \neq 0$.
In the following we distinguish between the case $\q{\theta}_1\in(0,1)$ and the case $\q{\theta}_1\notin (0,1)$.
We first show \cref{lemma:H:1:thesis} in the case
\begin{equation}\label{lemma:H:1:case1}
    \q{\theta}_1\in(0,1).
\end{equation} 
\Nobs that \cref{lemma:H:1:case1}, \cref{lemma:H:case4sx}, and  \cref{lemma:H:case4dx} prove that 
\begin{equation}\label{lemma:H:1:eq2}
    \cL^1_\infty(\theta)\geq\frac1{36}.
\end{equation}
This establishes \cref{lemma:H:1:thesis} in the case $ \q{\theta}_1\in(0,1).$
In the next step we demonstrate \cref{lemma:H:1:thesis} in the case
\begin{equation}\label{lemma:H:1:case2}
    \q{\theta}_1\notin(0,1).
\end{equation}
\Nobs that \cref{lemma:H:1:case2} ensures that there exist $a,b\in \R$ such that for all $x \in \R$ it holds that $\realization{1,\theta}_\infty(x)=a x + b$.
This and \cref{lemma:line} demonstrate that
\begin{equation}\label{lemma:H:1:eq3}
    \cL^1_\infty(\theta)\geq\frac1{16}.
\end{equation}
This establishes \cref{lemma:H:1:thesis} in the case $\q{\theta}_1\notin(0,1).$
\end{cproof}
\begin{lemma}\label{lemma:H:active:case0}
Assume \cref{H:setting:snn}
and let $\width\in \N\cap(1,\infty)$, $\theta\in (\cG^\width)^{-1}(\{0\})$ 
 satisfy $\{\vartheta\in (\cG^{\width-1})^{-1}(\{0\}) \colon {\realization{\width,\theta}_\infty}| _{[0,1]}={\realization{\width-1,\vartheta}_\infty}| _{[0,1]}\}=\emptyset$.
 Then it holds that
 \begin{equation}\label{lemma:H:active:case0:eq}
     \textstyle \prod_{k=1}^\width \v{\theta}_k\neq0.
 \end{equation}
\end{lemma}
\begin{cproof}{lemma:H:active:case0}
We prove \cref{lemma:H:active:case0:eq} by contradiction. Assume that $\prod_{k=1}^\width \v{\theta}_k=0$ and assume without loss of generality that $\v{\theta}_\width=0$. Throughout this proof let $\vartheta \in \R^{\fd_{\width-1}}$ satisfy for all $i \in \{1,2,\ldots,\width-1\}$ that 
\begin{equation}\label{lemma:H:active:case0:hp}
    \vartheta_i=\theta_i, \qquad     
    \vartheta_{\width-1+i}=\theta_{\width+i}, \qquad
    \vartheta_{2\width-2+i}=\theta_{2\width+i},
    \qandq     \vartheta_{3\width-2}=\theta_{3\width+1}.
\end{equation} 
\Nobs that \cref{lemma:H:active:case0:hp} assures that \begin{equation}\label{lemma:H:active:case0:th}
    {\realization{\width,\theta}_\infty}| _{[0,1]}={\realization{\width-1,\vartheta}_\infty}| _{[0,1]}.
\end{equation}
\Moreover that \cref{lemma:H:active:case0:hp} and \cref{prop:H:approximate:gradient} show that for all $i \in \{1,2,\ldots,3(\width-1)\}$ it holds that $\cG^{\width-1}_i(\vartheta)=\cG^{\width}_i(\theta)$ and $\cG^{\width-1}_{\fd_{\width-1}}(\vartheta)=\cG^{\width}_{\fd_\width}(\theta)$. This, \cref{lemma:H:active:case0:th}, and the assumption that $\theta\in (\cG^\width)^{-1}(\{0\})$ imply that $\vartheta \in \{v \in (\cG^{\width-1})^{-1}(\{0\}) \colon {\realization{\width,\theta}_\infty}| _{[0,1]}={\realization{\width-1,v}_\infty}|_{[0,1]}\}$ which is a contradiction.
\end{cproof}
\begin{lemma}\label{lemma:H:active}
Assume \cref{H:setting:snn}
and let $\width\in \N\cap(1,\infty)$.
 Then there exists $\varepsilon \in (0,\infty)$ such that 
 for all $\theta\in (\cG^\width)^{-1}(\{0\})$ 
 with $\{\vartheta\in (\cG^{\width-1})^{-1}(\{0\}) \colon {\realization{\width,\theta}_\infty}| _{[0,1]}={\realization{\width-1,\vartheta}_\infty}| _{[0,1]}\}=\emptyset$ it holds that
 $\cL^{\width}_\infty(\theta)\geq \varepsilon$.
\end{lemma}
\begin{cproof}{lemma:H:active}
Throughout this proof for every $\theta \in \R^{\fd_{\width}}$ let $R^{\theta}\subseteq\R^{\fd_{h-1}}$ 
satisfy $R^{\theta}= \{\vartheta\in (\cG^{\width-1})^{-1}(\{0\}) \colon {\realization{\width,\theta}_\infty}| _{[0,1]}={\realization{\width-1,\vartheta}_\infty}| _{[0,1]}\}$.
\Nobs that \cref{lemma:H:active:case0} ensures that for all $\theta \in  (\cG^\width)^{-1}(\{0\})$ with $R^{\theta}=\emptyset$ it holds that $\prod_{k=1}^\width \v{\theta}_k\neq0$.
\Moreover \cref{lemma:H:q0} demonstrates that for all $\theta \in (\cG^\width)^{-1}(\{0\})$ with $\prod_{k=1}^\width \v{\theta}_k\neq0$ and 
$\{k\in \{1,2,\ldots,\width\} \colon (0,1)\subseteq I_k^{\theta} \}\neq \emptyset$
 it holds that
\begin{equation}\label{lemma:H:active:1}
    \cL^{\width}_\infty(\theta)\geq \frac1{36}.
\end{equation}
\Moreover \cref{lemma:H:qq} and \cref{lemma:H:qqneg} assure that for all $\theta \in (\cG^\width)^{-1}(\{0\})$ with  $\prod_{k=1}^\width \v{\theta}_k\neq0$, $R^{\theta}=\emptyset$, and $\{k\in \{1,2,\ldots,\width\} \colon [\exists\, j\in \{1,2,\ldots,\width\}\backslash\{k\}\colon \q{\theta}_k=\q{\theta}_j \in (0,1)]\}\neq \emptyset$
it holds that 
\begin{equation}\label{lemma:H:active:2}
   \cL^{\width}_\infty(\theta)\geq \frac1{36}.
\end{equation}
\Moreover \cref{lemma:H:case4sx} and \cref{lemma:H:case4dx} prove that 
for all $\theta \in (\cG^\width)^{-1}(\{0\})$ with $ \prod_{k=1}^\width \v{\theta}_k\neq0<\q{\theta}_1 < \q{\theta}_2 < \ldots < \q{\theta}_\width<1$ and $\{\M^{\theta}_{0}, \M^{\theta}_{1}\}\supseteq\{\emptyset\}$
it holds that 
\begin{equation}\label{lemma:H:active:3}
    \cL^{\width}_\infty(\theta)\geq \frac1{36}.
\end{equation}
\Moreover \cref{lemma:H:updown} shows that there exists $\delta\in (0,\infty)$ such that for all $\theta \in(\cG^\width)^{-1}(\{0\})$ with $\prod_{k=1}^\width \v{\theta}_k\neq0<\q{\theta}_1<\q{\theta}_2 < \ldots < \q{\theta}_\width<1$, $\M^{\theta}_{0}\neq\emptyset\neq\M^{\theta}_{1}$,
and 
$\forall \, k \in \{0,1,\ldots,\width\}\backslash\{\m^{\theta}_{0,2}\} \colon \w{\theta}_k\w{\theta}_{k+1}<0\neq\alpha^{\theta}\beta^{\theta}$  we have that
 \begin{equation}\label{lemma:H:active:4}
   \cL^{\width}_\infty(\theta)\geq\delta.
\end{equation}
\Moreover \cref{lemma:H:N0:finN} establishes that for all $\theta \in (\cG^\width)^{-1}(\{0\})$ with $\prod_{k=1}^\width \v{\theta}_k\neq 0<\q{\theta}_1<\q{\theta}_2 < \ldots < \q{\theta}_\width<1$, $\M^{\theta}_{0}\neq\emptyset\neq\M^{\theta}_{1}$, and
$\{k\in \{0,1,\ldots,\width\}\backslash\{\m^{\theta}_{0,2}\} \colon \w{\theta}_k\w{\theta}_{k+1}>0 \}\neq \emptyset$ it holds that
 \begin{equation}\label{lemma:H:active:5}
   \cL^{\width}_\infty(\theta)\geq\frac{1}{384(1+\width)^2}.
\end{equation}
\Moreover \cref{lemma:H:N0:finN2} proves that for all $\theta \in (\cG^\width)^{-1}(\{0\})$ with $\prod_{k=1}^\width \v{\theta}_k\neq 0<\q{\theta}_1 <\q{\theta}_2< \ldots < \q{\theta}_\width<1$,  $\alpha^{\theta}\beta^{\theta}=0 $, and $\M^{\theta}_{0}\neq\emptyset\neq\M^{\theta}_{1}$ it holds that
 \begin{equation}\label{lemma:H:active:6}
  \cL^{\width}_\infty(\theta)\geq\frac{1}{384(1+\width)^2}.
\end{equation}
Next note that for every $\theta \in  (\cG^\width)^{-1}(\{0\})$ with $R^\theta=\emptyset$ there exists $\vartheta\in(\cG^\width)^{-1}(\{0\})$ which satisfies that $\realization{\width,\vartheta}_\infty=\realization{\width,\vartheta}_\infty$, $R^\vartheta=\emptyset$, and
$\q{\vartheta}_1 \leq \q{\vartheta}_2\leq    \ldots \leq \q{\vartheta}_\width$. 
Combining this, \cref{lemma:H:active:1}, \cref{lemma:H:active:2}, \cref{lemma:H:active:3}, \cref{lemma:H:active:4}, \cref{lemma:H:active:5}, and \cref{lemma:H:active:6}  
establishes that for all $\theta\in (\cG^\width)^{-1}(\{0\})$ with $R^{\theta} = \emptyset$ it holds that
\begin{equation}
    \cL^{\width}_\infty(\theta)\geq\min\left\{\frac{1}{384(1+\width)^2},\delta\right\}.
\end{equation}
\end{cproof}
\subsection{Blow up phenomena for GFs in the training of ANNs}\label{1div}
\begin{prop} \label{theo:intro:convergence}
Let $d, \width, \fd,  n \in \N$, $ \scra \in \R$, $\scrb \in ( \scra, \infty)$ satisfy $\fd = d\width + 2 \width + 1$, let $f \colon [a,b]^d \to \R$ be a function,
for every $i \in \cu{1,2, \ldots, n}$, $k \in \cu{0,1}$ let $\alpha_{i}^k \in \R^{n \times d}$,
let $\beta_{i }^k \in \R^n$, 
and
let $P_i^k \colon \R^d \to \R$ be a polynomial,
let $\dens \colon [\scra , \scrb ] ^d \to [0, \infty)$
satisfy for all $k \in \cu{0,1}$,
$x  \in [\scra , \scrb] ^d$ that
\begin{equation} 
\label{theo:intro:conv:eq0}
    k f ( x ) + ( 1 - k ) \dens ( x ) = \smallsum_{i=1}^n \br*{ P_i^k ( x ) \indicator{[0, \infty )^n} \rbr{ \alpha^k _i x + \beta^k_i } } ,
\end{equation}
let $\act_r \in C ( \R , \R )$, $r \in \N \cup \cu{ \infty } $, satisfy for all $x \in \R$ that $( \bigcup_{r \in \N} \cu{ \act_r } ) \subseteq C^1( \R , \R)$, $\act_\infty ( x ) = \max \cu{ x , 0 }$,
 $\sup_{r \in \N} \sup_{y \in [- \abs{x}, \abs{x} ] } \abs{ ( \act_r)'(y)} < \infty$, and
\begin{equation} \label{theo:intro:conv:eq1}
    \limsup\nolimits_{r \to \infty}  \rbr*{ \abs { \act_r ( x ) - \act _\infty ( x ) } + \abs { (\act_r)' ( x ) - \indicator{(0, \infty)} ( x ) } } = 0,
\end{equation}
let $\cL_r \colon \R^\fd \to \R$, $r \in \N \cup \cu{ \infty }$,
satisfy for all $r \in \N \cup \cu{ \infty }$, $\theta = (\theta_1, \ldots, \theta_\fd) \in \R^{\fd}$ that
\begin{multline} \label{theo:intro:conv:eq2}
   \cL_r ( \theta ) 
   = \int_{[\scra , \scrb]^d} \bigl( f ( x_1, \ldots, x_d )  \\ - \theta_{\fd} - \smallsum_{i=1}^\width \theta_{\width ( d + 1 ) + i } \br[\big]{ \act_r ( \theta_{\width d + i}  + \smallsum_{j=1}^d\theta_{(i-1)d + j } x_j ) } \bigr)^2  \dens ( x ) \,\d (x_1, \ldots, x_d ) ,
    \end{multline}
let $\cG  \colon \R^\fd \to \R^\fd$ satisfy for all
$\theta \in \cu{ \vartheta \in \R^\fd \colon ( ( \nabla \cL_r ) ( \vartheta ) ) _{r \in \N} 
\text{ is convergent} }$
that $\cG ( \theta ) = \lim_{r \to \infty} \allowbreak (\nabla \cL_r) ( \theta )$,
and let $\Theta \in C ( [ 0 , \infty ) , \R^\fd )$ satisfy $\liminf_{t \to  \infty } \norm{\Theta_t } < \infty$ and $\forall \, t \in [0, \infty ) \colon \Theta_t = \Theta_0 - \int_0^t \cG ( \Theta_s ) \, \d s$.
Then there exist $\vartheta \in \cG^{ - 1 } ( \cu{  0 } )$,
$ \fC,  \beta \in (0, \infty)$ which satisfy for all $t \in [ 0 , \infty )$ that
\begin{equation} \label{theo:intro:conv:eq3}
    \norm{\Theta_t - \vartheta} \leq  \fC ( 1 + t ) ^{- \beta } 
\qandq
  \abs{ \cL_\infty ( \Theta_t ) - \cL_\infty ( \vartheta ) } \leq \fC ( 1 + t ) ^{-1} .
\end{equation}
\end{prop}
\begin{cproof}{theo:intro:convergence}
The assertion is verified analogously to the proof of \cite[Theorem~1.3]{JentzenRiekert2021aa} 
(compare, e.g., \cite[Theorem~1.2]{Jentzen2}).
\end{cproof}
\begin{lemma}\label{lemma:H:bound}
Assume \cref{H:setting:snn}.
 Then it holds for all $\width\in \N$ that there exists $\varepsilon \in (0,\infty)$ such that for all $\theta\in (\cG^\width)^{-1}(\{0\})$ it holds that
 \begin{equation}\label{lemma:H:bound:thesis}
     \cL^\width_\infty(\theta)\geq \varepsilon.
 \end{equation}
\end{lemma}
\begin{cproof}{lemma:H:bound}
We prove \cref{lemma:H:bound:thesis} by induction. \Nobs that \cref{lemma:H:1}  assures that  for all $\theta \in (\cG^1)^{-1}(\{0\})$ it holds that \begin{equation}
    \cL^1_\infty(\theta)\geq \frac1{36}.
\end{equation}
For the induction step let $\width \in \N \cap (1,\infty)$ and assume that there exists $\epsilon \in (0,\infty)$ which satisfies for all $\vartheta\in (\cG^{\width-1})^{-1}(\{0\})$ that
 \begin{equation}\label{lemma:H:bound:induction}
     \cL^{\width-1}_\infty(\vartheta)\geq \epsilon.
 \end{equation}
\Nobs that \cref{lemma:H:bound:induction} shows that for all $\theta\in (\cG^\width)^{-1}(\{0\})$ 
 with
$\{\vartheta\in (\cG^{\width-1})^{-1}(\{0\}) \colon {\realization{\width,\theta}_\infty}| _{[0,1]}={\realization{\width-1,\vartheta}_\infty}| _{[0,1]}\} \neq \emptyset$
 there exists $\vartheta\in (\cG^{\width-1})^{-1}(\{0\})$ such that
  \begin{equation}\label{lemma:H:bound:hp}
     \cL_\infty^\width(\theta)=\cL_\infty^{\width-1}(\vartheta)\geq \epsilon.
 \end{equation}
\Nobs that \cref{lemma:H:active} demonstrates that there exists $\delta \in (0,\infty)$ which satisfies for all $\theta\in (\cG^\width)^{-1}(\{0\})$ with
$\{\vartheta\in (\cG^{\width-1})^{-1}(\{0\}) \colon {\realization{\width,\theta}_\infty}| _{[0,1]}={\realization{\width-1,\vartheta}_\infty}| _{[0,1]}\}= \emptyset$ that 
\begin{equation}\label{lemma:H:bound:hp3}
    \cL^\width_\infty(\theta)\geq \delta.
\end{equation}
 \Nobs that \cref{lemma:H:bound:hp} and \cref{lemma:H:bound:hp3} ensure that for all $\theta\in (\cG^\width)^{-1}(\{0\})$ it holds that $\cL^\width_\infty(\theta)\geq\min\{\epsilon,\delta\}$.
Induction thus establishes \cref{lemma:H:bound:thesis}.
\end{cproof}
\begin{theorem}\label{H:theorem:zerograd}
Let $\scra \in \R$, $\scrb \in (\scra,\infty)$, $ \width, \fd \in \N$ satisfy $\fd=3\width+1$,
let $\act_r \in C ( \R , \R )$, $r \in \N \cup \cu{ \infty } $, satisfy for all $x \in \R$ that $( \bigcup_{r \in \N} \cu{ \act_r } ) \subseteq C^1( \R , \R)$, $\act_\infty ( x ) = \max \cu{ x , 0 }$,
 $\sup_{r \in \N} \sup_{y \in [- \abs{x}, \abs{x} ] } \abs{ ( \act_r)'(y)} < \infty$, and
\begin{equation} 
    \limsup\nolimits_{r \to \infty}  \rbr*{ \abs { \act_r ( x ) - \act _\infty ( x ) } + \abs { (\act_r)' ( x ) - \indicator{(0, \infty)} ( x ) } } = 0,
\end{equation}
for every $r \in \N \cup \{\infty\}$, $\theta= (\theta_1, \ldots, \theta_\fd) \in \R^{\fd}$
let $\realization{\theta}_r \colon \R \to \R$
satisfy for all $x \in \R$ that
\begin{equation}
    \realization{\theta}_r (x) =  \theta_{\fd} + \smallsum_{i=1}^\width \theta_{2\width +i } \br[\big]{ \act_r ( \theta_{\width  + i}  +\theta_{i} x )} ,
\end{equation}
for every $r \in \N \cup \{\infty\}$ let $\cL_r \colon \R^\fd \to \R$
satisfy for all  $\theta  \in \R^{\fd}$ that
\begin{equation}
   \cL_r ( \theta ) 
   = \int_\scra^\scrb \bigl( \mathbbm{1}_{(\nicefrac{(\scra+\scrb)}2,\infty)}(x)  -\realization{\theta}_r (x)  \bigr)^2   \,\d x ,
   \end{equation}
and let $\cG  \colon \R^\fd \to \R^\fd$ satisfy for all
$\theta \in \cu{ \vartheta \in \R^\fd \colon ( ( \nabla \cL_r ) ( \vartheta ) ) _{r \in \N} \allowbreak
\text{ is convergent} }$
that $\cG ( \theta ) = \lim_{r \to \infty} \allowbreak (\nabla \cL_r) ( \theta )$.
 Then there exists $\varepsilon \in (0,\infty)$ such that for all $\theta\in \cG^{-1}(\{0\})$ it holds that $\cL_\infty(\theta)\geq\varepsilon$.
\end{theorem}
\begin{cproof}{H:theorem:zerograd}
Throughout this proof for every $\theta\in \R^{\fd}$ let $v^\theta \in \R^{\fd}$ satisfy for all $i \in \{1,2,\ldots,\width\}$ that
$v_i^{\theta}=\theta_i (\scrb-\scra)$, $v_{\width+i}^{\theta}=\theta_{\width+i}+\theta_{i}\scra$, $v_{2\width+i}^{\theta}=\theta_{2\width+i}$, and $v_{\fd}^{\theta}=\theta_{\fd}$ and for every $\theta \in \R^{\fd}$, $i\in \{1, 2,\ldots, \width \}$, $\scrc \in \R$, $\scrd \in (\scrc,\infty)$ 
let $I^{\theta,\scrc}_{i,\scrd} \subseteq \R$  satisfy  $I^{\theta,\scrc}_{i,\scrd} = \cu{ x \in [\scrc , \scrd ] \colon \theta_i x + \theta_{\width+i}  > 0 }$.
\Nobs that for all $\theta \in \R^{\fd}$ it holds that
\begin{equation}\label{H:theorem:zerograd:risk}
\begin{split}
     \cL_\infty(\theta) &=\int_\scra^\scrb \bigl( \mathbbm{1}_{(\nicefrac{(\scra+\scrb)}2,\infty)}(x) -\realization{\theta}_\infty (x)  \bigr)^2   \,\d x \\
    &=(\scrb-\scra) \int_0^1 \bigl( \mathbbm{1}_{(\nicefrac{1}2,\infty)}(y) -\realization{\theta}_\infty ((\scrb-\scra)y+\scra)  \bigr)^2   \,\d y \\
    &= (\scrb-\scra)\int_0^1 \bigl( \mathbbm{1}_{(\nicefrac{1}2,\infty)}(y) -\realization{v^\theta}_\infty (y)\bigr)^2   \,\d y.
\end{split}
\end{equation}
\Moreover \cite[Proposition 2.2]{Jentzen} establishes that for all $\theta\in \cG^{-1}(\{0\})$, $i \in \{1,2,\ldots,\width\}$ it holds that
\begin{equation} \label{H:theorem:zerograd:gradient}
\begin{split}
        0=\cG_{ i} ( \theta) &= 2 \theta_{2\width+i} \int_{I_{i,\scrb}^{\theta,\scra}} x  ( \realization{\theta}_\infty (x) - \mathbbm{1}_{(\nicefrac{(\scra+\scrb)}2,\infty)}(x)) \, \d x \\
        &= 2 v^\theta_{2\width+i} (\scrb-\scra) \int_{I_{i,1}^{v^\theta,0}} ((\scrb-\scra)y+\scra)  ( \realization{v^\theta}_\infty (y) - \mathbbm{1}_{(\nicefrac12,\infty)}(y)) \, \d y, \\
        0=\cG_{\width + i} ( \theta) &= 2 \theta_{2\width+i} \int_{I_{i,\scrb}^{\theta,\scra}} (\realization{\theta}_\infty (x) - \mathbbm{1}_{(\nicefrac{(\scra+\scrb)}2,\infty)}(x)) \,  \d x  \\
        &= 2 v^\theta_{2\width+i} (\scrb-\scra) \int_{I_{i,1}^{v^\theta,0}} (\realization{v^\theta}_\infty (y) - \mathbbm{1}_{(\nicefrac12,\infty)}(y)) \,  \d y ,\\
        0=\cG_{2\width + i} ( \theta) &= 2 \int_{\scra}^\scrb \br[\big]{\act_\infty \rbr{ \theta_{i} x + \theta_{\width+i} } } ( \realization{\theta}_\infty (x) - \mathbbm{1}_{(\nicefrac{(\scra+\scrb)}2,\infty)}(x) ) \, \d x\\
         &= 2  (\scrb-\scra)\int_{0}^1 \br[\big]{\act_\infty \rbr{ v^\theta_{i}  y+ v^\theta_{\width+i} } } ( \realization{v^\theta}_\infty (y) - \mathbbm{1}_{(\nicefrac{1}2,\infty)}(y) ) \, \d y, \\
        \text{and} \qquad 0=\cG_{\fd} ( \theta) &= 2 \int_{\scra}^\scrb (\realization{\theta}_\infty (x) - \mathbbm{1}_{(\nicefrac{(\scra+\scrb)}2,\infty)}(x) ) \, \d x\\
        &= 2  (\scrb-\scra)\int_{0}^1 (\realization{v^\theta}_\infty (y) - \mathbbm{1}_{(\nicefrac12,\infty)}(y) ) \, \d y.
        \end{split}
\end{equation}
Combining this and  \cref{H:theorem:zerograd:risk} with \cref{lemma:H:bound} demonstrates that there exists $\varepsilon \in (0,\infty)$ such that for all $\theta\in \cG^{-1}(\{0\})$ it holds that 
\begin{equation}\label{H:theorem:zerograd:eps}
         \cL_\infty(\theta) = (\scrb-\scra)\int_0^1 \bigl( \mathbbm{1}_{(\nicefrac{1}2,\infty)}(y) -\realization{v^\theta}_\infty (y ) \bigr) ^2  \,\d y \geq (\scrb-\scra) \varepsilon.
\end{equation}
\end{cproof}
\begin{theorem}\label{H:theorem:nozero}
Let $\scra \in \R$, $\scrb \in (\scra,\infty)$, $ \width, \fd \in \N$ satisfy $\fd=3\width+1$,
let $\act_r \in C ( \R , \R )$, $r \in \N \cup \cu{ \infty } $, satisfy for all $x \in \R$ that $( \bigcup_{r \in \N} \cu{ \act_r } ) \subseteq C^1( \R , \R)$, $\act_\infty ( x ) = \max \cu{ x , 0 }$,
 $\sup_{r \in \N} \sup_{y \in [- \abs{x}, \abs{x} ] } \abs{ ( \act_r)'(y)} < \infty$, and
\begin{equation} 
    \limsup\nolimits_{r \to \infty}  \rbr*{ \abs { \act_r ( x ) - \act _\infty ( x ) } + \abs { (\act_r)' ( x ) - \indicator{(0, \infty)} ( x ) } } = 0,
\end{equation}
for every $r \in \N \cup \{\infty\}$, $\theta= (\theta_1, \ldots, \theta_\fd) \in \R^{\fd}$
let $\realization{\theta}_r \colon \R \to \R$
satisfy for all $x \in \R$ that
\begin{equation}
    \realization{\theta}_r (x) =  \theta_{\fd} + \smallsum_{i=1}^\width \theta_{2\width +i } \br[\big]{ \act_r ( \theta_{\width  + i}  +\theta_{i} x )} ,
\end{equation}
for every $r \in \N \cup \{\infty\}$ let $\cL_r \colon \R^\fd \to \R$
satisfy for all  $\theta  \in \R^{\fd}$ that
\begin{equation}
   \cL_r ( \theta ) 
   = \int_\scra^\scrb \bigl( \mathbbm{1}_{(\nicefrac{(\scra+\scrb)}2,\infty)}(x)  -\realization{\theta}_r (x)  \bigr)^2   \,\d x ,
   \end{equation}
and let $\cG  \colon \R^\fd \to \R^\fd$ satisfy for all
$\theta \in \cu{ \vartheta \in \R^\fd \colon ( ( \nabla \cL_r ) ( \vartheta ) ) _{r \in \N} \allowbreak
\text{ is convergent} }$
that $\cG ( \theta ) = \lim_{r \to \infty} \allowbreak (\nabla \cL_r) ( \theta )$.
 Then there exists $\varepsilon \in (0,\infty)$ such that for all
$\Theta \in C([0, \infty) , \R^{\fd})$  with  $\forall \, t \in [0, \infty) \colon  \Theta_t = \Theta_0 - \int_0^t \cG ( \Theta_s ) \, \d s$ and $\cL_\infty(\Theta_0)<\varepsilon$
it holds that
$\liminf_{t \to\infty} \norm{\Theta_t } = \infty.$
\end{theorem}
\begin{cproof}{H:theorem:nozero}
\Nobs that \cref{H:theorem:zerograd} ensures that there exists $\varepsilon \in (0,\infty)$ which satisfies for all $\theta\in \cG^{-1}(\{0\})$ that 
\begin{equation}\label{cor:H:theorem:nozero:eps}
    \cL_\infty(\theta)\geq\varepsilon.
\end{equation}
Furthermore, observe that, e.g., \cite[Lemma 3.1]{Jentzen}  implies that 
for all
$\Theta \in C([0, \infty) , \R^{\fd})$  with  $\forall \, t \in [0, \infty) \colon  \Theta_t = \Theta_0 - \int_0^t \cG ( \Theta_s ) \, \d s$ it holds that
$[0, \infty) \ni t \mapsto \cL_\infty ( \Theta_t ) \in \R$ is non-increasing.
Combining this with \cref{cor:H:theorem:nozero:eps} and \cref{theo:intro:convergence} assures that for all
$\Theta \in C([0, \infty) , \R^{\fd})$  with  $\forall \, t \in [0, \infty) \colon  \Theta_t = \Theta_0 - \int_0^t \cG ( \Theta_s ) \, \d s$ and $\cL_\infty(\Theta_0)< \varepsilon$
it holds that
$\liminf_{t \to\infty} \norm{\Theta_t } = \infty.$
\end{cproof}
\subsection{Blow up phenomena for GFs in the training of ANNs with two hidden neurons}\label{subsectionestimate}
\begin{lemma} \label{lemma:2:zerogradient:bw}
Assume \cref{H:setting:snn} and let $\theta\in \R^{7}$ satisfy $\w{\theta}_1\w{\theta}_2>0\neq\v{\theta}_1\v{\theta}_2$,
$0<\q{\theta}_1 < \q{\theta}_2<1$, and $\nicefrac12 \in (\q{\theta}_1 , \q{\theta}_2)$. Then  \begin{equation}\label{lemma:2:zerogradient:bw:thesis}
    \cG^2(\theta)\neq0.
\end{equation}
\end{lemma}
\begin{cproof}{lemma:2:zerogradient:bw}
We prove \cref{lemma:2:zerogradient:bw:thesis} by contradiction. Assume that  $\cG^2(\theta)=0$.
In the following we distinguish between the case $\min\{\w{\theta}_1,\w{\theta}_2\}>0$ and the case $\max\{\w{\theta}_1,\w{\theta}_2\}<0$.
We first establish the contradiction in the case
\begin{equation}\label{lemma:2:zerogradient:bw:case1}
    \min\{\w{\theta}_1,\w{\theta}_2\}>0.
\end{equation}
\Nobs that \cref{lemma:2:zerogradient:bw:case1}, the fact that $\cG^2(\theta)=0$, and \cref{prop:H:approximate:gradient} imply that for all $i \in \{1,2\}$ it holds that
\begin{equation}
    \int_{\q{\theta}_{i}}^{\q{\theta}_{i+1}}   x(  \realization{2,\theta}_\infty (x) - \mathbbm{1}_{(\nicefrac12,\infty)}(x)) \, \d x=
    \int_{\q{\theta}_{i}}^{\q{\theta}_{i+1}}   (  \realization{2,\theta}_\infty (x) - \mathbbm{1}_{(\nicefrac12,\infty)}(x))\, \d x=0.
\end{equation}
This, \cref{lem:affine:integral:zero}, and \cref{cor:affine1:integral:zero2} demonstrate that for all $x \in [\q{\theta}_1,1]$ it holds that $\q{\theta}_2=\nicefrac34 -\nicefrac12 \q{\theta}_1$ and 
\begin{equation}
    \realization{2,\theta}_\infty(x)=\begin{cases}
    1 
    &  \colon x\in[\q{\theta}_2,1]\\
    -\frac{16x}{9(2\q{\theta}_1-1)}+\frac{10\q{\theta}_1+3}{9(2\q{\theta}_1-1)} 
    & \colon x\in[\q{\theta}_1,\q{\theta}_2).
    \end{cases}
\end{equation}
Combining this with the fact that $I^{\theta}_1=(\q{\theta}_1,1]$, the fact that $I^{\theta}_2=(\q{\theta}_2,1]$, and continuity of $\realization{2,\theta}_\infty$ shows that for all $x \in [0,\q{\theta}_1)$ it holds that
\begin{equation}
\realization{2,\theta}_\infty(x)=-\frac13.
\end{equation}
This and \cref{prop:H:approximate:gradient} imply that
\begin{equation}
     0=\int_{0}^{1}   (  \realization{2,\theta}_\infty (x) - \mathbbm{1}_{(\nicefrac12,\infty)}(x)) \, \d x=
      \int_{0}^{\q{\theta}_{1}}   (  \realization{2,\theta}_\infty (x)) \, \d x=
      \int_{0}^{\q{\theta}_{1}}   -\frac13 \, \d x= -\frac{\q{\theta}_{1}}3 .
\end{equation}
This is a contradiction.
In the next step we establish the contradiction in the case 
\begin{equation}\label{lemma:2:zerogradient:bw:case2}
    \max\{\w{\theta}_1,\w{\theta}_2\}<0.
\end{equation} 
\Nobs that \cref{lemma:2:zerogradient:bw:case2}, the fact that $\cG^2(\theta)=0$, and \cref{prop:H:approximate:gradient} imply that for all $i \in \{0,1\}$ it holds that
\begin{equation}
    \int_{\q{\theta}_{i}}^{\q{\theta}_{i+1}}   x(  \realization{2,\theta}_\infty (x) - \mathbbm{1}_{(\nicefrac12,\infty)}(x)) \, \d x=
    \int_{\q{\theta}_{i}}^{\q{\theta}_{i+1}}   (  \realization{2,\theta}_\infty (x) - \mathbbm{1}_{(\nicefrac12,\infty)}(x))\, \d x=0.
\end{equation}
This, \cref{lem:affine:integral:zero}, and \cref{cor:affine1:integral:zero} demonstrate that for all $x \in [0,\q{\theta}_2]$ it holds that $\q{\theta}_1=\nicefrac34 -\nicefrac12 \q{\theta}_2$ and 
\begin{equation}
    \realization{2,\theta}_\infty(x)=\begin{cases}
    0
    & \colon x\in[0,\q{\theta}_1]\\
    \frac{16x}{9(2\q{\theta}_2-1)}+\frac{4(2\q{\theta}_2-3)}{9(2\q{\theta}_2-1)}
    & \colon x\in(\q{\theta}_1,\q{\theta}_2].
    \end{cases}
\end{equation}
Combining this with the fact that $I^{\theta}_1=[0,\q{\theta}_1)$, the fact that $I^{\theta}_2=[0,\q{\theta}_2)$, and the fact that $\realization{2,\theta}_\infty$ is continuous shows that for all $x \in (\q{\theta}_2,1]$ it holds that
\begin{equation}
\realization{2,\theta}_\infty(x)=\frac43.
\end{equation}
This and \cref{prop:H:approximate:gradient} imply that
\begin{equation}
     0= \int_{0}^{1}   (  \realization{2,\theta}_\infty (x) - \mathbbm{1}_{(\nicefrac12,\infty)}(x)) \, \d x=
      \int_{\q{\theta}_{2}}^1   (  \realization{2,\theta}_\infty (x) - 1) \, \d x=
       \int_{\q{\theta}_{2}}^1 \frac13 \, \d x= \frac13 (1-\q{\theta}_{2}).
\end{equation}
This is a contradiction.
\end{cproof}
\begin{lemma} \label{lemma:2:zerogradient:bwmist}
Assume \cref{H:setting:snn} and let $\theta\in \R^{7}$ satisfy $\w{\theta}_1<0<\w{\theta}_2$,
$\v{\theta}_1\v{\theta}_2\neq0<\q{\theta}_1 < \q{\theta}_2<1$, and $\nicefrac12 \in (\q{\theta}_1 , \q{\theta}_2)$. Then  \begin{equation}\label{lemma:2:zerogradient:bwmist:thesis}
    \cG^2(\theta)\neq0.
\end{equation}
\end{lemma}
\begin{cproof}{lemma:2:zerogradient:bwmist}
We prove \cref{lemma:2:zerogradient:bwmist:thesis} by contradiction. Assume that $\cG^2(\theta)=0$.
This, the assumption that $\w{\theta}_1<0<\w{\theta}_2$, and \cref{prop:H:approximate:gradient} imply that for all $i \in \{0,2\}$ it holds that
\begin{equation}
    \int_{\q{\theta}_{i}}^{\q{\theta}_{i+1}}   x(  \realization{2,\theta}_\infty (x) - \mathbbm{1}_{(\nicefrac12,\infty)}(x)) \, \d x=
    \int_{\q{\theta}_{i}}^{\q{\theta}_{i+1}}   (  \realization{2,\theta}_\infty (x) - \mathbbm{1}_{(\nicefrac12,\infty)}(x))\, \d x=0.
\end{equation}
This and \cref{lem:affine:integral:zero} demonstrate that 
for all $x \in [0,1]\backslash(\q{\theta}_1,\q{\theta}_2)$ it holds that
\begin{equation}
    \realization{2,\theta}_\infty(x)=\begin{cases}
    0
    & \colon x\in[0,\q{\theta}_1]\\
    1
    & \colon x\in[\q{\theta}_2,1].
    \end{cases}
\end{equation}
Combining this with the fact that $I^{\theta}_1=[0,\q{\theta}_1)$, the fact that $I^{\theta}_2=(\q{\theta}_2,1]$, and the fact that $\realization{2,\theta}_\infty$ is continuous  shows that for all $x \in (\q{\theta}_1,\q{\theta}_2)$ it holds that
$\realization{2,\theta}_\infty(x)=0.$
This is a contradiction.
\end{cproof}
\begin{lemma} \label{lemma:2:zerogradient:bwfin}
Assume \cref{H:setting:snn} and let $\theta\in (\cG^2)^{-1}(\{0\})$ satisfy $\w{\theta}_2<0<\w{\theta}_1$,
$\v{\theta}_1\v{\theta}_2\neq0<\q{\theta}_1 < \q{\theta}_2<1$, and $\nicefrac12 \in (\q{\theta}_1 , \q{\theta}_2)$. Then  \begin{equation}\label{lemma:2:zerogradient:bwfin:thesis}
    \cL^2_\infty(\theta)\geq \frac1{864}.
\end{equation}
\end{lemma}
\begin{cproof}{lemma:2:zerogradient:bwfin}
\Nobs that \cref{lemma:H:updown:sx},  \cref{lemma:H:updown:dx}, and the fact that $\realization{\theta}_\infty$ is continuous  imply that 
\begin{equation}\label{lemma:2:zerogradient:bwfin:real}
    \realization{2,\theta}_\infty(x)=\begin{cases}
    \alpha^{\theta} x -\frac{ \alpha^{\theta} \q{\theta}_1}2 
    & \colon x\in[0,\q{\theta}_1]\\
     (\alpha^{\theta}+\beta^{\theta})(x-\q{\theta}_1)+\frac{ \alpha^{\theta} \q{\theta}_1}2 
    & \colon x\in(\q{\theta}_1,\q{\theta}_2]\\
     \beta^{\theta} x +1-\frac{\beta^{\theta}}2   (1+\q{\theta}_2)
    & \colon x\in(\q{\theta}_2,1]
    \end{cases}
\end{equation}
and \begin{equation}\label{lemma:2:zerogradient:bwfin:cont}
      (\alpha^{\theta}+\beta^{\theta})(\q{\theta}_2-\q{\theta}_1)+\frac{ \alpha^{\theta} \q{\theta}_1}2 =  1-\frac{\beta^{\theta}}2   (1-\q{\theta}_2).
\end{equation}
\Moreover the assumption that $\w{\theta}_2<0<\w{\theta}_1$ and \cref{prop:H:approximate:gradient} show that for all $i \in \{0,1\}$ it holds that
\begin{equation}
    \int_{\q{\theta}_{i}}^{\q{\theta}_{i+2}}   x(  \realization{2,\theta}_\infty (x) - \mathbbm{1}_{(\nicefrac12,\infty)}(x)) \, \d x=
    \int_{\q{\theta}_{i}}^{\q{\theta}_{i+2}}   (  \realization{2,\theta}_\infty (x) - \mathbbm{1}_{(\nicefrac12,\infty)}(x))\, \d x=0.
\end{equation}
Combining this and \cref{lemma:2:zerogradient:bwfin:real} ensures that 
\begin{equation}
\begin{split}
     0 &= \int_{0}^{\q{\theta}_{1}}   x(  \realization{2,\theta}_\infty (x) - \mathbbm{1}_{(\nicefrac12,\infty)}(x)) \, \d x-
     \int_{\q{\theta}_{2}}^{1}   x(  \realization{2,\theta}_\infty (x) - \mathbbm{1}_{(\nicefrac12,\infty)}(x)) \, \d x \\
     & =\frac{\alpha^{\theta}}{12} (\q{\theta}_1)^3-\frac{\beta^{\theta}}{12} (1-\q{\theta}_2)^3.
\end{split}
\end{equation}
This establishes that $\alpha^{\theta}\beta^{\theta}>0$. Combining this and \cref{lemma:2:zerogradient:bwfin:cont} proves that $\min\{\alpha^{\theta}, \beta^{\theta}\}>0$. This and \cref{lemma:2:zerogradient:bwfin:real} imply that for all $x\in[0,1]$ it holds that $\realization{2,\theta}_\infty(x)\leq (\alpha^{\theta}+\beta^{\theta}) x$. \Hence that
\begin{equation}\label{lemma:2:zerogradient:bwfin:ab}
    1<\realization{2,\theta}_\infty (1)\leq \alpha^{\theta}+\beta^{\theta}.
\end{equation}
This shows that $\max\{\alpha^{\theta},\beta^{\theta}\}\geq\nicefrac12$. In the following we distinguish between the case $\alpha^{\theta}\geq\nicefrac12$ and the case $\beta^{\theta}\geq\nicefrac12$.
We first prove \cref{lemma:2:zerogradient:bwfin:thesis} in the case \begin{equation}\label{lemma:2:zerogradient:bwfin:case1}
    \alpha^{\theta}\geq\frac12.
\end{equation} 
\Nobs that \cref{lemma:intmin},  \cref{lemma:2:zerogradient:bwfin:real}, \cref{lemma:2:zerogradient:bwfin:ab}, and \cref{lemma:2:zerogradient:bwfin:case1} demonstrate that
\begin{equation}\label{lemma:2:zerogradient:bwfin:est}
\begin{split}
    \cL^2_\infty(\theta) &\geq \int_{0}^{\q{\theta}_{1}}   (  \realization{2,\theta}_\infty (x) )^2 \, \d x
    +\int_{\q{\theta}_{1}}^{\frac12}   (  \realization{2,\theta}_\infty (x))^2 \, \d x\\
    &\geq\frac1{12}(\alpha^{\theta})^2 (\q{\theta}_1)^3 + \frac1{12}(\alpha^{\theta}+\beta^{\theta})^2\left(\frac12-\q{\theta}_1\right)^3\geq \frac 1{48} (\q{\theta}_1)^3+ \frac1{12}\left(\frac12-\q{\theta}_1\right)^3\\
    &\geq \frac{1}{48}\left(\frac13\right)^3+ \frac1{12}\left(\frac12-\frac13\right)^3=\frac{1}{864}.
\end{split}
\end{equation}
This proves \cref{lemma:2:zerogradient:bwfin:thesis} in the case $\alpha^{\theta}\geq\nicefrac12$.
Next we establish \cref{lemma:2:zerogradient:bwfin:thesis} in the case \begin{equation}\label{lemma:2:zerogradient:bwfin:case2}
    \beta^{\theta}\geq\frac12.
\end{equation}  
\Nobs that \cref{lemma:intmin},  \cref{lemma:2:zerogradient:bwfin:real}, \cref{lemma:2:zerogradient:bwfin:ab}, and \cref{lemma:2:zerogradient:bwfin:case2} assure that
\begin{equation}
\begin{split}
    \cL^2_\infty(\theta) &\geq \int_{\frac12}^{\q{\theta}_{2}}   (  \realization{2,\theta}_\infty (x) -1)^2 \, \d x
    +\int_{\q{\theta}_{2}}^{1}   (  \realization{2,\theta}_\infty (x)-1)^2 \, \d x\\
    &   \geq \frac1{12}(\alpha^{\theta}+\beta^{\theta})^2\left(\q{\theta}_2-\frac12\right)^3+\frac1{12}(\beta^{\theta})^2 (1-\q{\theta}_2)^3 \\
    & \geq \frac1{12}\left(\q{\theta}_2-\frac12\right)^3+\frac 1{48} (1-\q{\theta}_2)^3\geq\frac1{12}\left(\frac23-\frac12\right)^3+\frac 1{48} \left(1-\frac23\right)^3= \frac{1}{864}.
\end{split}
\end{equation}
This demonstrates \cref{lemma:2:zerogradient:bwfin:thesis} in the case $\beta^{\theta}\geq\nicefrac12$.
\end{cproof}
\begin{lemma}\label{lemma:2}
Assume \cref{H:setting:snn}.
 Then it holds for all $\theta\in (\cG^2)^{-1}(\{0\})$ that
 \begin{equation}\label{lemma:2:bound:thesis}
     \cL^2_\infty(\theta)\geq \frac{1}{864}.
 \end{equation}
\end{lemma}
\begin{cproof}{lemma:2}
\Nobs that for all $\theta\in (\cG^2)^{-1}(\{0\})$ with $\{k \in \{1,2\}\colon I^{\theta}_{k}= \emptyset\}\neq\emptyset$  there exists $\vartheta \in (\cG^{1})^{-1}(\{0\})$  such that for all $x \in [0,1]$ it holds that $\realization{2,\theta}_\infty(x)=\realization{1,\vartheta}_\infty(x)$. Combining this and \cref{lemma:H:1} shows that for all $\theta\in (\cG^2)^{-1}(\{0\})$ with $\{k \in \{1,2\}\colon I^{\theta}_{k}= \emptyset\}\neq\emptyset$  there exists $\vartheta \in (\cG^{1})^{-1}(\{0\})$  such that
 \begin{equation}\label{lemma:2:active:1}
     \cL_\infty^2(\theta)=\cL_\infty^{1}(\vartheta)\geq \frac1{36}.
 \end{equation}
 \Moreover for all $\theta\in (\cG^2)^{-1}(\{0\})$ with
 $\v{\theta}_1\v{\theta}_2=0$ there exists $\vartheta \in (\cG^{1})^{-1}(\{0\})$  such that for all $x \in [0,1]$ it holds that $\realization{2,\theta}_\infty(x)=\realization{1,\vartheta}_\infty(x)$. Combining this and \cref{lemma:H:1} establishes that for all $\theta\in (\cG^2)^{-1}(\{0\})$ with
 $\v{\theta}_1\v{\theta}_2=0$ there exists $\vartheta \in (\cG^{1})^{-1}(\{0\})$  such that
 \begin{equation}\label{lemma:2:active:2}
     \cL_\infty^2(\theta)=\cL_\infty^{1}(\vartheta)\geq \frac1{36}.
 \end{equation}
Furthermore, \nobs that \cref{lemma:H:q0} ensures that
for all $\theta\in (\cG^2)^{-1}(\{0\})$ with 
$\{k\in \{1,2\} \colon \allowbreak(0,1)\subseteq I_k^{2,\theta} \}\neq \emptyset$
it holds that
\begin{equation}\label{lemma:2:active:2a}
     \cL_\infty^2(\theta)\geq \frac1{36}.
 \end{equation}
Next, \nobs that \cref{lemma:H:qq} and \cref{lemma:H:qqneg} assure that for all $\theta \in (\cG^2)^{-1}(\{0\})$ with  $\v{\theta}_1\v{\theta}_2\neq0$, $\q{\theta}_1=\q{\theta}_2 \in (0,1)$ it holds that 
\begin{equation}\label{lemma:2:active:3}
   \cL^{2}_\infty(\theta)\geq \frac1{36}.
\end{equation}
In addition, \nobs that \cref{lemma:H:case4sx} and \cref{lemma:H:case4dx} demonstrate that 
for all $\theta \in (\cG^2)^{-1}(\{0\})$ with $ \v{\theta}_1\v{\theta}_2\neq0<\q{\theta}_1 < \q{\theta}_2<1$ and $\nicefrac12\notin(\q{\theta}_1,   \q{\theta}_2) $
it holds that 
\begin{equation}\label{lemma:2:active:4}
    \cL^{2}_\infty(\theta)\geq \frac1{36}.
\end{equation}
Moreover,
\nobs that \cref{lemma:2:zerogradient:bw}, \cref{lemma:2:zerogradient:bwmist}, and \cref{lemma:2:zerogradient:bwfin} prove that 
for all $\theta \in (\cG^2)^{-1}(\{0\})$ with $ \v{\theta}_1\v{\theta}_2\neq0<\q{\theta}_1 < \q{\theta}_2<1$ and $\nicefrac12\in(\q{\theta}_1,   \q{\theta}_2) $
it holds that 
\begin{equation}\label{lemma:2:active:5}
    \cL^{2}_\infty(\theta)\geq \frac1{864}.
\end{equation}
Combining this, \cref{lemma:2:active:1}, \cref{lemma:2:active:2}, \cref{lemma:2:active:2a},  \cref{lemma:2:active:3}, and \cref{lemma:2:active:4} implies that for all $\theta\in (\cG^2)^{-1}(\{0\})$ it holds that
 \begin{equation}
     \cL^2_\infty(\theta)\geq \frac{1}{864}.
 \end{equation}
\end{cproof}
\begin{lemma}\label{lemma:case2}
Let $\act_r \in C ( \R , \R )$, $r \in \N \cup \cu{ \infty } $, satisfy for all $x \in \R$ that $( \bigcup_{r \in \N} \cu{ \act_r } ) \subseteq C^1( \R , \R)$, $\act_\infty ( x ) = \max \cu{ x , 0 }$, $\sup_{r \in \N} \sup_{y \in [- \abs{x}, \abs{x} ] }  \abs{(\act_r)'(y)} < \infty$,
 and
\begin{equation}
    \limsup\nolimits_{r \to \infty}  \rbr*{ \abs { \act_r ( x ) - \act _\infty ( x ) } + \abs { (\act_r)' ( x ) - \indicator{(0, \infty)} ( x ) } } = 0,
\end{equation}
let $\cL_r \colon \R^7 \to \R$, $r \in \N \cup \cu{ \infty }$,
satisfy for all $r \in \N \cup \cu{ \infty }$, $\theta = (\theta_1, \ldots, \theta_7) \in \R^{7}$ that
\begin{equation}
    \cL_r ( \theta ) = \int_{0}^1 \rbr[\big]{  \mathbbm{1}_{(\nicefrac12,\infty)}(x) - \theta_{7} - \smallsum_{i=1}^2 \theta_{4 + i } \br[\big]{ \act_r ( \theta_{2  + i}  +\theta_{i } x ) } }^2  \, \d x,
\end{equation}
and let $\cG \colon \R^7 \to \R^7$ satisfy for all
$\theta \in \cu{ \vartheta \in \R^7 \colon ( ( \nabla \cL_r ) ( \vartheta ) ) _{r \in \N} \text{ is convergent} }$
that $\cG ( \theta ) = \lim_{r \to \infty} (\nabla \cL_r) ( \theta )$. 
Then it holds for all
$\Theta \in C([0, \infty) , \R^{7})$  with  $\forall \, t \in [0, \infty) \colon  \Theta_t = \Theta_0 - \int_0^t \cG ( \Theta_s ) \, \d s$ and $\cL_\infty(\Theta_0)<\nicefrac1{864}$ that
$\liminf_{t \to \infty} \norm{\Theta_t } = \infty.$
\end{lemma}
\begin{cproof}{lemma:case2}
\Nobs that  \cref{lemma:2} demonstrates that for all $\theta\in \cG^{-1}(\{0\})$ it holds that 
\begin{equation}\label{lemma:H:theorem:nozero:eps}
         \cL_\infty(\theta) \geq \frac1{864}.
\end{equation}
Furthermore, \nobs that, e.g., \cite[Lemma 3.1]{Jentzen}  implies that 
$[0, \infty) \ni t \mapsto \cL_\infty ( \Theta_t ) \in \R$ is non-increasing.
Combining this with \cref{lemma:H:theorem:nozero:eps} and \cref{theo:intro:convergence} assures that for all
$\Theta \in C([0, \infty) , \R^{7})$  with  $\forall \, t \in [0, \infty) \colon  \Theta_t = \Theta_0 - \int_0^t \cG ( \Theta_s ) \, \d s$ and $\cL_\infty(\Theta_0)<\nicefrac1{864}$
it holds that
$\liminf_{t \to \infty} \norm{\Theta_t } = \infty.$
\end{cproof}
\subsection{Upper bounds for GFs}
\label{subsec:upper_bounds}
\begin{prop}\label{prop:upperbound}
Let $\fd \in \N$, $\Theta \in C([0,\infty), \R^{\fd})$, let $\cL 
\colon \R^{\fd} \to [0,\infty)$ and $\cG \colon \R^{\fd}\to \R^{\fd}$ be measurable, and assume for all $t \in [0,\infty)$ that $\cL(\Theta_t)=\cL(\Theta_0)-\int_0^t\norm{\cG(\Theta_s)}^2 \, \d s$ and $\Theta_t = \Theta_0 - \int_0^t \cG ( \Theta_s ) \, \d s$.
Then it holds for all $t \in [0,\infty)$ that
\begin{equation}
    \norm{\Theta_t}\leq\norm{\Theta_0}+t^{\nicefrac12}|\cL(\Theta_0)-\cL(\Theta_t)|^{\nicefrac12}\leq  \norm{\Theta_0}+[t\cL(\Theta_0)]^{\nicefrac12}.
\end{equation}
\end{prop}
\begin{cproof}{prop:upperbound}
\Nobs that the assumption that for all $t \in [0,\infty)$ it holds that $\Theta_t = \Theta_0 - \int_0^t \cG ( \Theta_s ) \, \d s$, the triangle inequality, and the Cauchy-Schwarz inequality ensure for all $t \in [0,\infty)$ that
\begin{equation}
    \norm{\Theta_t} \leq \norm{\Theta_0} + \int_0^t \norm{\cG ( \Theta_s )} \, \d s\leq  \norm{\Theta_0} + t^{\nicefrac12} \left[\int_0^t \norm{\cG ( \Theta_s )}^2 \, \d s\right]^{\nicefrac12}.
\end{equation}
\Hence that for all $t \in [0,\infty)$ it holds that
\begin{equation}
     \norm{\Theta_t} \leq \norm{\Theta_0} + t^{\nicefrac12} |\cL(\Theta_0)-\cL(\Theta_t)|^{\nicefrac12}\leq  \norm{\Theta_0}+[t\cL(\Theta_0)]^{\nicefrac12}.
\end{equation}
\end{cproof}
\section{Non-existence of global minima of the risk and divergence of GFs and gradient descent (GD) for widely used activation functions}\label{generalframe}
%

In this section we establish, in the case of at least two neurons on the hidden layer, the non-existence of global minima employing various activation and target functions. Next we show the blow up of GFs under a specific asymptotic optimality assumption regarding the risk values. The key idea is to prove that there exists a sequence of ANN parameters such that the risk converges to zero and to consequently demonstrate that the set of global minima is empty. After proving this, \cref{theorem:globalminima:proof}, \cref{theorem:globalminima:proof2}, and \cref{newtechnique:risk} in \cref{subdivergence} assure the divergence of GFs under the assumption that the risk of GFs converges to the infimum of the risk while \cref{theorem:globalminima:proofb}, \cref{theorem:globalminima:proof2b}, and \cref{newtechnique:riskb} in \cref{subdivergence} prove the corresponding result in the discrete-time case.
Related results can be found in \cite[Proposition 3.6]{petersen}.
\cref{theorem:globalminima:proof}, \cref{theorem:globalminima:proof2}, and \cref{newtechnique:risk} are based on \cref{lemma:new}, demonstrated using compactness and continuity  properties, and the well-known
deterministic It\^{o}-type formula for continuously differentiable functions, see, e.g., \cite[Lemma 3.1]{CheriditoJentzenRiekert2021}.
\cref{theorem:globalminima:proofb}, \cref{theorem:globalminima:proof2b}, and \cref{newtechnique:riskb} follow from  \cref{lemma:new2}.
 
Choosing the square function as target function we establish the non-existence of global minima  in the case of softplus activation in \cref{H:noinf:soft:k2} in \cref{soft}, in the case of standard logistic, hyperbolic tangent, arctangent, and inverse square root unit activation in \cref{H:noinf:log2} in \cref{some}, in the case of exponential linear unit activation in \cref{H:noinf:elu} in \cref{elu}, and in the case of softsign activation in \cref{H:noinf:softsign} in \cref{softsign}.
The proofs of \cref{H:noinf:soft:k2} and  \cref{H:noinf:log2} employ properties of real analytic functions and are inspired by \cite[Theorem 3.3]{petersen}.
The proofs of \cref{H:noinf:elu} and \cref{H:noinf:softsign} instead use a comparison between target and realization function derivatives.

Employing an indicator function as target function we demonstrate the non-existence of global minima in the case of ReLU and leaky ReLU activation function in \cref{H:noinf} in \cref{ReLU}.
Its proof uses \cref{prop:relu:limzero}, \cref{lemma:leakyrelu:lessk},
 Lipschitz continuity results in  \cref{newtechnique:l2lowerbound}, and elementary properties of realization functions in \cref{newtechnique:l2lowerboundrealization}.

In the case of rectified power unit activation we instead establish in \cref{H:noinf:repu} in \cref{repu} the non-existence of global minima using as target function the rectified power unit itself with a smaller exponent. Ingredients employed in the proof are \cref{lemma:repuk:zero} and a continuity study inspired by \cite[Theorem 3.3]{petersen}.

Also notably,  we show the non-existence of global minima for every number of hidden neurons in \cref{H:noinf:soft} in \cref{soft} employing the ReLU function as target function and the softplus activation function and in \cref{H:noinf:log} in \cref{some} using the identity function as target function and standard logistic, arctangent, and inverse square root unit activation. The proofs of \cref{H:noinf:soft} and \cref{H:noinf:log} are inspired by \cite[Theorem 3.3]{petersen}.
\subsection{Mathematical description of ANNs}
\begin{setting} \label{H:setting:general}
Let $ \scra\in \R$, $\scrb \in (\scra,\infty)$, $\xii \in (0,\infty)$, $ \width, \fd \in \N$ satisfy $\fd = 3 \width + 1$,
let $\fw  = (( \w{\theta} _ {1}, \ldots, \w{\theta}_{\width}  ))_{ \theta \in \R^{\fd}} \colon \R^{\fd} \to \R^{\width}$,
$\fb =  (( \b{\theta} _ 1, \ldots, \b{\theta}_{\width} ) )_{ \theta \in \R^{\fd}} \colon \R^{\fd} \to \R^{\width}$,
$\fv = (( \v{\theta} _ 1, \ldots, \v{\theta}_{\width} ) )_{ \theta \in \R^{\fd}} \colon \R^{\fd} \to \R^{\width}$, and
$\fc = (\c{\theta})_{\theta \in \R^{\fd }} \colon \R^{\fd} \to \R$
 satisfy for all $\theta  = ( \theta_1 ,  \ldots, \theta_{\fd}) \in \R^{\fd}$, $j \in \{1, 2, \ldots, \width \}$ that $\w{\theta}_{ j} = \theta_{ j}$, $\b{\theta}_j = \theta_{\width + j}$, 
$\v{\theta}_j = \theta_{2 \width + j}$, and
$\c{\theta} = \theta_{\fd}$,
for every $k \in \Z $, $\gamma \in \R$
let $A_{k, \gamma} \colon \R \to \R$ satisfy for all $x \in \R$ that
\begin{equation}\label{H:setting:general:activation}
    A_{k, \gamma}(x)=
    \begin{cases}
           x(1+|x|)^{-1} & \colon k <-5 \\
           \arctan (x) & \colon k =-5 \\
           x (1+\xii x^2)^{-\nicefrac12}  & \colon k=-4 \\
            x \indicator{(0,\infty)}(x)+(\exp(x)-1)\indicator{(-\infty,0]}(x) & \colon k=-3 \\
           (\exp(x)-\exp(-x))(\exp(x)+\exp(-x))^{-1} & \colon k=-2  \\
          (1+\exp(-x))^{-1} & \colon k=-1 \\
          \ln (1+\exp(x))   & \colon k =0 \\
          (\max\{x,0\})^k+\min\{\gamma x,0\} & \colon k >0
    \end{cases} 
\end{equation}
and let $\act^r_{k,\gamma} \in C ( \R , \R )$, $r \in \N \cup \cu{ \infty } $,
satisfy for all $x \in \R$ that $( \bigcup_{r \in \N} \cu{ \act^r_{k,\gamma} } ) \subseteq C^1( \R , \R)$, $\act^\infty_{k,\gamma} ( x ) = A_{k, \gamma}(x)$,
 $\sup_{r \in \N} \sup_{y \in [- \abs{x}, \abs{x} ] }  \abs{ ( \act^r_{k,\gamma})'(y)}< \infty$, and
\begin{equation} 
    \limsup\nolimits_{r \to \infty} \Big(  \abs{ \act^r_{k,\gamma} ( x ) - \act_{k,\gamma}^\infty ( x ) }
    +  \bigl| (\act^r_{k,\gamma})' ( x ) - \lim\nolimits_{h \nearrow 0}\tfrac{\act_{k,\gamma}^\infty ( x +h)-\act_{k,\gamma}^\infty ( x )}{h} \bigr|\Big) = 0,
\end{equation}
for every  $\theta \in \R^{\fd}$, $r \in \N \cup \cu{ \infty } $, $k \in \Z $, $\gamma \in \R$
let $\realization{\theta,r}_{k,\gamma} \colon \R \to \R$
satisfy for all $x \in \R$ that
\begin{equation} \label{H:setting:general:realization}
    \realization{\theta,r}_{k,\gamma} (x) = \c{\theta} + \sum_{i=1}^{\width} \v{\theta}_i \br[\big]{ \act^r_{k,\gamma} ( \w{\theta}_i x + \b{\theta}_i )},
\end{equation}
let  $f \colon \R \to \R$ be measurable, for every $r \in \N \cup \cu{ \infty } $, $k \in \Z $, $\gamma \in \R$ let $\cL^r_{k,\gamma} \colon \R^{\fd} \to \R$
satisfy for all $\theta \in \R^{\fd}$ that
\begin{equation}
    \cL^r_{k,\gamma} ( \theta ) = \int_{\scra}^{\scrb} \rbr[\big]{ f(x) - \realization{\theta,r}_{k,\gamma} (x) }^2  \, \d x,
\end{equation}
and for every $k \in \Z $, $\gamma \in \R$ let $\cG_{k,\gamma} \colon \R^{\fd} \to \R^{\fd}$ satisfy for all
$\theta \in \cu{ \vartheta \in \R^{\fd} \colon ( ( \nabla \cL^r_{k,\gamma}) ( \vartheta ) ) _{r \in \N} \allowbreak \text{ is convergent} }$
that $\cG_{k,\gamma} ( \theta ) = \lim_{r \to \infty} (\nabla \cL_{k,\gamma}^{r}) ( \theta )$. 
\end{setting}
\subsection{
ANNs with ReLU and leaky ReLU activation}\label{ReLU}
\begin{prop}\label{prop:relu:limzero}
Assume \cref{H:setting:general}, assume $\width > 1$, assume for all $x \in \R$ that $f(x)=\indicator{(\nicefrac{(\scra+\scrb)}{2},\infty
)}(x)$, let $\gamma \in (-\infty,0]$, and let $(\theta_n)_{n \in \N}\subseteq \R^{\fd}$ satisfy for all $n \in \N$ that $\w{\theta_n}_1=\w{\theta_n}_2=(n+3)(2\scrb-2\scra)^{-1}$,   
$\b{\theta_n}_{1}=-(1+n)4^{-1}-\scra(n+3)(2\scrb-2\scra)^{-1}$,  
$\b{\theta_n}_{2}=-(5+n)4^{-1}-\scra(n+3)(2\scrb-2\scra)^{-1}$, 
$\v{\theta_n}_{1}=-\v{\theta_n}_{2} =(1-\gamma)^{-1}$, and  $|\c{\theta_n}|+ \sum_{j=3}^\width|\w{\theta_n}_{j}|+ |\b{\theta_n}_{j}|+|\v{\theta_n}_{j}| =0$. Then \begin{equation}
    \limsup\nolimits_{n \to \infty}\cL^\infty_{1,\gamma}(\theta_n)=0.
\end{equation}
\end{prop}
\begin{cproof}{prop:relu:limzero}
\Nobs that \cref{H:setting:general:realization} ensures that for all $n\in \N$, $x \in [\scra,\scrb]$ it holds that
\begin{equation}
\begin{split}
    \realization{\theta_n,\infty}_{1,\gamma} (x) &= \max \left\{ -\frac{1+n}{4}-\frac{\scra(n+3)}{2(\scrb-\scra)} + \frac{n+3}{2(\scrb-\scra)} x, 0 \right\} \\
    & \quad - \max\left\{ -\frac{5+n}{4}-\frac{\scra(n+3)}{2(\scrb-\scra)} + \frac{n+3}{2(\scrb-\scra)} x, 0 \right\}. 
\end{split}
\end{equation}
This implies that for all $n\in \N$ it holds that
\begin{equation}
\begin{split}
    \cL_{1,\gamma}^\infty (\theta_n) &= \int_{\scra}^{\scrb} (\mathbbm{1}_{(\nicefrac{(\scra+\scrb)}2,\infty)}(x) - \realization{\theta_n,\infty}_{1,\gamma} (x) )^2 \,\d x\\
    &= \int_{\scra+\frac{(\scrb-\scra)(1+n)}{2(n+3)}}^{\frac{\scra+\scrb}{2}} \left(-\frac{1+n}{4}-\frac{\scra(n+3)}{2(\scrb-\scra)} + \frac{n+3}{2(\scrb-\scra)} x \right)^2 \,\d x \\
    & \quad +\int_{\frac{\scra+\scrb}{2}}^{\scra+\frac{(\scrb-\scra)(5+n)}{2(n+3)}} \left(-\frac{5+n}{4}-\frac{\scra(n+3)}{2(\scrb-\scra)} + \frac{n+3}{2(\scrb-\scra)} x \right)^2 \,\d x \\
    &\leq\int_{\scra+\frac{(\scrb-\scra)(1+n)}{2(n+3)}}^{\frac{\scra+\scrb}{2}} \left(\frac{1}{2} \right)^2 \,\d x +\int_{\frac{\scra+\scrb}{2}}^{\scra+\frac{(\scrb-\scra)(5+n)}{2(n+3)}} \left(\frac{1}{2} \right)^2 \,\d x =  \frac{\scrb-\scra}{2(n+3)}.
\end{split}
\end{equation}
\Hence that
\begin{equation}
    \limsup\nolimits_{n \to \infty}\cL^\infty_{1,\gamma}(\theta_n)=0.
\end{equation}
\end{cproof}
\begin{prop}\label{lemma:leakyrelu:lessk}
Assume \cref{H:setting:general}, assume $\width>1$,
assume for all $x \in \R$ that $f(x)=\indicator{(\nicefrac{(\scra+\scrb)}{2},\infty
)}(x)$, let $\gamma \in (0,\infty) \backslash\{1\}$, and let $(\theta_n)_{n \in \N} \subseteq \R^{\fd} $ satisfy for all $n \in \N$ that
$\w{\theta_n}_1=\w{\theta_n}_2=(n+3)(2\scrb-2\scra)^{-1}$,
$\b{\theta_n}_{1}=-(1+n)4^{-1}-\scra(n+3)(2\scrb-2\scra)^{-1}$,  
$\b{\theta_n}_{2}=-(5+n)4^{-1}-\scra(n+3)(2\scrb-2\scra)^{-1}$,
$\v{\theta_n}_1=-\v{\theta_n}_2=(1-\gamma)^{-1}$,
$\c{\theta_n}=-\gamma(1-\gamma)^{-1}$, and 
$\sum_{j=3}^\width|\w{\theta_n}_{j}|+ |\b{\theta_n}_{j}|+|\v{\theta_n}_{j}| =0$.
Then  
\begin{equation}
    \limsup\nolimits_{n \to \infty} \cL^\infty_{1,\gamma}(\theta_n)=0.
\end{equation}
\end{prop}
\begin{cproof}{lemma:leakyrelu:lessk}
\Nobs that \cref{H:setting:general:realization} ensures that for all $n\in \N$, $x \in \R$ it holds that
\begin{equation}\label{lemma:leakyrelu:lessk:cases}
    \realization{\theta_n,\infty}_{1,\gamma}(x)=
    \begin{cases}
     0
     &
    \colon
    x \in [\scra,\scra+\frac{(\scrb-\scra)(1+n)}{2(n+3)}]
    \\
    -\frac{1+n}{4}-\frac{\scra(n+3)}{2(\scrb-\scra)} + \frac{n+3}{2(\scrb-\scra)} x.
     &
    \colon
    x \in (\scra+\frac{(\scrb-\scra)(1+n)}{2(n+3)},\scra+\frac{(\scrb-\scra)(5+n)}{2(n+3)}]
    \\
    1 
    &
    \colon
    x \in (\scra+\frac{(\scrb-\scra)(5+n)}{2(n+3)},\scrb].
    \end{cases}
\end{equation}
This implies that for all $x \in \R$ it holds that
\begin{equation}\label{lemma:leakyrelu:lessk:eq}
    \limsup\nolimits_{n \to \infty}|\realization{\theta_n,\infty}_{1,\gamma}(x)-\mathbbm{1}_{(\nicefrac{(\scra+\scrb)}2,\infty)}(x)|=0.
\end{equation}
\Moreover \cref{lemma:leakyrelu:lessk:cases} assures that for all $n\in \N$, $x \in [\scra,\scrb]$ it holds that $|\realization{\theta_n}_{1,\gamma}(x)-f(x)|\leq 1$. Combining this with \cref{lemma:leakyrelu:lessk:eq} and Lebesgue's dominated convergence theorem demonstrates that \begin{equation}
    \limsup\nolimits_{n \to \infty} \cL^\infty_{1,\gamma}(\theta_n)=0.
\end{equation}
\end{cproof}
\begin{lemma}\label{newtechnique:l2lowerbound}
Let $\scra \in \R$, $\scrb \in (\scra,\infty)$, let $f \colon [\scra,\scrb] \to \R$ satisfy for all $x \in [\scra,\scrb]$ that $f(x) = \mathbbm{1}_{(\nicefrac{(\scra+\scrb)}2,\infty)}(x)$, and let $L \in \R$, $g \in C( [\scra,\scrb], \R )$ satisfy for all $ x, y \in [\scra,\scrb] $ that $ | g(x) - g(y) | \leq L | x - y |$. 
Then 
\begin{equation}\label{newtechnique:l2lowerbound:thesis}
    \int_\scra^\scrb ( f(x) - g(x) )^2 \, \d x \geq \frac1{32\max\{L,\nicefrac1{(\scrb-\scra)}\}}.
\end{equation}
\end{lemma}
\begin{cproof}{newtechnique:l2lowerbound}
Throughout this proof let $c\in\R$ satisfy \begin{equation}\label{newtechnique:l2lowerbound:c}
    c=\max\left\{L,\nicefrac1{(\scrb-\scra)}\right\}.
\end{equation} 
\Nobs that \cref{newtechnique:l2lowerbound:c} assures that  for all $ x, y \in [\scra,\scrb] $ it holds that
\begin{equation}\label{eqC:newtechnique:l2lowerbound}
    | g(x) - g(y) | \leq c | x - y |.
\end{equation}
In the following we distinguish between the case $g(\nicefrac{(\scra+\scrb)}2)>\nicefrac12$ and the case $g(\nicefrac{(\scra+\scrb)}2)\leq \nicefrac12$.
We first prove \cref{newtechnique:l2lowerbound:thesis} in the case
\begin{equation}\label{eqbigger:newtechnique:l2lowerbound}
    g\left(\frac{\scra+\scrb}2\right)>\frac12.
\end{equation}
\Nobs that \cref{eqC:newtechnique:l2lowerbound} and \cref{eqbigger:newtechnique:l2lowerbound} imply that for all $x \in [\nicefrac{(\scra+\scrb)}2-\nicefrac{1}{2c},\nicefrac{(\scra+\scrb)}2]$ it holds that
\begin{equation}\label{eq:newtechnique:l2lowerbound}
    0 \leq \frac12 -c\left(\frac{\scra+\scrb}2-x\right)< g\left(\frac{\scra+\scrb}2\right)-c\left(\frac{\scra+\scrb}2-x\right) \leq g( x).
\end{equation}
This proves that
\begin{equation}
    \int_{\frac{\scra+\scrb}2 - \frac{1}{2c}}^{\frac{\scra+\scrb}2} |g(x)| \, \d x \geq
    \int_{\frac{\scra+\scrb}2 - \frac{1}{2c}}^{\frac{\scra+\scrb}2} \frac12(1-c(\scra+\scrb)) +cx \, \d x = \frac{1}{8c}.
\end{equation}
Combining this and the Cauchy-Schwarz inequality demonstrates that
\begin{equation}
\begin{split}
     \int_{\scra}^{\scrb} |f(x)-g(x)|^2 \, \d x & \geq \int_{\frac{\scra+\scrb}2 - \frac{1}{2c}}^{\frac{\scra+\scrb}2} |f(x)-g(x)|^2 \, \d x \\
     &\geq  2c \left(\int_{\frac{\scra+\scrb}2 - \frac{1}{2c}}^{\frac{\scra+\scrb}2} |g(x)| \, \d x \right)^2 = \frac{1}{32c}.
\end{split}
\end{equation}
This establishes \cref{newtechnique:l2lowerbound:thesis} in the case $g(\nicefrac{(\scra+\scrb)}2)>\nicefrac12$.
In the next step we prove \cref{newtechnique:l2lowerbound:thesis} in the case
\begin{equation}\label{eqlower:newtechnique:l2lowerbound}
    g\left(\frac{\scra+\scrb}2\right)\leq \frac12.
\end{equation}
\Nobs that \cref{eqC:newtechnique:l2lowerbound} and \cref{eqlower:newtechnique:l2lowerbound} imply that for all $x \in [\nicefrac{(\scra+\scrb)}2,\nicefrac{(\scra+\scrb)}2+\nicefrac{1}{2c}]$ it holds that
\begin{equation}
    g\big(x\big)\leq 
    c\left(x-\frac{\scra+\scrb}2\right)+ g\left(\frac{\scra+\scrb}2\right) \leq 
     c\left(x-\frac{\scra+\scrb}2\right) + \frac12 \leq
    1.
\end{equation}
This proves that
\begin{equation}
    \int_{\frac{\scra+\scrb}2}^{\frac{\scra+\scrb}2+\frac{1}{2c}} |f(x)-g(x)| \, \d x \geq
   \int_{\frac{\scra+\scrb}2}^{\frac{\scra+\scrb}2+\frac{1}{2c}} 1 - \frac12(1-c(\scra+\scrb)) - cx \, \d x = \frac{1}{8c}.
\end{equation}
Combining this and the Cauchy-Schwarz inequality demonstrates that
\begin{equation}
 \begin{split}
     \int_{\scra}^{\scrb} |f(x)-g(x)|^2 \, \d x & \geq \int_{\frac{\scra+\scrb}2}^{\frac{\scra+\scrb}2+\frac{1}{2c}} |f(x)-g(x)|^2 \, \d x \\
     &\geq 2c
     \left(\int_{\frac{\scra+\scrb}2}^{\frac{\scra+\scrb}2+\frac{1}{2c}} |1-g(x)| \, \d x \right)^2 = \frac{1}{32c}.
\end{split}
\end{equation}
This establishes \cref{newtechnique:l2lowerbound:thesis} in the case $g\left(\nicefrac{(\scra+\scrb)}2\right)\leq \nicefrac12$
\end{cproof}
\begin{lemma}\label{newtechnique:l2lowerboundrealization}
Assume \cref{H:setting:general}  and let $\theta \in \R^{ \fd}$,  $\gamma \in \R\backslash\{1\}$. Then
\begin{equation}
    \sup\nolimits_{ x, y \in [\scra,\scrb],\, x \neq y } ( | \realization{ \theta,\infty }_{1,\gamma}( x ) - \realization{ \theta,\infty }_{1,\gamma}( y ) |  | x - y |^{-1} ) \leq  (1+|\gamma|) \| \theta \|^2 .
\end{equation}
\end{lemma}
\begin{cproof}{newtechnique:l2lowerboundrealization}
\Nobs that the fact that $\realization{\theta,\infty}_{1,\gamma}$ is continuous and the fact that $\realization{\theta,\infty}_{1,\gamma}$ is piecewise affine linear imply for all $ x, y \in  [\scra,\scrb]$  that 
\begin{equation}
    | \realization{ \theta,\infty }_{1,\gamma}( x ) - \realization{ \theta,\infty }_{1,\gamma}( y ) | \leq \max\{1,\gamma\} \left(\smallsum_{i=1}^H |\w{\theta}_i\v{\theta}_i| \right) |x-y| \leq (1+|\gamma|) \| \theta \|^2 |x-y|.
\end{equation}
This demonstrates that 
\begin{equation}
   \sup\nolimits_{ x, y \in [\scra,\scrb],\, x \neq y } ( | \realization{ \theta,\infty }_{1,\gamma}( x ) - \realization{ \theta,\infty }_{1,\gamma}( y ) |  | x - y |^{-1} ) \leq  (1+|\gamma|) \| \theta \|^2 .
\end{equation}
\end{cproof}
\begin{lemma}\label{H:noinf}
Assume \cref{H:setting:general},  assume $\width> 1$, assume for all $x \in \R$ that $f(x)=\indicator{(\nicefrac{(\scra+\scrb)}{2},\infty
)}(x)$,
and let $\gamma \in \R \backslash\{1\}$. Then
\begin{equation}\label{H:noinf:leakyrelu:thesis}
   \big\{ \vartheta \in \R^{ \fd } \colon \cL^\infty_{1,\gamma}( \vartheta ) = \inf\nolimits_{ \theta \in \R^{ \fd } } \cL^\infty_{1,\gamma}( \theta ) \big\} = \emptyset.
\end{equation}
\end{lemma}
\begin{cproof}{H:noinf}
\Nobs that \cref{prop:relu:limzero} and \cref{lemma:leakyrelu:lessk} imply that $\inf_{ \theta \in \R^{ \fd } } \cL^\infty_{1,\gamma}( \theta )=0$.
\Moreover \cref{newtechnique:l2lowerbound} and \cref{newtechnique:l2lowerboundrealization} ensure for all $\theta \in \R^{ \fd}$ that
\begin{equation}
    \cL^\infty_{1,\gamma}(\theta)
    \geq \frac1{32\max\{(1+|\gamma|)\|\theta\|^2,\nicefrac1{(\scrb-\scra)}\}}.
\end{equation}
This implies that  $\{ \vartheta \in \R^{ \fd } \colon \cL^\infty_{1,\gamma}( \vartheta ) =0 \} = \emptyset$.
\end{cproof}
\subsection{ANNs with softplus activation }\label{soft}
\begin{lemma}\label{lemma:onesoft:lessk}
Assume \cref{H:setting:general} and let $r\in(0,\infty)$, $\theta \in \R^4 $ satisfy $\c{\theta}=0$, $\b{\theta}_1=-r(\scra+\scrb)2^{-1}$, $\w{\theta}_1=r$, and $\v{\theta}_1=\nicefrac1{r}$.
Then  it holds for all $x \in \R$ that $0\leq \realization{\theta,\infty}_{0,0}(x)- \max\{x-\nicefrac{(\scra+\scrb)}2, 0\}\leq \nicefrac{1}{r}$.
\end{lemma}
\begin{cproof}{lemma:onesoft:lessk}
\Nobs that \cref{H:setting:general:realization} ensures that for all $x \in \R$ it holds that 
\begin{equation}\label{lemma:onesoft:lessk:bound}
    \realization{\theta,\infty}_{0,0}(x)=\frac1{r} \ln\left(1+\exp\left(rx-r\frac{\scra+\scrb}2\right)\right).
\end{equation}
\Moreover for all $x\in [0, \infty)$ it holds that
\begin{equation}
    x\leq \ln(1+\exp(x))\leq x+1.
\end{equation}
Combining this and \cref{lemma:onesoft:lessk:bound} establishes for all $x \in [\nicefrac{(\scra+\scrb)}2, \infty)$  that
\begin{equation}\label{lemma:onesoft:lessk:part1}
    0\leq\realization{\theta,\infty}_{0,0}(x)-x+\frac{\scra+\scrb}2 \leq \frac{1}{r}.
\end{equation}
\Moreover for all $x\in (-\infty,0]$ it holds that
\begin{equation}
    0\leq \ln(1+\exp(rx))\leq \exp(rx) \leq 1.
\end{equation}
Combining this and \cref{lemma:onesoft:lessk:bound} implies for all $x \in (-\infty,\nicefrac{(\scra+\scrb)}2]$  that $0\leq\realization{\theta,\infty}_{0,0}(x)\leq \nicefrac{1}{r}$.
This and \cref{lemma:onesoft:lessk:part1} show that
for all $x \in \R$ it holds that $0\leq \realization{\theta,\infty}_{0,0}(x)- \max\{x-\nicefrac{(\scra+\scrb)}2, 0\}\leq \nicefrac{1}{r}$.
\end{cproof}
\begin{prop}\label{prop:soft:lessk}
Assume \cref{H:setting:general}, assume $\width >1$, assume for all $x \in \R$ that $f(x)=x^2$,
and let $(\theta_n)_{n \in \N}\subseteq \R^{\fd}$ satisfy for all $n \in \N$ that
$\w{\theta_n}_1=-\w{\theta_n}_2=\nicefrac1n$,   $\b{\theta_n}_{1}=\b{\theta_n}_{2}=0$, $\v{\theta_n}_{1}=\v{\theta_n}_{2}=4n^2$, 
$\c{\theta_n}=-8n^2\ln(2)$,  and $\sum_{j=3}^\width|\w{\theta_n}_{j}|+ |\b{\theta_n}_{j}|+|\v{\theta_n}_{j}| =0.$
Then  
\begin{equation}
    \limsup\nolimits_{n \to \infty} \cL^\infty_{0,0}(\theta_n)=0.
\end{equation}
\end{prop}
\begin{cproof}{prop:soft:lessk}
\Nobs that \cref{H:setting:general:realization} ensures that for all $x\in \R$, $n \in \N$ it holds that
\begin{equation}
    \realization{\theta_n,\infty}_{0,0} (x) = 4n^2 \ln\left(1+ \exp\left(\frac xn\right)\right)+ 4n^2 \ln\left(1+\exp\left(-\frac xn \right)\right) -8n^2\ln(2).
\end{equation}
This and the fact that there exists $c \in (0,\infty)$ such that for all $x \in [-1,1]$ it holds that $|\ln(1+\exp(x))-\ln(2)- \nicefrac x2-\nicefrac {x^2}{8}|\leq c |x^3| $ assure that there exist $c,M\in (0,\infty)$ such that for all $x\in [\scra,\scrb]$, $n\geq M$ it holds that 
\begin{equation}\label{prop:soft:lessk:eq}
    |\realization{\theta_n,\infty}_{0,0}(x)-x^2 |\leq c \left|\frac 1{n^3}\right|.
\end{equation}
This and Lebesgue's dominated convergence theorem demonstrate that 
\begin{equation}
    \limsup\nolimits_{n \to \infty}\cL^\infty_{0,0}(\theta_n)=\int_{\scra}^{\scrb}\limsup\nolimits_{n \to \infty} \rbr[\big]{ x^2 - \realization{\theta_n,\infty}_{0,0}(x) }^2  \, \d x=0.
\end{equation}
\end{cproof}
%
%
%
%
%
%
%
\begin{lemma}\label{H:noinf:soft}
Assume \cref{H:setting:general} and assume for all $x \in \R$ that $f(x)=\max\{x-\nicefrac{(\scra+\scrb)}2,0\}$. Then
\begin{equation}\label{H:noinf:soft:thesis}
    \big\{ \vartheta \in \R^{ \fd } \colon \cL^\infty_{0,0}( \vartheta ) = \inf\nolimits_{ \theta \in \R^{ \fd } } \cL^\infty_{0,0}( \theta ) \big\} = \emptyset.
\end{equation}
\end{lemma}
\begin{cproof}{H:noinf:soft}
\Nobs that \cref{lemma:onesoft:lessk} proves that for every $r\in (0,\infty)$ there exists $\vartheta_r \in \R^{\fd}$ which satisfies for all $x\in \R$ that
$|\realization{\vartheta_r,\infty}_{0,0}(x)- \max\{x-\nicefrac{(\scra+\scrb)}2, 0\}|\leq \nicefrac1r$.
This implies that for every $r\in (0,\infty)$ there exists $\vartheta_r \in \R^{\fd}$ such that $\cL^\infty_{0,0}(\vartheta_r)\leq \nicefrac{(\scrb-\scra)}{r^2}$.
\Hence that $\inf_{\theta \in \R^{\fd}} \cL^\infty_{0,0}( \theta) = 0$.
\Moreover for all $\theta \in \R^{\fd}$ it holds that $\realization{\theta,\infty}_{0,0}\in C^\infty(\R,\R).$
We prove \cref{H:noinf:soft:thesis} by contradiction. Assume that there exists $\vartheta \in \R^{ \fd }$ which satisfies that \begin{equation}\label{H:noinf:soft:ab}
    \cL^\infty_{0,0}(\vartheta)=0.
\end{equation} 
\Nobs that \cref{H:noinf:soft:ab} ensures that for all $x \in [\scra, \scrb]$ it holds that $\realization{\vartheta,\infty}_{0,0}(x)=\max\{x-\nicefrac{(\scra+\scrb)}2,0\}$. This demonstrates that $\realization{\vartheta,\infty}_{0,0}\in C([\scra, \scrb],\R)\backslash C^1([\scra, \scrb],\R)$ which is a contradiction. 
\end{cproof}
%
%
%
%
\begin{lemma}\label{H:noinf:soft:pos}
Assume \cref{H:setting:general}, assume $\width >1$, and assume for all $x \in \R$ that $f(x)=x$. Then
\begin{equation}
    \big\{ \vartheta \in \R^{ \fd } \colon \cL^\infty_{0,0}( \vartheta ) = \inf\nolimits_{ \theta \in \R^{ \fd } } \cL^\infty_{0,0}( \theta ) \big\} \neq \emptyset.
\end{equation}
\end{lemma}
\begin{cproof}{H:noinf:soft:pos}
Let $\vartheta \in \R^{\fd}$ satisfy  for all $i \in \{1,2,\ldots, \width\} \backslash \{1,2 \}$ that $\w{\vartheta}_1=-\w{\vartheta}_2=\v{\vartheta}_1=-\v{\vartheta}_2=1$
and $\b{\vartheta}_1=\b{\vartheta}_2=\w{\vartheta}_i=\v{\vartheta}_i=\b{\vartheta}_i=\c{\vartheta}=0$. This proves that for all $x \in \R$ it holds that 
\begin{equation}
     \realization{\vartheta,\infty}_{0,0}(x)=\ln (1+\exp(x))- \ln (1+ \exp(-x))=\ln\left(\frac{1+\exp(x)}{1+\exp(-x)}\right)=\ln(\exp(x))=x.
\end{equation}
\Hence that $\cL^\infty_{0,0}(\vartheta)=0$. This and the fact that for all $\theta \in \R^{\fd}$ it holds that $\cL^\infty_{0,0}(\theta)\geq0$ demonstrates that $\vartheta \in   \{ v \in \R^{ \fd } \colon \cL^\infty_{0,0}(v) = \inf_{ \theta \in \R^{ \fd } } \cL^\infty_{0,0}( \theta ) \}$.
\end{cproof}
%
%
%
%
\begin{lemma}\label{H:noinf:soft:k2}
Assume \cref{H:setting:general}, assume $\width>1$, and assume for all $x \in \R$ that $f(x)=x^2$. Then
\begin{equation}\label{H:noinf:soft:pos:thesis}
    \big\{ \vartheta \in \R^{ \fd } \colon \cL^\infty_{0,0}( \vartheta ) = \inf\nolimits_{ \theta \in \R^{ \fd } } \cL^\infty_{0,0}( \theta ) \big\} = \emptyset.
\end{equation}
\end{lemma}
\begin{cproof}{H:noinf:soft:k2}
\Nobs that $f$ is real analytic.
\Moreover for all $\theta \in \R^{ \fd }$ it holds that $\realization{\theta,\infty}_{0,0}$ is real analytic.
\Moreover \cref{prop:soft:lessk} ensures that $\inf_{ \theta \in \R^{ \fd } } \cL^\infty_{0,0}( \theta )=0$. We prove \cref{H:noinf:soft:pos:thesis} by contradiction. Assume that there exists $\vartheta \in \R^{ \fd }$ such that 
\begin{equation}\label{H:noinf:soft:k2:ab}
    \cL^\infty_{0,0}(\vartheta)=0.
\end{equation}
\Nobs that \cref{H:noinf:soft:k2:ab} establishes that for all $x \in [\scra,\scrb]$ it holds that $f(x) = \realization{\vartheta,\infty}_{0,0}(x)$.
Combining this with the fact that $f$ and $\realization{\vartheta,\infty}_{0,0}$ are real analytic implies for all $x \in \R$ that
\begin{equation}\label{H:noinf:soft:k2:eq}
    f(x)=\realization{\vartheta,\infty}_{0,0}(x).
\end{equation}
\Nobs that for all $x \in \R$ it holds that
\begin{equation}
    |(\realization{\vartheta,\infty}_{0,0})'(x)|=\left|\sum_{k=1}^\width \frac{\v{\vartheta}_k\w{\vartheta}_k\exp(\w{\vartheta}_k x + \b{\vartheta}_k)}{1+\exp(\w{\vartheta}_k x + \b{\vartheta}_k)}\right|\leq  \sum_{k=1}^\width |\v{\vartheta}_k\w{\vartheta}_k|.
\end{equation}
This, the fact that $f'$ is unbounded, and \cref{H:noinf:soft:k2:eq} show the contradiction.
\end{cproof}
\subsection{ANNs with standard logistic, hyperbolic tangent, arctangent, and inverse square root unit activation}\label{some}
\begin{prop}\label{lemma:log:lessk}
Assume \cref{H:setting:general}, assume for all $x \in \R$ that $f(x)=x$, and let $(\theta_n)_{n \in \N}\subseteq \R^{\fd}$ satisfy for all $n \in \N$ that $\w{\theta_n}_1=\nicefrac1n$,
$\b{\theta_n}_1=0$,
$\v{\theta_n}_1=4n$,
$\c{\theta_n}=-n$, and
$\sum_{j=2}^\width|\w{\theta_n}_{j}|+ |\b{\theta_n}_{j}|+|\v{\theta_n}_{j}| =0.$
Then  
\begin{equation}
    \limsup\nolimits_{n \to \infty} \cL^\infty_{-1,0}(\theta_n)=0.
\end{equation}
\end{prop}
\begin{cproof}{lemma:log:lessk}
\Nobs that \cref{H:setting:general:realization} ensures that for all $x \in \R$, $n \in \N$ it holds  that 
\begin{equation}\label{lemma:log:lessk:bound}
    \realization{\theta_n,\infty}_{-1,0}(x)= a_{-1,0}\left(\frac xn\right) 4n
    -2n.
\end{equation}
This and the fact that there exists $c \in (0,\infty)$ such that for all $x \in [-1,1]$ it holds that $|a_{-1,0}(x)-\nicefrac12-\nicefrac x4 +\nicefrac{x^3}{48}| \leq c |x^4|$ demonstrate that for all $x\in \R$ it holds that \begin{equation}\label{lemma:log:lessk:eq}
     \limsup\nolimits_{n \to \infty}|\realization{\theta_n,\infty}_{-1,0}(x)-x|=0.
\end{equation}
\Moreover the mean-value theorem demonstrates that for all  $x \in \R$  there exists $\tilde{x} \in [\min\{0,x\}, \max\{0,x\}]$ which satisfies that $a_{-1,0} (x)- \nicefrac12=  x \, a_{-1,0}'(\tilde{x})$. This and the fact that $a_{-1,0}'$ is continuous imply that there exists $L\in \R$ such that for all $x \in [\scra,\scrb]$, $n \in \N$ 
it holds that
\begin{equation}
\begin{split}
     |\realization{\theta_n,\infty}_{-1,0}(x)-x|
    &=|4n(a_{-1,0} (\nicefrac xn)- \nicefrac12)-x|
    =|x| |4a_{-1,0}'(\tilde{x})-1|\\
    &\leq\max\{|\scra|,|\scrb|\} |4a_{-1,0}'(\tilde{x})-1|\leq L.
\end{split}
\end{equation}
Combining this, \cref{lemma:log:lessk:eq}, and Lebesgue's dominated convergence theorem proves that 
\begin{equation}
    \limsup\nolimits_{n \to \infty}\cL^\infty_{-1,0} ( \theta_n ) =  \int_{\scra}^\scrb \limsup\nolimits_{n \to \infty}\rbr[\big]{ \realization{\theta_n,\infty}_{-1,0} (x) - x  }^2  \, \d x=0.
\end{equation}
\end{cproof}
\begin{prop}\label{lemma:log:lessk2}
Assume \cref{H:setting:general}, assume $\width>1$, assume for all $x \in \R$ that $f(x)=x^2$,  and let $(\theta_n)_{n \in \N}\subseteq \R^{\fd}$ satisfy for all $n \in \N$ that $\w{\theta_n}_1=-\w{\theta_n}_2=-\nicefrac1n$,
$\b{\theta_n}_1=\b{\theta_n}_2=-1$,
$\v{\theta_n}_1=\v{\theta_n}_2=n^2 (1+e)^3 (e(e-1))^{-1}$,
$\c{\theta_n}=-2n^2 (1+e)^2 (e(e-1))^{-1}$, and
$\sum_{j=3}^\width|\w{\theta_n}_{j}|+ |\b{\theta_n}_{j}|+|\v{\theta_n}_{j}| =0.$
Then  
\begin{equation}
    \limsup\nolimits_{n \to \infty} \cL^\infty_{-1,0}(\theta_n)=0.
\end{equation}
\end{prop}
\begin{cproof}{lemma:log:lessk2}
\Nobs that \cref{H:setting:general:realization} ensures that for all $x \in \R$, $n \in \N$ it holds  that 
\begin{equation}\label{lemma:log:lessk2:bound}
    \realization{\theta_n,\infty}_{-1,0}(x)=
    n^2\frac{(1+e)^3}{e(e-1)}\bigg(\frac{1}{1+e^{-\frac xn +1}}+\frac{1}{1+e^{\frac xn +1}}\bigg)-2n^2\frac{(1+e)^2}{e(e-1)}.
\end{equation}
This and the fact that there exists $c\in (0,\infty)$ such that for all $x\in[-1,1]$ it holds that
\begin{equation}
      \bigg|\frac1{1+e^{x+1}} -\frac{1}{1 + e} + \frac{e\, x}{(1 + e)^2} - \frac{e(e - 1)  x^2}{2 (1 + e)^3}  \bigg| \leq c|x^3|
\end{equation}
 assure that there exist $c,M\in (0,\infty)$ such that for all $x\in [\scra, \scrb]$, $n\geq M$ it holds that 
\begin{equation}\label{lemma:log:lessk2:eq}
    |\realization{\theta_n,\infty}_{-1,0}(x)-x^2|\leq  c \left|\frac 1{n}\right|.
\end{equation}
This and Lebesgue's dominated convergence theorem demonstrate that 
\begin{equation}
    \limsup\nolimits_{n \to \infty}\cL^\infty_{-1,0}(\theta_n)=\int_{\scra}^{\scrb}\limsup\nolimits_{n \to \infty} \rbr[\big]{ x^2 - \realization{\theta_n,\infty}_{-1,0} (x) }^2  \, \d x=0.
\end{equation}
\end{cproof}
\begin{prop}\label{lemma:tanh:lessk}
Assume \cref{H:setting:general}, assume $\width>1$, assume for all $x \in \R$ that $f(x)=x^2$,  and let $(\theta_n)_{n \in \N}\subseteq \R^{\fd}$ satisfy for all $n \in \N$ that $\w{\theta_n}_1=-\w{\theta_n}_2=\nicefrac1n$,
$\b{\theta_n}_1=\b{\theta_n}_2=-1$,
$\v{\theta_n}_1=\v{\theta_n}_2=n^2 (1+e^2)^3 (8 e^2 (e^2-1))^{-1}$,
$\c{\theta_n}=n^2 (1+e^2)^2(4 e^2)^{-1}$, and
$\sum_{j=3}^\width|\w{\theta_n}_{j}|+ |\b{\theta_n}_{j}|+|\v{\theta_n}_{j}| =0.$
Then  
\begin{equation}
    \limsup\nolimits_{n \to \infty} \cL^\infty_{-2,0}(\theta_n)=0.
\end{equation}
\end{prop}
\begin{cproof}{lemma:tanh:lessk}
\Nobs that \cref{H:setting:general:realization} ensures that for all $x \in \R$, $n \in \N$ it holds  that 
\begin{equation}\label{lemma:tanh:lessk:bound}
    \realization{\theta_n,\infty}_{-2,0}(x)=\frac{n^2 (1+e^2)^3}{8 e^2 (e^2-1)} \bigg(\frac{e^{\frac xn -1}-e^{-\frac xn +1}}{e^{\frac xn -1}+e^{-\frac xn +1}} + \frac{e^{-\frac xn -1}-e^{\frac xn +1}}{e^{-\frac xn -1}+e^{\frac xn +1}}\bigg)+ \frac{n^2 (1+e^2)^2}{4 e^2}.
\end{equation}
This and the fact that there exists $c\in (0,\infty)$ such that for all $x\in[-1,1]$ it holds that
\begin{equation}
    \bigg|\frac{e^{x -1}-e^{-x +1}}{e^{x -1}+e^{-x +1}}-
    \frac{1 - e^2}{1 + e^2} - \frac{4 e^2 x}{(1 + e^2)^2} - \frac{4 e^2 (e^2 - 1) x^2}{(1 + e^2)^3}\bigg| \leq c|x^3|
\end{equation}
 assure that there exist $c,M\in (0,\infty)$ such that for all  $x\in [\scra,\scrb]$, $n\geq M$ it holds that 
\begin{equation}\label{lemma:tanh:lessk:eq}
    |\realization{\theta_n,\infty}_{-2,0}(x)-x^2| \leq c \left|\frac1{n}\right|.
\end{equation}
This and Lebesgue's dominated convergence theorem demonstrate that 
\begin{equation}
    \limsup\nolimits_{n \to \infty}\cL^\infty_{-2,0}(\theta_n)=\int_{-1}^{1}\limsup\nolimits_{n \to \infty} \rbr[\big]{ x^2 - \realization{\theta_n,\infty}_{-2,0} (x) }^2  \, \d x=0.
\end{equation}
\end{cproof}
\begin{prop}\label{lemma:isru:zero}
Assume \cref{H:setting:general}, assume for all $x \in \R$ that $f(x)=x$,  and let $(\theta_n)_{n \in \N}\subseteq \R^{\fd}$ satisfy for all $n \in \N$ that $\w{\theta_n}_1=\nicefrac1n$,
$\v{\theta_n}_1=n$, and
$|\b{\theta_n}_1|+|\c{\theta_n}|+\sum_{j=2}^\width|\w{\theta_n}_{j}|+ |\b{\theta_n}_{j}|+|\v{\theta_n}_{j}| =0.$
Then  
\begin{equation}
    \limsup\nolimits_{n \to \infty} \cL^\infty_{-4,0}(\theta_n)=0.
\end{equation}
\end{prop}
\begin{cproof}{lemma:isru:zero}
\Nobs that \cref{H:setting:general:realization} ensures that for all $x \in \R$, $n \in \N$ it holds  that 
\begin{equation}\label{lemma:isru:lessk:bound}
    \realization{\theta_n,\infty}_{-4,0}(x)=\frac{x}{\sqrt{1+\xii(\frac{x}{n})^2}}.
\end{equation}
This shows that for all $x \in \R$ it holds that 
\begin{equation}\label{lemma:isru:zero:eq}
    \limsup\nolimits_{n\to \infty}|\realization{\theta_n,\infty}_{-4,0}(x)-x|=0.
\end{equation}
\Moreover \cref{lemma:isru:lessk:bound} implies that for all $x\in \R$, $n \in \N$ it holds that
\begin{equation}
  |\realization{\theta_n,\infty}_{-4,0}(x)| \leq x.
\end{equation}
This, \cref{lemma:isru:zero:eq}, and Lebesgue's dominated convergence theorem demonstrate that 
\begin{equation}
    \limsup\nolimits_{n \to \infty}\cL^\infty_{-4,0}(\theta_n)=\int_{\scra}^{\scrb}\limsup\nolimits_{n \to \infty} \rbr[\big]{ x - \realization{\theta_n,\infty}_{-4,0} (x) }^2  \, \d x=0.
\end{equation}
\end{cproof}
\begin{prop}\label{lemma:isru:zero2}
Assume \cref{H:setting:general}, assume  $\width>1$ and $\xii<3$, assume for all $x \in \R$ that $f(x)=x^2$, and let $(\theta_n)_{n \in \N}\subseteq \R^{\fd}$ satisfy for all $n \in \N$ that $\w{\theta_n}_1=-\w{\theta_n}_2=\nicefrac1n$,
$\v{\theta_n}_1=\v{\theta_n}_2= - (\xii+4)^{\frac52} n^2 (48 \xii)^{-1}$,
$\b{\theta_n}_1=\b{\theta_n}_2=\nicefrac12$,
$\c{\theta_n}= (\xii+4)^{\frac52} n^2 (48 \xii \sqrt{1+\nicefrac{\xii}{4}})^{-1} $, and
$\sum_{j=3}^\width|\w{\theta_n}_{j}|+ |\b{\theta_n}_{j}|+|\v{\theta_n}_{j}| =0.$
Then  
\begin{equation}
    \limsup\nolimits_{n \to \infty} \cL^\infty_{-4,0}(\theta_n)=0.
\end{equation}
\end{prop}
\begin{cproof}{lemma:isru:zero2}
\Nobs that \cref{H:setting:general:realization} ensures that for all $x \in \R$, $n \in \N$ it holds  that 
\begin{equation}\label{lemma:isru:zero2:realization}
    \realization{\theta_n,\infty}_{-4,0}(x)=
   - \frac{(\xii+4)^{\frac52} n^2}{48 \xii}\left( \frac{\frac xn +\frac12}{\sqrt{1+\xii(\frac{x}{n}+\frac12)^2}}+ 
   \frac{-\frac xn +\frac12}{\sqrt{1+\xii(-\frac{x}{n}+\frac12)^2}}  - \frac{1}{\sqrt{1+\frac{\xii}{4}}}\right).
\end{equation}
This and the fact that there exists $c\in (0,\infty)$ such that for all $x\in[-1,1]$ it holds that
\begin{equation}
    \bigg|\frac{x+\frac12}{\sqrt{(1+\xii (x+\frac12)^2)}} -
    \frac1{(\xii + 4)^{\frac12}} - \frac{8 x}{(\xii + 4)^{\frac32}} + \frac{24 \xii x^2}{(\xii + 4)^{\frac52}}  \bigg| \leq c|x^3|
\end{equation} 
assure that there exist $c,M\in (0,\infty)$ such that for all  $x\in [\scra,\scrb]$, $n\geq M$ it holds that 
\begin{equation}\label{lemma:isru:zero2:eq}
    |\realization{\theta_n,\infty}_{-4,0}(x)-x^2|\leq c\left|\frac1{n}\right|.
\end{equation}
This and Lebesgue's dominated convergence theorem demonstrate that 
\begin{equation}
    \limsup\nolimits_{n \to \infty}\cL^\infty_{-4,0}(\theta_n)=\int_{\scra}^{\scrb}\limsup\nolimits_{n \to \infty} \rbr[\big]{ x^2 - \realization{\theta_n,\infty}_{-4,0} (x) }^2  \, \d x=0.
\end{equation}
\end{cproof}
\begin{prop}\label{lemma:arctan:zero}
Assume \cref{H:setting:general}, assume for all $x \in \R$ that $f(x)=x$, and let $(\theta_n)_{n \in \N}\subseteq \R^{\fd}$ satisfy for all $n \in \N$ that $\w{\theta_n}_1=\nicefrac1n$,
$\v{\theta_n}_1=n$, and
$|\b{\theta_n}_1|+|\c{\theta_n}|+\sum_{j=2}^\width|\w{\theta_n}_{j}|+ |\b{\theta_n}_{j}|+|\v{\theta_n}_{j}| =0.$
Then  
\begin{equation}
    \limsup\nolimits_{n \to \infty} \cL^\infty_{-5,0}(\theta_n)=0.
\end{equation}
\end{prop}
\begin{cproof}{lemma:arctan:zero}
\Nobs that \cref{H:setting:general:realization} ensures that for all $x \in \R$, $n \in \N$ it holds  that 
\begin{equation}\label{lemma:arctan:lessk:bound}
    \realization{\theta_n,\infty}_{-5,0}(x)=n \arctan\Big(\frac{x}{n}\Big).
\end{equation}
This shows that for all $x \in \R$ it holds that 
\begin{equation}\label{lemma:arctan:zero:eq}
    \limsup\nolimits_{n\to \infty}|\realization{\theta_n,\infty}_{-5,0}(x)-x|=0.
\end{equation}
\Moreover \cref{lemma:arctan:lessk:bound} implies that for all $x\in \R$, $n \in \N$ it holds that
\begin{equation}
  |\realization{\theta_n,\infty}_{-5,0}(x)| \leq |x|.
\end{equation}
This, \cref{lemma:arctan:zero:eq}, and Lebesgue's dominated convergence theorem demonstrate that 
\begin{equation}
    \limsup\nolimits_{n \to \infty}\cL^\infty_{-5,0}(\theta_n)=\int_{\scra}^{\scrb}\limsup\nolimits_{n \to \infty} \rbr[\big]{ x - \realization{\theta_n,\infty}_{-5,0} (x) }^2  \, \d x=0.
\end{equation}
\end{cproof}
\begin{prop}\label{lemma:arctan:zero2}
Assume \cref{H:setting:general}, assume $\width>1$,  assume for all $x \in \R$ that $f(x)=x^2$, and let $(\theta_n)_{n \in \N}\subseteq \R^{\fd}$ satisfy for all $n \in \N$ that $\w{\theta_n}_1=-\w{\theta_n}_2=\nicefrac1n$,
$\v{\theta_n}_1=\v{\theta_n}_2=-2n^2$,
$\b{\theta_n}_1=\b{\theta_n}_2=1$,
$\c{\theta_n}= n^2 \pi $, and
$\sum_{j=3}^\width|\w{\theta_n}_{j}|+ |\b{\theta_n}_{j}|+|\v{\theta_n}_{j}| =0.$
Then  
\begin{equation}
    \limsup\nolimits_{n \to \infty} \cL^\infty_{-5,0}(\theta_n)=0.
\end{equation}
\end{prop}
\begin{cproof}{lemma:arctan:zero2}
\Nobs that \cref{H:setting:general:realization} ensures that for all $x \in \R$, $n \in \N$ it holds  that 
\begin{equation}\label{lemma:arctan:zero2:realization}
    \realization{\theta_n,\infty}_{-5,0}(x)=
   -2 n^2 \arctan\Big(\frac xn+1\Big)-2 n^2\arctan\Big(-\frac xn +1\Big)+ \pi n^2.
\end{equation}
This and the fact that there exists $c\in (0,\infty)$ such that for all $x\in[-1,1]$ it holds that
\begin{equation}
    \bigg| \arctan(x+1)-\frac{\pi}{4}-\frac x2 + \frac{x^2}{4} \bigg|\leq c |x^3|
\end{equation}
 assure that there exist $c,M \in (0,\infty)$ such that for all  $x\in [\scra,\scrb]$, $n \geq M$ it holds that 
\begin{equation}\label{lemma:arctan:zero2:eq}
  | \realization{\theta_n,\infty}_{-5,0}(x)-x^2|\leq c \left|\frac1{n}\right|.
\end{equation}
 This and Lebesgue's dominated convergence theorem demonstrate that 
\begin{equation}
    \limsup\nolimits_{n \to \infty}\cL^\infty_{-5,0}(\theta_n)=\int_{\scra}^{\scrb}\limsup\nolimits_{n \to \infty} \rbr[\big]{ x^2 - \realization{\theta_n,\infty}_{-5,0} (x) }^2  \, \d x=0.
\end{equation}
\end{cproof}
%
%
%
\begin{lemma}\label{H:noinf:sup}
Assume \cref{H:setting:general} and let $\theta\in\R^{\fd}$. Then there exists $c \in \R$ such that for all $k \in\{-1,-2,-4,-5\}$, $x \in \R$ it holds that
\begin{equation}
    |\realization{\theta,\infty}_{k,0}(x)|\leq c.
\end{equation}
\end{lemma}
\begin{cproof}{H:noinf:sup}
\Nobs that \cref{H:setting:general:activation} ensures that for all $j \in \{-1,-2,-4,-5\}$, $x\in\R$ it holds that
\begin{equation}
    |a_{j,0}(x)|\leq\begin{cases}
    \frac{\pi}{2} & \colon k=-5 \\
    \frac{|x|}{\sqrt{\xii x^2}}\leq\frac1{\sqrt\xii}  & \colon  k=-4\\
    \frac{\exp(x)+\exp(-x)}{\exp(x)+\exp(-x)}\leq1 & \colon  k=-2\\
    1 &  \colon k=-1.
    \end{cases}
\end{equation}
\Hence that 
for all $j \in \{-1,-2,-4,-5\}$, $x\in\R$ it holds that
\begin{equation}
    |\realization{\theta,\infty}_{j,0}(x)|\leq\begin{cases}
    |\c{\theta}|+\sum_{i=1}^\width \frac{\pi}{2}|\v{\theta}_i|\leq  \frac{\pi}{2} \norm{\theta}^2 & \colon  k=-5 \\
    |\c{\theta}|+\sum_{i=1}^\width \frac{1}{\sqrt{\xii}}|\v{\theta}_i|\leq  \frac{1}{\sqrt{\xii}} \norm{\theta}^2  &  \colon k=-4\\
    |\c{\theta}|+\sum_{i=1}^\width|\v{\theta}_i|\leq   \norm{\theta}^2 & \colon  k=-2\\
       |\c{\theta}|+\sum_{i=1}^\width|\v{\theta}_i|\leq   \norm{\theta}^2  & \colon  k=-1.
    \end{cases}
\end{equation}
This implies that for all $j \in \{-1,-2,-4,-5\}$, $x\in\R$ it holds that \begin{equation}
    |\realization{\theta,\infty}_{j,0}(x)|\leq \max\{\frac{\pi}{2} \norm{\theta}^2,\frac{1}{\sqrt{\xii}} \norm{\theta}^2\}.
\end{equation}
\end{cproof}
\begin{lemma}\label{H:noinf:log}
Assume \cref{H:setting:general}, let $k \in\{-1,-4,-5\}$, and assume for all $x \in \R$ that $f(x)=x$. Then
\begin{equation}\label{H:noinf:log:thesis}
   \big\{ \vartheta \in \R^{ \fd } \colon \cL^\infty_{k,0}( \vartheta ) = \inf\nolimits_{ \theta \in \R^{ \fd } } \cL^\infty_{k,0}( \theta ) \big\} = \emptyset.
\end{equation}
\end{lemma}
\begin{cproof}{H:noinf:log}
\Nobs that $f$ is real analytic.
\Moreover for all $\theta \in \R^{ \fd }$ it holds that $\realization{\theta,\infty}_{k,0}$ is real analytic.
\Moreover \cref{lemma:log:lessk}, \cref{lemma:isru:zero2}, and \cref{lemma:arctan:zero} ensure that for all $j \in \{-1,-4,-5\}$ it holds that $\inf_{ \theta \in \R^{ \fd } } \cL^\infty_{j,0}( \theta )=0$. We prove \cref{H:noinf:log:thesis} by contradiction. Assume that there exists $\vartheta \in \R^{ \fd }$ which satisfies that
\begin{equation}\label{H:noinf:log:ab}
    \cL^\infty_{k,0}(\vartheta)=0.
\end{equation} 
\Nobs that \cref{H:noinf:log:ab}  establishes that for all $x \in [\scra,\scrb]$ it holds that $
    f(x) = \realization{\vartheta,\infty}_{k,0}(x).$
Combining this with the fact that $f$ and $\realization{\vartheta,\infty}_{k,0}$ are real analytic implies for all $x \in \R$ that 
\begin{equation}\label{H:noinf:log:eq}
    f(x)=\realization{\vartheta,\infty}_{k,0}(x).
\end{equation}
\Nobs that \cref{H:noinf:sup} assures that there exists $c \in \R$ such that for all $j \in\{-1,-4,-5\}$, $x \in \R$ it holds that
\begin{equation}
    |\realization{\vartheta,\infty}_{j,0}(x)|\leq c.
\end{equation}
Combining this, the fact that $f$ is unbounded, and \cref{H:noinf:log:eq} shows the contradiction.
\end{cproof}
%
%
%
%
\begin{lemma}\label{H:noinf:log2}
Assume \cref{H:setting:general}, let $k \in \{-1,-2,-4,-5\}$, assume $\width>1$, assume $\xii <3$, and assume for all $x \in \R$ that $f(x)=x^2$. Then
\begin{equation}\label{H:noinf:log2:thesis}
    \big\{ \vartheta \in \R^{ \fd } \colon \cL^\infty_{k,0}( \vartheta ) = \inf\nolimits_{ \theta \in \R^{ \fd } } \cL^\infty_{k,0}( \theta ) \big\} = \emptyset.
\end{equation}
\end{lemma}
\begin{cproof}{H:noinf:log2}
\Nobs that $f$ is real analytic.
\Moreover for all $\theta \in \R^{ \fd }$ it holds that $\realization{\theta,\infty}_{k,0}$ is real analytic.
\Moreover \cref{lemma:log:lessk2}, \cref{lemma:tanh:lessk}, \cref{lemma:arctan:zero2}, and \cref{lemma:isru:zero2}  ensure that for all $j \in \{-1,-2,-4,-5\}$ it holds that $\inf_{ \theta \in \R^{ \fd } } \cL^\infty_{j,0}( \theta )=0$. We prove \cref{H:noinf:log2:thesis} by contradiction. Assume that there exists $\vartheta \in \R^{ \fd }$ such that
\begin{equation}\label{H:noinf:log2:ab}
    \cL^\infty_{k,0}(\vartheta)=0.
\end{equation} 
\Nobs that \cref{H:noinf:log2:ab} establishes that for all $x \in [\scra,\scrb]$ it holds that $f(x)=\realization{\vartheta,\infty}_{k,0}(x)$.
Combining this with the fact that $f$ and $\realization{\vartheta,\infty}_{k,0}$ are real analytic implies for all $x \in \R$ that \begin{equation}\label{H:noinf:log2:eq}
    f(x) = \realization{\vartheta,\infty}_{k,0}(x).
\end{equation}
\Moreover \cref{H:noinf:sup} assures that there exists $c \in \R$ such that for all $j \in\{-1,-2,-4,-5\}$, $x \in \R$ it holds that
\begin{equation}
    |\realization{\vartheta,\infty}_{j,0}(x)|\leq c.
\end{equation}
This shows that for all $j \in\{-1,-2,-4,-5\}$, $x \in \R$ it holds that $\sup_{x\in \R} |\realization{\vartheta,\infty}_{j,0}(x)|<\infty.$ 
Combining this with \cref{H:noinf:log2:eq} assures that
\begin{equation}
    \infty=\sup\nolimits_{x\in \R}|x^2|=\sup\nolimits_{x\in \R}|f(x)|=\sup\nolimits_{x\in \R}|\realization{\vartheta,\infty}_{k,0}(x)|<\infty.
\end{equation}
This contradiction establishes \cref{H:noinf:log2:thesis}.
\end{cproof}
\subsection{ANNs with rectified power unit activation}\label{repu}
\begin{prop}\label{lemma:repuk:zero}
Assume \cref{H:setting:general}, assume $\width>1$, let $k \in \N\backslash\{1\}$ satisfy for all $x \in \R$ that $f(x)=(\max\{x- \nicefrac{(\scra+\scrb)}2,0\})^{k-1}$, and let $(\theta_n)_{n \in \N} \subseteq \R^{\fd} $ satisfy for all $n \in \N$ that $\w{\theta_n}_1=\w{\theta_n}_2=1$,
$\b{\theta_n}_1=\nicefrac1n- \nicefrac{(\scra+\scrb)}2$,
$\b{\theta_n}_2=- \nicefrac{(\scra+\scrb)}2$,
$\v{\theta_n}_1=-\v{\theta_n}_2=\nicefrac nk$, and
$|\c{\theta_n}|+\sum_{j=3}^\width|\w{\theta_n}_{j}|+ |\b{\theta_n}_{j}|+|\v{\theta_n}_{j}| =0.$
Then  
\begin{equation}
    \limsup\nolimits_{n \to \infty} \cL^\infty_{k,0}(\theta_n)=0.
\end{equation}
\end{prop}
\begin{cproof}{lemma:repuk:zero}
\Nobs that \cref{H:setting:general:realization} ensures that for all $x \in \R$, $n \in \N$ it holds that
\begin{equation}\label{lemma:repuk:zero:realization}
\begin{split}
    \realization{\theta_n,\infty}_{k,0}(x)&=\frac nk(\max\{x + \nicefrac 1n - \nicefrac{(\scra+\scrb)}2,0\})^k-\frac nk(\max\{x- \nicefrac{(\scra+\scrb)}2,0\})^k\\
    &=
    \begin{cases}
    0
     &
    \colon
    x \in \left(-\infty ,-\frac1n+\frac{\scra+\scrb}2\right)\\
    \frac1k \sum_{i=0}^{k}{k \choose i} \left(x-\frac{\scra+\scrb}2\right)^{k-i}(\frac1{n})^{i-1}
     &
    \colon
    x \in \left[-\frac 1n+\frac{\scra+\scrb}2 ,\frac{\scra+\scrb}2\right)\\
    \left(x-\frac{\scra+\scrb}2\right)^{k-1}+\frac nk \sum_{i=2}^{k}{k \choose i} \left(x-\frac{\scra+\scrb}2\right)^{k-i}(\frac1{n})^{i}
     &
    \colon
    x \in \left[\frac{\scra+\scrb}2,\infty\right).
\end{cases}
\end{split}
\end{equation}
This implies that for all $x \in \R$ it holds that 
\begin{equation}\label{lemma:repuk:zero:eq}
    \limsup\nolimits_{n\to \infty}|\realization{\theta_n,\infty}_{k,0}(x)-(\max \{x-\nicefrac{(\scra+\scrb)}2,0\})^{k-1}|=0.
\end{equation}
\Moreover \cref{lemma:repuk:zero:realization} proves that
for all $x \in \R$, $n \in \N$ it holds that
\begin{equation}
    |\realization{\theta_n,\infty}_{k,0}(x)|\leq 
     \begin{cases}
     0
     &
    \colon
    x \in \left(-\infty ,-\frac1n+\frac{\scra+\scrb}2\right)
    \\
\frac 1k 
     &
  \colon
    x \in \left[-\frac 1n+\frac{\scra+\scrb}2 ,\frac{\scra+\scrb}2\right)\\
    \left(x-\frac{\scra+\scrb}2\right)^{k-1}+\frac 1k \sum_{i=2}^{k}{k \choose i} \left(x-\frac{\scra+\scrb}2\right)^{k-i}
     &
     \colon
    x \in \left[\frac{\scra+\scrb}2,\infty\right).
\end{cases}
\end{equation}
This, \cref{lemma:repuk:zero:eq}, and Lebesgue's dominated convergence theorem demonstrate that \begin{equation}
    \limsup\nolimits_{n \to \infty} \cL^\infty_{k,0}(\theta_n)=\int_{\scra}^{\scrb}\limsup\nolimits_{n \to \infty} \rbr[\big]{ (\max \{x-\nicefrac{(\scra+\scrb)}2,0\})^{k-1} - \realization{\theta_n,\infty}_{k,0} (x) }^2  \, \d x=0.
\end{equation}
\end{cproof}
%
%
%
%
\begin{lemma}\label{H:noinf:repu}
Assume \cref{H:setting:general}, assume $\width>1$, and let $k \in \N\backslash\{1\}$ satisfy for all $x \in \R$ that $f(x)=(\max\{x-\nicefrac{(\scra+\scrb)}2,0\})^{k-1}$. Then
\begin{equation}\label{H:noinf:repu:thesis}
    \big\{ \vartheta \in \R^{ \fd } \colon \cL^\infty_{k,0}( \vartheta ) = \inf\nolimits_{ \theta \in \R^{ \fd } } \cL^\infty_{k,0}( \theta ) \big\} = \emptyset.
\end{equation}
\end{lemma}
\begin{cproof}{H:noinf:repu}
\Nobs that the fact that $k \in \N\backslash\{1\}$ ensures that for all $\theta \in \R^{\fd}$ it holds that $\realization{\theta,\infty}_{k,0}\in C^{k-1}([\scra,\scrb],\R)$. Furthermore, \nobs that \cref{lemma:repuk:zero} establishes that  $\inf_{ \theta \in \R^{ \fd } } \cL^\infty_{k,0}( \theta )=0$.
We prove \cref{H:noinf:repu:thesis} by contradiction. Assume that there exists $\vartheta \in \R^{ \fd }$ such that 
\begin{equation}\label{H:noinf:repu:ab}
    \cL^\infty_{k,0}(\vartheta)=0.
\end{equation} 
\Nobs that \cref{H:noinf:repu:ab} implies that for all $x \in [\scra,\scrb]$ it holds that
\begin{equation}
    \realization{\vartheta,\infty}_{k,0}(x)=(\max\{x-\nicefrac{(\scra+\scrb)}2,0\})^{k-1}.
\end{equation}
This establishes that $\realization{\vartheta,\infty}_{k,0}\in C^{k-2}([\scra,\scrb],\R)\backslash C^{k-1}([\scra,\scrb],\R)$ which is a contradiction.
\end{cproof}
\subsection{ANNs with exponential linear unit activation}\label{elu}
\begin{prop}\label{prop:elu:inf}
Assume \cref{H:setting:general}, assume $\width>1$, assume for all $x \in \R$ that $f(x)=x^2$, and let $(\theta_n)_{n \in \N}\subseteq \R^{\fd}$ satisfy for all $n \in \N$ that $\w{\theta_n}_1=-\w{\theta_n}_2=-\nicefrac1n$,
$\b{\theta_n}_1=\b{\theta_n}_2=-2|\scrb|$,
$\v{\theta_n}_1=\v{\theta_n}_2=n^2e^{2|\scrb|}$,
$\c{\theta_n}=-2n^2(1-e^{2|\scrb|})$, and
$\sum_{j=3}^\width|\w{\theta_n}_{j}|+ |\b{\theta_n}_{j}|+|\v{\theta_n}_{j}| =0.$
Then  
\begin{equation}
    \limsup\nolimits_{n \to \infty} \cL^\infty_{-3,0}(\theta_n)=0.
\end{equation}
\end{prop}
\begin{cproof}{prop:elu:inf}
\Nobs that \cref{H:setting:general:realization} ensures that for all $x \in [\scra,\scrb]$, $n \in \N$ it holds that
\begin{equation}\label{lemma:elu:zero:realization}
\begin{split}
    \realization{\theta_n,\infty}_{-3,0}(x)& =e^{2|\scrb|}n^2( e^{-2|\scrb|-\frac xn}-1+
    e^{-2|\scrb|+\frac xn}-1)-2n^2(1-e^{2|\scrb|})\\
    &=n^2 ( e^{-\frac xn}+
    e^{\frac xn}-2).
\end{split}
\end{equation}
This and the fact that there exists $c\in (0,\infty)$ such that for all $x\in[-1,1]$ it holds that
\begin{equation}
    \bigg| e^{x}-1-x-\frac{x^2}2 \bigg| \leq c|x^3|
\end{equation}
 assure that there exist $c,M\in (0,\infty)$ such that for all $x\in [\scra,\scrb]$, $n\geq M$ it holds that 
\begin{equation}\label{lemma:elu:zero:eq}
    |\realization{\theta_n,\infty}_{-3,0}(x)-x^2|\leq c\left|\frac1{n}\right|.
\end{equation}
This
and Lebesgue's dominated convergence theorem demonstrate that
\begin{equation}
    \limsup\nolimits_{n \to \infty}\cL^\infty_{-3,0}(\theta_n)=\int_{\scra}^{\scrb}\limsup \nolimits_{n \to \infty} \rbr[\big]{ x^2 - \realization{\theta_n,\infty}_{-3,0} (x) }^2  \, \d x=0.
\end{equation}
\end{cproof}
%
%
%
%
\begin{lemma}\label{H:noinf:elu}
Assume \cref{H:setting:general}, assume $\width>1$, and assume for all $x \in \R$ that $f(x)=x^2$. Then
\begin{equation}\label{H:noinf:elu:thesis}
    \big\{ \vartheta \in \R^{ \fd } \colon \cL^\infty_{-3,0}( \vartheta ) = \inf\nolimits_{ \theta \in \R^{ \fd } } \cL^\infty_{-3,0}( \theta ) \big\} = \emptyset.
\end{equation}
\end{lemma}
\begin{cproof}{H:noinf:elu}
\Nobs that \cref{prop:elu:inf} implies that $\inf_{ \theta \in \R^{ \fd } } \cL^\infty_{-3,0}( \theta )=0$.
We prove \cref{H:noinf:elu:thesis} by contradiction.
We thus assume that there exists $\vartheta \in \R^{\fd}$ which satisfies that
 \begin{equation}\label{H:noinf:elu:ab}
    \cL^\infty_{-3,0}(\vartheta)=0.
\end{equation}
\Nobs that \cref{H:noinf:elu:ab} implies that for all $x \in [\scra,\scrb]$ it holds that
\begin{equation}\label{H:noinf:elu:realization}
    \realization{\vartheta,\infty}_{-3,0}(x)=x^2.
\end{equation}
\Hence $\realization{\vartheta,\infty}_{-3,0}\in C^\infty([\scra,\scrb],\R)$. This demonstrates that for all  $i \in \{1,2,\ldots,\width\}$ it holds that 
$\{x\in (\scra,\scrb) \colon \w{\vartheta}_i x + \b{\vartheta}_i\leq 0 \}\in\{(\scra,\scrb),\emptyset \}$.
For every $i \in \{1,2,\ldots,\width\}$ let $Q^i\subseteq \R$ satisfy $Q^i=\{x\in (\scra,\scrb) \colon \w{\vartheta}_i x + \b{\vartheta}_i\leq 0 \}$ and let $ S \subseteq\N$ satisfy $ S=\{i\in \{1,2,\ldots, \width\} \colon Q^i=(\scra,\scrb) \}$.
In the following we distinguish between the case $|S|=0$ and the case $|S|>0$.
We first establish the contradiction in the case \begin{equation}\label{H:noinf:elu:case1}
    |S|=0.
\end{equation}
\Nobs that \cref{H:noinf:elu:case1} ensures that there exist $\alpha, \beta \in \R$ such that for all $x \in [\scra,\scrb]$ it holds that $\realization{\vartheta,\infty}_{-3,0}(x)=\alpha x+ \beta$ which is a contradiction.
In the next step we establish the contradiction in the case \begin{equation}\label{H:noinf:elu:case2}
    |S|>0.
\end{equation}
For every $i \in S$ let $y_i\in \R$ satisfy  $y_i=\v{\vartheta}_i \exp\left(\w{\vartheta}_i(\nicefrac{\scra+\scrb}{2})+\b{\vartheta}_i\right)$, let $k \in \N$ satisfy $|S|=k$, and assume without loss of generality that $S=\{1,2,\ldots, k \}$.
\Nobs that \cref{H:noinf:elu:realization} proves that for all $n\in \N\cap(2,\infty)$, $x\in[\scra,\scrb]$ it holds that 
\begin{equation}\label{H:noinf:elu:der}
      0=(\realization{\vartheta,\infty}_{-3,0})^{(n)}(x)= \smallsum_{i \in S} \v{\vartheta}_i (\w{\vartheta}_i)^n\exp(\w{\vartheta}_i x + \b{\vartheta}_i) \qandq 2=(\realization{\vartheta,\infty}_{-3,0})^{(2)}(x).
\end{equation}
This implies that
\begin{equation} \label{H:noinf:elu:mat}
    \begin{split}
        (\w{\vartheta}_1)^2 y_1+\ldots+ (\w{\vartheta}_k)^2 y_k & =2,\\
        (\w{\vartheta}_1)^3 y_1+\ldots+ (\w{\vartheta}_k)^3 y_k & =0, \\
       \ldots & \\
    \andq    (\w{\vartheta}_1)^{k+2} y_1+\ldots+ (\w{\vartheta}_k)^{k+2} y_k & =0.
    \end{split}
\end{equation}
\Hence there exists $\eta=(\eta_1,\ldots,\eta_k) \in \R^k$ such that for all $j \in S$ it holds that
\begin{equation}
    \smallsum_{i=1}^k \eta_i (\w{\vartheta}_j)^{2+i}=(\w{\vartheta}_j)^{2}.
\end{equation}
This and \cref{H:noinf:elu:mat} show that $0=(\w{\vartheta}_1)^2 y_1+\ldots+ (\w{\vartheta}_k)^2 y_k  =2$ which is a contradiction.
\end{cproof}
\subsection{ANNs with softsign activation}\label{softsign}
\begin{prop}\label{lemma:softsign:zero2}
Assume \cref{H:setting:general}, assume $\width>1$, assume for all $x \in \R$ that $f(x)=x^2$, let $\scrc \in \R $ satisfy $\scrc=\max\{|\scra|,|\scrb|\}$,
and let $(\theta_n)_{n \in \N}\subseteq \R^{\fd}$ satisfy for all $n \in \N$ that $\w{\theta_n}_1=-\w{\theta_n}_2=\nicefrac1n$,
$\v{\theta_n}_1=\v{\theta_n}_2=\nicefrac{-(\scrc+1)^3n^2}{2}$,
$\b{\theta_n}_1=\b{\theta_n}_2=\scrc$,
$\c{\theta_n}=  \scrc(\scrc+1)^2 n^2 $, and 
$\sum_{j=3}^\width|\w{\theta_n}_{j}|+ |\b{\theta_n}_{j}|+|\v{\theta_n}_{j}| =0.$
Then  
\begin{equation}
    \limsup\nolimits_{n \to \infty} \cL^\infty_{-6,0}(\theta_n)=0.
\end{equation}
\end{prop}
\begin{cproof}{lemma:softsign:zero2}
\Nobs that \cref{H:setting:general:realization} ensures that for all $x \in \R$, $n \in \N$ it holds  that 
\begin{equation}\label{lemma:softsign:zero2:realization}
    \realization{\theta_n,\infty}_{-6,0}(x)=
   -\frac{(\scrc+1)^3n^2}{2}\left( \frac{\frac xn +\scrc}{1+|\frac xn +\scrc|}
   + \frac{-\frac xn +\scrc}{1+|-\frac xn +\scrc|}\right)
   + \scrc(\scrc+1)^2 n^2.
\end{equation}
This and the fact that there exists $c\in (0,\infty)$ such that for all $x\in[-1,1]$ it holds that
\begin{equation}
    \bigg| \frac{x +\scrc}{1+x +\scrc}-\frac{\scrc}{\scrc+1}-
    \frac{x}{(\scrc+1)^2}+\frac{x^2}{(\scrc+1)^3}\bigg| \leq c
    |x^3|
\end{equation}
 assure that there exist $c,M\in (0,\infty)$ such that for all $x\in [\scra,\scrb]$, $n\geq M$ it holds that 
\begin{equation}\label{lemma:softsign:zero2:eq}
    |\realization{\theta_n,\infty}_{-6,0}(x)-x^2|\leq c\left|\frac1{n}\right|.
\end{equation}
This and Lebesgue's dominated convergence theorem demonstrate that 
\begin{equation}
    \limsup\nolimits_{n \to \infty}\cL^\infty_{-6,0}(\theta_n)=\int_{\scra}^{\scrb}\limsup\nolimits_{n \to \infty} \rbr[\big]{ x^2 - \realization{\theta_n,\infty}_{-6,0} (x) }^2  \, \d x=0.
\end{equation}
\end{cproof}
%
%
%
%
\begin{lemma}\label{H:noinf:softsign}
Assume \cref{H:setting:general}, assume $\width>1$, and assume for all $x \in \R$ that $f(x)=x^2$. Then
\begin{equation}\label{H:noinf:softsign:thesis}
    \big\{ \vartheta \in \R^{ \fd } \colon \cL^\infty_{-6,0}( \vartheta ) = \inf\nolimits_{ \theta \in \R^{ \fd } } \cL^\infty_{-6,0}( \theta ) \big\} = \emptyset.
\end{equation}
\end{lemma}
\begin{cproof}{H:noinf:softsign}
\Nobs that \cref{lemma:softsign:zero2} implies that $\inf_{ \theta \in \R^{ \fd } } \cL^\infty_{-6,0}( \theta )=0$.
We prove \cref{H:noinf:softsign:thesis} by contradiction.
Assume that there exists $\vartheta \in \R^{\fd}$ which satisfies
\begin{equation}\label{H:noinf:softsign:ab}
    \cL^\infty_{-6,0}(\vartheta)=0.
\end{equation}
\Nobs that \cref{H:noinf:softsign:ab} implies that for all $x \in [\scra,\scrb]$ it holds that
\begin{equation}\label{H:noinf:softsign:realization}
    \realization{\vartheta,\infty}_{-6,0}(x)=x^2.
\end{equation}
Therefore we obtain that $\realization{\vartheta,\infty}_{-6,0}\in C^\infty([\scra,\scrb],\R)$. This demonstrates that for all  $i \in \{1,2,\ldots,\width\}$ it holds that 
$\{x\in (\scra,\scrb) \colon \w{\vartheta}_i x + \b{\vartheta}_i = 0 \}=\emptyset$.
This implies that for all $i \in \{1,2,\ldots, \width\}$ there exists $k_i\in \{-1,1\}$ which satisfies for all $x \in [\scra,\scrb]$ that 
\begin{equation}
    \realization{\vartheta,\infty}_{-6,0}(x)=\smallsum_{i=1}^{\width}\v{\vartheta}_i\frac{\w{\vartheta}_ix+\b{\vartheta}_i}{1+k_i(\w{\vartheta}_ix+\b{\vartheta}_i)}.
\end{equation}
Combining this with \cref{H:noinf:softsign:realization} proves that for all $n\in \N\cap(2,\infty)$, $x\in[\scra,\scrb]$ it holds that 
\begin{equation}\label{H:noinf:softsign:der}
    0=(\realization{\vartheta,\infty}_{-6,0})^{(n)}(x)= \smallsum_{i=1}^{\width} \v{\vartheta}_i \frac{n!\,\w{\vartheta}_i(-k_i \w{\vartheta}_i)^{n-1}}{(1+ k_i(\w{\vartheta}_i x + \b{\vartheta}_i))^{n+1}} \qandq  2=(\realization{\vartheta,\infty}_{-6,0})^{(2)}(x) .
\end{equation}
Let $S \subseteq\N$ satisfy $S=\{i\in \{1,2,\ldots, \width\} \colon \w{\vartheta}_i\neq0\}$ and for every $i \in S$ let $\xii_i \in  \R$ satisfy 
\begin{equation}
    \xii_i=\frac{-k_i \w{\vartheta}_i}{1+ k_i(\w{\vartheta}_i \frac{\scra+\scrb}{2} + \b{\vartheta}_i)}.
\end{equation}
\Nobs that \cref{H:noinf:softsign:realization} assures that $S\neq\emptyset$.
Let $k\in \N$ satisfy $|S|=k$ and assume without loss of generality that $S=\{1,2,\ldots,k\}$.
Combining this with \cref{H:noinf:softsign:der} shows that
\begin{equation}\label{H:noinf:softsign:mat}
    \begin{split}
        (\eta_1)^3 \frac{\v{\vartheta}_1}{\w{\vartheta}_1}+\ldots+ (\eta_k)^3 \frac{\v{\vartheta}_k}{\w{\vartheta}_k} & =2,\\
        (\eta_1)^4 \frac{\v{\vartheta}_1}{\w{\vartheta}_1}+\ldots+ (\eta_k)^4 \frac{\v{\vartheta}_k}{\w{\vartheta}_k} & =0,\\
       \ldots & \\
    \andq    (\eta_1)^{k+3} \frac{\v{\vartheta}_1}{\w{\vartheta}_1}+\ldots+ (\eta_k)^{k+3} \frac{\v{\vartheta}_k}{\w{\vartheta}_k} & =0.
    \end{split}
\end{equation}
\Hence that there exists $\scrc=(\scrc_1,\ldots,\scrc_k) \in \R^k$ such that for all $j \in S$ it holds that
\begin{equation}
    \smallsum_{i=1}^k \scrc_i (\eta_j)^{3+i}=(\eta_j)^{3}.
\end{equation}
This and \cref{H:noinf:softsign:mat} show that $0=(\eta_1)^3 \frac{\v{\vartheta}_1}{\w{\vartheta}_1}+\ldots+ (\eta_k)^3 \frac{\v{\vartheta}_k}{\w{\vartheta}_k} =2$ which is a contradiction.
\end{cproof}
\subsection{Divergence of GFs} \label{subdivergence}
\begin{lemma}\label{lemma:new}
Let $ \fd \in \N$, $\Theta \in C([0,\infty), \R^{\fd})$, $\cL \in C(\R^{\fd  },\R)$ satisfy $\{ \vartheta \in \R^{ \fd } \colon \cL( \vartheta ) = \inf_{ \theta \in \R^{ \fd } } \cL( \theta ) \} \allowbreak= \emptyset$ and $\liminf_{ t \to \infty } \cL( \Theta_t ) = \inf_{ \theta \in  \R^{ \fd } } \cL( \theta )$, let $\cG \colon \R^{\fd}\to \R^{\fd}$ be measurable, and assume for all $t \in [0,\infty)$ that $\cL(\Theta_t)=\cL(\Theta_0)-\int_0^t\norm{\cG(\Theta_s)}^2 \, \d s$.
Then
 \begin{equation}\label{lemma:new:thesis}
     \liminf\nolimits_{ t \to \infty } \| \Theta_t \| = \infty.
 \end{equation}
\end{lemma}
\begin{cproof}{lemma:new}
\Nobs that the assumption that for all $t \in [0,\infty)$ it holds that $\cL(\Theta_t)=\cL(\Theta_0)-\int_0^t\norm{\cG(\Theta_s)} ^{2}\, \d s$ assures that $[0, \infty) \ni t \mapsto \cL ( \Theta_t ) \in \R$ is non-increasing. This demonstrates that \begin{equation}\label{lemma:new:liminf}
    \limsup\nolimits_{ t \to \infty } \cL( \Theta_t ) = \liminf\nolimits_{ t \to \infty } \cL( \Theta_t ) = \inf\nolimits_{ \theta \in  \R^{ \fd } } \cL( \theta ).
\end{equation}
We prove \cref{lemma:new:thesis} by contradiction.
We thus assume that $\liminf\nolimits_{ t \to \infty } \| \Theta_t \| < \infty.$
Therefore, by compactness, there exist $\vartheta \in \R^{ \fd }$
and $\tau_n \in [0, \infty )$, $n \in \N$, which satisfy $\liminf_{n \to \infty} \tau_n = \infty$ and
\begin{equation}
    \limsup\nolimits_{ n \to \infty } \| \Theta_{\tau_n} - \vartheta \| = 0.
\end{equation}
Hence, continuity of $\cL$ shows that
$\limsup_{ n \to \infty } | \cL( \Theta_{\tau_n} ) - \cL( \vartheta ) | = 0.$
Combining this with \cref{lemma:new:liminf}  proves that
\begin{equation}
    \cL( \vartheta ) = \inf\nolimits_{ \theta \in \R^{ \fd }  } \cL( \theta ).
\end{equation}
This implies that
$\vartheta \in \{ \theta \in \R^{ \fd } \colon \cL( \theta ) = \inf_{ v \in \R^{ \fd } } \cL( v ) \}$ which is a contradiction.
\end{cproof}
\begin{cor}\label{theorem:globalminima:proof}
Assume \cref{H:setting:general}, assume $\width>1$, assume $\xii<3$, let $k \in \Z \backslash \N$, assume   for all $x \in \R$ that $f(x)=x^2$, and let $\Theta \in C([0,\infty), \R^{\fd})$ satisfy $\liminf_{ t \to \infty } \cL^\infty_{k,0}( \Theta_t ) = \inf_{ \theta \in  \R^{ \fd } } \cL^\infty_{k,0}( \theta )$
and $\forall \, t \in [0, \infty ) \colon \Theta_t = \Theta_0 - \int_0^t \cG_{k,0} ( \Theta_s ) \, \d s$. 
Then $\liminf_{ t \to \infty } \| \Theta_t \| = \infty$.
\end{cor}
\begin{cproof}{theorem:globalminima:proof}
\Nobs that, e.g., \cite[Lemma 3.1]{CheriditoJentzenRiekert2021} implies that for all $t \in [0,\infty)$ it holds that 
\begin{equation}\label{theorem:globalminima:proof:eq}
    \cL^\infty_{k,0}(\Theta_t)=\cL^\infty_{k,0}(\Theta_0)-\int_0^t\norm{\cG_{k,0}(\Theta_s)}^2 \, \d s.
\end{equation}
\Moreover \cref{H:noinf:soft:k2}, \cref{H:noinf:log2}, \cref{H:noinf:elu}, and \cref{H:noinf:softsign} assure that for all $j \in  \Z \backslash \N$ it holds that
 $\{ \vartheta \in \R^{ \fd } \colon \cL^\infty_{j,0}( \vartheta ) = \inf\nolimits_{ \theta \in \R^{ \fd } } \cL^\infty_{j,0}( \theta ) \} = \emptyset.$
 Combining this and \cref{theorem:globalminima:proof:eq} with \cref{lemma:new} proves that 
 $\liminf_{ t \to \infty } \| \Theta_t \| = \infty$.
\end{cproof}
\begin{cor}\label{theorem:globalminima:proof2}
Assume \cref{H:setting:general}, assume $\width>1$, let $k \in \N \backslash \{1\}$ satisfy for all $x \in \R$ that $f(x)=(\max\{x-\nicefrac{(\scra+\scrb)}2,0\})^{k-1}$, and let $\Theta \in C([0,\infty), \R^{\fd})$ satisfy $\liminf_{ t \to \infty } \cL^\infty_{k,0}( \Theta_t ) = \inf_{ \theta \in  \R^{\fd } } \cL^\infty_{k,0}( \theta )$
and $\forall \, t \in [0, \infty ) \colon \Theta_t = \Theta_0 - \int_0^t \cG_{k,0} ( \Theta_s ) \, \d s$.
Then $\liminf_{ t \to \infty } \| \Theta_t \| = \infty$.
\end{cor}
\begin{cproof}{theorem:globalminima:proof2}
\Nobs that, e.g., \cite[Lemma 3.1]{CheriditoJentzenRiekert2021} assures that for all $t \in [0,\infty)$ it holds that 
\begin{equation}\label{theorem:globalminima:proof2:eq}
    \cL^\infty_{k,0}(\Theta_t)=\cL^\infty_{k,0}(\Theta_0)-\int_0^t\norm{\cG_{k,0}(\Theta_s)}^2 \, \d s.
\end{equation}
\Moreover \cref{H:noinf:repu} shows that 
 $\{ \vartheta \in \R^{ \fd } \colon \cL^\infty_{k,0}( \vartheta ) = \inf\nolimits_{ \theta \in \R^{ \fd } } \cL^\infty_{k,0}( \theta ) \} = \emptyset.$
 Combining this and \cref{theorem:globalminima:proof2:eq} with \cref{lemma:new} demonstrates that 
 $\liminf_{ t \to \infty } \| \Theta_t \| = \infty$.
\end{cproof}
\begin{cor}\label{newtechnique:risk}
Assume \cref{H:setting:general},  assume $\width> 1$, assume for all $x \in [\scra,\scrb]$ that $f(x)=\indicator{(\nicefrac{(\scra+\scrb)}{2},\infty
)}(x)$, let $\gamma \in \R \backslash\{1\}$, and let  $\Theta \in C([0,\infty), \R^{\fd})$ satisfy $\liminf_{ t \to \infty } \cL^\infty_{1,\gamma}( \Theta_t ) =\allowbreak \inf_{ \theta \in  \R^{ \fd } } \cL^\infty_{1,\gamma}( \theta )$ and $\forall \, t \in [0, \infty ) \colon \Theta_t = \Theta_0 - \int_0^t \cG_{1,\gamma} ( \Theta_s ) \, \d s$.
Then $\liminf_{ t \to \infty } \| \Theta_t \| = \infty$.
\end{cor}
\begin{cproof}{newtechnique:risk}
\Nobs that the assumption that for all $t \in [0,\infty)$ it holds that $\Theta_t = \Theta_0 - \int_0^t \cG_{1,\gamma} ( \Theta_s ) \, \d s$
assures that for all $t \in [0,\infty)$ it holds that 
\begin{equation}\label{newtechnique:risk:eq}
    \cL^\infty_{1,\gamma}(\Theta_t)=\cL^\infty_{1,\gamma}(\Theta_0)-\int_0^t\norm{\cG_{1,\gamma}(\Theta_s)}^2 \, \d s
\end{equation}
(cf., e.g., Cheridito et al. \cite[Lemma 3.5]{CheriditoJentzenRiekert2021}).
\Moreover \cref{H:noinf} shows that  $\left\{ \vartheta \in \R^{ \fd } \colon \cL^\infty_{1,\gamma}( \vartheta ) = \inf\nolimits_{ \theta \in \R^{ \fd } } \cL^\infty_{1,\gamma}( \theta ) \right\} = \emptyset.$
 Combining this and \cref{newtechnique:risk:eq} with \cref{lemma:new} demonstrates that 
 $\liminf_{ t \to \infty } \| \Theta_t \| = \infty$.
\end{cproof}
\subsection{Divergence of GD}\label{subdivergenceGD}
\begin{lemma}\label{lemma:new2}
Let $ \fd \in \N$, $\cL \in C(\R^{\fd  },\R)$ satisfy $\{ \vartheta \in \R^{ \fd } \colon \cL( \vartheta ) = \inf_{ \theta \in \R^{ \fd } } \cL( \theta ) \} = \emptyset$ and let $\Theta=(\Theta_n)_{n \in \N_0} \colon \N_0 \to \R^{ \fd }$ satisfy  $\limsup_{ n \to \infty } \cL( \Theta_n ) = \inf_{ \theta \in \R^{ \fd } } \cL( \theta )$. Then
\begin{equation}\label{lemma:new2:thesis}
    \liminf\nolimits_{ n \to \infty } \| \Theta_n \| = \infty.
\end{equation}
\end{lemma}
\begin{cproof}{lemma:new2}
We prove \cref{lemma:new2:thesis} by contradiction.
We assume that $\liminf_{ n \to \infty } \| \Theta_n \| < \infty.$
Therefore, by compactness, there exist $\vartheta \in \R^{ \fd }$
and a strictly increasing $n \colon \N \to \N$ which satisfies that
\begin{equation}
  \limsup\nolimits_{ k \to \infty } \| \Theta_{n(k)} - \vartheta \| = 0.
\end{equation}
Hence, continuity of $\cL$ shows that
$\limsup_{ k \to \infty } | \cL( \Theta_{n(k)} ) - \cL( \vartheta ) | = 0.$
Combining this with the assumption that $\limsup_{ n \to \infty } \cL( \Theta_n ) = \inf_{ \theta \in \R^{ \fd } } \cL( \theta )$ proves that
\begin{equation}
    \cL( \vartheta ) = \inf\nolimits_{ \theta \in \R^{ \fd }  } \cL( \theta ).
\end{equation}
This implies that
$\vartheta \in \{ \theta \in \R^{ \fd } \colon \cL( \theta ) = \inf_{ v \in \R^{ \fd } } \cL( v ) \}$ which is a contradiction.
\end{cproof}
\begin{cor}\label{theorem:globalminima:proofb}
Assume \cref{H:setting:general}, assume $\width>1$, assume $\xii<3$, let $k \in \Z \backslash \N$, assume   for all $x \in \R$ that $f(x)=x^2$, and let $\Theta \colon \N_0 \to \R^{\fd}$ satisfy $\limsup_{ n \to \infty } \cL^\infty_{k,0}( \Theta_n ) = \inf_{ \theta \in  \R^{ \fd } } \cL^\infty_{k,0}( \theta )$.
Then $\liminf_{ n \to \infty } \| \Theta_n \| = \infty$.
\end{cor}
\begin{cproof}{theorem:globalminima:proofb}
\Nobs that \cref{H:noinf:soft:k2}, \cref{H:noinf:log2}, \cref{H:noinf:elu}, and \cref{H:noinf:softsign} assure that for all $j \in  \Z \backslash \N$ it holds that
 $\{ \vartheta \in \R^{ \fd } \colon \cL^\infty_{j,0}( \vartheta ) = \inf\nolimits_{ \theta \in \R^{ \fd } } \cL^\infty_{j,0}( \theta ) \} = \emptyset.$
 Combining this with \cref{lemma:new2} proves that 
 $\liminf_{ n \to \infty } \| \Theta_n \| = \infty$.
\end{cproof}
\begin{cor}\label{cor:globalminima:proofb}
Assume \cref{H:setting:general}, assume $\width>1$, assume $\xii<3$, let $k \in \Z \backslash \N$, assume   for all $x \in \R$ that $f(x)=x^2$, and let $\Theta \colon \N_0 \to \R^{\fd}$ satisfy $\liminf_{ n \to \infty } \cL^\infty_{k,0}( \Theta_n ) = \inf_{ \theta \in  \R^{ \fd } } \cL^\infty_{k,0}( \theta )$.
Then $\limsup_{ n \to \infty } \| \Theta_n \| = \infty$.
\end{cor}
\begin{cproof}{cor:globalminima:proofb}
\Nobs that the assumption that $\liminf_{ n \to \infty } \cL^\infty_{k,0}( \Theta_n ) = \inf_{ \theta \in  \R^{ \fd } } \cL^\infty_{k,0}( \theta )$ assures that there exists $n \colon \N \to \N $ which satisfies that 
\begin{equation}
    \lim\nolimits_{ j \to \infty } \cL^\infty_{k,0}( \Theta_{n(j)} ) = \inf\nolimits_{ \theta \in  \R^{ \fd } } \cL^\infty_{k,0}( \theta ).
\end{equation}
This and \cref{theorem:globalminima:proofb} imply that $\liminf_{ j \to \infty } \| \Theta_{n(j)} \| = \infty$. \Hence that $\limsup_{ j \to \infty } \allowbreak \| \Theta_j \| = \infty$.
\end{cproof}
\begin{cor}\label{theorem:globalminima:proof2b}
Assume \cref{H:setting:general}, assume $\width>1$, let $k \in \N \backslash \{1\}$ satisfy for all $x \in \R$ that $f(x)=(\max\{x-\nicefrac{(\scra+\scrb)}2,0\})^{k-1}$, and let $\Theta \colon \N_0 \to \R^{\fd}$ satisfy $\limsup_{ n \to \infty } \cL^\infty_{k,0}( \Theta_n ) = \inf_{ \theta \in  \R^{ \fd } } \cL^\infty_{k,0}( \theta )$.
Then $\liminf_{ n \to \infty } \| \Theta_n \| = \infty$.
\end{cor}
\begin{cproof}{theorem:globalminima:proof2b}
\Nobs that \cref{H:noinf:repu} demonstrates that 
 $\{ \vartheta \in \R^{ \fd } \colon \cL^\infty_{k,0}( \vartheta ) =\allowbreak \inf\nolimits_{ \theta \in \R^{ \fd } } \cL^\infty_{k,0}( \theta ) \} \allowbreak= \emptyset.$
 Combining this with \cref{lemma:new2} shows that 
 $\liminf_{ n \to \infty } \| \Theta_n \| = \infty$.
\end{cproof}
\begin{cor}\label{cor:globalminima:proof2b}
Assume \cref{H:setting:general}, assume $\width>1$, let $k \in \N \backslash \{1\}$ satisfy for all $x \in \R$ that $f(x)=(\max\{x-\nicefrac{(\scra+\scrb)}2,0\})^{k-1}$, and let $\Theta \colon \N_0 \to \R^{\fd}$ satisfy $\liminf_{ n \to \infty } \cL^\infty_{k,0}( \Theta_n ) = \inf_{ \theta \in  \R^{ \fd } } \cL^\infty_{k,0}( \theta )$.
Then $\limsup_{ n \to \infty } \| \Theta_n \| = \infty$.
\end{cor}
\begin{cproof}{cor:globalminima:proof2b}
\Nobs that the assumption that $\liminf_{ n \to \infty } \cL^\infty_{k,0}( \Theta_n ) = \inf_{ \theta \in  \R^{ \fd } } \cL^\infty_{k,0}( \theta )$ assures that there exists $n \colon \N \to \N $ which satisfies that 
\begin{equation}
    \lim\nolimits_{ j \to \infty } \cL^\infty_{k,0}( \Theta_{n(j)} ) = \inf\nolimits_{ \theta \in  \R^{ \fd } } \cL^\infty_{k,0}( \theta ).
\end{equation}
This and \cref{theorem:globalminima:proof2b} imply that $\liminf_{ j \to \infty } \| \Theta_{n(j)} \| = \infty$. \Hence that $\limsup_{ j \to \infty } \allowbreak \| \Theta_j \| = \infty$.
\end{cproof}
\begin{cor}\label{newtechnique:riskb}
Assume \cref{H:setting:general},  assume $\width> 1$, assume for all $x \in [\scra,\scrb]$ that $f(x)=\indicator{(\nicefrac{(\scra+\scrb)}{2},\infty
)}(x)$, let $\gamma \in \R \backslash\{1\}$, and let $\Theta \colon \N_0 \to \R^{\fd}$ satisfy $\limsup_{ n \to \infty } \cL^\infty_{1,\gamma}( \Theta_n ) = \inf_{ \theta \in  \R^{ \fd } } \cL^\infty_{1,\gamma}( \theta )$.
Then $\liminf_{ n \to \infty } \| \Theta_n \| = \infty$.
\end{cor}
\begin{cproof}{newtechnique:riskb}
    \Nobs that \cref{H:noinf} assures that  $\left\{ \vartheta \in \R^{ \fd } \colon \cL^\infty_{1,\gamma}( \vartheta ) = \inf\nolimits_{ \theta \in \R^{ \fd } } \cL^\infty_{1,\gamma}( \theta ) \right\} \allowbreak= \emptyset.$ This and \cref{lemma:new2} proves that $\liminf_{ n \to \infty } \| \Theta_n \| = \infty$.
\end{cproof}
\begin{cor}\label{cor:newtechnique:riskb}
Assume \cref{H:setting:general},  assume $\width> 1$, assume for all $x \in [\scra,\scrb]$ that $f(x)=\indicator{(\nicefrac{(\scra+\scrb)}{2},\infty
)}(x)$, let $\gamma \in \R \backslash\{1\}$, and let $\Theta \colon \N_0 \to \R^{\fd}$ satisfy $\liminf_{ n \to \infty } \cL^\infty_{1,\gamma}( \Theta_n ) = \inf_{ \theta \in  \R^{ \fd } } \cL^\infty_{1,\gamma}( \theta )$.
Then $\limsup_{ n \to \infty } \| \Theta_n \| = \infty$.
\end{cor}
\begin{cproof}{cor:newtechnique:riskb}
\Nobs that the assumption that $\liminf_{ n \to \infty } \cL^\infty_{1,\gamma}( \Theta_n ) = \inf_{ \theta \in  \R^{ \fd } } \cL^\infty_{1,\gamma}( \theta )$ assures that there exists $n \colon \N \to \N $ which satisfies that 
\begin{equation}
    \lim\nolimits_{ j \to \infty } \cL^\infty_{1,\gamma}( \Theta_{n(j)} ) = \inf\nolimits_{ \theta \in  \R^{ \fd } } \cL^\infty_{1,\gamma}( \theta ).
\end{equation}
This and \cref{newtechnique:riskb} imply that $\liminf_{ j \to \infty } \| \Theta_{n(j)} \| = \infty$. \Hence that $\limsup_{ j \to \infty } \allowbreak \| \Theta_j \| = \infty$.
\end{cproof}
\section{Blow up phenomena for data driven supervised learning problems}\label{discretemeasure}
In this section we analyze the existence of global minima in the case where the risk is defined using a discrete measure, the activation function is the standard logistic function, and the hidden layer is made up of one neuron.
In \cref{lemma_discrete:gen:non_existence} in \cref{sub3} and \cref{lemma_discrete:gen:non_existence2} in \cref{sub3},
assuming to have three non-strictly increasing or decreasing and non-constant data points $\fy_1,\fy_2,\fy_3\in \R$, we prove the non-existence of global minima of the risk function.
The proofs of \cref{lemma_discrete:gen:non_existence} and \cref{lemma_discrete:gen:non_existence2} are based 
on \cref{theorem_discrete:gen:inf} and on \cref{cor_discrete:gen:constant}.
In \cref{theorem_discrete:gen:inf} we find an upper bound for the infimum of the risk assuming that the data points do not coincide, $\max\{|\fy_1-\fy_2|,|\fy_3-\fy_2|\}>0$,
and are not non-strictly increasing or decreasing, $0\leq(\fy_1-\fy_2)(\fy_3-\fy_2)$.
In \cref{cor_discrete:gen:constant} we provide a lower bound for the risk in the case where the realization function is constant. The proof of \cref{cor_discrete:gen:constant} employs the elementary result for the first derivative of the realization function in
\cref{lemma_discrete:derivative}.

In \cref{theorem_discrete:existence} in \cref{sub1} and \cref{lemma_discrete:existence} in \cref{sub2} we establish the existence of global minima of the risk function in the case of two data points and in the case of three data points.
\subsection{Mathematical description of ANNs}
\begin{setting}\label{setting:discrete}
Let
 $\fw, \fb, \fv,  \fc \in C(\R^4, \R )$ 
 satisfy for all $\theta  = ( \theta_1 ,  \ldots, \theta_{4}) \in \R^{4}$ that $\w{\theta} 
 = \theta_{1}$, $\b{\theta} = \theta_{2}$,
$\v{\theta} = \theta_{3}$, and $\c{\theta} = \theta_{4}$,
let $A \colon \R \to \R$ satisfy for all $x \in \R$ that
\begin{equation} \label{dis:setting:relu:elu}
    A(x)=\frac{1}{1+\exp(-x)},
\end{equation}
for every $\theta \in \R^{4}$
let $\realization{\theta} \colon \R \to \R$
satisfy for all $x \in \R$ that
\begin{equation} \label{dis:setting:relu:realization}
    \realization{\theta} (x) = \c{\theta} + \v{\theta} \br[\big]{ A ( \w{\theta} x + \b{\theta} )},
\end{equation}
let $M \in \N $, $\fx=(\fx_1, \ldots, \fx_M)\in \R^{M}$, $\fy=(\fy_1, \ldots, \fy_M)\in \R^{M}$,
$\cL \in C( \R^{4}, \R)$
satisfy for all $\theta \in \R^{4}$ that
 \begin{equation}
    \cL ( \theta ) = \tfrac{1}{M}\smallsum_{i=1}^M \rbr[\big]{ \realization{\theta} (\fx_i) -\fy_i) } ^2,
\end{equation}
and let $ \operatorname{sgn} \colon \R \to \R$  satisfy for all $x \in \R$ that
\begin{equation}
     \operatorname{sgn}(x)= \begin{cases}
    1 &\colon x\geq 0 \\
    -1 &\colon x< 0.
    \end{cases}
\end{equation}
\end{setting}
\subsection{Existence of global minima for two data points}\label{sub1}
\begin{lemma}\label{theorem_discrete:existence}
Assume \cref{setting:discrete}, assume $M=2$, and assume $\fx_1<\fx_2$. Then there exists $\theta \in \R^4$ such that $\cL(\theta)=0$.
\end{lemma}
\begin{cproof}{theorem_discrete:existence}
Throughout this proof let $\vartheta\in \R^4$ satisfy 
\begin{equation}
\begin{split}
    \w{\vartheta} &=1, \qquad  \v{\vartheta} =(\fy_2-\fy_1)\Big(\frac{1}{1+\exp(-\fx_2)}-\frac{1}{1+\exp(-\fx_1)}\Big)^{-1}, \\
    \b{\vartheta} &=0, \qandq \c{\vartheta}=\fy_1-\frac{\v{\vartheta}}{1+\exp(-\fx_1)}.
\end{split}
\end{equation}
This implies that 
\begin{equation}
\begin{split}
    &\realization{\vartheta}(\fx_1)=\fy_1-\frac{\v{\vartheta}}{1+\exp(-\fx_1)}+\frac{\v{\vartheta}}{1+\exp(-\fx_1)}=\fy_1  \\
    \andq &\realization{\vartheta}(\fx_2)=\fy_1-\frac{\v{\vartheta}}{1+\exp(-\fx_1)}+\frac{\v{\vartheta}}{1+\exp(-\fx_2)}=\fy_2.
\end{split}
\end{equation}
\Hence that
 $\cL(\vartheta)=0$.
\end{cproof}
\subsection{Existence of global minima for three data points}\label{sub2}
\begin{lemma}\label{lemma_discrete:existence}
Assume \cref{setting:discrete} and assume $M=3$,  $\fx_1<\fx_2<\fx_3$, and $\min\{\fy_1,\fy_3\}< \fy_2<\max\{\fy_1,\fy_3\}$. Then there exists $\theta \in \R^4$ such that $\cL(\theta)=0$.
\end{lemma}
\begin{cproof}{lemma_discrete:existence}
Throughout this proof let $f \colon \R^2 \to \R$ satisfy for all $(x_1,x_2)\in \R^2$ that 
\begin{equation}\label{lemma_discrete:existence:f}
    f(x_1,x_2)=\frac{\big(\exp(-\fx_1 x_1)-\exp(-\fx_3x_1 )\big)\big(1+\exp(-\fx_2x_1-x_2)\big)}{\big(\exp(-\fx_1x_1 )-\exp(-\fx_2x_1 )\big)\big(1+\exp(-\fx_3x_1-x_2)\big)}.
\end{equation}
\Nobs that \cref{lemma_discrete:existence:f} assures that
\begin{equation}
    \liminf\nolimits_{x_1\to \infty}f(x_1,-\fx_3 x_1)=\infty \qandq \limsup\nolimits_{x_1\to -\infty}|f(x_1,0)-1|=0.
\end{equation}
Combining this with intermediate value theorem implies that for all $ y \in (1,\infty)$ there exist
$x_1,x_2 \in \R$ such that $ f(x_1,x_2)=y.$ 
\Nobs that $(\fy_3-\fy_1)(\fy_2-\fy_1)^{-1}>1$. 
Throughout this proof let $\vartheta\in \R^4$ satisfy 
\begin{equation}
\begin{split}
    (\w{\vartheta},\b{\vartheta})&=f^{-1}\Big(\frac{\fy_3-\fy_1}{\fy_2-\fy_1}\Big), \qquad 
     \c{\vartheta}=\fy_1-\frac{\v{\vartheta}}{1+\exp(-\w{\vartheta}\fx_1-\b{\vartheta})},
    \\ \andq
    \v{\vartheta} &=(\fy_2-\fy_1)
     \Big(\frac{(1+\exp(-\w{\vartheta}\fx_2-\b{\vartheta}))(1+\exp(-\w{\vartheta}\fx_1-\b{\vartheta})}{\exp(-\w{\vartheta}\fx_1-\b{\vartheta})-\exp(-\w{\vartheta}\fx_2-\b{\vartheta})}   \Big).
\end{split}
\end{equation}
This shows that 
\begin{equation}
\begin{split}
    \realization{\vartheta}(\fx_1) &=\fy_1-\frac{\v{\vartheta}}{1+\exp(-\w{\vartheta}\fx_1-\b{\vartheta})}+\frac{\v{\vartheta}}{1+\exp(-\w{\vartheta}\fx_1-\b{\vartheta})} =\fy_1,\\ 
    \realization{\vartheta}(\fx_2) &=\fy_1-(\fy_2-\fy_1)
    \frac{1+\exp(-\w{\vartheta}\fx_2-\b{\vartheta})}{\exp(-\w{\vartheta}\fx_1-\b{\vartheta})-\exp(-\w{\vartheta}\fx_2-\b{\vartheta})}   \\
     & \quad +(\fy_2-\fy_1)
     \frac{1+\exp(-\w{\vartheta}\fx_2-\b{\vartheta})}{\exp(-\w{\vartheta}\fx_1-\b{\vartheta})-\exp(-\w{\vartheta}\fx_2-\b{\vartheta})}  
    =\fy_2, \qand\\
    \realization{\vartheta}(\fx_3) &=\fy_1-(\fy_2-\fy_1)
     \frac{1+\exp(-\w{\vartheta}\fx_2-\b{\vartheta})}{\exp(-\w{\vartheta}\fx_1-\b{\vartheta})-\exp(-\w{\vartheta}\fx_2-\b{\vartheta})}    +(\fy_2-\fy_1) \\
     & \quad \frac{(1+\exp(-\w{\vartheta}\fx_2-\b{\vartheta}))(1+\exp(-\w{\vartheta}\fx_1-\b{\vartheta})}{(\exp(-\w{\vartheta}\fx_1-\b{\vartheta})-\exp(-\w{\vartheta}\fx_2-\b{\vartheta}))(1+\exp(-\w{\vartheta}\fx_3-\b{\vartheta}))}  \\
     &=\fy_1+(\fy_2-\fy_1)\frac{\fy_3-\fy_1}{\fy_2-\fy_1}=\fy_3.
\end{split}
\end{equation} 
\Hence that
 $\cL(\vartheta)=0$.
\end{cproof}
\subsection{Non-existence of global minima for three data points}\label{sub3}
\begin{prop}\label{lemma_discrete:derivative}
Assume \cref{setting:discrete}. 
Then it holds for all  $x \in \R$, $\theta\in \R^{4}$ with $\w{\theta}\v{\theta}\neq0$ that 
\begin{equation}
    \operatorname{sgn}(\w{\theta}\v{\theta})(\realization{\theta})'(x)>0.
\end{equation}
\end{prop}
\begin{cproof}{lemma_discrete:derivative}
 \Nobs that \cref{dis:setting:relu:realization} ensures that for all  $x \in \R$, $\theta\in \R^{4}$ it holds that 
 \begin{equation}
     (\realization{\theta})'(x)=(\w{\theta}\v{\theta})\frac{\exp(-\w{\theta}x-\b{\theta})}{(1+\exp(-\w{\theta}x-\b{\theta}))^2}.
 \end{equation}
 This implies that for all  $x \in \R$, $\theta\in \R^{4}$ with $\w{\theta}\v{\theta}\neq0$ it holds that 
\begin{equation}
    \operatorname{sgn}(\w{\theta}\v{\theta})(\realization{\theta})'(x)>0.
\end{equation}
\end{cproof}
\begin{prop}\label{theorem_discrete:gen:inf}
Assume \cref{setting:discrete}, assume $M=3$, $\fx_1<\fx_2<\fx_3$, and $-\max\{|\fy_1-\fy_2|,|\fy_3-\fy_2|\}<0\leq(\fy_1-\fy_2)(\fy_3-\fy_2)$,
let $I=\{i\in\{1,3\}\colon |\fy_i-\fy_2|=\max\{|\fy_1-\fy_2|,|\fy_3-\fy_2|\}\}$, let $j,k\in \N$ satisfy $j=\min I$, $k\in \{1,3\}\backslash\{j\}$, 
and let
$(\theta_n)_{n \in \N} \subseteq \R^{4}$ satisfy for all $n \in \N$ that 
$\w{\theta_n}=(2-j)n$, $\b{\theta_n}=(j-2)n\fx_j$,
$\c{\theta_n}=2\fy_j-(\nicefrac{(\fy_2+\fy_k)}{2})$, and
$\v{\theta_n}=\fy_2+\fy_k-2\fy_j$. Then 
\begin{equation}
    \limsup_{n \to \infty}\Big|3\cL(\theta_n)-2 \Big(\frac{\fy_2-\fy_k}{2} \Big)^2\Big|=0.
\end{equation}
\end{prop}
\begin{cproof}{theorem_discrete:gen:inf}
\Nobs that \cref{dis:setting:relu:realization} ensures thst for all $x \in \R$, $n\in \N$ it holds that
\begin{equation}
    \realization{\theta_n}(x)=\frac{\fy_2+\fy_k-2\fy_j}{1+\exp(n(j-2)(x-\fx_j))}+2\fy_j-\frac{\fy_2+\fy_k}{2}.
\end{equation}
This implies that for all $n \in \N$ it holds that
\begin{equation}
\begin{split}
    3\cL(\theta_n)&=(\realization{\theta_n}(\fx_j)-\fy_j)^2+(\realization{\theta_n}(\fx_2)-\fy_2)^2+(\realization{\theta_n}(\fx_k)-\fy_k)^2 \\
    &=0 + \Big(\frac{\fy_2+\fy_k-2\fy_j}{1+\exp(n(j-2)(\fx_2-\fx_j))}+2\fy_j-\frac{3\fy_2+\fy_k}{2}\Big)^2\\
    &\quad +\Big(\frac{\fy_2+\fy_k-2\fy_j}{1+\exp(n(j-2)(\fx_k-\fx_j))}+2\fy_j-\frac{\fy_2+3\fy_k}{2}\Big)^2.
\end{split}
\end{equation} 
\Hence that 
\begin{equation}
\begin{split}
    &\limsup_{n \to \infty}\Big|3\cL(\theta_n)-2 \Big(\frac{\fy_2-\fy_k}{2} \Big)^2\Big| \\
    & \quad =
   \Big| \Big(\fy_2+\fy_k-\frac{3\fy_2+\fy_k}{2}\Big)^2+\Big(\fy_2+\fy_k-\frac{\fy_2+3\fy_k}{2}\Big)^2 -2 \Big(\frac{\fy_2-\fy_k}{2} \Big)^2 \Big| \\
     & \quad =
    \Big|\Big(\frac{\fy_k-\fy_2}{2}\Big)^2+\Big(\frac{\fy_2-\fy_k}{2}\Big)^2-2 \Big(\frac{\fy_2-\fy_k}{2} \Big)^2 \Big|=0.
\end{split}
\end{equation}
\end{cproof}
\begin{prop}\label{cor_discrete:gen:constant}
Assume \cref{setting:discrete}, assume $M=3$,  $\fx_1<\fx_2<\fx_3$, and $\max\{|\fy_1-\fy_2|,|\fy_3-\fy_2|\}>0$, let $k \in \N$ satisfy $|\fy_k-\fy_2|=\min\{|\fy_1-\fy_2|,|\fy_3-\fy_2|\}$, and
let $\vartheta \in \R^4$ satisfy $\w{\vartheta}\v{\vartheta}=0$.
Then 
\begin{equation}
  3\cL(\vartheta)>2  \Big(\frac{\fy_2-\fy_k}{2} \Big)^2.
\end{equation}
\end{prop}
\begin{cproof}{cor_discrete:gen:constant}
\Nobs that the assumption that $\w{\vartheta}\v{\vartheta}=0$ assures that there exists $r \in \R$ which satisfies for all $x\in \R$ that $\realization{\vartheta}(x)=r$.
This implies that
\begin{equation}
    3\cL(\vartheta)=(r-\fy_1)^2+(r-\fy_2)^2+(r-\fy_3)^2.
\end{equation}
Assume without loss of generality that $k=1$. 
This and the assumption that $\max\{|\fy_1-\fy_2|,|\fy_3-\fy_2|\}>0$ assure that $|\fy_3-\fy_2|>0$.
\Hence that in the case $r=\fy_3$ it holds that
\begin{equation}\label{cor_discrete:gen:constant:eq}
    3\cL(\vartheta)=(\fy_3-\fy_1)^2+(\fy_3-\fy_2)^2\geq(\fy_3-\fy_1)^2+(\fy_1-\fy_2)^2>\frac12(\fy_2-\fy_1)^2
\end{equation}
and in the case $r \neq \fy_3$ it holds that
\begin{equation}
    3\cL(\vartheta)=(r-\fy_1)^2+(r-\fy_2)^2+(r-\fy_3)^2>(r-\fy_1)^2+(r-\fy_2)^2\geq \frac12 (\fy_2-\fy_1)^2.
\end{equation}
Combining this and \cref{cor_discrete:gen:constant:eq} shows that
\begin{equation}
      3\cL(\vartheta)>2  \Big(\frac{\fy_2-\fy_1}{2} \Big)^2.
\end{equation}
\end{cproof}
\begin{lemma}\label{lemma_discrete:gen:non_existence}
Assume \cref{setting:discrete} and assume $M=3$,  $\fx_1<\fx_2<\fx_3$, $\max\{|\fy_1-\fy_2|,|\fy_3-\fy_2|\}>0$, and $\min\{\fy_1,\fy_3\}\geq \fy_2$. Then 
\begin{equation}\label{lemma_discrete:gen:non_existence:thesis}
    \big\{ \vartheta \in \R^{4 } \colon \cL( \vartheta ) = \inf\nolimits_{ \theta \in \R^{ 4} } \cL( \theta ) \big\} = \emptyset.
\end{equation}
\end{lemma}
\begin{cproof}{lemma_discrete:gen:non_existence}
We prove \cref{lemma_discrete:gen:non_existence:thesis} by contradiction. We thus assume that there exists $\vartheta \in \R^4$ such that $\cL( \vartheta ) = \inf\nolimits_{ \theta \in \R^{ 4} } \cL( \theta ) $ and let $a_1,a_2,a_3 \in \R$ satisfy for all $n \in \{1,2,3\}$ that
$\realization{\vartheta}(\fx_n)=a_n$. \Nobs that  \cref{theorem_discrete:gen:inf} implies that 
\begin{equation}\label{lemma_discrete:gen:non_existence:inf}
    \cL(\vartheta)=\inf\nolimits_{ \theta \in \R^{ 4} } \cL( \theta ) \leq \frac23 \Big(\frac{\fy_2-\min\{\fy_1,\fy_3\}}{2} \Big)^2.
\end{equation}
This and \cref{cor_discrete:gen:constant} show that $\w{\vartheta}\v{\vartheta}\neq0$.
Combining this with \cref{lemma_discrete:derivative} demonstrates that for all $x \in \R $ it holds that $(\realization{\vartheta})'(x)\neq 0$.
In the following we distinguish between the case $\min_{x \in [\fx_1,\fx_3]}(\realization{\vartheta})'(x)>0$ and the case $\max_{x \in [\fx_1,\fx_3]}(\realization{\vartheta})'(x)<0$. We first establish the contradiction in the case  \begin{equation}\label{lemma_discrete:gen:non_existence:case1}
    \min\nolimits_{x \in [\fx_1,\fx_3]}(\realization{\vartheta})'(x)>0.
\end{equation}
\Nobs that \cref{lemma_discrete:gen:non_existence:case1} assures that $a_1<a_2<a_3$. Combining this and the assumption that $\min\{\fy_1,\fy_3\}\geq \fy_2$ proves that
\begin{equation}
    (a_1-\fy_1)^2+(a_2-\fy_2)^2>
    \begin{cases} (a_1-\fy_1)^2+(a_1-\fy_2)^2 
    & \colon a_1\geq \fy_2 \\
    (\fy_1-\fy_2)^2+ (a_2-\fy_2)^2
    & \colon a_1<\fy_2.
    \end{cases}
\end{equation}
This implies that \begin{equation}
    3\cL(\vartheta)\geq(a_1-\fy_1)^2+(a_2-\fy_2)^2>2 \Big(\frac{\fy_2-\fy_1}{2} \Big)^2\geq2\Big(\frac{\fy_2-\min\{\fy_1,\fy_3\}}{2} \Big)^2.
\end{equation}
Combining this with \cref{lemma_discrete:gen:non_existence:inf} shows that 
$\cL( \vartheta ) >\inf\nolimits_{ \theta \in \R^{ 4} } \cL( \theta )$ which is a contradiction.
In the next step we establish the contradiction in the case  
\begin{equation}\label{lemma_discrete:gen:non_existence:case2}
    \max\nolimits_{x \in [\fx_1,\fx_3]}(\realization{\vartheta})'(x)<0.
\end{equation}
\Nobs that \cref{lemma_discrete:gen:non_existence:case2} assures that $a_1>a_2>a_3$. Combining this and the assumption that $\min\{\fy_1,\fy_3\}\geq \fy_2$ proves that
\begin{equation}
    (a_2-\fy_2)^2+(a_3-\fy_3)^2>
    \begin{cases} 
    (a_3-\fy_2)^2+(a_3-\fy_3)^2 
    & \colon a_3\geq \fy_2 \\
    (a_2-\fy_2)^2 +(\fy_2-\fy_3)^2
    & \colon a_3<\fy_2.
    \end{cases}
\end{equation}
This implies that \begin{equation}
    3\cL(\vartheta)\geq (a_2-\fy_2)^2+(a_3-\fy_3)^2>2 \Big(\frac{\fy_2-\fy_3}{2} \Big)^2\geq 2\Big(\frac{\fy_2-\min\{\fy_1,\fy_3\}}{2} \Big)^2.
\end{equation}
Combining this with \cref{lemma_discrete:gen:non_existence:inf} shows that 
$\cL( \vartheta ) > \inf\nolimits_{ \theta \in \R^{ 4} } \cL( \theta )$ which is a contradiction.
\end{cproof}
\begin{lemma}\label{lemma_discrete:gen:non_existence2}
Assume \cref{setting:discrete} and assume $M=3$,  $\fx_1<\fx_2<\fx_3$, $\max\{|\fy_1-\fy_2|,|\fy_3-\fy_2|\}>0$, $\max\{\fy_1,\fy_3\}\leq \fy_2$. Then 
\begin{equation}\label{lemma_discrete:gen:non_existence2:thesis}
    \big\{ \vartheta \in \R^{4 } \colon \cL( \vartheta ) = \inf\nolimits_{ \theta \in \R^{ 4} } \cL( \theta ) \big\} = \emptyset.
\end{equation}
\end{lemma}
\begin{cproof}{lemma_discrete:gen:non_existence2}
We prove \cref{lemma_discrete:gen:non_existence2:thesis} by contradiction. Assume that there exists $\vartheta \in \R^4$ such that $\cL( \vartheta ) = \inf\nolimits_{ \theta \in \R^{ 4} } \cL( \theta ) $ and let $a_1,a_2,a_3 \in \R$ satisfy for all $n \in \{1,2,3\}$ that
$\realization{\vartheta}(\fx_n)=a_n$. \Nobs that  \cref{theorem_discrete:gen:inf} implies that 
\begin{equation}\label{lemma_discrete:gen:non_existence2:inf}
    \cL(\vartheta)=\inf\nolimits_{ \theta \in \R^{ 4} } \cL( \theta ) \leq \frac23 \Big(\frac{\fy_2-\max\{\fy_1,\fy_3\}}{2} \Big)^2.
\end{equation}
This and \cref{cor_discrete:gen:constant} show that $\w{\vartheta}\v{\vartheta}\neq0$.
Combining this with \cref{lemma_discrete:derivative} demonstrates that for all $x \in \R $ it holds that $(\realization{\vartheta})'(x)\neq 0$.
In the following we distinguish between the case $\min_{x \in [\fx_1,\fx_3]}(\realization{\vartheta})'(x)>0$ and the case $\max_{x \in [\fx_1,\fx_3]}(\realization{\vartheta})'(x)<0$. We first establish the contradiction in the case  
\begin{equation}\label{lemma_discrete:gen:non_existence2:case1}
    \min\nolimits_{x \in [\fx_1,\fx_3]}(\realization{\vartheta})'(x)>0.
\end{equation}
\Nobs that \cref{lemma_discrete:gen:non_existence2:case1} assures that $a_1<a_2<a_3$. Combining this and the assumption that $\max\{\fy_1,\fy_3\}\leq \fy_2$ proves that
\begin{equation}
    (a_2-\fy_2)^2+(a_3-\fy_3)^2>
    \begin{cases} (a_2-\fy_2)^2+(\fy_3-\fy_2)^2 
    & \colon a_3\geq \fy_2 \\
    (a_3-\fy_2)^2+ (a_3-\fy_3)^2
    & \colon a_3<\fy_2.
    \end{cases}
\end{equation}
This implies that \begin{equation}
    3\cL(\vartheta)\geq(a_2-\fy_2)^2+(a_3-\fy_3)^2>2 \Big(\frac{\fy_2-\fy_3}{2} \Big)^2\geq2\Big(\frac{\fy_2-\max\{\fy_1,\fy_3\}}{2} \Big)^2.
\end{equation}
Combining this with \cref{lemma_discrete:gen:non_existence2:inf} shows that 
$\cL( \vartheta ) >\inf\nolimits_{ \theta \in \R^{ 4} } \cL( \theta )$ which is a contradiction.
In the next step we establish the contradiction in the case  
\begin{equation}\label{lemma_discrete:gen:non_existence2:case2}
    \max\nolimits_{x \in [\fx_1,\fx_3]}(\realization{\vartheta})'(x)<0.
\end{equation}
\Nobs that \cref{lemma_discrete:gen:non_existence2:case2} assures that $a_1>a_2>a_3$. Combining this and the assumption that $\max\{\fy_1,\fy_3\}\leq \fy_2$ proves that
\begin{equation}
    (a_1-\fy_1)^2+(a_2-\fy_2)^2>
    \begin{cases} 
    (a_1-\fy_1)^2+(a_1-\fy_2)^2 
    & \colon a_1\leq \fy_2 \\
    (\fy_2-\fy_1)^2+ (a_2-\fy_2)^2
    & \colon a_1>\fy_2.
    \end{cases}
\end{equation}
This implies that \begin{equation}
    3\cL(\vartheta)\geq (a_1-\fy_1)^2+(a_2-\fy_2)^2>2 \Big(\frac{\fy_2-\fy_1}{2} \Big)^2\geq 2\Big(\frac{\fy_2-\max\{\fy_1,\fy_3\}}{2} \Big)^2.
\end{equation}
Combining this with \cref{lemma_discrete:gen:non_existence2:inf} shows that 
$\cL( \vartheta ) > \inf\nolimits_{ \theta \in \R^{ 4} } \cL( \theta )$ which is a contradiction.
\end{cproof}

\subsection*{Acknowledgements}
Adrian Riekert is gratefully acknowledged for several useful comments. The second author acknowledges funding by the Deutsche Forschungsgemeinschaft (DFG, German Research Foundation) under Germany’s Excellence Strategy EXC 2044-390685587, Mathematics Muenster: Dynamics-Geometry-Structure. 

\input{Bibliography.bbl}
\end{document}

%% file: 0_commands.tex
\DeclareMathAlphabet{\mathpzc}{OT1}{pzc}{m}{it}

\DeclareFontEncoding{LS1}{}{}
\DeclareFontSubstitution{LS1}{stix}{m}{n}
\DeclareMathAlphabet{\mathscr}{LS1}{stixscr}{m}{n}

\newcommand{\R}{\mathbb{R}}
\newcommand{\N}{\mathbb{N}}

\newcommand{\Z}{\mathbb{Z}}

\newcommand{\bbA}{\mathbb{A}}

\newcommand{\dens}{\mathfrak{p}}

\newcommand{\w}[1]{\mathfrak{w}^{#1}}
\renewcommand{\b}[1]{\mathfrak{b}^{#1}}
\renewcommand{\v}[1]{\mathfrak{v}^{#1}}
\newcommand{\q}[1]{\mathfrak{q}^{#1}}
\renewcommand{\d}[1]{\mathfrak{d}^{#1}}
\renewcommand{\c}[1]{\mathfrak{c}^{#1}}

\newcommand{\smallsum}{\textstyle\sum}

\newcommand{\m}{m}
\newcommand{\xii}{\xi}
\newcommand{\M}{M}

\newcommand{\cG}{\mathcal{G}}

\newcommand{\cL}{\mathcal{L}}

\newcommand{\bfa}{\mathbf{a}}
\newcommand{\bfb}{\mathbf{b}}

\newcommand{\bfm}{\mathbf{m}}
\newcommand{\bfn}{\mathbf{n}}

\newcommand{\bfq}{\mathbf{q}}

\newcommand{\act}{\mathbb{A}}

\newcommand{\scra}{\mathscr{a}}
\newcommand{\scrb}{\mathscr{b}}
\newcommand{\scrc}{\mathscr{c}}
\newcommand{\scrd}{\mathscr{d}}

\newcommand{\fC}{\mathfrak{C}}

\newcommand{\fa}{\mathfrak{a}}
\newcommand{\fb}{\mathfrak{b}}
\newcommand{\fc}{\mathfrak{c}}
\newcommand{\fd}{\mathfrak{d}}

\newcommand{\fv}{\mathfrak{v}}
\newcommand{\fw}{\mathfrak{w}}
\newcommand{\fx}{\mathscr{x}}
\newcommand{\fy}{\mathscr{y}}

\renewcommand{\emptyset}{\varnothing}

\DeclarePairedDelimiter{\norm}{\lVert}{\rVert}
\DeclarePairedDelimiter{\abs}{\lvert}{\rvert}
\DeclarePairedDelimiter{\rbr}{(}{)}
\DeclarePairedDelimiter{\br}{[}{]}
\DeclarePairedDelimiter{\cu}{\{}{\}}

\renewcommand{\d}{ \mathrm{d}}

\newcommand{\qand}{\qquad\text{and}}
\newcommand{\qandq}{\qquad\text{and}\qquad}
\newcommand{\andq}{\text{and}\qquad}

\newcommand{\indicator}[1]{\mathbbm{1}_{\smash{#1}}}

\newcommand{\realization}[1] {\mathscr{N} ^{ #1  }}

\newcommand{\width}{h}

\ExplSyntaxOn

\ExplSyntaxOn

\bool_new:N \g_noteobserve

\NewDocumentCommand{\setnote}{}{
  \bool_gset_true:N \g_noteobserve
}

\NewDocumentCommand{\setobserve}{}{
  \bool_gset_false:N \g_noteobserve
}

\NewDocumentCommand{\nobs}{ o }{
  \IfValueT{#1}{
    \str_if_eq:noTF {note} {#1} {
      \bool_gset_true:N \g_noteobserve
    } {
      \str_if_eq:noTF {Note} {#1} {
        \bool_gset_true:N \g_noteobserve
      } {
        \bool_gset_false:N \g_noteobserve
      }
    }
  }
  \bool_if:nTF { \g_noteobserve } {
    \bool_gset_false:N \g_noteobserve
    note
  } {
    \bool_gset_true:N \g_noteobserve
    observe
  }
  \IfValueF{#1}{~}
}

\NewDocumentCommand{\Nobs}{ o }{
  \IfValueT{#1}{
    \str_if_eq:noTF {note} {#1} {
      \bool_gset_true:N \g_noteobserve
    } {
      \str_if_eq:noTF {Note} {#1} {
        \bool_gset_true:N \g_noteobserve
      } {
        \bool_gset_false:N \g_noteobserve
      }
    }
  }
  \bool_if:nTF { \g_noteobserve } {
    \bool_gset_false:N \g_noteobserve
    Note
  } {
    \bool_gset_true:N \g_noteobserve
    Observe
  }
  \IfValueF{#1}{~}
}

\int_new:N \g_furthermore

\NewDocumentCommand{\Moreover}{ o o }{
  \IfValueT{#1}{
    \str_case:nn {#1} {
      {Furthermore} {\int_set:Nn {\g_furthermore} {0}}
      {Moreover} {\int_set:Nn {\g_furthermore} {1}}
      {In~addition} {\int_set:Nn {\g_furthermore} {2}}
      {note} {\bool_gset_true:N \g_noteobserve}
      {observe} {\bool_gset_false:N \g_noteobserve}
    }
    \IfValueT{#2}{
      \str_case:nn {#2} {
        {Furthermore} {\int_set:Nn {\g_furthermore} {0}}
        {Moreover} {\int_set:Nn {\g_furthermore} {1}}
        {In~addition} {\int_set:Nn {\g_furthermore} {2}}
        {note} {\bool_gset_true:N \g_noteobserve}
        {observe} {\bool_gset_false:N \g_noteobserve}
      }
    }
  }
  \int_case:nn { \int_mod:nn {\g_furthermore} {3} } {
    { 0 } { Furthermore,~\nobs that}
    { 1 } { Moreover,~\nobs that}
    { 2 } { In~addition,~\nobs that}
  }
  \int_incr:N \g_furthermore
  \IfValueF{#1}{~}
}

\bool_new:N \g_hencetherefore

\NewDocumentCommand{\hence}{}{
  \bool_if:nTF { \g_hencetherefore } {
    \bool_gset_false:N \g_hencetherefore
    hence~
  } {
    \bool_gset_true:N \g_hencetherefore
    therefore~
  }
}

\NewDocumentCommand{\Hence}{}{
  \bool_if:nTF { \g_hencetherefore } {
    \bool_gset_false:N \g_hencetherefore
    Hence,~we~obtain~
  } {
    \bool_gset_true:N \g_hencetherefore
    Therefore,~we~obtain~
  }
}

\seq_new:N \g_cflist_loaded
\seq_new:N \g_cflist_pending

\NewDocumentCommand{\cfadd}{ m }
{
  \seq_if_in:NnF \g_cflist_loaded { #1 } {
    \seq_if_in:NnF \g_cflist_pending { #1 } {
      \seq_gput_right:Nn \g_cflist_pending { #1 }
    }
  }
}

\NewDocumentCommand{\cfconsiderloaded}{ m }{
  \seq_gput_right:Nn \g_cflist_loaded {#1}
}

\NewDocumentCommand{\cfremove}{ m }
{
  \seq_gremove_all:Nn \g_cflist_pending { #1 }
}

\NewDocumentCommand{\cfload}{ o }
{
  \seq_if_empty:NTF \g_cflist_pending {\unskip} {
    (cf.\ \cref{\seq_use:Nn \g_cflist_pending {,}})\IfValueTF{#1}{#1~}{\unskip}
    \seq_gconcat:NNN \g_cflist_loaded \g_cflist_loaded \g_cflist_pending
    \seq_gclear:N \g_cflist_pending
  }
}

\NewDocumentCommand{\cfclear} {} {
  \seq_gclear:N \g_cflist_loaded
  \seq_gclear:N \g_cflist_pending
}

\NewDocumentCommand{\cfout}{ o }
{
  \seq_if_empty:NTF \g_cflist_pending {\unskip} {
    (cf.\ \cref{\seq_use:Nn \g_cflist_pending {,}})\IfValueTF{#1}{#1~}{\unskip}
    \seq_gclear:N \g_cflist_pending
  }
}

\NewDocumentCommand{\ifnocf} { m } {
  \seq_if_empty:NT \g_cflist_pending { #1 }
}

\ExplSyntaxOff

\NewDocumentEnvironment{cproof}{m}
{\begin{proof}[Proof of \cref{#1}]}%
{\noindent The proof of~\cref{#1} is thus complete.
\end{proof}}

\NewDocumentEnvironment{cproof2}{m}
{\begin{proof}[Proof of \cref{#1}]}%
{\noindent This completes the proof of~\cref{#1}.
\end{proof}}